\documentclass{amsproc}
\usepackage{amscd, amssymb}
\ifx\fivrm\undefined
  \font\fivrm=cmr5\relax
\fi

\font\fiverm=cmr5 
\catcode`@=11 \catcode`!=11

\@ifundefined{normalshape}{%
     \def\pdollar{{\ifdim \fontdimen\@ne\font >\z@ \sl \fi\char36}}}{%
        \def\pdollar{\text{\ifdim \fontdimen\@ne\font >\z@
              \sl \else \normalshape \fi\char36}}}

\expandafter\ifx\csname fiverm\endcsname\relax
  \let\fiverm\fivrm
\fi
  
\let\!latexendpicture=\endpicture 
\let\!latexframe=\frame
\let\!latexlinethickness=\linethickness
\let\!latexmultiput=\multiput
\let\!latexput=\put
 
\def\@picture(#1,#2)(#3,#4){%
  \@picht #2\unitlength
  \setbox\@picbox\hbox to #1\unitlength\bgroup 
  \let\endpicture=\!latexendpicture
  \let\frame=\!latexframe
  \let\linethickness=\!latexlinethickness
  \let\multiput=\!latexmultiput
  \let\put=\!latexput
  \hskip -#3\unitlength \lower #4\unitlength \hbox\bgroup}

\catcode`@=12 \catcode`!=12

\catcode`!=11 
 
\def\PiC{P\kern-.12em\lower.5ex\hbox{I}\kern-.075emC}
\def\PiCTeX{\PiC\kern-.11em\TeX}

\def\!ifnextchar#1#2#3{%
  \let\!testchar=#1%
  \def\!first{#2}%
  \def\!second{#3}%
  \futurelet\!nextchar\!testnext}
\def\!testnext{%
  \ifx \!nextchar \!spacetoken 
    \let\!next=\!skipspacetestagain
  \else
    \ifx \!nextchar \!testchar
      \let\!next=\!first
    \else 
      \let\!next=\!second 
    \fi 
  \fi
  \!next}
\def\\{\!skipspacetestagain} 
  \expandafter\def\\ {\futurelet\!nextchar\!testnext} 
\def\\{\let\!spacetoken= } \\  

\def\!tfor#1:=#2\do#3{%
  \edef\!fortemp{#2}%
  \ifx\!fortemp\!empty 
    \else
    \!tforloop#2\!nil\!nil\!!#1{#3}%
  \fi}
\def\!tforloop#1#2\!!#3#4{%
  \def#3{#1}%
  \ifx #3\!nnil
    \let\!nextwhile=\!fornoop
  \else
    #4\relax
    \let\!nextwhile=\!tforloop
  \fi 
  \!nextwhile#2\!!#3{#4}}

\def\!etfor#1:=#2\do#3{%
  \def\!!tfor{\!tfor#1:=}%
  \edef\!!!tfor{#2}%
  \expandafter\!!tfor\!!!tfor\do{#3}}

\def\!cfor#1:=#2\do#3{%
  \edef\!fortemp{#2}%
  \ifx\!fortemp\!empty 
  \else
    \!cforloop#2,\!nil,\!nil\!!#1{#3}%
  \fi}
\def\!cforloop#1,#2\!!#3#4{%
  \def#3{#1}%
  \ifx #3\!nnil
    \let\!nextwhile=\!fornoop 
  \else
    #4\relax
    \let\!nextwhile=\!cforloop
  \fi
  \!nextwhile#2\!!#3{#4}}

\def\!ecfor#1:=#2\do#3{%
  \def\!!cfor{\!cfor#1:=}%
  \edef\!!!cfor{#2}%
  \expandafter\!!cfor\!!!cfor\do{#3}}

\def\!empty{}
\def\!nnil{\!nil}
\def\!fornoop#1\!!#2#3{}

\def\!ifempty#1#2#3{%
  \edef\!emptyarg{#1}%
  \ifx\!emptyarg\!empty
    #2%
  \else
    #3%
  \fi}
 
\def\!getnext#1\from#2{%
  \expandafter\!gnext#2\!#1#2}%
\def\!gnext\\#1#2\!#3#4{%
  \def#3{#1}%
  \def#4{#2\\{#1}}%
  \ignorespaces}

\def\!getnextvalueof#1\from#2{%
  \expandafter\!gnextv#2\!#1#2}%
\def\!gnextv\\#1#2\!#3#4{%
  #3=#1%
  \def#4{#2\\{#1}}%
  \ignorespaces}

\def\!copylist#1\to#2{%
  \expandafter\!!copylist#1\!#2}
\def\!!copylist#1\!#2{%
  \def#2{#1}\ignorespaces}

\def\!wlet#1=#2{%
  \let#1=#2 
  \wlog{\string#1=\string#2}}
 
\def\!listaddon#1#2{%
  \expandafter\!!listaddon#2\!{#1}#2}
\def\!!listaddon#1\!#2#3{%
  \def#3{#1\\#2}}
 
\def\!rightappend#1\withCS#2\to#3{\expandafter\!!rightappend#3\!#2{#1}#3}
\def\!!rightappend#1\!#2#3#4{\def#4{#1#2{#3}}}

\def\!leftappend#1\withCS#2\to#3{\expandafter\!!leftappend#3\!#2{#1}#3}
\def\!!leftappend#1\!#2#3#4{\def#4{#2{#3}#1}}

\def\!lop#1\to#2{\expandafter\!!lop#1\!#1#2}
\def\!!lop\\#1#2\!#3#4{\def#4{#1}\def#3{#2}}

\def\!loop#1\repeat{\def\!body{#1}\!iterate}
\def\!iterate{\!body\let\!next=\!iterate\else\let\!next=\relax\fi\!next}
 
\def\!!loop#1\repeat{\def\!!body{#1}\!!iterate}
\def\!!iterate{\!!body\let\!!next=\!!iterate\else\let\!!next=\relax\fi\!!next}
 
\def\!removept#1#2{\edef#2{\expandafter\!!removePT\the#1}}
{\catcode`p=12 \catcode`t=12 \gdef\!!removePT#1pt{#1}}

\def\placevalueinpts of <#1> in #2 {%
  \!removept{#1}{#2}}
 
\def\!mlap#1{\hbox to 0pt{\hss#1\hss}}
\def\!vmlap#1{\vbox to 0pt{\vss#1\vss}}
 
\def\!not#1{%
  #1\relax
    \!switchfalse
  \else
    \!switchtrue
  \fi
  \if!switch
  \ignorespaces}

\let\!!!wlog=\wlog              
\def\wlog#1{}    

\newdimen\headingtoplotskip     
\newdimen\linethickness         
\newdimen\longticklength        
\newdimen\plotsymbolspacing     
\newdimen\shortticklength       
\newdimen\stackleading          
\newdimen\tickstovaluesleading  
\newdimen\totalarclength        
\newdimen\valuestolabelleading  

\newbox\!boxA                   
\newbox\!boxB                   
\newbox\!picbox                 
\newbox\!plotsymbol             
\newbox\!putobject              
\newbox\!shadesymbol            

\newcount\!countA               
\newcount\!countB               
\newcount\!countC               
\newcount\!countD               
\newcount\!countE               
\newcount\!countF               
\newcount\!countG               
\newcount\!fiftypt              
\newcount\!intervalno           
\newcount\!npoints              
\newcount\!nsegments            
\newcount\!ntemp                
\newcount\!parity               
\newcount\!scalefactor          
\newcount\!tfs                  
\newcount\!tickcase             

\newdimen\!Xleft                
\newdimen\!Xright               
\newdimen\!Xsave                
\newdimen\!Ybot                 
\newdimen\!Ysave                
\newdimen\!Ytop                 
\newdimen\!angle                
\newdimen\!arclength            
\newdimen\!areabloc             
\newdimen\!arealloc             
\newdimen\!arearloc             
\newdimen\!areatloc             
\newdimen\!bshrinkage           
\newdimen\!checkbot             
\newdimen\!checkleft            
\newdimen\!checkright           
\newdimen\!checktop             
\newdimen\!dimenA               
\newdimen\!dimenB               
\newdimen\!dimenC               
\newdimen\!dimenD               
\newdimen\!dimenE               
\newdimen\!dimenF               
\newdimen\!dimenG               
\newdimen\!dimenH               
\newdimen\!dimenI               
\newdimen\!distacross           
\newdimen\!downlength           
\newdimen\!dp                   
\newdimen\!dshade               
\newdimen\!dxpos                
\newdimen\!dxprime              
\newdimen\!dypos                
\newdimen\!dyprime              
\newdimen\!ht                   
\newdimen\!leaderlength         
\newdimen\!lshrinkage           
\newdimen\!midarclength         
\newdimen\!offset               
\newdimen\!plotheadingoffset    
\newdimen\!plotsymbolxshift     
\newdimen\!plotsymbolyshift     
\newdimen\!plotxorigin          
\newdimen\!plotyorigin          
\newdimen\!rootten              
\newdimen\!rshrinkage           
\newdimen\!shadesymbolxshift    
\newdimen\!shadesymbolyshift    
\newdimen\!tenAa                
\newdimen\!tenAc                
\newdimen\!tenAe                
\newdimen\!tshrinkage           
\newdimen\!uplength             
\newdimen\!wd                   
\newdimen\!wmax                 
\newdimen\!wmin                 
\newdimen\!xB                   
\newdimen\!xC                   
\newdimen\!xE                   
\newdimen\!xM                   
\newdimen\!xS                   
\newdimen\!xaxislength          
\newdimen\!xdiff                
\newdimen\!xleft                
\newdimen\!xloc                 
\newdimen\!xorigin              
\newdimen\!xpivot               

\newdimen\!xpos                 
\newdimen\!xprime               
\newdimen\!xright               
\newdimen\!xshade               
\newdimen\!xshift               
\newdimen\!xtemp                
\newdimen\!xunit                
\newdimen\!xxE                  
\newdimen\!xxM                  
\newdimen\!xxS                  
\newdimen\!xxloc                
\newdimen\!yB                   
\newdimen\!yC                   
\newdimen\!yE                   
\newdimen\!yM                   
\newdimen\!yS                   
\newdimen\!yaxislength          
\newdimen\!ybot                 
\newdimen\!ydiff                
\newdimen\!yloc                 
\newdimen\!yorigin              
\newdimen\!ypivot               
\newdimen\!ypos                 
\newdimen\!yprime               
\newdimen\!yshade               
\newdimen\!yshift               
\newdimen\!ytemp                
\newdimen\!ytop                 
\newdimen\!yunit                
\newdimen\!yyE                  
\newdimen\!yyM                  
\newdimen\!yyS                  
\newdimen\!yyloc                
\newdimen\!zpt                  

\newif\if!axisvisible           
\newif\if!gridlinestoo          
\newif\if!keepPO                
\newif\if!placeaxislabel        
\newif\if!switch                
\newif\if!xswitch               

\newtoks\!axisLaBeL             
\newtoks\!keywordtoks           

\newwrite\!replotfile           

\newhelp\!keywordhelp{The keyword mentioned in the error message in unknown. 
Replace NEW KEYWORD in the indicated response by the keyword that 
should have been specified.}    

\!wlet\!!origin=\!xM                   
\!wlet\!!unit=\!uplength               
\!wlet\!Lresiduallength=\!dimenG       
\!wlet\!Rresiduallength=\!dimenF       
\!wlet\!axisLength=\!distacross        
\!wlet\!axisend=\!ydiff                
\!wlet\!axisstart=\!xdiff              
\!wlet\!axisxlevel=\!arclength         
\!wlet\!axisylevel=\!downlength        
\!wlet\!beta=\!dimenE                  
\!wlet\!gamma=\!dimenF                 
\!wlet\!shadexorigin=\!plotxorigin     
\!wlet\!shadeyorigin=\!plotyorigin     
\!wlet\!ticklength=\!xS                
\!wlet\!ticklocation=\!xE              
\!wlet\!ticklocationincr=\!yE          
\!wlet\!tickwidth=\!yS                 
\!wlet\!totalleaderlength=\!dimenE     
\!wlet\!xone=\!xprime                  
\!wlet\!xtwo=\!dxprime                 
\!wlet\!ySsave=\!yM                    
\!wlet\!ybB=\!yB                       
\!wlet\!ybC=\!yC                       
\!wlet\!ybE=\!yE                       
\!wlet\!ybM=\!yM                       
\!wlet\!ybS=\!yS                       
\!wlet\!ybpos=\!yyloc                  
\!wlet\!yone=\!yprime                  
\!wlet\!ytB=\!xB                       
\!wlet\!ytC=\!xC                       
\!wlet\!ytE=\!downlength               
\!wlet\!ytM=\!arclength                
\!wlet\!ytS=\!distacross               
\!wlet\!ytpos=\!xxloc                  
\!wlet\!ytwo=\!dyprime                 

\!zpt=0pt                              
\!xunit=1pt
\!yunit=1pt
\!arearloc=\!xunit
\!areatloc=\!yunit
\!dshade=5pt
\!leaderlength=24in
\!tfs=256                              
\!wmax=5.3pt                           
\!wmin=2.7pt                           
\!xaxislength=\!xunit
\!xpivot=\!zpt
\!yaxislength=\!yunit 
\!ypivot=\!zpt
  \!dimenA=50pt \!fiftypt=\!dimenA     

\!rootten=3.162278pt                   
\!tenAa=8.690286pt                     
\!tenAc=2.773839pt                     
\!tenAe=2.543275pt                     

\def\!cosrotationangle{1}      
\def\!sinrotationangle{0}      
\def\!xpivotcoord{0}           
\def\!xref{0}                  
\def\!xshadesave{0}            
\def\!ypivotcoord{0}           
\def\!yref{0}                  
\def\!yshadesave{0}            
\def\!zero{0}                  

\let\wlog=\!!!wlog
\def\normalgraphs{%
  \longticklength=.4\baselineskip
  \shortticklength=.25\baselineskip
  \tickstovaluesleading=.25\baselineskip
  \valuestolabelleading=.8\baselineskip
  \linethickness=.4pt
  \stackleading=.17\baselineskip
  \headingtoplotskip=1.5\baselineskip
  \visibleaxes
  \ticksout
  \nogridlines
  \unloggedticks}
\def\setplotarea x from #1 to #2, y from #3 to #4 {%
  \!arealloc=\!M{#1}\!xunit \advance \!arealloc -\!xorigin
  \!areabloc=\!M{#3}\!yunit \advance \!areabloc -\!yorigin
  \!arearloc=\!M{#2}\!xunit \advance \!arearloc -\!xorigin
  \!areatloc=\!M{#4}\!yunit \advance \!areatloc -\!yorigin
  \!initinboundscheck
  \!xaxislength=\!arearloc  \advance\!xaxislength -\!arealloc
  \!yaxislength=\!areatloc  \advance\!yaxislength -\!areabloc
  \!plotheadingoffset=\!zpt
  \!dimenput {{\setbox0=\hbox{}\wd0=\!xaxislength\ht0=\!yaxislength\box0}}
     [bl] (\!arealloc,\!areabloc)}
\def\visibleaxes{%
  \def\!axisvisibility{\!axisvisibletrue}}

\def\!fixkeyword#1{%
  \errhelp=\!keywordhelp
  \errmessage{Unrecognized keyword `#1': \the\!keywordtoks{NEW KEYWORD}'}}

\!keywordtoks={enter `i\fixkeyword}

\def\fixkeyword#1{%
  \!nextkeyword#1 }

\def\axis {%
  \def\!nextkeyword##1 {%
    \expandafter\ifx\csname !axis##1\endcsname \relax
      \def\!next{\!fixkeyword{##1}}%
    \else
      \def\!next{\csname !axis##1\endcsname}%
    \fi
    \!next}%
  \!offset=\!zpt
  \!axisvisibility
  \!placeaxislabelfalse
  \!nextkeyword}

\def\!axisbottom{%
  \!axisylevel=\!areabloc
  \def\!tickxsign{0}%
  \def\!tickysign{-}%
  \def\!axissetup{\!axisxsetup}%
  \def\!axislabeltbrl{t}%
  \!nextkeyword}

\def\!axistop{%
  \!axisylevel=\!areatloc
  \def\!tickxsign{0}%
  \def\!tickysign{+}%
  \def\!axissetup{\!axisxsetup}%
  \def\!axislabeltbrl{b}%
  \!nextkeyword}

\def\!axisleft{%
  \!axisxlevel=\!arealloc
  \def\!tickxsign{-}%
  \def\!tickysign{0}%
  \def\!axissetup{\!axisysetup}%
  \def\!axislabeltbrl{r}%
  \!nextkeyword}

\def\!axisright{%
  \!axisxlevel=\!arearloc
  \def\!tickxsign{+}%
  \def\!tickysign{0}%
  \def\!axissetup{\!axisysetup}%
  \def\!axislabeltbrl{l}%
  \!nextkeyword}

\def\!axisshiftedto#1=#2 {%
  \if 0\!tickxsign
    \!axisylevel=\!M{#2}\!yunit
    \advance\!axisylevel -\!yorigin
  \else
    \!axisxlevel=\!M{#2}\!xunit
    \advance\!axisxlevel -\!xorigin
  \fi
  \!nextkeyword}

\def\!axisvisible{%
  \!axisvisibletrue  
  \!nextkeyword}

\def\!axisinvisible{%
  \!axisvisiblefalse
  \!nextkeyword}

\def\!axislabel#1 {%
  \!axisLaBeL={#1}%
  \!placeaxislabeltrue
  \!nextkeyword}

\expandafter\def\csname !axis/\endcsname{%
  \!axissetup 
  \if!placeaxislabel
    \!placeaxislabel
  \fi
  \if +\!tickysign 
    \!dimenA=\!axisylevel
    \advance\!dimenA \!offset 
    \advance\!dimenA -\!areatloc 
    \ifdim \!dimenA>\!plotheadingoffset
      \!plotheadingoffset=\!dimenA 
    \fi
  \fi}

\def\grid #1 #2 {%
  \!countA=#1\advance\!countA 1
  \axis bottom invisible ticks length <\!zpt> andacross quantity {\!countA} /
  \!countA=#2\advance\!countA 1
  \axis left   invisible ticks length <\!zpt> andacross quantity {\!countA} / }

\def\plotheading#1 {%
  \advance\!plotheadingoffset \headingtoplotskip
  \!dimenput {#1} [B] <.5\!xaxislength,\!plotheadingoffset>
    (\!arealloc,\!areatloc)}

\def\!axisxsetup{%
  \!axisxlevel=\!arealloc
  \!axisstart=\!arealloc
  \!axisend=\!arearloc
  \!axisLength=\!xaxislength
  \!!origin=\!xorigin
  \!!unit=\!xunit
  \!xswitchtrue
  \if!axisvisible 
    \!makeaxis
  \fi}

\def\!axisysetup{%
  \!axisylevel=\!areabloc
  \!axisstart=\!areabloc
  \!axisend=\!areatloc
  \!axisLength=\!yaxislength
  \!!origin=\!yorigin
  \!!unit=\!yunit
  \!xswitchfalse
  \if!axisvisible
    \!makeaxis
  \fi}

\def\!makeaxis{%
  \setbox\!boxA=\hbox{
    \beginpicture
      \!setdimenmode
      \setcoordinatesystem point at {\!zpt} {\!zpt}   
      \putrule from {\!zpt} {\!zpt} to
        {\!tickysign\!tickysign\!axisLength} 
        {\!tickxsign\!tickxsign\!axisLength}
    \endpicturesave <\!Xsave,\!Ysave>}%
    \wd\!boxA=\!zpt
    \!placetick\!axisstart}

\def\!placeaxislabel{%
  \advance\!offset \valuestolabelleading
  \if!xswitch
    \!dimenput {\the\!axisLaBeL} [\!axislabeltbrl]
      <.5\!axisLength,\!tickysign\!offset> (\!axisxlevel,\!axisylevel)
    \advance\!offset \!dp  
    \advance\!offset \!ht  
  \else
    \!dimenput {\the\!axisLaBeL} [\!axislabeltbrl]
      <\!tickxsign\!offset,.5\!axisLength> (\!axisxlevel,\!axisylevel)
  \fi
  \!axisLaBeL={}}

\def\arrow <#1> [#2,#3]{%
  \!ifnextchar<{\!arrow{#1}{#2}{#3}}{\!arrow{#1}{#2}{#3}<\!zpt,\!zpt> }}

\def\!arrow#1#2#3<#4,#5> from #6 #7 to #8 #9 {%
  \!xloc=\!M{#8}\!xunit   
  \!yloc=\!M{#9}\!yunit
  \!dxpos=\!xloc  \!dimenA=\!M{#6}\!xunit  \advance \!dxpos -\!dimenA
  \!dypos=\!yloc  \!dimenA=\!M{#7}\!yunit  \advance \!dypos -\!dimenA
  \let\!MAH=\!M
  \!setdimenmode
  \!xshift=#4\relax  \!yshift=#5\relax
  \!reverserotateonly\!xshift\!yshift
  \advance\!xshift\!xloc  \advance\!yshift\!yloc
  \!xS=-\!dxpos  \advance\!xS\!xshift
  \!yS=-\!dypos  \advance\!yS\!yshift
  \!start (\!xS,\!yS)
  \!ljoin (\!xshift,\!yshift)
  \!Pythag\!dxpos\!dypos\!arclength
  \!divide\!dxpos\!arclength\!dxpos  
  \!dxpos=32\!dxpos  \!removept\!dxpos\!!cos
  \!divide\!dypos\!arclength\!dypos  
  \!dypos=32\!dypos  \!removept\!dypos\!!sin
  \!halfhead{#1}{#2}{#3}
  \!halfhead{#1}{-#2}{-#3}
  \let\!M=\!MAH
  \ignorespaces}
  \def\!halfhead#1#2#3{%
    \!dimenC=-#1%
    \divide \!dimenC 2 
    \!dimenD=#2\!dimenC
    \!rotate(\!dimenC,\!dimenD)by(\!!cos,\!!sin)to(\!xM,\!yM)
    \!dimenC=-#1
    \!dimenD=#3\!dimenC
    \!dimenD=.5\!dimenD
    \!rotate(\!dimenC,\!dimenD)by(\!!cos,\!!sin)to(\!xE,\!yE)
    \!start (\!xshift,\!yshift)
    \advance\!xM\!xshift  \advance\!yM\!yshift
    \advance\!xE\!xshift  \advance\!yE\!yshift
    \!qjoin (\!xM,\!yM) (\!xE,\!yE) 
    \ignorespaces}

\def\betweenarrows #1#2 from #3 #4 to #5 #6 {%
  \!xloc=\!M{#3}\!xunit  \!xxloc=\!M{#5}\!xunit%
  \!yloc=\!M{#4}\!yunit  \!yyloc=\!M{#6}\!yunit%
  \!dxpos=\!xxloc  \advance\!dxpos by -\!xloc
  \!dypos=\!yyloc  \advance\!dypos by -\!yloc
  \advance\!xloc .5\!dxpos
  \advance\!yloc .5\!dypos
  \let\!MBA=\!M
  \!setdimenmode
  \ifdim\!dypos=\!zpt
    \ifdim\!dxpos<\!zpt \!dxpos=-\!dxpos \fi
    \put {\!lrarrows{\!dxpos}{#1}}#2{} at {\!xloc} {\!yloc}
  \else
    \ifdim\!dxpos=\!zpt
      \ifdim\!dypos<\!zpt \!dypos=-\!zpt \fi
      \put {\!udarrows{\!dypos}{#1}}#2{} at {\!xloc} {\!yloc}
    \fi
  \fi
  \let\!M=\!MBA
  \ignorespaces}

\def\!lrarrows#1#2{
  {\setbox\!boxA=\hbox{$\mkern-2mu\mathord-\mkern-2mu$}%
   \setbox\!boxB=\hbox{$\leftarrow$}\!dimenE=\ht\!boxB
   \setbox\!boxB=\hbox{}\ht\!boxB=2\!dimenE
   \hbox to #1{$\mathord\leftarrow\mkern-6mu
     \cleaders\copy\!boxA\hfil
     \mkern-6mu\mathord-$%
     \kern.4em $\vcenter{\box\!boxB}$$\vcenter{\hbox{#2}}$\kern.4em
     $\mathord-\mkern-6mu
     \cleaders\copy\!boxA\hfil
     \mkern-6mu\mathord\rightarrow$}}}

\def\!udarrows#1#2{
  {\setbox\!boxB=\hbox{#2}%
   \setbox\!boxA=\hbox to \wd\!boxB{\hss$\vert$\hss}%
   \!dimenE=\ht\!boxA \advance\!dimenE \dp\!boxA \divide\!dimenE 2
   \vbox to #1{\offinterlineskip
      \vskip .05556\!dimenE
      \hbox to \wd\!boxB{\hss$\mkern.4mu\uparrow$\hss}\vskip-\!dimenE
      \cleaders\copy\!boxA\vfil
      \vskip-\!dimenE\copy\!boxA
      \vskip\!dimenE\copy\!boxB\vskip.4em
      \copy\!boxA\vskip-\!dimenE
      \cleaders\copy\!boxA\vfil
      \vskip-\!dimenE \hbox to \wd\!boxB{\hss$\mkern.4mu\downarrow$\hss}
      \vskip .05556\!dimenE}}}

\def\putbar#1breadth <#2> from #3 #4 to #5 #6 {%
  \!xloc=\!M{#3}\!xunit  \!xxloc=\!M{#5}\!xunit%
  \!yloc=\!M{#4}\!yunit  \!yyloc=\!M{#6}\!yunit%
  \!dypos=\!yyloc  \advance\!dypos by -\!yloc
  \!dimenI=#2  
  \ifdim \!dimenI=\!zpt 
    \putrule#1from {#3} {#4} to {#5} {#6} 
  \else 
    \let\!MBar=\!M
    \!setdimenmode 
    \divide\!dimenI 2
    \ifdim \!dypos=\!zpt             
      \advance \!yloc -\!dimenI 
      \advance \!yyloc \!dimenI
    \else
      \advance \!xloc -\!dimenI 
      \advance \!xxloc \!dimenI
    \fi
    \putrectangle#1corners at {\!xloc} {\!yloc} and {\!xxloc} {\!yyloc}
    \let\!M=\!MBar 
  \fi
  \ignorespaces}

\def\setbars#1breadth <#2> baseline at #3 = #4 {%
  \edef\!barshift{#1}%
  \edef\!barbreadth{#2}%
  \edef\!barorientation{#3}%
  \edef\!barbaseline{#4}%
  \def\!bardobaselabel{\!bardoendlabel}%
  \def\!bardoendlabel{\!barfinish}%
  \let\!drawcurve=\!barcurve
  \!setbars}
\def\!setbars{%
  \futurelet\!nextchar\!!setbars}
\def\!!setbars{%
  \if b\!nextchar
    \def\!!!setbars{\!setbarsbget}%
  \else 
    \if e\!nextchar
      \def\!!!setbars{\!setbarseget}%
    \else
      \def\!!!setbars{\relax}%
    \fi
  \fi
  \!!!setbars}
\def\!setbarsbget baselabels (#1) {%
  \def\!barbaselabelorientation{#1}%
  \def\!bardobaselabel{\!!bardobaselabel}%
  \!setbars}
\def\!setbarseget endlabels (#1) {%
  \edef\!barendlabelorientation{#1}%
  \def\!bardoendlabel{\!!bardoendlabel}%
  \!setbars}

\def\!barcurve #1 #2 {%
  \if y\!barorientation
    \def\!basexarg{#1}%
    \def\!baseyarg{\!barbaseline}%
  \else
    \def\!basexarg{\!barbaseline}%
    \def\!baseyarg{#2}%
  \fi
  \expandafter\putbar\!barshift breadth <\!barbreadth> from {\!basexarg}
    {\!baseyarg} to {#1} {#2}
  \def\!endxarg{#1}%
  \def\!endyarg{#2}%
  \!bardobaselabel}

\def\!!bardobaselabel "#1" {%
  \put {#1}\!barbaselabelorientation{} at {\!basexarg} {\!baseyarg}
  \!bardoendlabel}
 
\def\!!bardoendlabel "#1" {%
  \put {#1}\!barendlabelorientation{} at {\!endxarg} {\!endyarg}
  \!barfinish}

\def\!barfinish{%
  \!ifnextchar/{\!finish}{\!barcurve}}

\def\putrectangle{%
  \!ifnextchar<{\!putrectangle}{\!putrectangle<\!zpt,\!zpt> }}
\def\!putrectangle<#1,#2> corners at #3 #4 and #5 #6 {%
  \!xone=\!M{#3}\!xunit  \!xtwo=\!M{#5}\!xunit%
  \!yone=\!M{#4}\!yunit  \!ytwo=\!M{#6}\!yunit%
  \ifdim \!xtwo<\!xone
    \!dimenI=\!xone  \!xone=\!xtwo  \!xtwo=\!dimenI
  \fi
  \ifdim \!ytwo<\!yone
    \!dimenI=\!yone  \!yone=\!ytwo  \!ytwo=\!dimenI
  \fi
  \!dimenI=#1\relax  \advance\!xone\!dimenI  \advance\!xtwo\!dimenI
  \!dimenI=#2\relax  \advance\!yone\!dimenI  \advance\!ytwo\!dimenI
  \let\!MRect=\!M
  \!setdimenmode
  \!shaderectangle
  \!dimenI=.5\linethickness
  \advance \!xone  -\!dimenI
  \advance \!xtwo   \!dimenI
  \putrule from {\!xone} {\!yone} to {\!xtwo} {\!yone} 
  \putrule from {\!xone} {\!ytwo} to {\!xtwo} {\!ytwo} 
  \advance \!xone   \!dimenI
  \advance \!xtwo  -\!dimenI%
  \advance \!yone  -\!dimenI
  \advance \!ytwo   \!dimenI
  \putrule from {\!xone} {\!yone} to {\!xone} {\!ytwo} 
  \putrule from {\!xtwo} {\!yone} to {\!xtwo} {\!ytwo} 
  \let\!M=\!MRect
  \ignorespaces}
 
\def\shaderectanglesoff{%
  \def\!shaderectangle{}%
  \ignorespaces}

\shaderectanglesoff
 
\def\!!shaderectangle{%
  \!dimenA=\!xtwo  \advance \!dimenA -\!xone
  \!dimenB=\!ytwo  \advance \!dimenB -\!yone
  \ifdim \!dimenA<\!dimenB
    \!startvshade (\!xone,\!yone,\!ytwo)
    \!lshade      (\!xtwo,\!yone,\!ytwo)
  \else
    \!starthshade (\!yone,\!xone,\!xtwo)
    \!lshade      (\!ytwo,\!xone,\!xtwo)
  \fi
  \ignorespaces}
  
\def\frame{%
  \!ifnextchar<{\!frame}{\!frame<\!zpt> }}
\long\def\!frame<#1> #2{%
  \beginpicture
    \setcoordinatesystem units <1pt,1pt> point at 0 0 
    \put {#2} [Bl] at 0 0 
    \!dimenA=#1\relax
    \!dimenB=\!wd \advance \!dimenB \!dimenA
    \!dimenC=\!ht \advance \!dimenC \!dimenA
    \!dimenD=\!dp \advance \!dimenD \!dimenA
    \let\!MFr=\!M
    \!setdimenmode
    \putrectangle corners at {-\!dimenA} {-\!dimenD} and {\!dimenB} {\!dimenC}
    \!setcoordmode
    \let\!M=\!MFr
  \endpicture
  \ignorespaces}
 
\def\rectangle <#1> <#2> {%
  \setbox0=\hbox{}\wd0=#1\ht0=#2\frame {\box0}}

\def\plot{%
  \!ifnextchar"{\!plotfromfile}{\!drawcurve}}
\def\!plotfromfile"#1"{%
  \expandafter\!drawcurve \input #1 /}

\def\setquadratic{%
  \let\!drawcurve=\!qcurve
  \let\!!Shade=\!!qShade
  \let\!!!Shade=\!!!qShade}

\def\setlinear{%
  \let\!drawcurve=\!lcurve
  \let\!!Shade=\!!lShade
  \let\!!!Shade=\!!!lShade}

\def\sethistograms{%
  \let\!drawcurve=\!hcurve}

\def\!qcurve #1 #2 {%
  \!start (#1,#2)
  \!Qjoin}
\def\!Qjoin#1 #2 #3 #4 {%
  \!qjoin (#1,#2) (#3,#4)             
  \!ifnextchar/{\!finish}{\!Qjoin}}

\def\!lcurve #1 #2 {%
  \!start (#1,#2)
  \!Ljoin}
\def\!Ljoin#1 #2 {%
  \!ljoin (#1,#2)                    
  \!ifnextchar/{\!finish}{\!Ljoin}}

\def\!finish/{\ignorespaces}

\def\!hcurve #1 #2 {%
  \edef\!hxS{#1}%
  \edef\!hyS{#2}%
  \!hjoin}
\def\!hjoin#1 #2 {%
  \putrectangle corners at {\!hxS} {\!hyS} and {#1} {#2}
  \edef\!hxS{#1}%
  \!ifnextchar/{\!finish}{\!hjoin}}

\def\vshade #1 #2 #3 {%
  \!startvshade (#1,#2,#3)
  \!Shadewhat}

\def\hshade #1 #2 #3 {%
  \!starthshade (#1,#2,#3)
  \!Shadewhat}

\def\!Shadewhat{%
  \futurelet\!nextchar\!Shade}
\def\!Shade{%
  \if <\!nextchar
    \def\!nextShade{\!!Shade}%
  \else
    \if /\!nextchar
      \def\!nextShade{\!finish}%
    \else
      \def\!nextShade{\!!!Shade}%
    \fi
  \fi
  \!nextShade}
\def\!!lShade<#1> #2 #3 #4 {%
  \!lshade <#1> (#2,#3,#4)                 
  \!Shadewhat}
\def\!!!lShade#1 #2 #3 {%
  \!lshade (#1,#2,#3)
  \!Shadewhat} 
\def\!!qShade<#1> #2 #3 #4 #5 #6 #7 {%
  \!qshade <#1> (#2,#3,#4) (#5,#6,#7)      
  \!Shadewhat}
\def\!!!qShade#1 #2 #3 #4 #5 #6 {%
  \!qshade (#1,#2,#3) (#4,#5,#6)
  \!Shadewhat} 

\setlinear

\def\setdashpattern <#1>{%
  \def\!Flist{}\def\!Blist{}\def\!UDlist{}%
  \!countA=0
  \!ecfor\!item:=#1\do{%
    \!dimenA=\!item\relax
    \expandafter\!rightappend\the\!dimenA\withCS{\\}\to\!UDlist%
    \advance\!countA  1
    \ifodd\!countA
      \expandafter\!rightappend\the\!dimenA\withCS{\!Rule}\to\!Flist%
      \expandafter\!leftappend\the\!dimenA\withCS{\!Rule}\to\!Blist%
    \else 
      \expandafter\!rightappend\the\!dimenA\withCS{\!Skip}\to\!Flist%
      \expandafter\!leftappend\the\!dimenA\withCS{\!Skip}\to\!Blist%
    \fi}%
  \!leaderlength=\!zpt
  \def\!Rule##1{\advance\!leaderlength  ##1}%
  \def\!Skip##1{\advance\!leaderlength  ##1}%
  \!Flist%
  \ifdim\!leaderlength>\!zpt 
  \else
    \def\!Flist{\!Skip{24in}}\def\!Blist{\!Skip{24in}}\ignorespaces
    \def\!UDlist{\\{\!zpt}\\{24in}}\ignorespaces
    \!leaderlength=24in
  \fi
  \!dashingon}

\def\!dashingon{%
  \def\!advancedashing{\!!advancedashing}%
  \def\!drawlinearsegment{\!lineardashed}%
  \def\!puthline{\!putdashedhline}%
  \def\!putvline{\!putdashedvline}%
  \ignorespaces}%
\def\!dashingoff{%
  \def\!advancedashing{\relax}%
  \def\!drawlinearsegment{\!linearsolid}%
  \def\!puthline{\!putsolidhline}%
  \def\!putvline{\!putsolidvline}%
  \ignorespaces}

\def\setdots{%
  \!ifnextchar<{\!setdots}{\!setdots<5pt>}}
\def\!setdots<#1>{%
  \!dimenB=#1\advance\!dimenB -\plotsymbolspacing
  \ifdim\!dimenB<\!zpt
    \!dimenB=\!zpt
  \fi
\setdashpattern <\plotsymbolspacing,\!dimenB>}
 
\def\setdotsnear <#1> for <#2>{%
  \!dimenB=#2\relax  \advance\!dimenB -.05pt  
  \!dimenC=#1\relax  \!countA=\!dimenC 
  \!dimenD=\!dimenB  \advance\!dimenD .5\!dimenC  \!countB=\!dimenD
  \divide \!countB  \!countA
  \ifnum 1>\!countB 
    \!countB=1
  \fi
  \divide\!dimenB  \!countB
  \setdots <\!dimenB>}
 
\def\setdashes{%
  \!ifnextchar<{\!setdashes}{\!setdashes<5pt>}}
\def\!setdashes<#1>{\setdashpattern <#1,#1>}
 
\def\setdashesnear <#1> for <#2>{%
  \!dimenB=#2\relax  
  \!dimenC=#1\relax  \!countA=\!dimenC 
  \!dimenD=\!dimenB  \advance\!dimenD .5\!dimenC  \!countB=\!dimenD
  \divide \!countB  \!countA
  \ifodd \!countB 
  \else 
    \advance \!countB  1
  \fi
  \divide\!dimenB  \!countB
  \setdashes <\!dimenB>}
 
\def\setsolid{%
  \def\!Flist{\!Rule{24in}}\def\!Blist{\!Rule{24in}}%
  \def\!UDlist{\\{24in}\\{\!zpt}}%
  \!dashingoff}  
\setsolid

\def\!divide#1#2#3{%
  \!dimenB=#1
  \!dimenC=#2
  \!dimenD=\!dimenB
  \divide \!dimenD \!dimenC
  \!dimenA=\!dimenD
  \multiply\!dimenD \!dimenC
  \advance\!dimenB -\!dimenD
  \!dimenD=\!dimenC
    \ifdim\!dimenD<\!zpt \!dimenD=-\!dimenD 
  \fi
  \ifdim\!dimenD<64pt
    \!divstep[\!tfs]\!divstep[\!tfs]%
  \else 
    \!!divide
  \fi
  #3=\!dimenA\ignorespaces}

\def\!!divide{%
  \ifdim\!dimenD<256pt
    \!divstep[64]\!divstep[32]\!divstep[32]%
  \else 
    \!divstep[8]\!divstep[8]\!divstep[8]\!divstep[8]\!divstep[8]%
    \!dimenA=2\!dimenA
  \fi}

\def\!divstep[#1]{
  \!dimenB=#1\!dimenB
  \!dimenD=\!dimenB
    \divide \!dimenD by \!dimenC
  \!dimenA=#1\!dimenA
    \advance\!dimenA by \!dimenD%
  \multiply\!dimenD by \!dimenC
    \advance\!dimenB by -\!dimenD}
 
\def\Divide <#1> by <#2> forming <#3> {%
  \!divide{#1}{#2}{#3}}

\def\circulararc{%
  \ellipticalarc axes ratio 1:1 }

\def\ellipticalarc axes ratio #1:#2 #3 degrees from #4 #5 center at #6 #7 {%
  \!angle=#3pt\relax
  \ifdim\!angle>\!zpt 
    \def\!sign{}
  \else 
    \def\!sign{-}\!angle=-\!angle
  \fi
  \!xxloc=\!M{#6}\!xunit
  \!yyloc=\!M{#7}\!yunit     
  \!xxS=\!M{#4}\!xunit
  \!yyS=\!M{#5}\!yunit
  \advance\!xxS -\!xxloc
  \advance\!yyS -\!yyloc
  \!divide\!xxS{#1pt}\!xxS 
  \!divide\!yyS{#2pt}\!yyS 
  \let\!MC=\!M
  \!setdimenmode
  \!xS=#1\!xxS  \advance\!xS\!xxloc
  \!yS=#2\!yyS  \advance\!yS\!yyloc
  \!start (\!xS,\!yS)%
  \!loop\ifdim\!angle>14.9999pt
    \!rotate(\!xxS,\!yyS)by(\!cos,\!sign\!sin)to(\!xxM,\!yyM) 
    \!rotate(\!xxM,\!yyM)by(\!cos,\!sign\!sin)to(\!xxE,\!yyE)
    \!xM=#1\!xxM  \advance\!xM\!xxloc  \!yM=#2\!yyM  \advance\!yM\!yyloc
    \!xE=#1\!xxE  \advance\!xE\!xxloc  \!yE=#2\!yyE  \advance\!yE\!yyloc
    \!qjoin (\!xM,\!yM) (\!xE,\!yE)
    \!xxS=\!xxE  \!yyS=\!yyE 
    \advance \!angle -15pt
  \repeat
  \ifdim\!angle>\!zpt
    \!angle=100.53096\!angle
    \divide \!angle 360 
    \!sinandcos\!angle\!!sin\!!cos
    \!rotate(\!xxS,\!yyS)by(\!!cos,\!sign\!!sin)to(\!xxM,\!yyM) 
    \!rotate(\!xxM,\!yyM)by(\!!cos,\!sign\!!sin)to(\!xxE,\!yyE)
    \!xM=#1\!xxM  \advance\!xM\!xxloc  \!yM=#2\!yyM  \advance\!yM\!yyloc
    \!xE=#1\!xxE  \advance\!xE\!xxloc  \!yE=#2\!yyE  \advance\!yE\!yyloc
    \!qjoin (\!xM,\!yM) (\!xE,\!yE)
  \fi
  \let\!M=\!MC
  \ignorespaces}

\def\!rotate(#1,#2)by(#3,#4)to(#5,#6){%
  \!dimenA=#3#1\advance \!dimenA -#4#2
  \!dimenB=#3#2\advance \!dimenB  #4#1
  \divide \!dimenA 32  \divide \!dimenB 32 
  #5=\!dimenA  #6=\!dimenB
  \ignorespaces}
\def\!sin{4.17684}
\def\!cos{31.72624}

\def\!sinandcos#1#2#3{%
 \!dimenD=#1
 \!dimenA=\!dimenD
 \!dimenB=32pt
 \!removept\!dimenD\!value
 \!dimenC=\!dimenD
 \!dimenC=\!value\!dimenC \divide\!dimenC by 64 
 \advance\!dimenB by -\!dimenC
 \!dimenC=\!value\!dimenC \divide\!dimenC by 96 
 \advance\!dimenA by -\!dimenC
 \!dimenC=\!value\!dimenC \divide\!dimenC by 128 
 \advance\!dimenB by \!dimenC%
 \!removept\!dimenA#2
 \!removept\!dimenB#3
 \ignorespaces}

\def\putrule#1from #2 #3 to #4 #5 {%
  \!xloc=\!M{#2}\!xunit  \!xxloc=\!M{#4}\!xunit%
  \!yloc=\!M{#3}\!yunit  \!yyloc=\!M{#5}\!yunit%
  \!dxpos=\!xxloc  \advance\!dxpos by -\!xloc
  \!dypos=\!yyloc  \advance\!dypos by -\!yloc
  \ifdim\!dypos=\!zpt
    \def\!!Line{\!puthline{#1}}\ignorespaces
  \else
    \ifdim\!dxpos=\!zpt
      \def\!!Line{\!putvline{#1}}\ignorespaces
    \else 
       \def\!!Line{}
    \fi
  \fi
  \let\!ML=\!M
  \!setdimenmode
  \!!Line%
  \let\!M=\!ML
  \ignorespaces}

\def\!putsolidhline#1{%
  \ifdim\!dxpos>\!zpt 
    \put{\!hline\!dxpos}#1[l] at {\!xloc} {\!yloc}
  \else 
    \put{\!hline{-\!dxpos}}#1[l] at {\!xxloc} {\!yyloc}
  \fi
  \ignorespaces}
 
\def\!putsolidvline#1{%
  \ifdim\!dypos>\!zpt 
    \put{\!vline\!dypos}#1[b] at {\!xloc} {\!yloc}
  \else 
    \put{\!vline{-\!dypos}}#1[b] at {\!xxloc} {\!yyloc}
  \fi
  \ignorespaces}
 
\def\!hline#1{\hbox to #1{\leaders \hrule height\linethickness\hfill}}
\def\!vline#1{\vbox to #1{\leaders \vrule width\linethickness\vfill}}

\def\!putdashedhline#1{%
  \ifdim\!dxpos>\!zpt 
    \!DLsetup\!Flist\!dxpos
    \put{\hbox to \!totalleaderlength{\!hleaders}\!hpartialpattern\!Rtrunc}
      #1[l] at {\!xloc} {\!yloc} 
  \else 
    \!DLsetup\!Blist{-\!dxpos}
    \put{\!hpartialpattern\!Ltrunc\hbox to \!totalleaderlength{\!hleaders}}
      #1[r] at {\!xloc} {\!yloc} 
  \fi
  \ignorespaces}
 
\def\!putdashedvline#1{%
  \!dypos=-\!dypos
  \ifdim\!dypos>\!zpt 
    \!DLsetup\!Flist\!dypos 
    \put{\vbox{\vbox to \!totalleaderlength{\!vleaders}
      \!vpartialpattern\!Rtrunc}}#1[t] at {\!xloc} {\!yloc} 
  \else 
    \!DLsetup\!Blist{-\!dypos}
    \put{\vbox{\!vpartialpattern\!Ltrunc
      \vbox to \!totalleaderlength{\!vleaders}}}#1[b] at {\!xloc} {\!yloc} 
  \fi
  \ignorespaces}

\def\!DLsetup#1#2{
  \let\!RSlist=#1
  \!countB=#2
  \!countA=\!leaderlength
  \divide\!countB by \!countA
  \!totalleaderlength=\!countB\!leaderlength
  \!Rresiduallength=#2%
  \advance \!Rresiduallength by -\!totalleaderlength
  \!Lresiduallength=\!leaderlength
  \advance \!Lresiduallength by -\!Rresiduallength
  \ignorespaces}
 
\def\!hleaders{%
  \def\!Rule##1{\vrule height\linethickness width##1}%
  \def\!Skip##1{\hskip##1}%
  \leaders\hbox{\!RSlist}\hfill}
 
\def\!hpartialpattern#1{%
  \!dimenA=\!zpt \!dimenB=\!zpt 
  \def\!Rule##1{#1{##1}\vrule height\linethickness width\!dimenD}%
  \def\!Skip##1{#1{##1}\hskip\!dimenD}%
  \!RSlist}
 
\def\!vleaders{%
  \def\!Rule##1{\hrule width\linethickness height##1}%
  \def\!Skip##1{\vskip##1}%
  \leaders\vbox{\!RSlist}\vfill}
 
\def\!vpartialpattern#1{%
  \!dimenA=\!zpt \!dimenB=\!zpt 
  \def\!Rule##1{#1{##1}\hrule width\linethickness height\!dimenD}%
  \def\!Skip##1{#1{##1}\vskip\!dimenD}%
  \!RSlist}
 
\def\!Rtrunc#1{\!trunc{#1}>\!Rresiduallength}
\def\!Ltrunc#1{\!trunc{#1}<\!Lresiduallength}
 
\def\!trunc#1#2#3{%
  \!dimenA=\!dimenB         
  \advance\!dimenB by #1%
  \!dimenD=\!dimenB  \ifdim\!dimenD#2#3\!dimenD=#3\fi
  \!dimenC=\!dimenA  \ifdim\!dimenC#2#3\!dimenC=#3\fi
  \advance \!dimenD by -\!dimenC}

\def\!start (#1,#2){%
  \!plotxorigin=\!xorigin  \advance \!plotxorigin by \!plotsymbolxshift
  \!plotyorigin=\!yorigin  \advance \!plotyorigin by \!plotsymbolyshift
  \!xS=\!M{#1}\!xunit \!yS=\!M{#2}\!yunit
  \!rotateaboutpivot\!xS\!yS
  \!copylist\!UDlist\to\!!UDlist
  \!getnextvalueof\!downlength\from\!!UDlist
  \!distacross=\!zpt
  \!intervalno=0 
  \global\totalarclength=\!zpt
  \ignorespaces}

\def\!ljoin (#1,#2){%
  \advance\!intervalno by 1
  \!xE=\!M{#1}\!xunit \!yE=\!M{#2}\!yunit
  \!rotateaboutpivot\!xE\!yE
  \!xdiff=\!xE \advance \!xdiff by -\!xS
  \!ydiff=\!yE \advance \!ydiff by -\!yS
  \!Pythag\!xdiff\!ydiff\!arclength
  \global\advance \totalarclength by \!arclength%
  \!drawlinearsegment
  \!xS=\!xE \!yS=\!yE
  \ignorespaces}

\def\!linearsolid{%
  \!npoints=\!arclength
  \!countA=\plotsymbolspacing
  \divide\!npoints by \!countA
  \ifnum \!npoints<1 
    \!npoints=1 
  \fi
  \divide\!xdiff by \!npoints
  \divide\!ydiff by \!npoints
  \!xpos=\!xS \!ypos=\!yS
  \loop\ifnum\!npoints>-1
    \!plotifinbounds
    \advance \!xpos by \!xdiff
    \advance \!ypos by \!ydiff
    \advance \!npoints by -1
  \repeat
  \ignorespaces}

\def\!lineardashed{%
  \ifdim\!distacross>\!arclength
    \advance \!distacross by -\!arclength  
  \else
    \loop\ifdim\!distacross<\!arclength
      \!divide\!distacross\!arclength\!dimenA
      \!removept\!dimenA\!t
      \!xpos=\!t\!xdiff \advance \!xpos by \!xS
      \!ypos=\!t\!ydiff \advance \!ypos by \!yS
      \!plotifinbounds
      \advance\!distacross by \plotsymbolspacing
      \!advancedashing
    \repeat  
    \advance \!distacross by -\!arclength
  \fi
  \ignorespaces}

\def\!!advancedashing{%
  \advance\!downlength by -\plotsymbolspacing
  \ifdim \!downlength>\!zpt
  \else
    \advance\!distacross by \!downlength
    \!getnextvalueof\!uplength\from\!!UDlist
    \advance\!distacross by \!uplength
    \!getnextvalueof\!downlength\from\!!UDlist
  \fi}

\def\inboundscheckoff{%
  \def\!plotifinbounds{\!plot(\!xpos,\!ypos)}%
  \def\!initinboundscheck{\relax}\ignorespaces}
 
\inboundscheckoff
 
\def\!!plotifinbounds{%
  \ifdim \!xpos<\!checkleft
  \else
    \ifdim \!xpos>\!checkright
    \else
      \ifdim \!ypos<\!checkbot
      \else
         \ifdim \!ypos>\!checktop
         \else
           \!plot(\!xpos,\!ypos)
         \fi 
      \fi
    \fi
  \fi}

\def\!!initinboundscheck{%
  \!checkleft=\!arealloc     \advance\!checkleft by \!xorigin
  \!checkright=\!arearloc    \advance\!checkright by \!xorigin
  \!checkbot=\!areabloc      \advance\!checkbot by \!yorigin
  \!checktop=\!areatloc      \advance\!checktop by \!yorigin}

\def\!logten#1#2{%
  \expandafter\!!logten#1\!nil
  \!removept\!dimenF#2%
  \ignorespaces}

\def\!!logten#1#2\!nil{%
  \if -#1%
    \!dimenF=\!zpt
    \def\!next{\ignorespaces}%
  \else
    \if +#1%
      \def\!next{\!!logten#2\!nil}%
    \else
      \if .#1%
        \def\!next{\!!logten0.#2\!nil}%
      \else
        \def\!next{\!!!logten#1#2..\!nil}%
      \fi
    \fi
  \fi
  \!next}

\def\!!!logten#1#2.#3.#4\!nil{%
  \!dimenF=1pt 
  \if 0#1%
    \!!logshift#3pt 
  \else 
    \!logshift#2/
    \!dimenE=#1.#2#3pt 
  \fi 
  \ifdim \!dimenE<\!rootten
    \multiply \!dimenE 10 
    \advance  \!dimenF -1pt
  \fi
  \!dimenG=\!dimenE
    \advance\!dimenG 10pt
  \advance\!dimenE -10pt 
  \multiply\!dimenE 10 
  \!divide\!dimenE\!dimenG\!dimenE
  \!removept\!dimenE\!t
  \!dimenG=\!t\!dimenE
  \!removept\!dimenG\!tt
  \!dimenH=\!tt\!tenAe
    \divide\!dimenH 100
  \advance\!dimenH \!tenAc
  \!dimenH=\!tt\!dimenH
    \divide\!dimenH 100   
  \advance\!dimenH \!tenAa
  \!dimenH=\!t\!dimenH
    \divide\!dimenH 100 
  \advance\!dimenF \!dimenH}

\def\!logshift#1{%
  \if #1/%
    \def\!next{\ignorespaces}%
  \else
    \advance\!dimenF 1pt 
    \def\!next{\!logshift}%
  \fi 
  \!next}
 
 \def\!!logshift#1{%
   \advance\!dimenF -1pt
   \if 0#1%
     \def\!next{\!!logshift}%
   \else
     \if p#1%
       \!dimenF=1pt
       \def\!next{\!dimenE=1p}%
     \else
       \def\!next{\!dimenE=#1.}%
     \fi
   \fi
   \!next}

\def\beginpicture{%
  \setbox\!picbox=\hbox\bgroup%
  \!xleft=\maxdimen  
  \!xright=-\maxdimen
  \!ybot=\maxdimen
  \!ytop=-\maxdimen}
 
\def\endpicture{%
  \ifdim\!xleft=\maxdimen
    \!xleft=\!zpt \!xright=\!zpt \!ybot=\!zpt \!ytop=\!zpt 
  \fi
  \global\!Xleft=\!xleft \global\!Xright=\!xright
  \global\!Ybot=\!ybot \global\!Ytop=\!ytop
  \egroup%
  \ht\!picbox=\!Ytop  \dp\!picbox=-\!Ybot
  \ifdim\!Ybot>\!zpt
  \else 
    \ifdim\!Ytop<\!zpt
      \!Ybot=\!Ytop
    \else
      \!Ybot=\!zpt
    \fi
  \fi
  \hbox{\kern-\!Xleft\lower\!Ybot\box\!picbox\kern\!Xright}}
 
\def\endpicturesave <#1,#2>{%
  \endpicture \global #1=\!Xleft \global #2=\!Ybot \ignorespaces}

\def\setcoordinatesystem{%
  \!ifnextchar{u}{\!getlengths }
    {\!getlengths units <\!xunit,\!yunit>}}
\def\!getlengths units <#1,#2>{%
  \!xunit=#1\relax
  \!yunit=#2\relax
  \!ifcoordmode 
    \let\!SCnext=\!SCccheckforRP
  \else
    \let\!SCnext=\!SCdcheckforRP
  \fi
  \!SCnext}
\def\!SCccheckforRP{%
  \!ifnextchar{p}{\!cgetreference }
    {\!cgetreference point at {\!xref} {\!yref} }}
\def\!cgetreference point at #1 #2 {%
  \edef\!xref{#1}\edef\!yref{#2}%
  \!xorigin=\!xref\!xunit  \!yorigin=\!yref\!yunit  
  \!initinboundscheck 
  \ignorespaces}
\def\!SCdcheckforRP{%
  \!ifnextchar{p}{\!dgetreference}%
    {\ignorespaces}}
\def\!dgetreference point at #1 #2 {%
  \!xorigin=#1\relax  \!yorigin=#2\relax
  \ignorespaces}

\long\def\put#1#2 at #3 #4 {%
  \!setputobject{#1}{#2}%
  \!xpos=\!M{#3}\!xunit  \!ypos=\!M{#4}\!yunit  
  \!rotateaboutpivot\!xpos\!ypos%
  \advance\!xpos -\!xorigin  \advance\!xpos -\!xshift
  \advance\!ypos -\!yorigin  \advance\!ypos -\!yshift
  \kern\!xpos\raise\!ypos\box\!putobject\kern-\!xpos%
  \!doaccounting\ignorespaces}
 
\long\def\multiput #1#2 at {%
  \!setputobject{#1}{#2}%
  \!ifnextchar"{\!putfromfile}{\!multiput}}
\def\!putfromfile"#1"{%
  \expandafter\!multiput \input #1 /}
\def\!multiput{%
  \futurelet\!nextchar\!!multiput}
\def\!!multiput{%
  \if *\!nextchar
    \def\!nextput{\!alsoby}%
  \else
    \if /\!nextchar
      \def\!nextput{\!finishmultiput}%
    \else
      \def\!nextput{\!alsoat}%
    \fi
  \fi
  \!nextput}
\def\!finishmultiput/{%
  \setbox\!putobject=\hbox{}%
  \ignorespaces}
 
\def\!alsoat#1 #2 {%
  \!xpos=\!M{#1}\!xunit  \!ypos=\!M{#2}\!yunit  
  \!rotateaboutpivot\!xpos\!ypos%
  \advance\!xpos -\!xorigin  \advance\!xpos -\!xshift
  \advance\!ypos -\!yorigin  \advance\!ypos -\!yshift
  \kern\!xpos\raise\!ypos\copy\!putobject\kern-\!xpos%
  \!doaccounting
  \!multiput}
 
\def\!alsoby*#1 #2 #3 {%
  \!dxpos=\!M{#2}\!xunit \!dypos=\!M{#3}\!yunit 
  \!rotateonly\!dxpos\!dypos
  \!ntemp=#1%
  \!!loop\ifnum\!ntemp>0
    \advance\!xpos by \!dxpos  \advance\!ypos by \!dypos
    \kern\!xpos\raise\!ypos\copy\!putobject\kern-\!xpos%
    \advance\!ntemp by -1
  \repeat
  \!doaccounting 
  \!multiput}
 
\def\accountingon{\def\!doaccounting{\!!doaccounting}\ignorespaces}

\accountingon
\def\!!doaccounting{%
  \!xtemp=\!xpos  
  \!ytemp=\!ypos
  \ifdim\!xtemp<\!xleft 
     \!xleft=\!xtemp 
  \fi
  \advance\!xtemp by  \!wd 
  \ifdim\!xright<\!xtemp 
    \!xright=\!xtemp
  \fi
  \advance\!ytemp by -\!dp
  \ifdim\!ytemp<\!ybot  
    \!ybot=\!ytemp
  \fi
  \advance\!ytemp by  \!dp
  \advance\!ytemp by  \!ht 
  \ifdim\!ytemp>\!ytop  
    \!ytop=\!ytemp  
  \fi}
 
\long\def\!setputobject#1#2{%
  \setbox\!putobject=\hbox{#1}%
  \!ht=\ht\!putobject  \!dp=\dp\!putobject  \!wd=\wd\!putobject
  \wd\!putobject=\!zpt
  \!xshift=.5\!wd   \!yshift=.5\!ht   \advance\!yshift by -.5\!dp
  \edef\!putorientation{#2}%
  \expandafter\!SPOreadA\!putorientation[]\!nil%
  \expandafter\!SPOreadB\!putorientation<\!zpt,\!zpt>\!nil\ignorespaces}
 
\def\!SPOreadA#1[#2]#3\!nil{\!etfor\!orientation:=#2\do\!SPOreviseshift}
 
\def\!SPOreadB#1<#2,#3>#4\!nil{\advance\!xshift by -#2\advance\!yshift by -#3}
 
\def\!SPOreviseshift{%
  \if l\!orientation 
    \!xshift=\!zpt
  \else 
    \if r\!orientation 
      \!xshift=\!wd
    \else 
      \if b\!orientation
        \!yshift=-\!dp
      \else 
        \if B\!orientation 
          \!yshift=\!zpt
        \else 
          \if t\!orientation 
            \!yshift=\!ht
          \fi 
        \fi
      \fi
    \fi
  \fi}

\long\def\!dimenput#1#2(#3,#4){%
  \!setputobject{#1}{#2}%
  \!xpos=#3\advance\!xpos by -\!xshift
  \!ypos=#4\advance\!ypos by -\!yshift
  \kern\!xpos\raise\!ypos\box\!putobject\kern-\!xpos%
  \!doaccounting\ignorespaces}

\def\!setdimenmode{%
  \let\!M=\!M!!\ignorespaces}
\def\!setcoordmode{%
  \let\!M=\!M!\ignorespaces}
\def\!ifcoordmode{%
  \ifx \!M \!M!}
\def\!ifdimenmode{%
  \ifx \!M \!M!!}
\def\!M!#1#2{#1#2} 
\def\!M!!#1#2{#1}
\!setcoordmode
\let\setdimensionmode=\!setdimenmode
\let\setcoordinatemode=\!setcoordmode

\def\!stack[#1]{%
  \let\!lglue=\hfill \let\!rglue=\hfill
  \expandafter\let\csname !#1glue\endcsname=\relax
  \!ifnextchar<{\!!stack}{\!!stack<\stackleading>}}
\def\!!stack<#1>#2{%
  \vbox{\def\!valueslist{}\!ecfor\!value:=#2\do{%
    \expandafter\!rightappend\!value\withCS{\\}\to\!valueslist}%
    \!lop\!valueslist\to\!value
    \let\\=\cr\lineskiplimit=\maxdimen\lineskip=#1%
    \baselineskip=-1000pt\halign{\!lglue##\!rglue\cr \!value\!valueslist\cr}}%
  \ignorespaces}

\def\!lines[#1]#2{%
  \let\!lglue=\hfill \let\!rglue=\hfill
  \expandafter\let\csname !#1glue\endcsname=\relax
  \vbox{\halign{\!lglue##\!rglue\cr #2\crcr}}%
  \ignorespaces}

\def\!Lines[#1]#2{%
  \let\!lglue=\hfill \let\!rglue=\hfill
  \expandafter\let\csname !#1glue\endcsname=\relax
  \vtop{\halign{\!lglue##\!rglue\cr #2\crcr}}%
  \ignorespaces}

\def\setplotsymbol(#1#2){%
  \!setputobject{#1}{#2}
  \setbox\!plotsymbol=\box\!putobject%
  \!plotsymbolxshift=\!xshift 
  \!plotsymbolyshift=\!yshift 
  \ignorespaces}
 
\setplotsymbol({\fiverm .})

\def\!!plot(#1,#2){%
  \!dimenA=-\!plotxorigin \advance \!dimenA by #1
  \!dimenB=-\!plotyorigin \advance \!dimenB by #2
  \kern\!dimenA\raise\!dimenB\copy\!plotsymbol\kern-\!dimenA%
  \ignorespaces}
 
\def\!!!plot(#1,#2){%
  \!dimenA=-\!plotxorigin \advance \!dimenA by #1

  \!dimenB=-\!plotyorigin \advance \!dimenB by #2
  \kern\!dimenA\raise\!dimenB\copy\!plotsymbol\kern-\!dimenA%
  \!countE=\!dimenA
  \!countF=\!dimenB
  \immediate\write\!replotfile{\the\!countE,\the\!countF.}%
  \ignorespaces}

\def\savelinesandcurves on "#1" {%
  \immediate\closeout\!replotfile
  \immediate\openout\!replotfile=#1%
  \let\!plot=\!!!plot}

\def\dontsavelinesandcurves {%
  \let\!plot=\!!plot}
\dontsavelinesandcurves

{\catcode`\%=11\xdef\!Commentsignal{
\def\writesavefile#1 {%
  \immediate\write\!replotfile{\!Commentsignal #1}%
  \ignorespaces}

\def\replot"#1" {%
  \expandafter\!replot\input #1 /}
\def\!replot#1,#2. {%
  \!dimenA=#1sp
  \kern\!dimenA\raise#2sp\copy\!plotsymbol\kern-\!dimenA
  \futurelet\!nextchar\!!replot}
\def\!!replot{%
  \if /\!nextchar 
    \def\!next{\!finish}%
  \else
    \def\!next{\!replot}%
  \fi
  \!next}

\def\!Pythag#1#2#3{%
  \!dimenE=#1\relax                                     
  \ifdim\!dimenE<\!zpt 
    \!dimenE=-\!dimenE 
  \fi
  \!dimenF=#2\relax
  \ifdim\!dimenF<\!zpt 
    \!dimenF=-\!dimenF 
  \fi
  \advance \!dimenF by \!dimenE
  \ifdim\!dimenF=\!zpt 
    \!dimenG=\!zpt
  \else 
    \!divide{8\!dimenE}\!dimenF\!dimenE
    \advance\!dimenE by -4pt
      \!dimenE=2\!dimenE
    \!removept\!dimenE\!!t
    \!dimenE=\!!t\!dimenE
    \advance\!dimenE by 64pt
    \divide \!dimenE by 2
    \!dimenH=7pt
    \!!Pythag\!!Pythag\!!Pythag
    \!removept\!dimenH\!!t
    \!dimenG=\!!t\!dimenF
    \divide\!dimenG by 8
  \fi
  #3=\!dimenG
  \ignorespaces}

\def\!!Pythag{
  \!divide\!dimenE\!dimenH\!dimenI
  \advance\!dimenH by \!dimenI
    \divide\!dimenH by 2}

\def\placehypotenuse for <#1> and <#2> in <#3> {%
  \!Pythag{#1}{#2}{#3}}

\def\!qjoin (#1,#2) (#3,#4){%
  \advance\!intervalno by 1
  \!ifcoordmode
    \edef\!xmidpt{#1}\edef\!ymidpt{#2}%
  \else
    \!dimenA=#1\relax \edef\!xmidpt{\the\!dimenA}%
    \!dimenA=#2\relax \edef\!xmidpt{\the\!dimenA}%
  \fi
  \!xM=\!M{#1}\!xunit  \!yM=\!M{#2}\!yunit   \!rotateaboutpivot\!xM\!yM
  \!xE=\!M{#3}\!xunit  \!yE=\!M{#4}\!yunit   \!rotateaboutpivot\!xE\!yE
  \!dimenA=\!xM  \advance \!dimenA by -\!xS
  \!dimenB=\!xE  \advance \!dimenB by -\!xM
  \!xB=3\!dimenA \advance \!xB by -\!dimenB
  \!xC=2\!dimenB \advance \!xC by -2\!dimenA
  \!dimenA=\!yM  \advance \!dimenA by -\!yS%
  \!dimenB=\!yE  \advance \!dimenB by -\!yM%
  \!yB=3\!dimenA \advance \!yB by -\!dimenB%
  \!yC=2\!dimenB \advance \!yC by -2\!dimenA%
  \!xprime=\!xB  \!yprime=\!yB
  \!dxprime=.5\!xC  \!dyprime=.5\!yC
  \!getf \!midarclength=\!dimenA
  \!getf \advance \!midarclength by 4\!dimenA
  \!getf \advance \!midarclength by \!dimenA
  \divide \!midarclength by 12
  \!arclength=\!dimenA
  \!getf \advance \!arclength by 4\!dimenA
  \!getf \advance \!arclength by \!dimenA
  \divide \!arclength by 12
  \advance \!arclength by \!midarclength
  \global\advance \totalarclength by \!arclength
  \ifdim\!distacross>\!arclength 
    \advance \!distacross by -\!arclength
  \else
    \!initinverseinterp
    \loop\ifdim\!distacross<\!arclength
      \!inverseinterp
      \!xpos=\!t\!xC \advance\!xpos by \!xB
        \!xpos=\!t\!xpos \advance \!xpos by \!xS
      \!ypos=\!t\!yC \advance\!ypos by \!yB
        \!ypos=\!t\!ypos \advance \!ypos by \!yS
      \!plotifinbounds
      \advance\!distacross \plotsymbolspacing
      \!advancedashing
    \repeat  
    \advance \!distacross by -\!arclength
  \fi
  \!xS=\!xE
  \!yS=\!yE
  \ignorespaces}

\def\!getf{\!Pythag\!xprime\!yprime\!dimenA%
  \advance\!xprime by \!dxprime
  \advance\!yprime by \!dyprime}

\def\!initinverseinterp{%
  \ifdim\!arclength>\!zpt
    \!divide{8\!midarclength}\!arclength\!dimenE
    \ifdim\!dimenE<\!wmin \!setinverselinear
    \else 
      \ifdim\!dimenE>\!wmax \!setinverselinear
      \else
        \def\!inverseinterp{\!inversequad}\ignorespaces
         \!removept\!dimenE\!Ew
         \!dimenF=-\!Ew\!dimenE
         \advance\!dimenF by 32pt
         \!dimenG=8pt 
         \advance\!dimenG by -\!dimenE
         \!dimenG=\!Ew\!dimenG
         \!divide\!dimenF\!dimenG\!beta
         \!gamma=1pt
         \advance \!gamma by -\!beta
      \fi
    \fi
  \fi
  \ignorespaces}

\def\!inversequad{%
  \!divide\!distacross\!arclength\!dimenG
  \!removept\!dimenG\!v
  \!dimenG=\!v\!gamma
  \advance\!dimenG by \!beta
  \!dimenG=\!v\!dimenG
  \!removept\!dimenG\!t}

\def\!setinverselinear{%
  \def\!inverseinterp{\!inverselinear}%
  \divide\!dimenE by 8 \!removept\!dimenE\!t
  \!countC=\!intervalno \multiply \!countC 2
  \!countB=\!countC     \advance \!countB -1
  \!countA=\!countB     \advance \!countA -1
  \wlog{\the\!countB th point (\!xmidpt,\!ymidpt) being plotted 
    doesn't lie in the}%
  \wlog{ middle third of the arc between the \the\!countA th 
    and \the\!countC th points:}%
  \wlog{ [arc length \the\!countA\space to \the\!countB]/[arc length 
    \the \!countA\space to \the\!countC]=\!t.}%
  \ignorespaces}
 
\def\!inverselinear{%
  \!divide\!distacross\!arclength\!dimenG
  \!removept\!dimenG\!t}

\def\startrotation{%
  \let\!rotateaboutpivot=\!!rotateaboutpivot
  \let\!rotateonly=\!!rotateonly
  \!ifnextchar{b}{\!getsincos }%
    {\!getsincos by {\!cosrotationangle} {\!sinrotationangle} }}
\def\!getsincos by #1 #2 {%
  \edef\!cosrotationangle{#1}%
  \edef\!sinrotationangle{#2}%
  \!ifcoordmode 
    \let\!ROnext=\!ccheckforpivot
  \else
    \let\!ROnext=\!dcheckforpivot
  \fi
  \!ROnext}
\def\!ccheckforpivot{%
  \!ifnextchar{a}{\!cgetpivot}%
    {\!cgetpivot about {\!xpivotcoord} {\!ypivotcoord} }}
\def\!cgetpivot about #1 #2 {%
  \edef\!xpivotcoord{#1}%
  \edef\!ypivotcoord{#2}%
  \!xpivot=#1\!xunit  \!ypivot=#2\!yunit
  \ignorespaces}
\def\!dcheckforpivot{%
  \!ifnextchar{a}{\!dgetpivot}{\ignorespaces}}
\def\!dgetpivot about #1 #2 {%
  \!xpivot=#1\relax  \!ypivot=#2\relax
  \ignorespaces}

\def\stoprotation{%
  \let\!rotateaboutpivot=\!!!rotateaboutpivot
  \let\!rotateonly=\!!!rotateonly
  \ignorespaces}
 
\def\!!rotateaboutpivot#1#2{%
  \!dimenA=#1\relax  \advance\!dimenA -\!xpivot
  \!dimenB=#2\relax  \advance\!dimenB -\!ypivot
  \!dimenC=\!cosrotationangle\!dimenA
    \advance \!dimenC -\!sinrotationangle\!dimenB
  \!dimenD=\!cosrotationangle\!dimenB
    \advance \!dimenD  \!sinrotationangle\!dimenA
  \advance\!dimenC \!xpivot  \advance\!dimenD \!ypivot
  #1=\!dimenC  #2=\!dimenD
  \ignorespaces}

\def\!!rotateonly#1#2{%
  \!dimenA=#1\relax  \!dimenB=#2\relax 
  \!dimenC=\!cosrotationangle\!dimenA
    \advance \!dimenC -\!rotsign\!sinrotationangle\!dimenB
  \!dimenD=\!cosrotationangle\!dimenB
    \advance \!dimenD  \!rotsign\!sinrotationangle\!dimenA
  #1=\!dimenC  #2=\!dimenD
  \ignorespaces}
\def\!rotsign{}
\def\!!!rotateaboutpivot#1#2{\relax}
\def\!!!rotateonly#1#2{\relax}
\stoprotation

\def\!reverserotateonly#1#2{%
  \def\!rotsign{-}%
  \!rotateonly{#1}{#2}%
  \def\!rotsign{}%
  \ignorespaces}

\def\!getspan span <#1>{%
  \!dshade=#1\relax
  \!ifcoordmode 
    \let\!GRnext=\!GRccheckforAP
  \else
    \let\!GRnext=\!GRdcheckforAP
  \fi
  \!GRnext}
\def\!GRccheckforAP{%
  \!ifnextchar{p}{\!cgetanchor }
    {\!cgetanchor point at {\!xshadesave} {\!yshadesave} }}
\def\!cgetanchor point at #1 #2 {%
  \edef\!xshadesave{#1}\edef\!yshadesave{#2}%
  \!xshade=\!xshadesave\!xunit  \!yshade=\!yshadesave\!yunit
  \ignorespaces}
\def\!GRdcheckforAP{%
  \!ifnextchar{p}{\!dgetanchor}%
    {\ignorespaces}}
\def\!dgetanchor point at #1 #2 {%
  \!xshade=#1\relax  \!yshade=#2\relax
  \ignorespaces}

\def\setshadesymbol{%
  \!ifnextchar<{\!setshadesymbol}{\!setshadesymbol<,,,> }}

\def\!setshadesymbol <#1,#2,#3,#4> (#5#6){%
  \!setputobject{#5}{#6}%
  \setbox\!shadesymbol=\box\!putobject%
  \!shadesymbolxshift=\!xshift \!shadesymbolyshift=\!yshift
  \!dimenA=\!xshift \advance\!dimenA \!smidge
  \!override\!dimenA{#1}\!lshrinkage%
  \!dimenA=\!wd \advance \!dimenA -\!xshift
    \advance\!dimenA \!smidge
    \!override\!dimenA{#2}\!rshrinkage
  \!dimenA=\!dp \advance \!dimenA \!yshift
    \advance\!dimenA \!smidge
    \!override\!dimenA{#3}\!bshrinkage
  \!dimenA=\!ht \advance \!dimenA -\!yshift
    \advance\!dimenA \!smidge
    \!override\!dimenA{#4}\!tshrinkage
  \ignorespaces}
\def\!smidge{-.2pt}%

\def\!override#1#2#3{%
  \edef\!!override{#2}%
  \ifx \!!override\empty
    #3=#1\relax
  \else
    \if z\!!override
      #3=\!zpt
    \else
      \ifx \!!override\!blankz
        #3=\!zpt
      \else
        #3=#2\relax
      \fi
    \fi
  \fi
  \ignorespaces}
\def\!blankz{ z}

\setshadesymbol ({\fiverm .})

\def\!startvshade#1(#2,#3,#4){%
  \let\!!xunit=\!xunit%
  \let\!!yunit=\!yunit%
  \let\!!xshade=\!xshade%
  \let\!!yshade=\!yshade%
  \def\!getshrinkages{\!vgetshrinkages}%
  \let\!setshadelocation=\!vsetshadelocation%
  \!xS=\!M{#2}\!!xunit
  \!ybS=\!M{#3}\!!yunit
  \!ytS=\!M{#4}\!!yunit
  \!shadexorigin=\!xorigin  \advance \!shadexorigin \!shadesymbolxshift
  \!shadeyorigin=\!yorigin  \advance \!shadeyorigin \!shadesymbolyshift
  \ignorespaces}
 
\def\!starthshade#1(#2,#3,#4){%
  \let\!!xunit=\!yunit%
  \let\!!yunit=\!xunit%
  \let\!!xshade=\!yshade%
  \let\!!yshade=\!xshade%
  \def\!getshrinkages{\!hgetshrinkages}%
  \let\!setshadelocation=\!hsetshadelocation%
  \!xS=\!M{#2}\!!xunit
  \!ybS=\!M{#3}\!!yunit
  \!ytS=\!M{#4}\!!yunit
  \!shadexorigin=\!xorigin  \advance \!shadexorigin \!shadesymbolxshift
  \!shadeyorigin=\!yorigin  \advance \!shadeyorigin \!shadesymbolyshift
  \ignorespaces}

\def\!lattice#1#2#3#4#5{%
  \!dimenA=#1
  \!dimenB=#2
  \!countB=\!dimenB
  \!dimenC=#3
  \advance\!dimenC -\!dimenA
  \!countA=\!dimenC
  \divide\!countA \!countB
  \ifdim\!dimenC>\!zpt
    \!dimenD=\!countA\!dimenB
    \ifdim\!dimenD<\!dimenC
      \advance\!countA 1 
    \fi
  \fi
  \!dimenC=\!countA\!dimenB
    \advance\!dimenC \!dimenA
  #4=\!countA
  #5=\!dimenC
  \ignorespaces}

\def\!qshade#1(#2,#3,#4)#5(#6,#7,#8){%
  \!xM=\!M{#2}\!!xunit
  \!ybM=\!M{#3}\!!yunit
  \!ytM=\!M{#4}\!!yunit
  \!xE=\!M{#6}\!!xunit
  \!ybE=\!M{#7}\!!yunit
  \!ytE=\!M{#8}\!!yunit
  \!getcoeffs\!xS\!ybS\!xM\!ybM\!xE\!ybE\!ybB\!ybC
  \!getcoeffs\!xS\!ytS\!xM\!ytM\!xE\!ytE\!ytB\!ytC
  \def\!getylimits{\!qgetylimits}%
  \!shade{#1}\ignorespaces}
 
\def\!lshade#1(#2,#3,#4){%
  \!xE=\!M{#2}\!!xunit
  \!ybE=\!M{#3}\!!yunit
  \!ytE=\!M{#4}\!!yunit
  \!dimenE=\!xE  \advance \!dimenE -\!xS
  \!dimenC=\!ytE \advance \!dimenC -\!ytS
  \!divide\!dimenC\!dimenE\!ytB
  \!dimenC=\!ybE \advance \!dimenC -\!ybS
  \!divide\!dimenC\!dimenE\!ybB
  \def\!getylimits{\!lgetylimits}%
  \!shade{#1}\ignorespaces}
 
\def\!getcoeffs#1#2#3#4#5#6#7#8{%
  \!dimenC=#4\advance \!dimenC -#2
  \!dimenE=#3\advance \!dimenE -#1
  \!divide\!dimenC\!dimenE\!dimenF
  \!dimenC=#6\advance \!dimenC -#4
  \!dimenH=#5\advance \!dimenH -#3
  \!divide\!dimenC\!dimenH\!dimenG
  \advance\!dimenG -\!dimenF
  \advance \!dimenH \!dimenE
  \!divide\!dimenG\!dimenH#8
  \!removept#8\!t
  #7=-\!t\!dimenE
  \advance #7\!dimenF
  \ignorespaces}

\def\!shade#1{%
  \!getshrinkages#1<,,,>\!nil
  \advance \!dimenE \!xS
  \!lattice\!!xshade\!dshade\!dimenE
    \!parity\!xpos
  \!dimenF=-\!dimenF
    \advance\!dimenF \!xE
  \!loop\!not{\ifdim\!xpos>\!dimenF}
    \!shadecolumn%
    \advance\!xpos \!dshade
    \advance\!parity 1
  \repeat
  \!xS=\!xE
  \!ybS=\!ybE
  \!ytS=\!ytE
  \ignorespaces}

\def\!vgetshrinkages#1<#2,#3,#4,#5>#6\!nil{%
  \!override\!lshrinkage{#2}\!dimenE
  \!override\!rshrinkage{#3}\!dimenF
  \!override\!bshrinkage{#4}\!dimenG
  \!override\!tshrinkage{#5}\!dimenH
  \ignorespaces}
\def\!hgetshrinkages#1<#2,#3,#4,#5>#6\!nil{%
  \!override\!lshrinkage{#2}\!dimenG
  \!override\!rshrinkage{#3}\!dimenH
  \!override\!bshrinkage{#4}\!dimenE
  \!override\!tshrinkage{#5}\!dimenF
  \ignorespaces}

\def\!shadecolumn{%
  \!dxpos=\!xpos
  \advance\!dxpos -\!xS
  \!removept\!dxpos\!dx
  \!getylimits
  \advance\!ytpos -\!dimenH
  \advance\!ybpos \!dimenG
  \!yloc=\!!yshade
  \ifodd\!parity 
     \advance\!yloc \!dshade
  \fi
  \!lattice\!yloc{2\!dshade}\!ybpos%
    \!countA\!ypos
  \!dimenA=-\!shadexorigin \advance \!dimenA \!xpos
  \loop\!not{\ifdim\!ypos>\!ytpos}
    \!setshadelocation
    \!rotateaboutpivot\!xloc\!yloc%
    \!dimenA=-\!shadexorigin \advance \!dimenA \!xloc
    \!dimenB=-\!shadeyorigin \advance \!dimenB \!yloc
    \kern\!dimenA \raise\!dimenB\copy\!shadesymbol \kern-\!dimenA
    \advance\!ypos 2\!dshade
  \repeat
  \ignorespaces}
 
\def\!qgetylimits{%
  \!dimenA=\!dx\!ytC              
  \advance\!dimenA \!ytB
  \!ytpos=\!dx\!dimenA
  \advance\!ytpos \!ytS
  \!dimenA=\!dx\!ybC              
  \advance\!dimenA \!ybB
  \!ybpos=\!dx\!dimenA
  \advance\!ybpos \!ybS}
 
\def\!lgetylimits{%
  \!ytpos=\!dx\!ytB
  \advance\!ytpos \!ytS
  \!ybpos=\!dx\!ybB
  \advance\!ybpos \!ybS}
 
\def\!vsetshadelocation{
  \!xloc=\!xpos
  \!yloc=\!ypos}
\def\!hsetshadelocation{
  \!xloc=\!ypos
  \!yloc=\!xpos}

\def\!axisticks {%
  \def\!nextkeyword##1 {%
    \expandafter\ifx\csname !ticks##1\endcsname \relax
      \def\!next{\!fixkeyword{##1}}%
    \else
      \def\!next{\csname !ticks##1\endcsname}%
    \fi
    \!next}%
  \!axissetup
    \def\!axissetup{\relax}%
  \edef\!ticksinoutsign{\!ticksinoutSign}%
  \!ticklength=\longticklength
  \!tickwidth=\linethickness
  \!gridlinestatus
  \!setticktransform
  \!maketick
  \!tickcase=0
  \def\!LTlist{}%
  \!nextkeyword}

\def\ticksout{%
  \def\!ticksinoutSign{+}}

\ticksout

\def\nogridlines{%
  \def\!gridlinestatus{\!gridlinestoofalse}}
\nogridlines

\def\loggedticks{%
  \def\!setticktransform{\let\!ticktransform=\!logten}}
\def\unloggedticks{%
  \def\!setticktransform{\let\!ticktransform=\!donothing}}
\def\!donothing#1#2{\def#2{#1}}
\unloggedticks

\expandafter\def\csname !ticks/\endcsname{%
  \!not {\ifx \!LTlist\empty}
    \!placetickvalues
  \fi
  \def\!tickvalueslist{}%
  \def\!LTlist{}%
  \expandafter\csname !axis/\endcsname}

\def\!maketick{%
  \setbox\!boxA=\hbox{%
    \beginpicture
      \!setdimenmode
      \setcoordinatesystem point at {\!zpt} {\!zpt}   
      \linethickness=\!tickwidth
      \ifdim\!ticklength>\!zpt
        \putrule from {\!zpt} {\!zpt} to
          {\!ticksinoutsign\!tickxsign\!ticklength}
          {\!ticksinoutsign\!tickysign\!ticklength}
      \fi
      \if!gridlinestoo
        \putrule from {\!zpt} {\!zpt} to
          {-\!tickxsign\!xaxislength} {-\!tickysign\!yaxislength}
      \fi
    \endpicturesave <\!Xsave,\!Ysave>}%
    \wd\!boxA=\!zpt}
  
\def\!ticksin{%
  \def\!ticksinoutsign{-}%
  \!maketick
  \!nextkeyword}

\def\!ticksout{%
  \def\!ticksinoutsign{+}%
  \!maketick
  \!nextkeyword}

\def\!tickslength<#1> {%
  \!ticklength=#1\relax
  \!maketick
  \!nextkeyword}

\def\!tickslong{%
  \!tickslength<\longticklength> }

\def\!ticksshort{%
  \!tickslength<\shortticklength> }

\def\!tickswidth<#1> {%
  \!tickwidth=#1\relax
  \!maketick
  \!nextkeyword}

\def\!ticksandacross{%
  \!gridlinestootrue
  \!maketick
  \!nextkeyword}

\def\!ticksbutnotacross{%
  \!gridlinestoofalse
  \!maketick
  \!nextkeyword}

\def\!tickslogged{%
  \let\!ticktransform=\!logten
  \!nextkeyword}

\def\!ticksunlogged{%
  \let\!ticktransform=\!donothing
  \!nextkeyword}

\def\!ticksunlabeled{%
  \!tickcase=0
  \!nextkeyword}

\def\!ticksnumbered{%
  \!tickcase=1
  \!nextkeyword}

\def\!tickswithvalues#1/ {%
  \edef\!tickvalueslist{#1! /}%
  \!tickcase=2
  \!nextkeyword}

\def\!ticksquantity#1 {%
  \ifnum #1>1
    \!updatetickoffset
    \!countA=#1\relax
    \advance \!countA -1
    \!ticklocationincr=\!axisLength
      \divide \!ticklocationincr \!countA
    \!ticklocation=\!axisstart
    \loop \!not{\ifdim \!ticklocation>\!axisend}
      \!placetick\!ticklocation
      \ifcase\!tickcase
          \relax 
        \or
          \relax 
        \or
          \expandafter\!gettickvaluefrom\!tickvalueslist
          \edef\!tickfield{{\the\!ticklocation}{\!value}}%
          \expandafter\!listaddon\expandafter{\!tickfield}\!LTlist%
      \fi
      \advance \!ticklocation \!ticklocationincr
    \repeat
  \fi
  \!nextkeyword}

\def\!ticksat#1 {%
  \!updatetickoffset
  \edef\!Loc{#1}%
  \if /\!Loc
    \def\next{\!nextkeyword}%
  \else
    \!ticksincommon
    \def\next{\!ticksat}%
  \fi
  \next}    
      
\def\!ticksfrom#1 to #2 by #3 {%
  \!updatetickoffset
  \edef\!arg{#3}%
  \expandafter\!separate\!arg\!nil
  \!scalefactor=1
  \expandafter\!countfigures\!arg/
  \edef\!arg{#1}%
  \!scaleup\!arg by\!scalefactor to\!countE
  \edef\!arg{#2}%
  \!scaleup\!arg by\!scalefactor to\!countF
  \edef\!arg{#3}%
  \!scaleup\!arg by\!scalefactor to\!countG
  \loop \!not{\ifnum\!countE>\!countF}
    \ifnum\!scalefactor=1
      \edef\!Loc{\the\!countE}%
    \else
      \!scaledown\!countE by\!scalefactor to\!Loc
    \fi
    \!ticksincommon
    \advance \!countE \!countG
  \repeat
  \!nextkeyword}

\def\!updatetickoffset{%
  \!dimenA=\!ticksinoutsign\!ticklength
  \ifdim \!dimenA>\!offset
    \!offset=\!dimenA
  \fi}

\def\!placetick#1{%
  \if!xswitch
    \!xpos=#1\relax
    \!ypos=\!axisylevel
  \else
    \!xpos=\!axisxlevel
    \!ypos=#1\relax
  \fi
  \advance\!xpos \!Xsave
  \advance\!ypos \!Ysave
  \kern\!xpos\raise\!ypos\copy\!boxA\kern-\!xpos
  \ignorespaces}

\def\!gettickvaluefrom#1 #2 /{%
  \edef\!value{#1}%
  \edef\!tickvalueslist{#2 /}%
  \ifx \!tickvalueslist\!endtickvaluelist
    \!tickcase=0
  \fi}
\def\!endtickvaluelist{! /}

\def\!ticksincommon{%
  \!ticktransform\!Loc\!t
  \!ticklocation=\!t\!!unit
  \advance\!ticklocation -\!!origin
  \!placetick\!ticklocation
  \ifcase\!tickcase
    \relax 
  \or 
    \ifdim\!ticklocation<-\!!origin
      \edef\!Loc{$\!Loc$}%
    \fi
    \edef\!tickfield{{\the\!ticklocation}{\!Loc}}%
    \expandafter\!listaddon\expandafter{\!tickfield}\!LTlist%
  \or 
    \expandafter\!gettickvaluefrom\!tickvalueslist
    \edef\!tickfield{{\the\!ticklocation}{\!value}}%
    \expandafter\!listaddon\expandafter{\!tickfield}\!LTlist%
  \fi}

\def\!separate#1\!nil{%
  \!ifnextchar{-}{\!!separate}{\!!!separate}#1\!nil}
\def\!!separate-#1\!nil{%
  \def\!sign{-}%
  \!!!!separate#1..\!nil}
\def\!!!separate#1\!nil{%
  \def\!sign{+}%
  \!!!!separate#1..\!nil}
\def\!!!!separate#1.#2.#3\!nil{%
  \def\!arg{#1}%
  \ifx\!arg\!empty
    \!countA=0
  \else
    \!countA=\!arg

  \fi
  \def\!arg{#2}%
  \ifx\!arg\!empty
    \!countB=0
  \else
    \!countB=\!arg
  \fi}
 
\def\!countfigures#1{%
  \if #1/%
    \def\!next{\ignorespaces}%
  \else
    \multiply\!scalefactor 10
    \def\!next{\!countfigures}%
  \fi
  \!next}

\def\!scaleup#1by#2to#3{%
  \expandafter\!separate#1\!nil
  \multiply\!countA #2\relax
  \advance\!countA \!countB
  \if -\!sign
    \!countA=-\!countA
  \fi
  #3=\!countA
  \ignorespaces}

\def\!scaledown#1by#2to#3{%
  \!countA=#1\relax
  \ifnum \!countA<0 
    \def\!sign{-}
    \!countA=-\!countA
  \else
    \def\!sign{}%
  \fi
  \!countB=\!countA
  \divide\!countB #2\relax
  \!countC=\!countB
    \multiply\!countC #2\relax
  \advance \!countA -\!countC
  \edef#3{\!sign\the\!countB.}
  \!countC=\!countA 
  \ifnum\!countC=0 
    \!countC=1
  \fi
  \multiply\!countC 10
  \!loop \ifnum #2>\!countC
    \edef#3{#3\!zero}%
    \multiply\!countC 10
  \repeat
  \edef#3{#3\the\!countA}
  \ignorespaces}

\def\!placetickvalues{%
  \advance\!offset \tickstovaluesleading
  \if!xswitch
    \setbox\!boxA=\hbox{%
      \def\\##1##2{%
        \!dimenput {##2} [B] (##1,\!axisylevel)}%
      \beginpicture 
        \!LTlist
      \endpicturesave <\!Xsave,\!Ysave>}%
    \!dimenA=\!axisylevel
      \advance\!dimenA -\!Ysave
      \advance\!dimenA \!tickysign\!offset
      \if -\!tickysign
        \advance\!dimenA -\ht\!boxA
      \else
        \advance\!dimenA  \dp\!boxA
      \fi
    \advance\!offset \ht\!boxA 
      \advance\!offset \dp\!boxA
    \!dimenput {\box\!boxA} [Bl] <\!Xsave,\!Ysave> (\!zpt,\!dimenA)
  \else
    \setbox\!boxA=\hbox{%
      \def\\##1##2{%
        \!dimenput {##2} [r] (\!axisxlevel,##1)}%
      \beginpicture 
        \!LTlist
      \endpicturesave <\!Xsave,\!Ysave>}%
    \!dimenA=\!axisxlevel
      \advance\!dimenA -\!Xsave
      \advance\!dimenA \!tickxsign\!offset
      \if -\!tickxsign
        \advance\!dimenA -\wd\!boxA
      \fi
    \advance\!offset \wd\!boxA
    \!dimenput {\box\!boxA} [Bl] <\!Xsave,\!Ysave> (\!dimenA,\!zpt)
  \fi}

\normalgraphs
\catcode`!=12 

\catcode`@=11 \catcode`!=11
  
\let\!pictexendpicture=\endpicture 
\let\!pictexframe=\frame
\let\!pictexlinethickness=\linethickness
\let\!pictexmultiput=\multiput
\let\!pictexput=\put

\def\beginpicture{%
  \setbox\!picbox=\hbox\bgroup%
  \let\endpicture=\!pictexendpicture
  \let\frame=\!pictexframe
  \let\linethickness=\!pictexlinethickness
  \let\multiput=\!pictexmultiput
  \let\put=\!pictexput
  \let\input=\@@input   
  \!xleft=\maxdimen  
  \!xright=-\maxdimen
  \!ybot=\maxdimen
  \!ytop=-\maxdimen}

\let\frame=\!latexframe

\let\pictexframe=\!pictexframe

\let\linethickness=\!latexlinethickness
\let\pictexlinethickness=\!pictexlinethickness

\let\\=\@normalcr
\catcode`@=12 \catcode`!=12

\begingroup\makeatletter
\def\x#1#2#3#4#5#6#7\relax{\def\x{#1#2#3#4#5#6}}
\expandafter\x\fmtname xxxxxx\relax \def\y{splain}
\ifx\x\y   
\gdef\SetFigFont#1#2#3{%
  \ifnum #1<17 \tiny\else \ifnum #1<20 \small\else
  \ifnum #1<24 \normalsize\else \ifnum #1<29 \large\else
  \ifnum #1<34 \Large\else \ifnum #1<41 \LARGE\else
     \huge\fi\fi\fi\fi\fi\fi
  \csname #3\endcsname}
\else
\gdef\SetFigFont#1#2#3{\begingroup
  \count@#1\relax \ifnum 25<\count@ \count@25 \fi
  \def\x{\endgroup\@setsize\SetFigFont{#2pt}}%
  \expandafter\x
    \csname \romannumeral\the\count@ pt\expandafter\endcsname
    \csname @\romannumeral\the\count@ pt\endcsname
  \csname #3\endcsname}
\fi
\endgroup

\newtheorem{pr}{Proposition}
\newtheorem{lm}{Lemma}
\newtheorem{tm}{Theorem}

\newcommand{\proj}{\mathbf P}
\newcommand{\grass}{\mathbf G}
\newcommand{\barr}{\overline}
\newcommand{\rarr}{\rightarrow}
\newcommand{\oh}{{\mathcal{O}}}
\newcommand{\com}{\mathbb{C}}
\newcommand{\tl}{\tilde}

\newcommand{\map}{\barr{\mathcal{M}}}
\newcommand{\smap}{\barr{M}}
\newcommand{\M}{\barr{M}}
\newcommand{\st}{_{g,n}(X,\beta)}
\newcommand{\sto}{_{0,n}(X,\beta)}
\newcommand{\stp}{_{g,n}(\proj^r,d)}
\newcommand{\stpo}{_{0,n}(\proj^r,d)}
\newcommand{\stv}{_{g,n}(\proj^r,d,\barr{t})}
\newcommand{\stvo}{_{0,n}(\proj^r,d,\barr{t})}
\newcommand{\mk}{\barr{M}_{0,m}}
\newcommand{\mku} {\barr{U}_{0,m}}

\newcommand{\HH}{\mathcal{H}}
\newcommand{\LL}{\mathcal{L}}
\newcommand{\GG}{\mathcal{G}}
\newcommand{\goth}{\mathfrak{S}}
\newcommand{\Q}{\mathbb Q}
\newcommand{\Z}{\mathbb Z}

\newcommand{\cal}{\mathcal}
\newcommand{\smallcup}{\mathbin{\text{\scriptsize$\cup$}}}
\newcommand{\smallcap}{\mathbin{\text{\scriptsize$\cap$}}}
\newcommand{\smallstm}{\mathbin{\text{\scriptsize$\setminus$}}}

\newcommand{\eqq}{\stackrel {\sim}{=}}

\newcommand{\bpf}{\noindent {\em Proof.} }
\newcommand{\epf}{\qed \vspace{+10pt}}

\begin{document}


\title{Notes On Stable Maps And Quantum Cohomology}

\author{W. Fulton}
\address{Department of Mathematics, University of Chicago,
Chicago, Illinois, 60637}
\email{fulton@math.uchicago.edu}
\thanks{The first author was supported in part by NSF
Grant DMS 9307922.}
\author{R. Pandharipande}
\address{Department of Mathematics, University of Chicago,
Chicago, Illinois, 60637}
\email{rahul@math.uchicago.edu}
\thanks{The second author was supported in part by
an NSF Post-Doctoral Fellowship.}
\subjclass{Primary 14N10, 14H10;
Secondary 14E99}
\date{14 June 1996}
\dedicatory{Dedicated to the memory of Claude Itzykson}
\maketitle

\setcounter{section}{-1}
\tableofcontents

\section{{\bf Introduction}}
\subsection{Overview}
The aim of these notes is to describe an exciting chapter 
in the recent development of quantum cohomology.  
Guided by ideas from physics (see [W]), a remarkable
structure on the solutions of certain rational enumerative geometry
problems has been found: the solutions are coefficients
in the multiplication table 
of a quantum cohomology ring. Associativity of the
ring yields non-trivial relations among the enumerative solutions.
In many cases, these relations suffice to solve the enumerative
problem. For example, let $N_d$ be the
number of degree $d$, rational plane curves passing
through $3d-1$ general points in $\proj^2$. 
Since there is a unique line passing through $2$ points,
$N_1=1$. 
The quantum cohomology ring of $\proj^2$
yields the following beautiful associativity 
relation determining all 
$N_d$ for $d\geq 2$:
$$N_d= \sum_{d_1+d_2=d, \ d_1, d_2 > 0}
N_{d_1} N_{d_2} \bigg( d_1^2 d_2^2 \binom{3d-4}{3d_1-2} -
d_1^3 d_2 \binom{3d-4}{3d_1-1} \bigg).$$
Similar enumerative  formulas
are valid on other homogeneous varieties.
Viewed from classical enumerative geometry, the quantum ring
structure is a complete surprise.
  
The path to quantum
cohomology presented here follows the work of Kontsevich
and Manin. The approach is algebro-geometric and involves
the construction and geometry of a natural 
compactification of the moduli space of maps.
The large and exciting conjectural parts of the
subject of quantum cohomology are avoided here.
We focus on a part of the story where the 
proofs are complete. We also make many assumptions that are not 
strictly necessary, but which simplify the presentation.

It should be emphasized that this is in no way a survey of quantum 
cohomology, or any attempt at evaluating various approaches.
In particular, the symplectic point of view is not covered
(see [R-T]). 
Another 
algebro-geometric approach, 
using a different compactification, can be found 
in [L-T 1].  

These notes  are based on a 
jointly taught course at the University of 
Chicago  in which our main efforts were aimed at 
understanding the papers of Kontsevich and Manin.  
We thank R. Donagi for instigating this course.
Thanks are due to D. Abramovich, P. Belorousski, I.
Ciocan-Fontanine, C. Faber, T. Graber, S. Kleiman,
A. Kresch, C. Procesi, K. Ranestad,
H. Tamvakis, J. Thomsen, E. Tj\o{}tta, and A. Vistoli
for comments and suggestions.
A seminar course at the Mittag-Leffler
Institute has led to many improvements. Some related
preprints that have appeared since the Santa Cruz
conference are pointed out in footnotes.

\subsection{Notation}
\label{nota}
In this exposition, for simplicity, we consider
only homology classes of even dimension. To avoid
doubling  indices, we set, for a complete variety $X$,
$$A_d X=H_{2d}(X, \mathbb{Z}), 
\ \ A^d X = H^{2d}(X, \mathbb{Z}).$$
When $X$ is nonsingular of dimension $n$, identify
$A^d X$ with $A_{n-d} X$ by the Poincar\'e duality
isomorphism
$$A^d X \stackrel{\sim}{\rarr} A_{n-d} X, \ \ c \mapsto
c \smallcap [X].$$
In particular, a closed subvariety $V$ of $X$ of
pure codimension $d$ determines classes in $A_{n-d} X$ and $A^d X$
via the duality isomorphism. Both of these classes are
denoted by $[V]$.
For homogeneous varieties, which are our main concern,
the Chow groups coincide with the topological groups. 
Hence $A_d X$ and $A^d X$ can be taken to be the
Chow homology and cohomology groups for homogeneous 
varieties (see [F]).

If  $X$  is complete, and  $c$  
is a class in the ring  $A^*X = \bigoplus A^dX$,  and  $\beta$  is a 
class in  $A_kX$,  we denote by  $\int_\beta c$  the degree of the 
class of the zero cycle obtained by evaluating  $c_k$  on  $\beta$,  
where  $c_k$  is the component of  $c$  in  $A^kX$.  When  $V$  is a 
closed, pure dimensional
 subvariety of  $X$,  we write  $\int_V c$  instead of  
$\int_{[V]} c$.   We use cup 
product notation  $\smallcup$  for the product in  $A^*X$. 

We 
will be concerned only with varieties over $\com$
since the relevant moduli spaces have 
not yet been constructed in positive characteristic.  
Let  $[n]$  denote the finite set of integers $\{1,2, \ldots,n\}$.

\subsection{{Compactifications of moduli spaces }}
\label{rmod}
An important feature of quantum cohomology is the use of 
intersection theory in a space of maps of curves into a variety, 
rather than in the variety itself.  To carry this out, a good 
compactification of such a space is required.
At least when  $X$  is sufficiently 
positive, Kontsevich has constructed such a compactification. 
We start, in this section, by 
reviewing some related moduli spaces with similar properties.
Kontsevich's space of stable maps will be introduced
in section \ref{kontdef}.

Algebraic geometers by now have become quite comfortable working 
with the moduli space  $M_g$  of projective nonsingular curves of 
genus  $g$,  and its compactification  $\M_g$,  whose points 
correspond to projective, connected, nodal curves of arithmetic genus  
$g$,  satisfying a stability condition (due to
Deligne and Mumford) that guarantees the curve 
has only a finite automorphism group.  These 
moduli spaces are irreducible 
varieties of dimension  $3g - 3$  if  $g \geq 2$,  smooth if regarded 
as (Deligne-Mumford) stacks, and with orbifold singularities if 
regarded as ordinary coarse moduli spaces.

Some related spaces have become increasingly important.  The 
moduli space  $M_{g,n}$  parametrizes projective nonsingular curves  
$C$  of genus  $g$  together with  $n$  distinct marked points  $p_1, 
\ldots , p_n$  on  $C$.   $M_{g,n}$ has a compactification  $\M_{g,n}$  
whose points correspond to projective, connected, nodal curves  $C$,  
together with  $n$  distinct, nonsingular, marked
points, again with a 
stability condition equivalent to the finiteness of automorphism 
groups.  $\M_{g,1}$ is often called
the universal curve over  $\M_g$ (although, 
as coarse moduli spaces, this is
a slight abuse of language).

The first remarkable feature of the space $\M_{g,n}$ 
is that it compactifies  
$M_{g,n}$  without ever allowing the points to come together.  When 
points on a smooth curve approach each other, in fact, the curve 
sprouts off one or more components, each isomorphic to the 
projective line, and the points distribute themselves at smooth 
points on these new components, in a way that reflects the relative 
rates of approach.

\vspace{-0pt}
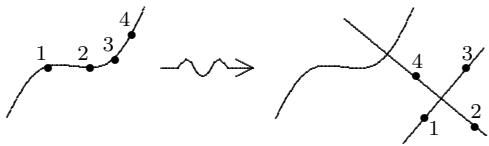
\begin{center}
\font\thinlinefont=cmr5
\begingroup\makeatletter\ifx\SetFigFont\undefined
\def\x#1#2#3#4#5#6#7\relax{\def\x{#1#2#3#4#5#6}}%
\expandafter\x\fmtname xxxxxx\relax \def\y{splain}%
\ifx\x\y   
\gdef\SetFigFont#1#2#3{%
  \ifnum #1<17\tiny\else \ifnum #1<20\small\else
  \ifnum #1<24\normalsize\else \ifnum #1<29\large\else
  \ifnum #1<34\Large\else \ifnum #1<41\LARGE\else
     \huge\fi\fi\fi\fi\fi\fi
  \csname #3\endcsname}%
\else
\gdef\SetFigFont#1#2#3{\begingroup
  \count@#1\relax \ifnum 25<\count@\count@25\fi
  \def\x{\endgroup\@setsize\SetFigFont{#2pt}}%
  \expandafter\x
    \csname \romannumeral\the\count@ pt\expandafter\endcsname
    \csname @\romannumeral\the\count@ pt\endcsname
  \csname #3\endcsname}%
\fi
\fi\endgroup
\mbox{\beginpicture
\setcoordinatesystem units <0.70000cm,0.70000cm>
\unitlength=0.70000cm
\linethickness=1pt
\setplotsymbol ({\makebox(0,0)[l]{\tencirc\symbol{'160}}})
\setshadesymbol ({\thinlinefont .})
\setlinear
%
%
\linethickness= 0.500pt
\setplotsymbol ({\thinlinefont .})
\put{\makebox(0,0)[l]{\circle*{ 0.157}}} at  3.888 23.495
%
%
\linethickness= 0.500pt
\setplotsymbol ({\thinlinefont .})
\put{\makebox(0,0)[l]{\circle*{ 0.157}}} at  3.571 23.019
%
%
\linethickness= 0.500pt
\setplotsymbol ({\thinlinefont .})
\put{\makebox(0,0)[l]{\circle*{ 0.157}}} at  3.095 22.860
%
%
\linethickness= 0.500pt
\setplotsymbol ({\thinlinefont .})
\put{\makebox(0,0)[l]{\circle*{ 0.157}}} at  2.301 22.860
%
%
\linethickness= 0.500pt
\setplotsymbol ({\thinlinefont .})
\put{\makebox(0,0)[l]{\circle*{ 0.157}}} at  9.286 22.701
%
%
\linethickness= 0.500pt
\setplotsymbol ({\thinlinefont .})
\put{\makebox(0,0)[l]{\circle*{ 0.157}}} at 10.397 21.749
%
%
\linethickness= 0.500pt
\setplotsymbol ({\thinlinefont .})
\put{\makebox(0,0)[l]{\circle*{ 0.157}}} at 10.238 22.860
%
%
\linethickness= 0.500pt
\setplotsymbol ({\thinlinefont .})
\put{\makebox(0,0)[l]{\circle*{ 0.157}}} at  9.445 21.907
%
%
\linethickness= 0.500pt
\setplotsymbol ({\thinlinefont .})
\plot  6.191 22.860  5.874 22.701 /
%
%
\linethickness= 0.500pt
\setplotsymbol ({\thinlinefont .})
\plot  6.193 22.862  5.870 23.021 /
%
%
\linethickness= 0.500pt
\setplotsymbol ({\thinlinefont .})
\plot  7.925 23.796 10.624 21.573 /
%
%
\linethickness= 0.500pt
\setplotsymbol ({\thinlinefont .})
\plot 10.573 23.233  9.049 21.431 /
\linethickness= 0.500pt
\setplotsymbol ({\thinlinefont .})
%
%
\plot  6.623 21.907      6.677 22.013
         6.728 22.111
         6.777 22.201
         6.823 22.284
         6.868 22.360
         6.911 22.430
         6.993 22.551
         7.070 22.650
         7.146 22.730
         7.222 22.793
         7.300 22.843
         7.425 22.886
         7.496 22.896
         7.571 22.900
         7.651 22.898
         7.733 22.893
         7.816 22.886
         7.901 22.878
         7.986 22.869
         8.070 22.862
         8.152 22.858
         8.231 22.857
         8.307 22.861
         8.378 22.872
         8.503 22.917
         8.591 22.975
         8.675 23.049
         8.759 23.142
         8.844 23.256
         8.887 23.322
         8.931 23.395
         8.977 23.475
         9.025 23.563
         9.074 23.658
         9.125 23.761
         9.179 23.873
         9.235 23.995
        /
\linethickness= 0.500pt
\setplotsymbol ({\thinlinefont .})
%
%
\plot  4.445 22.860      4.549 22.860
         4.604 22.860
         4.683 22.860
         4.763 22.860
         4.763 22.860
         4.830 22.948
         4.921 23.019
         5.020 22.969
         5.080 22.860
         5.140 22.751
         5.239 22.701
         5.337 22.751
         5.397 22.860
         5.458 22.969
         5.556 23.019
         5.647 22.948
         5.715 22.860
         5.796 22.854
         5.874 22.860
         5.953 22.860
         6.032 22.860
         6.088 22.860
         6.191 22.860
        /
\linethickness= 0.500pt
\setplotsymbol ({\thinlinefont .})
%
%
\plot  1.518 21.922      1.572 22.028
         1.623 22.126
         1.671 22.216
         1.718 22.299
         1.763 22.375
         1.805 22.444
         1.887 22.565
         1.965 22.665
         2.041 22.744
         2.117 22.808
         2.195 22.858
         2.319 22.901
         2.390 22.911
         2.466 22.914
         2.545 22.913
         2.627 22.908
         2.711 22.901
         2.796 22.892
         2.880 22.884
         2.964 22.877
         3.046 22.873
         3.126 22.872
         3.201 22.876
         3.272 22.887
         3.397 22.932
         3.485 22.990
         3.570 23.064
         3.654 23.156
         3.738 23.271
         3.782 23.337
         3.826 23.410
         3.872 23.490
         3.919 23.577
         3.968 23.673
         4.020 23.776
         4.073 23.888
         4.130 24.009
        /
%
%
\put{\SetFigFont{8}{9.6}{rm}4} [lB] at  3.651 23.654
%
%
\put{\SetFigFont{8}{9.6}{rm}2} [lB] at  2.857 23.019
%
%
\put{\SetFigFont{8}{9.6}{rm}1} [lB] at  2.064 23.019
%
%
\put{\SetFigFont{8}{9.6}{rm}3} [lB] at  3.334 23.178
%
%
\put{\SetFigFont{8}{9.6}{rm}1} [lB] at  9.525 21.590
%
%
\put{\SetFigFont{8}{9.6}{rm}2} [lB] at 10.319 21.907
%
%
\put{\SetFigFont{8}{9.6}{rm}4} [lB] at  9.207 22.860
%
%
\put{\SetFigFont{8}{9.6}{rm}3} [lB] at 10.160 23.019
\linethickness=0pt
\putrectangle corners at  1.501 24.026 and 10.649 21.406
\endpicture}
\end{center}
\vspace{-0pt}

The spaces $\M_{g,n}$ again 
are smooth stacks, or orbifold coarse moduli spaces, of 
dimension  $3g - 3 + n$,  as long as this integer is nonnegative.  The 
case of genus zero plays a prominent role in our story.  In this case, 
  $\M_{0,n}$  is a fine moduli space and a nonsingular variety.  A 
point in  $\M_{0,n}$  corresponds to a curve which is a tree of 
projective lines meeting transversally, with  $n$ distinct,
nonsingular, marked  
points;  the stability condition is that each component must have at 
least three special points, which are either the marked points or 
the nodes where the component meets the other components.

For  $n = 3$,  of course,  $M_{0,3} = \M_{0,3}$  is a point.  For  $n = 
4$,  $M_{0,4}$  parametrizes  $4$  distinct marked points on a 
projective line.  Since, up to isomorphism, one can fix the first three 
of these points, say to be  $0$,  $1$,  and  $\infty$,  $M_{0,4}$ is 
isomorphic to  $\proj^1 \smallstm \{0,1,\infty\}$.  It is not 
hard to guess what  $\M_{0,4}$  must be.  In fact, the three added 
points are represented by the following three marked curves:

\vspace{-0pt}
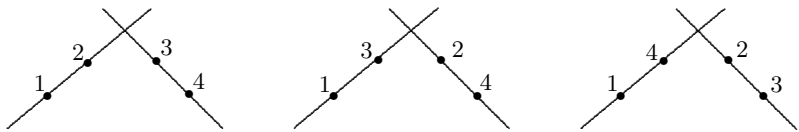
\begin{center}

\font\thinlinefont=cmr5
\begingroup\makeatletter\ifx\SetFigFont\undefined
\def\x#1#2#3#4#5#6#7\relax{\def\x{#1#2#3#4#5#6}}%
\expandafter\x\fmtname xxxxxx\relax \def\y{splain}%
\ifx\x\y   
\gdef\SetFigFont#1#2#3{%
  \ifnum #1<17\tiny\else \ifnum #1<20\small\else
  \ifnum #1<24\normalsize\else \ifnum #1<29\large\else
  \ifnum #1<34\Large\else \ifnum #1<41\LARGE\else
     \huge\fi\fi\fi\fi\fi\fi
  \csname #3\endcsname}%
\else
\gdef\SetFigFont#1#2#3{\begingroup
  \count@#1\relax \ifnum 25<\count@\count@25\fi
  \def\x{\endgroup\@setsize\SetFigFont{#2pt}}%
  \expandafter\x
    \csname \romannumeral\the\count@ pt\expandafter\endcsname
    \csname @\romannumeral\the\count@ pt\endcsname
  \csname #3\endcsname}%
\fi
\fi\endgroup
\mbox{\beginpicture
\setcoordinatesystem units <0.50000cm,0.50000cm>
\unitlength=0.50000cm
\linethickness=1pt
\setplotsymbol ({\makebox(0,0)[l]{\tencirc\symbol{'160}}})
\setshadesymbol ({\thinlinefont .})
\setlinear
%
%
\linethickness= 0.500pt
\setplotsymbol ({\thinlinefont .})
\put{\makebox(0,0)[l]{\circle*{ 0.182}}} at  2.328 22.464
%
%
\linethickness= 0.500pt
\setplotsymbol ({\thinlinefont .})
\put{\makebox(0,0)[l]{\circle*{ 0.182}}} at  3.397 23.336
%
%
\linethickness= 0.500pt
\setplotsymbol ({\thinlinefont .})
\put{\makebox(0,0)[l]{\circle*{ 0.182}}} at  5.222 23.368
%
%
\linethickness= 0.500pt
\setplotsymbol ({\thinlinefont .})
\put{\makebox(0,0)[l]{\circle*{ 0.182}}} at  6.096 22.511
%
%
\linethickness= 0.500pt
\setplotsymbol ({\thinlinefont .})
\put{\makebox(0,0)[l]{\circle*{ 0.182}}} at  9.938 22.447
%
%
\linethickness= 0.500pt
\setplotsymbol ({\thinlinefont .})
\put{\makebox(0,0)[l]{\circle*{ 0.182}}} at 11.127 23.417
%
%
\linethickness= 0.500pt
\setplotsymbol ({\thinlinefont .})
\put{\makebox(0,0)[l]{\circle*{ 0.182}}} at 12.795 23.417
%
%
\linethickness= 0.500pt
\setplotsymbol ({\thinlinefont .})
\put{\makebox(0,0)[l]{\circle*{ 0.182}}} at 13.760 22.447
%
%
\linethickness= 0.500pt
\setplotsymbol ({\thinlinefont .})
\put{\makebox(0,0)[l]{\circle*{ 0.182}}} at 17.570 22.464
%
%
\linethickness= 0.500pt
\setplotsymbol ({\thinlinefont .})
\put{\makebox(0,0)[l]{\circle*{ 0.182}}} at 18.713 23.385
%
%
\linethickness= 0.500pt
\setplotsymbol ({\thinlinefont .})
\put{\makebox(0,0)[l]{\circle*{ 0.182}}} at 20.428 23.400
%
%
\linethickness= 0.500pt
\setplotsymbol ({\thinlinefont .})
\put{\makebox(0,0)[l]{\circle*{ 0.182}}} at 21.364 22.447
%
%
\linethickness= 0.500pt
\setplotsymbol ({\thinlinefont .})
\plot  1.270 21.590  5.080 24.765 /
%
%
\linethickness= 0.500pt
\setplotsymbol ({\thinlinefont .})
\plot  3.810 24.765  6.985 21.590 /
%
%
\linethickness= 0.500pt
\setplotsymbol ({\thinlinefont .})
\plot  8.901 21.607 12.711 24.782 /
%
%
\linethickness= 0.500pt
\setplotsymbol ({\thinlinefont .})
\plot 11.424 24.765 14.599 21.590 /
%
%
\linethickness= 0.500pt
\setplotsymbol ({\thinlinefont .})
\plot 16.521 21.590 20.331 24.765 /
%
%
\linethickness= 0.500pt
\setplotsymbol ({\thinlinefont .})
\plot 19.061 24.750 22.236 21.575 /
%
%
\put{\SetFigFont{9}{10.8}{rm}1} [lB] at  9.542 22.543
%
%
\put{\SetFigFont{9}{10.8}{rm}3} [lB] at 10.653 23.448
%
%
\put{\SetFigFont{9}{10.8}{rm}2} [lB] at 13.075 23.480
%
%
\put{\SetFigFont{9}{10.8}{rm}4} [lB] at 13.822 22.608
%
%
\put{\SetFigFont{9}{10.8}{rm}1} [lB] at 17.187 22.591
%
%
\put{\SetFigFont{9}{10.8}{rm}4} [lB] at 18.235 23.400
%
%
\put{\SetFigFont{9}{10.8}{rm}2} [lB] at 20.633 23.480
%
%
\put{\SetFigFont{9}{10.8}{rm}3} [lB] at 21.552 22.559
%
%
\put{\SetFigFont{9}{10.8}{rm}1} [lB] at  1.969 22.591
%
%
\put{\SetFigFont{9}{10.8}{rm}2} [lB] at  2.985 23.417
%
%
\put{\SetFigFont{9}{10.8}{rm}3} [lB] at  5.332 23.544
%
%
\put{\SetFigFont{9}{10.8}{rm}4} [lB] at  6.191 22.623
\linethickness=0pt
\putrectangle corners at  1.245 24.807 and 22.261 21.550
\endpicture}

\end{center}
\vspace{-0pt}

In general, the closures of the loci of trees of a given combinatorial 
type are smooth subvarieties of  $\M_{0,n}$,  and all such loci meet 
transversally.  There is a divisor  $D(A | B)$  in  $\M_{0,n}$  for 
each partition of  $[n]$  into two disjoint sets  $A$  and  $B$,  each 
with at least two elements.  A generic point of  $D(A | B)$  is 
represented by two lines meeting transversally, with points labeled 
by  $A$  and  $B$  on each:

\vspace{-0pt}
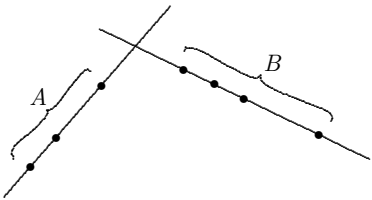
\begin{center}

\font\thinlinefont=cmr5
\begingroup\makeatletter\ifx\SetFigFont\undefined
\def\x#1#2#3#4#5#6#7\relax{\def\x{#1#2#3#4#5#6}}%
\expandafter\x\fmtname xxxxxx\relax \def\y{splain}%
\ifx\x\y   
\gdef\SetFigFont#1#2#3{%
  \ifnum #1<17\tiny\else \ifnum #1<20\small\else
  \ifnum #1<24\normalsize\else \ifnum #1<29\large\else
  \ifnum #1<34\Large\else \ifnum #1<41\LARGE\else
     \huge\fi\fi\fi\fi\fi\fi
  \csname #3\endcsname}%
\else
\gdef\SetFigFont#1#2#3{\begingroup
  \count@#1\relax \ifnum 25<\count@\count@25\fi
  \def\x{\endgroup\@setsize\SetFigFont{#2pt}}%
  \expandafter\x
    \csname \romannumeral\the\count@ pt\expandafter\endcsname
    \csname @\romannumeral\the\count@ pt\endcsname
  \csname #3\endcsname}%
\fi
\fi\endgroup
\mbox{\beginpicture
\setcoordinatesystem units <0.50000cm,0.50000cm>
\unitlength=0.50000cm
\linethickness=1pt
\setplotsymbol ({\makebox(0,0)[l]{\tencirc\symbol{'160}}})
\setshadesymbol ({\thinlinefont .})
\setlinear
%
%
\linethickness= 0.500pt
\setplotsymbol ({\thinlinefont .})
\put{\makebox(0,0)[l]{\circle*{ 0.178}}} at 10.863 20.489
%
%
\linethickness= 0.500pt
\setplotsymbol ({\thinlinefont .})
\put{\makebox(0,0)[l]{\circle*{ 0.178}}} at 11.553 21.241
%
%
\linethickness= 0.500pt
\setplotsymbol ({\thinlinefont .})
\put{\makebox(0,0)[l]{\circle*{ 0.178}}} at 12.759 22.638
%
%
\linethickness= 0.500pt
\setplotsymbol ({\thinlinefont .})
\put{\makebox(0,0)[l]{\circle*{ 0.178}}} at 14.935 23.067
%
%
\linethickness= 0.500pt
\setplotsymbol ({\thinlinefont .})
\put{\makebox(0,0)[l]{\circle*{ 0.178}}} at 15.761 22.686
%
%
\linethickness= 0.500pt
\setplotsymbol ({\thinlinefont .})
\put{\makebox(0,0)[l]{\circle*{ 0.178}}} at 16.538 22.289
%
%
\linethickness= 0.500pt
\setplotsymbol ({\thinlinefont .})
\put{\makebox(0,0)[l]{\circle*{ 0.178}}} at 18.538 21.336
%
%
\linethickness= 0.500pt
\setplotsymbol ({\thinlinefont .})
\plot 10.173 19.685 14.586 24.750 /
%
%
\linethickness= 0.500pt
\setplotsymbol ({\thinlinefont .})
\plot 12.713 24.147 19.888 20.686 /
\linethickness= 0.500pt
\setplotsymbol ({\thinlinefont .})
%
%
\plot 10.422 20.688     10.422 20.758
        10.374 20.815
        10.394 20.925
        10.445 21.025
        10.520 21.114
        10.610 21.195
        10.705 21.268
        10.797 21.335
        10.878 21.396
        10.939 21.452
        11.030 21.547
        11.148 21.664
        11.260 21.790
        11.333 21.907
        11.335 21.958
        11.314 21.994
        11.350 21.990
        11.419 22.005
        11.543 22.111
        11.607 22.183
        11.670 22.263
        11.731 22.343
        11.788 22.420
        11.883 22.543
        11.974 22.657
        12.030 22.729
        12.088 22.804
        12.149 22.878
        12.209 22.946
        12.321 23.048
        12.393 23.074
        12.474 23.021
        /
\linethickness= 0.500pt
\setplotsymbol ({\thinlinefont .})
%
%
\plot 14.967 23.552     14.997 23.654
        15.109 23.639
        15.170 23.645
        15.276 23.616
        15.396 23.567
        15.460 23.536
        15.525 23.503
        15.591 23.468
        15.657 23.432
        15.723 23.395
        15.787 23.357
        15.910 23.286
        16.020 23.224
        16.112 23.175
        16.178 23.142
        16.258 23.101
        16.348 23.054
        16.443 23.005
        16.539 22.958
        16.631 22.916
        16.716 22.882
        16.789 22.860
        16.857 22.851
        16.927 22.858
        16.988 22.913
        16.988 22.839
        17.028 22.767
        17.081 22.705
        17.180 22.636
        17.306 22.567
        17.376 22.532
        17.450 22.499
        17.526 22.465
        17.603 22.433
        17.680 22.401
        17.757 22.371
        17.832 22.341
        17.903 22.313
        17.971 22.287
        18.034 22.262
        18.140 22.217
        18.209 22.188
        18.295 22.156
        18.392 22.122
        18.493 22.084
        18.594 22.043
        18.689 22.000
        18.771 21.954
        18.836 21.905
        18.900 21.789
        18.862 21.709
        /
%
%
\put{\SetFigFont{9}{10.8}{it}A} [lB] at 10.833 22.132
%
%
\put{\SetFigFont{9}{10.8}{it}B} [lB] at 17.062 23.042
\linethickness=0pt
\putrectangle corners at 10.147 24.776 and 19.914 19.660
\endpicture}

\end{center}
\vspace{-0pt}

        It is convenient to allow labeling by finite sets other than  
$[n]$; we write  $\M_{g,A}$  for the corresponding moduli space 
where  $A$  is the labeling set.  
Let $B\subset A$ (if $g=0$, then let $|B|\geq 3$).
It is a fundamental fact that there 
is a morphism  $\M_{g,A} \rarr \M_{g,B}$ 
which ``forgets'' the points marked in  $A 
\smallstm B$.  On the open locus  $M_{g,n}$  this map is the 
obvious one, but it is more subtle on the boundary: removing some 
points may make a component unstable, and such a component must 
be collapsed.  For example, the map from  $\M_{0,5}$  to  $\M_{0,4}$  
forgetting the point labeled  $5$  sends

\vspace{-0pt}
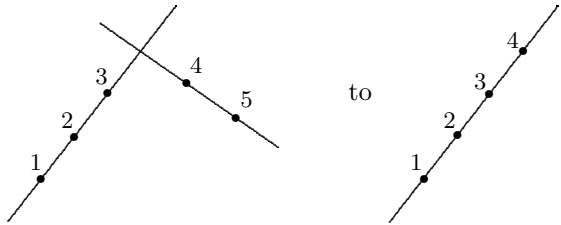
\begin{center}

\font\thinlinefont=cmr5
\begingroup\makeatletter\ifx\SetFigFont\undefined
\def\x#1#2#3#4#5#6#7\relax{\def\x{#1#2#3#4#5#6}}%
\expandafter\x\fmtname xxxxxx\relax \def\y{splain}%
\ifx\x\y   
\gdef\SetFigFont#1#2#3{%
  \ifnum #1<17\tiny\else \ifnum #1<20\small\else
  \ifnum #1<24\normalsize\else \ifnum #1<29\large\else
  \ifnum #1<34\Large\else \ifnum #1<41\LARGE\else
     \huge\fi\fi\fi\fi\fi\fi
  \csname #3\endcsname}%
\else
\gdef\SetFigFont#1#2#3{\begingroup
  \count@#1\relax \ifnum 25<\count@\count@25\fi
  \def\x{\endgroup\@setsize\SetFigFont{#2pt}}%
  \expandafter\x
    \csname \romannumeral\the\count@ pt\expandafter\endcsname
    \csname @\romannumeral\the\count@ pt\endcsname
  \csname #3\endcsname}%
\fi
\fi\endgroup
\mbox{\beginpicture
\setcoordinatesystem units <0.50000cm,0.50000cm>
\unitlength=0.50000cm
\linethickness=1pt
\setplotsymbol ({\makebox(0,0)[l]{\tencirc\symbol{'160}}})
\setshadesymbol ({\thinlinefont .})
\setlinear
%
%
\linethickness= 0.500pt
\setplotsymbol ({\thinlinefont .})
\put{\makebox(0,0)[l]{\circle*{ 0.207}}} at  9.751 20.170
%
%
\linethickness= 0.500pt
\setplotsymbol ({\thinlinefont .})
\put{\makebox(0,0)[l]{\circle*{ 0.207}}} at 10.636 21.294
%
%
\linethickness= 0.500pt
\setplotsymbol ({\thinlinefont .})
\put{\makebox(0,0)[l]{\circle*{ 0.207}}} at 11.525 22.454
%
%
\linethickness= 0.500pt
\setplotsymbol ({\thinlinefont .})
\put{\makebox(0,0)[l]{\circle*{ 0.207}}} at 13.621 22.708
%
%
\linethickness= 0.500pt
\setplotsymbol ({\thinlinefont .})
\put{\makebox(0,0)[l]{\circle*{ 0.207}}} at 14.937 21.787
%
%
\linethickness= 0.500pt
\setplotsymbol ({\thinlinefont .})
\put{\makebox(0,0)[l]{\circle*{ 0.207}}} at 19.939 20.146
%
%
\linethickness= 0.500pt
\setplotsymbol ({\thinlinefont .})
\put{\makebox(0,0)[l]{\circle*{ 0.207}}} at 20.828 21.328
%
%
\linethickness= 0.500pt
\setplotsymbol ({\thinlinefont .})
\put{\makebox(0,0)[l]{\circle*{ 0.207}}} at 21.670 22.432
%
%
\linethickness= 0.500pt
\setplotsymbol ({\thinlinefont .})
\put{\makebox(0,0)[l]{\circle*{ 0.207}}} at 22.574 23.559
%
%
\linethickness= 0.500pt
\setplotsymbol ({\thinlinefont .})
\plot  8.888 19.035 13.348 24.765 /
\putrule from 13.348 24.765 to 13.348 24.782
%
%
\linethickness= 0.500pt
\setplotsymbol ({\thinlinefont .})
\plot 11.333 24.304 16.063 21.004 /
%
%
\linethickness= 0.500pt
\setplotsymbol ({\thinlinefont .})
\plot 19.033 19.018 23.493 24.748 /
\putrule from 23.493 24.748 to 23.493 24.765
%
%
\put{\SetFigFont{9}{10.8}{rm}1} [lB] at 19.571 20.400
%
%
\put{\SetFigFont{9}{10.8}{rm}2} [lB] at 20.460 21.512
%
%
\put{\SetFigFont{9}{10.8}{rm}3} [lB] at 21.270 22.545
%
%
\put{\SetFigFont{9}{10.8}{rm}4} [lB] at 22.142 23.654
%
%
\put{\SetFigFont{9}{10.8}{rm}3} [lB] at 11.206 22.686
%
%
\put{\SetFigFont{9}{10.8}{rm}2} [lB] at 10.285 21.529
%
%
\put{\SetFigFont{9}{10.8}{rm}1} [lB] at  9.459 20.400
%
%
\put{\SetFigFont{9}{10.8}{rm}4} [lB] at 13.714 22.972
%
%
\put{\SetFigFont{9}{10.8}{rm}5} [lB] at 15.064 22.081
%
%
\put{\SetFigFont{10}{12.0}{rm}to} [lB] at 17.937 22.272
\linethickness=0pt
\putrectangle corners at  8.862 24.807 and 23.518 18.993
\endpicture}

\end{center}
\vspace{-0pt}

\noindent and

\vspace{-0pt}
\begin{center}

\font\thinlinefont=cmr5
\begingroup\makeatletter\ifx\SetFigFont\undefined
\def\x#1#2#3#4#5#6#7\relax{\def\x{#1#2#3#4#5#6}}%
\expandafter\x\fmtname xxxxxx\relax \def\y{splain}%
\ifx\x\y   
\gdef\SetFigFont#1#2#3{%
  \ifnum #1<17\tiny\else \ifnum #1<20\small\else
  \ifnum #1<24\normalsize\else \ifnum #1<29\large\else
  \ifnum #1<34\Large\else \ifnum #1<41\LARGE\else
     \huge\fi\fi\fi\fi\fi\fi
  \csname #3\endcsname}%
\else
\gdef\SetFigFont#1#2#3{\begingroup
  \count@#1\relax \ifnum 25<\count@\count@25\fi
  \def\x{\endgroup\@setsize\SetFigFont{#2pt}}%
  \expandafter\x
    \csname \romannumeral\the\count@ pt\expandafter\endcsname
    \csname @\romannumeral\the\count@ pt\endcsname
  \csname #3\endcsname}%
\fi
\fi\endgroup
\mbox
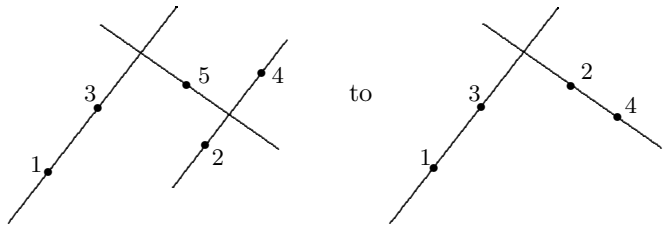
{\beginpicture
\setcoordinatesystem units <0.50000cm,0.50000cm>
\unitlength=0.50000cm
\linethickness=1pt
\setplotsymbol ({\makebox(0,0)[l]{\tencirc\symbol{'160}}})
\setshadesymbol ({\thinlinefont .})
\setlinear
%
%
\linethickness= 0.500pt
\setplotsymbol ({\thinlinefont .})
\put{\makebox(0,0)[l]{\circle*{ 0.207}}} at 13.621 22.708
%
%
\linethickness= 0.500pt
\setplotsymbol ({\thinlinefont .})
\put{\makebox(0,0)[l]{\circle*{ 0.207}}} at 14.110 21.099
%
%
\linethickness= 0.500pt
\setplotsymbol ({\thinlinefont .})
\put{\makebox(0,0)[l]{\circle*{ 0.207}}} at 15.608 23.027
%
%
\linethickness= 0.500pt
\setplotsymbol ({\thinlinefont .})
\put{\makebox(0,0)[l]{\circle*{ 0.207}}} at  9.936 20.384
%
%
\linethickness= 0.500pt
\setplotsymbol ({\thinlinefont .})
\put{\makebox(0,0)[l]{\circle*{ 0.207}}} at 11.252 22.081
%
%
\linethickness= 0.500pt
\setplotsymbol ({\thinlinefont .})
\put{\makebox(0,0)[l]{\circle*{ 0.207}}} at 23.834 22.686
%
%
\linethickness= 0.500pt
\setplotsymbol ({\thinlinefont .})
\put{\makebox(0,0)[l]{\circle*{ 0.207}}} at 25.074 21.838
%
%
\linethickness= 0.500pt
\setplotsymbol ({\thinlinefont .})
\put{\makebox(0,0)[l]{\circle*{ 0.207}}} at 21.448 22.126
%
%
\linethickness= 0.500pt
\setplotsymbol ({\thinlinefont .})
\put{\makebox(0,0)[l]{\circle*{ 0.207}}} at 20.185 20.496
%
%
\linethickness= 0.500pt
\setplotsymbol ({\thinlinefont .})
\plot  8.888 19.035 13.348 24.765 /
\putrule from 13.348 24.765 to 13.348 24.782
%
%
\linethickness= 0.500pt
\setplotsymbol ({\thinlinefont .})
\plot 11.333 24.304 16.063 21.004 /
%
%
\linethickness= 0.500pt
\setplotsymbol ({\thinlinefont .})
\plot 19.033 19.018 23.493 24.748 /
\putrule from 23.493 24.748 to 23.493 24.765
%
%
\linethickness= 0.500pt
\setplotsymbol ({\thinlinefont .})
\plot 16.284 23.893 13.252 19.971 /
%
%
\linethickness= 0.500pt
\setplotsymbol ({\thinlinefont .})
\plot 21.507 24.335 26.238 21.035 /
%
%
\put{\SetFigFont{9}{10.8}{rm}1} [lB] at  9.459 20.400
%
%
\put{\SetFigFont{10}{12.0}{rm}to} [lB] at 17.937 22.272
%
%
\put{\SetFigFont{9}{10.8}{rm}2} [lB] at 14.285 20.561
%
%
\put{\SetFigFont{9}{10.8}{rm}4} [lB] at 15.888 22.765
%
%
\put{\SetFigFont{9}{10.8}{rm}5} [lB] at 13.921 22.748
%
%
\put{\SetFigFont{9}{10.8}{rm}3} [lB] at 10.888 22.274
%
%
\put{\SetFigFont{9}{10.8}{rm}2} [lB] at 24.096 22.862
%
%
\put{\SetFigFont{9}{10.8}{rm}4} [lB] at 25.271 21.971
%
%
\put{\SetFigFont{9}{10.8}{rm}3} [lB] at 21.112 22.274
%
%
\put{\SetFigFont{9}{10.8}{rm}1} [lB] at 19.808 20.591
\linethickness=0pt
\putrectangle corners at  8.862 24.807 and 26.264 18.993
\endpicture}

\end{center}
\vspace{-0pt}

\noindent
The algebra that shows this is a morphism is carried out in [Kn].  

In particular, for any subset  $\{i,j,k,l\}$  of four integers in  $[n]$,  
we have a morphism from  $\M_{0,n}$  to  $\M_{0,\{i,j,k,l\}}$.  The 
inverse image of the point  $P(i,j \mid k,l)$     

\vspace{-0pt}
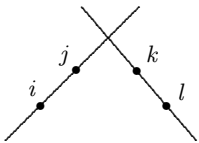
\begin{center}
  
\font\thinlinefont=cmr5
\begingroup\makeatletter\ifx\SetFigFont\undefined
\def\x#1#2#3#4#5#6#7\relax{\def\x{#1#2#3#4#5#6}}%
\expandafter\x\fmtname xxxxxx\relax \def\y{splain}%
\ifx\x\y   
\gdef\SetFigFont#1#2#3{%
  \ifnum #1<17\tiny\else \ifnum #1<20\small\else
  \ifnum #1<24\normalsize\else \ifnum #1<29\large\else
  \ifnum #1<34\Large\else \ifnum #1<41\LARGE\else
     \huge\fi\fi\fi\fi\fi\fi
  \csname #3\endcsname}%
\else
\gdef\SetFigFont#1#2#3{\begingroup
  \count@#1\relax \ifnum 25<\count@\count@25\fi
  \def\x{\endgroup\@setsize\SetFigFont{#2pt}}%
  \expandafter\x
    \csname \romannumeral\the\count@ pt\expandafter\endcsname
    \csname @\romannumeral\the\count@ pt\endcsname
  \csname #3\endcsname}%
\fi
\fi\endgroup
\mbox{\beginpicture
\setcoordinatesystem units <0.40000cm,0.40000cm>
\unitlength=0.40000cm
\linethickness=1pt
\setplotsymbol ({\makebox(0,0)[l]{\tencirc\symbol{'160}}})
\setshadesymbol ({\thinlinefont .})
\setlinear
%
%
\linethickness= 0.500pt
\setplotsymbol ({\thinlinefont .})
\put{\makebox(0,0)[l]{\circle*{ 0.224}}} at  6.905 21.368
%
%
\linethickness= 0.500pt
\setplotsymbol ({\thinlinefont .})
\put{\makebox(0,0)[l]{\circle*{ 0.224}}} at  7.889 20.210
%
%
\linethickness= 0.500pt
\setplotsymbol ({\thinlinefont .})
\put{\makebox(0,0)[l]{\circle*{ 0.224}}} at  4.887 21.406
%
%
\linethickness= 0.500pt
\setplotsymbol ({\thinlinefont .})
\put{\makebox(0,0)[l]{\circle*{ 0.224}}} at  3.700 20.201
%
%
\linethickness= 0.500pt
\setplotsymbol ({\thinlinefont .})
\plot  5.074 23.493  8.884 19.048 /
%
%
\linethickness= 0.500pt
\setplotsymbol ({\thinlinefont .})
\plot  2.540 19.050  6.985 23.495 /
%
%
\put{\SetFigFont{9}{10.8}{it}l} [lB] at  8.240 20.398
%
%
\put{\SetFigFont{9}{10.8}{it}i} [lB] at  3.270 20.527
%
%
\put{\SetFigFont{9}{10.8}{it}j} [lB] at  4.350 21.685
%
%
\put{\SetFigFont{9}{10.8}{it}k} [lB] at  7.205 21.668
\linethickness=0pt
\putrectangle corners at  2.515 23.520 and  8.909 19.022
\endpicture}

\end{center}
\vspace{-0pt}

\noindent
is a divisor on  $\M_{0,n}$.  This inverse image 
is a multiplicity-free sum of divisors  $D(A 
| B)$:  the sum is taken over all partitions 
$A\cup B= [n]$ satisfying   $i,j \in A$    
and $k,l\in B$. 
The fact that the three boundary
points  in   
$\M_{0,\{i,j,k,l\}} \cong \proj^1$  are linearly equivalent implies 
their inverse images in  $\M_{0,n}$  are linearly equivalent as well. 
Hence,  the fundamental relation is obtained:
\begin{equation}
\label{divisor}
\sum_{i,j \in A \  k,l \in B} D(A | B) \, = \; 
 \sum_{i,k \in A \  j,l \in B} D(A | B) \; = \; 
\sum_{i,l \in A \  j,k \in B} D(A | B)  
\end{equation}
in  $A^1(\M_{0,n})$.  Keel [Ke] has shown that the classes of these 
divisors $D(A | B)$  generate the Chow ring, and that the 
relations (\ref{divisor}), together with the (geometrically obvious) 
relations  $D(A | B) {\mathbf \cdot} D(A' | B') = 0$  
if there are no inclusions among the sets  $A$, $B$, $A'$, 
$B'$,  give a complete set of relations.

\subsection{The space of stable maps}
\label{kontdef}
In the remainder of the introduction, the basic
ideas and constructions in quantum cohomology are introduced.
The goal here is to give a precise overview with
no proofs. The ideas introduced here are
covered carefully (with proofs) in the main
sections of these notes.

Let $X$  be  a smooth projective variety, and let  $\beta$  be an 
element in $A_1 X$.  Let  $M_{g,n}(X,\beta)$  be
the set of isomorphism classes of pointed maps  
$(C, p_1, \ldots, p_n, \mu)$ where $C$ is a 
projective nonsingular curve 
of genus   $g$, the markings
$p_1, \ldots , p_n$   are distinct points of  $C$, and
$\mu$ is  a morphism from  $C$ to $ X$  
satisfying $\mu_*([C]) = 
\beta$.
$(C, p_1, \ldots , p_n, \mu)$  is {\em isomorphic}
 to  $(C', 
p_1', \ldots, p_n', \mu' )$  if there is a
scheme 
isomorphism  $\tau : C \rarr C'$  taking  $p_i$  to  
$p_i'$,  with  $\mu'\circ\tau = \mu$.  Of course, if  $\beta \neq 
0$,  $M_{g,n}(X,\beta)$  is empty unless  $\beta$  is the class of a 
curve in  $X$.  There are also 
other problems.  For example, if  $g = 0$,  
which will be the case of interest to us, 
$M_{g,n}(X,\beta)$ is empty if  $\beta \neq 
0$  and  $X$  contains no rational curves.  To obtain a well-behaved 
space, one needs to make strong assumptions on  $X$.  In 
general, there is a compactification
$$
M_{g,n}(X,\beta)  \subset  \M_{g,n}(X,\beta) ,
$$
whose points correspond to stable maps
$(C, p_1, \ldots , p_n, \mu)$ where  $C$  a projective, 
connected, nodal curve of arithmetic genus  $g$,  
the markings $p_1, \ldots , p_n$ are  
distinct nonsingular points of  $C$,  and  $\mu$  is a morphism from  
$C$  such that  $\mu_*([C]) = \beta$.
Again, the stability condition (due to Kontsevich) 
is equivalent to finiteness
of automorphisms of the map. Alternatively, 
$(C, p_1, \ldots, p_n, \mu)$ is a {\em stable} map if both of the
following conditions hold for every irreducible
component $E\subset C$:
\begin{enumerate}
\item[(1)]
If $E\eqq \proj^1$ and $E$ is  mapped to a point by  $\mu$, then  
$E$ must contain at least three special points (either marked points 
or points where  $E$  meets the other components of  $C$). 
\item[(2)]
If $E$ has arithmetic genus 1 and $E$ is mapped to
a point by $\mu$, then $E$ must contain at least one special point. 
\end{enumerate}
Condition (2) is relevant only  in case $g=1$, $n=0$,
and $\beta=0$ (in other cases, (2) is automatically
satisfied). From conditions (1) and (2), it follows that
$\M_{1,0}(X,0)= \emptyset$. Thus, in practice, (1) is the
important condition.

When  $X$  is a point, so  $\beta = 0$,  one recovers the
pointed moduli space of curves  
$\M_{g,n} \eqq \M_{g,n}(\text{point}, 0)$.  
When  $X\eqq \proj^r$  is a projective space, we write  
$\M_{g,n}(\proj^r,d)$  in place of  $\M_{g,n}(\proj^r,
\,d[\text{line}])$.
 
The simplest example is $\M_{0,0}(\proj^r,1)$, which
is the Grassmannian $\grass(\proj^1, \proj^r)$.
If $n\geq 1$, $\M_{0,n}(\proj^r,1)$ is a locally
trivial fibration over $\grass(\proj^1, \proj^r)$
with the configuration space $\proj^1[n]$ of [F-M] as
the fiber.
Let us look at 
the space  $\M_{0,0}(\proj^2,2)$.  An open set in this space is the 
space of nonsingular conics, since to each such conic  $D$  there is 
an isomorphism  $\proj^1 \stackrel {\sim} {\rarr}   D \subset 
\proj^2$,  unique up to equivalence.  Singular conics  $D$  that are 
the unions of two lines are similarly the isomorphic image  $C  
\stackrel {\sim} {\rarr}  D \subset \proj^2$,  where  $C$  is the union 
of two projective lines meeting transversally at a point.  This gives:

\vspace{-0pt}
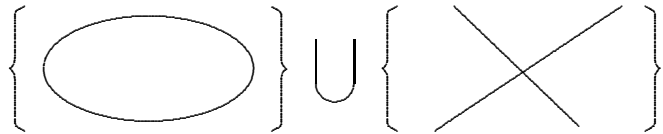
\begin{center}
  
\font\thinlinefont=cmr5
\begingroup\makeatletter\ifx\SetFigFont\undefined
\def\x#1#2#3#4#5#6#7\relax{\def\x{#1#2#3#4#5#6}}%
\expandafter\x\fmtname xxxxxx\relax \def\y{splain}%
\ifx\x\y   
\gdef\SetFigFont#1#2#3{%
  \ifnum #1<17\tiny\else \ifnum #1<20\small\else
  \ifnum #1<24\normalsize\else \ifnum #1<29\large\else
  \ifnum #1<34\Large\else \ifnum #1<41\LARGE\else
     \huge\fi\fi\fi\fi\fi\fi
  \csname #3\endcsname}%
\else
\gdef\SetFigFont#1#2#3{\begingroup
  \count@#1\relax \ifnum 25<\count@\count@25\fi
  \def\x{\endgroup\@setsize\SetFigFont{#2pt}}%
  \expandafter\x
    \csname \romannumeral\the\count@ pt\expandafter\endcsname
    \csname @\romannumeral\the\count@ pt\endcsname
  \csname #3\endcsname}%
\fi
\fi\endgroup
\mbox{\beginpicture
\setcoordinatesystem units <0.40000cm,0.40000cm>
\unitlength=0.40000cm
\linethickness=1pt
\setplotsymbol ({\makebox(0,0)[l]{\tencirc\symbol{'160}}})
\setshadesymbol ({\thinlinefont .})
\setlinear
\linethickness= 0.500pt
\setplotsymbol ({\thinlinefont .})
%
%
%
\plot    2.223 24.765  2.064 24.686
         1.945 24.606
         1.905 24.527
         1.905 24.448
         1.905 24.368
         1.905 24.304
         1.905 24.190
         1.905 24.114
         1.905 24.026
         1.905 23.925
         1.905 23.812
         1.905 23.697
         1.905 23.589
         1.905 23.489
         1.905 23.396
         1.905 23.310
         1.905 23.232
         1.905 23.161
         1.905 23.098
         1.900 22.989
         1.885 22.900
         1.826 22.781
         /
\plot  1.826 22.781  1.746 22.701 /
\linethickness= 0.500pt
\setplotsymbol ({\thinlinefont .})
%
%
%
\plot    2.223 20.637  2.064 20.717
         1.945 20.796
         1.905 20.876
         1.905 20.955
         1.905 21.034
         1.905 21.099
         1.905 21.213
         1.905 21.289
         1.905 21.377
         1.905 21.477
         1.905 21.590
         1.905 21.705
         1.905 21.813
         1.905 21.914
         1.905 22.007
         1.905 22.092
         1.905 22.170
         1.905 22.241
         1.905 22.304
         1.900 22.414
         1.885 22.503
         1.826 22.622
         /
\plot  1.826 22.622  1.746 22.701 /
\linethickness= 0.500pt
\setplotsymbol ({\thinlinefont .})
%
%
%
\plot   10.446 24.750 10.604 24.671
        10.724 24.591
        10.763 24.512
        10.763 24.433
        10.763 24.353
        10.763 24.289
        10.763 24.175
        10.763 24.099
        10.763 24.011
        10.763 23.911
        10.763 23.798
        10.763 23.682
        10.763 23.574
        10.763 23.474
        10.763 23.381
        10.763 23.295
        10.763 23.217
        10.763 23.147
        10.763 23.083
        10.768 22.974
        10.783 22.885
        10.843 22.766
         /
\plot 10.843 22.766 10.922 22.686 /
\linethickness= 0.500pt
\setplotsymbol ({\thinlinefont .})
%
%
%
\plot   10.446 20.623 10.604 20.702
        10.724 20.781
        10.763 20.861
        10.763 20.940
        10.763 21.020
        10.763 21.084
        10.763 21.198
        10.763 21.274
        10.763 21.362
        10.763 21.462
        10.763 21.575
        10.763 21.691
        10.763 21.798
        10.763 21.899
        10.763 21.992
        10.763 22.077
        10.763 22.156
        10.763 22.226
        10.763 22.290
        10.768 22.399
        10.783 22.488
        10.843 22.607
         /
\plot 10.843 22.607 10.922 22.686 /
\linethickness= 0.500pt
\setplotsymbol ({\thinlinefont .})
%
%
%
\plot   14.605 24.765 14.446 24.686
        14.327 24.606
        14.287 24.527
        14.287 24.448
        14.287 24.368
        14.287 24.304
        14.287 24.190
        14.287 24.114
        14.287 24.026
        14.287 23.925
        14.287 23.812
        14.287 23.697
        14.287 23.589
        14.287 23.489
        14.287 23.396
        14.287 23.310
        14.287 23.232
        14.287 23.161
        14.287 23.098
        14.283 22.989
        14.268 22.900
        14.208 22.781
         /
\plot 14.208 22.781 14.129 22.701 /
\linethickness= 0.500pt
\setplotsymbol ({\thinlinefont .})
%
%
%
\plot   14.605 20.637 14.446 20.717
        14.327 20.796
        14.287 20.876
        14.287 20.955
        14.287 21.034
        14.287 21.099
        14.287 21.213
        14.287 21.289
        14.287 21.377
        14.287 21.477
        14.287 21.590
        14.287 21.705
        14.287 21.813
        14.287 21.914
        14.287 22.007
        14.287 22.092
        14.287 22.170
        14.287 22.241
        14.287 22.304
        14.283 22.414
        14.268 22.503
        14.208 22.622
         /
\plot 14.208 22.622 14.129 22.701 /
%
%
\linethickness= 0.500pt
\setplotsymbol ({\thinlinefont .})
\circulararc 180.000 degrees from 11.906 22.225 center at 12.541 22.225
%
%
\linethickness= 0.500pt
\setplotsymbol ({\thinlinefont .})
\putrule from 11.906 22.225 to 11.906 23.654
\putrule from 11.906 23.654 to 11.906 23.654
%
%
\linethickness= 0.500pt
\setplotsymbol ({\thinlinefont .})
\putrule from 13.176 22.225 to 13.176 23.654
\putrule from 13.176 23.654 to 13.176 23.654
\linethickness= 0.500pt
\setplotsymbol ({\thinlinefont .})
%
%
%
\plot   22.860 24.765 23.019 24.686
        23.138 24.606
        23.178 24.527
        23.178 24.448
        23.178 24.368
        23.178 24.304
        23.178 24.190
        23.178 24.114
        23.178 24.026
        23.178 23.925
        23.178 23.812
        23.178 23.697
        23.178 23.589
        23.178 23.489
        23.178 23.396
        23.178 23.310
        23.178 23.232
        23.178 23.161
        23.178 23.098
        23.182 22.989
        23.197 22.900
        23.257 22.781
         /
\plot 23.257 22.781 23.336 22.701 /
\linethickness= 0.500pt
\setplotsymbol ({\thinlinefont .})
%
%
%
\plot   22.860 20.637 23.019 20.717
        23.138 20.796
        23.178 20.876
        23.178 20.955
        23.178 21.034
        23.178 21.099
        23.178 21.213
        23.178 21.289
        23.178 21.377
        23.178 21.477
        23.178 21.590
        23.178 21.705
        23.178 21.813
        23.178 21.914
        23.178 22.007
        23.178 22.092
        23.178 22.170
        23.178 22.241
        23.178 22.304
        23.182 22.414
        23.197 22.503
        23.257 22.622
         /
\plot 23.257 22.622 23.336 22.701 /
%
%
\linethickness= 0.500pt
\setplotsymbol ({\thinlinefont .})
\ellipticalarc axes ratio  3.493:1.746  360 degrees 
        from  9.842 22.701 center at  6.350 22.701
%
%
\linethickness= 0.500pt
\setplotsymbol ({\thinlinefont .})
\plot 15.875 20.637 22.066 24.765 /
%
%
\linethickness= 0.500pt
\setplotsymbol ({\thinlinefont .})
\plot 16.510 24.765 20.637 20.796 /
\linethickness=0pt
\putrectangle corners at  1.729 24.790 and 23.353 20.606
\endpicture}

\end{center}
\vspace{-0pt}

\noindent We also have maps 
from the same  $C$  to  $\proj^2$ sending each 
line in the domain onto the same line in $\proj^2$.
To determine this map up to isomorphism, 
however, the point that is the image of the intersection of the two 
lines must be specified, so the data for this is a line in  $\proj^2$  
together with a point on it.  Finally, there are maps from  $\proj^1$  
to a line in the plane that are branched coverings of degree two onto 
a line in the plane.  These are determined by specifying the line 
together with the two branch points.  The added points consist of:

\vspace{-0pt}
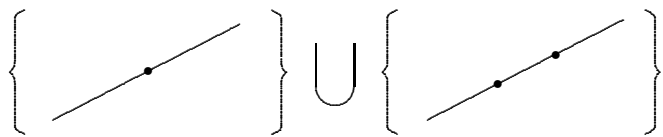
\begin{center}
  
\font\thinlinefont=cmr5
\begingroup\makeatletter\ifx\SetFigFont\undefined
\def\x#1#2#3#4#5#6#7\relax{\def\x{#1#2#3#4#5#6}}%
\expandafter\x\fmtname xxxxxx\relax \def\y{splain}%
\ifx\x\y   
\gdef\SetFigFont#1#2#3{%
  \ifnum #1<17\tiny\else \ifnum #1<20\small\else
  \ifnum #1<24\normalsize\else \ifnum #1<29\large\else
  \ifnum #1<34\Large\else \ifnum #1<41\LARGE\else
     \huge\fi\fi\fi\fi\fi\fi
  \csname #3\endcsname}%
\else
\gdef\SetFigFont#1#2#3{\begingroup
  \count@#1\relax \ifnum 25<\count@\count@25\fi
  \def\x{\endgroup\@setsize\SetFigFont{#2pt}}%
  \expandafter\x
    \csname \romannumeral\the\count@ pt\expandafter\endcsname
    \csname @\romannumeral\the\count@ pt\endcsname
  \csname #3\endcsname}%
\fi
\fi\endgroup
\mbox{\beginpicture
\setcoordinatesystem units <0.40000cm,0.40000cm>
\unitlength=0.40000cm
\linethickness=1pt
\setplotsymbol ({\makebox(0,0)[l]{\tencirc\symbol{'160}}})
\setshadesymbol ({\thinlinefont .})
\setlinear
\linethickness= 0.500pt
\setplotsymbol ({\thinlinefont .})
%
%
%
\plot    2.223 24.765  2.064 24.686
         1.945 24.606
         1.905 24.527
         1.905 24.448
         1.905 24.368
         1.905 24.304
         1.905 24.190
         1.905 24.114
         1.905 24.026
         1.905 23.925
         1.905 23.812
         1.905 23.697
         1.905 23.589
         1.905 23.489
         1.905 23.396
         1.905 23.310
         1.905 23.232
         1.905 23.161
         1.905 23.098
         1.900 22.989
         1.885 22.900
         1.826 22.781
         /
\plot  1.826 22.781  1.746 22.701 /
\linethickness= 0.500pt
\setplotsymbol ({\thinlinefont .})
%
%
%
\plot    2.223 20.637  2.064 20.717
         1.945 20.796
         1.905 20.876
         1.905 20.955
         1.905 21.034
         1.905 21.099
         1.905 21.213
         1.905 21.289
         1.905 21.377
         1.905 21.477
         1.905 21.590
         1.905 21.705
         1.905 21.813
         1.905 21.914
         1.905 22.007
         1.905 22.092
         1.905 22.170
         1.905 22.241
         1.905 22.304
         1.900 22.414
         1.885 22.503
         1.826 22.622
         /
\plot  1.826 22.622  1.746 22.701 /
\linethickness= 0.500pt
\setplotsymbol ({\thinlinefont .})
%
%
%
\plot   10.446 24.750 10.604 24.671
        10.724 24.591
        10.763 24.512
        10.763 24.433
        10.763 24.353
        10.763 24.289
        10.763 24.175
        10.763 24.099
        10.763 24.011
        10.763 23.911
        10.763 23.798
        10.763 23.682
        10.763 23.574
        10.763 23.474
        10.763 23.381
        10.763 23.295
        10.763 23.217
        10.763 23.147
        10.763 23.083
        10.768 22.974
        10.783 22.885
        10.843 22.766
         /
\plot 10.843 22.766 10.922 22.686 /
\linethickness= 0.500pt
\setplotsymbol ({\thinlinefont .})
%
%
%
\plot   10.446 20.623 10.604 20.702
        10.724 20.781
        10.763 20.861
        10.763 20.940
        10.763 21.020
        10.763 21.084
        10.763 21.198
        10.763 21.274
        10.763 21.362
        10.763 21.462
        10.763 21.575
        10.763 21.691
        10.763 21.798
        10.763 21.899
        10.763 21.992
        10.763 22.077
        10.763 22.156
        10.763 22.226
        10.763 22.290
        10.768 22.399
        10.783 22.488
        10.843 22.607
         /
\plot 10.843 22.607 10.922 22.686 /
\linethickness= 0.500pt
\setplotsymbol ({\thinlinefont .})
%
%
%
\plot   14.605 24.765 14.446 24.686
        14.327 24.606
        14.287 24.527
        14.287 24.448
        14.287 24.368
        14.287 24.304
        14.287 24.190
        14.287 24.114
        14.287 24.026
        14.287 23.925
        14.287 23.812
        14.287 23.697
        14.287 23.589
        14.287 23.489
        14.287 23.396
        14.287 23.310
        14.287 23.232
        14.287 23.161
        14.287 23.098
        14.283 22.989
        14.268 22.900
        14.208 22.781
         /
\plot 14.208 22.781 14.129 22.701 /
\linethickness= 0.500pt
\setplotsymbol ({\thinlinefont .})
%
%
%
\plot   14.605 20.637 14.446 20.717
        14.327 20.796
        14.287 20.876
        14.287 20.955
        14.287 21.034
        14.287 21.099
        14.287 21.213
        14.287 21.289
        14.287 21.377
        14.287 21.477
        14.287 21.590
        14.287 21.705
        14.287 21.813
        14.287 21.914
        14.287 22.007
        14.287 22.092
        14.287 22.170
        14.287 22.241
        14.287 22.304
        14.283 22.414
        14.268 22.503
        14.208 22.622
         /
\plot 14.208 22.622 14.129 22.701 /
%
%
\linethickness= 0.500pt
\setplotsymbol ({\thinlinefont .})
\circulararc 180.000 degrees from 11.906 22.225 center at 12.541 22.225
%
%
\linethickness= 0.500pt
\setplotsymbol ({\thinlinefont .})
\putrule from 11.906 22.225 to 11.906 23.654
\putrule from 11.906 23.654 to 11.906 23.654
%
%
\linethickness= 0.500pt
\setplotsymbol ({\thinlinefont .})
\putrule from 13.176 22.225 to 13.176 23.654
\putrule from 13.176 23.654 to 13.176 23.654
\linethickness= 0.500pt
\setplotsymbol ({\thinlinefont .})
%
%
%
\plot   22.860 24.765 23.019 24.686
        23.138 24.606
        23.178 24.527
        23.178 24.448
        23.178 24.368
        23.178 24.304
        23.178 24.190
        23.178 24.114
        23.178 24.026
        23.178 23.925
        23.178 23.812
        23.178 23.697
        23.178 23.589
        23.178 23.489
        23.178 23.396
        23.178 23.310
        23.178 23.232
        23.178 23.161
        23.178 23.098
        23.182 22.989
        23.197 22.900
        23.257 22.781
         /
\plot 23.257 22.781 23.336 22.701 /
\linethickness= 0.500pt
\setplotsymbol ({\thinlinefont .})
%
%
%
\plot   22.860 20.637 23.019 20.717
        23.138 20.796
        23.178 20.876
        23.178 20.955
        23.178 21.034
        23.178 21.099
        23.178 21.213
        23.178 21.289
        23.178 21.377
        23.178 21.477
        23.178 21.590
        23.178 21.705
        23.178 21.813
        23.178 21.914
        23.178 22.007
        23.178 22.092
        23.178 22.170
        23.178 22.241
        23.178 22.304
        23.182 22.414
        23.197 22.503
        23.257 22.622
         /
\plot 23.257 22.622 23.336 22.701 /
%
%
\linethickness= 0.500pt
\setplotsymbol ({\thinlinefont .})
\put{\makebox(0,0)[l]{\circle*{ 0.288}}} at  6.331 22.739
%
%
\linethickness= 0.500pt
\setplotsymbol ({\thinlinefont .})
\put{\makebox(0,0)[l]{\circle*{ 0.288}}} at 17.949 22.297
%
%
\linethickness= 0.500pt
\setplotsymbol ({\thinlinefont .})
\put{\makebox(0,0)[l]{\circle*{ 0.288}}} at 19.873 23.288
%
%
\linethickness= 0.500pt
\setplotsymbol ({\thinlinefont .})
\plot  3.175 21.114  9.366 24.289 /
%
%
\linethickness= 0.500pt
\setplotsymbol ({\thinlinefont .})
\plot 15.621 21.114 22.130 24.448 /
\linethickness=0pt
\putrectangle corners at  1.729 24.782 and 23.353 20.606
\endpicture}

\end{center}
\vspace{-0pt}

\noindent
Thus, we recover the classical space of complete conics -- but in 
quite a different realization from the usual one.  The same 
discussion is valid when  $\proj^2$  is replaced by  $\proj^r$,  but 
this time the space is not the classical space of complete conics in 
space. The classical space specifies a plane together with a 
complete conic contained in the plane;  Kontsevich's space has blown down 
all the planes containing a given line.

Let $X$ be a complete nonsingular variety with
tangent bundle $T_X$. 
$X$ is {\em convex} if,
for every 
morphism  $\mu: \proj^1 \rarr X$, 
\begin{equation}
\label{cconnv}
H^1(\proj^1,\,\mu^*(T_X)) = 0.
\end{equation}
Convexity is a very restrictive condition on $X$. 
The main examples of convex varieties are 
homogeneous spaces 
$X=G/P$, where $G$ is a Lie group
and $P$ is a parabolic subgroup.
Since $G$ acts transitively on $X$, $T_X$ is
generated by global sections. Hence,
$\mu^*(T_X)$ is globally generated for
every morphism of $\proj^1$, and the vanishing
(\ref{cconnv}) is obtained.
Projective spaces, 
Grassmannians, smooth quadrics, flag varieties,
and products of such varieties are all homogeneous.
It is for homogeneous spaces
that the theory of
quantum cohomology takes its simplest form. The development
of quantum cohomology in sections 7--10 is carried out
only in the homogeneous case. Other examples
of convex varieties include abelian varieties and projective bundles over
curves of positive genus.

The genus 0 moduli
space of stable maps is well-behaved in case $X$ is convex.  
In this case,
$\M_{0,n}(X,\beta)$  exists as a projective nonsingular 
stack or orbifold coarse moduli space, containing  $M_{0,n}(X,\beta)$  
as a dense open subset.  
When $\M_{0,n}(X,\beta)$ is nonempty,
its dimension is given by
$$
\dim\M_{0,n}(X,\beta)  =  \dim X  +  \int_\beta c_1(T_X)  +  n - 3 .
$$
Here,  $c_1(T_X)$  is the first Chern class of the tangent 
bundle to  $X$.  We assume always that the right side of this 
equation is nonnegative.  In the stack or orbifold sense, this is a 
compactification with normal crossing divisors.  When  $X$  is 
projective space,  $\M_{0,n}(X,d)$  is irreducible.  
These assertions are Theorems 1--3 in these notes and are established 
in sections 1--6.

We will also write   $\M_{0,A}(X,\beta)$  when the index set is a set  
$A$  instead of  $[n]$.  These varieties also have 
forgetful morphisms   $\M_{0,A}(X,\beta) \rarr  \M_{0,B}(X,\beta) $ 
when  $B$  is a subset of  $A$.  In addition, if
$|A|\geq 3$, there are morphisms  
$\M_{0,A}(X,\beta) \rarr  \M_{0,A}$  that forget the map  $\mu$.  
In both these cases, as before, one must collapse components that 
become unstable.

When $X$ is convex, 
the spaces $\M_{0,n}(X, \beta)$ have fundamental boundary
divisors analogous to the divisors $D(A |B)$ on $\M_{0,n}$.
Let $n\geq 4$.
Let $A\cup B$ be a partition of $[n]$. Let $\beta_1 + \beta_2=\beta$
be a sum in $A_1 X$.
There is a divisor in $\M_{0,n}(X,\beta)$ determined by:
\begin{equation}
\label{bddeff}
D(A,B; \beta_1, \beta_2) = \M_{0, A\cup\{\bullet\}}(X, \beta_1)
                             \times_X
                             \M_{0, B\cup \{\bullet\}}(X, \beta_2),
\end{equation}
$$ D(A,B; \beta_1, \beta_2) \subset \M_{0,n}(X,\beta).$$
A moduli point in $D(A,B, \beta_1, \beta_2)$ corresponds
to a map with a reducible domain $C= C_1 \cup C_2$
where $\mu_*([C_1])= \beta_1$ and $\mu_*([C_2])=\beta_2$.
The points labeled by $A$ lie on $C_1$ and points
labeled by $B$ lie on $C_2$. The curves $C_1$ and $C_2$
are connected at the points labeled $\bullet$.

\vspace{-0pt}
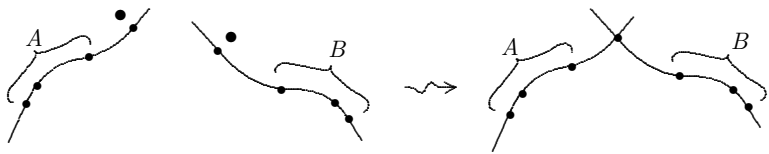
\begin{center}
  
\font\thinlinefont=cmr5
\begingroup\makeatletter\ifx\SetFigFont\undefined
\def\x#1#2#3#4#5#6#7\relax{\def\x{#1#2#3#4#5#6}}%
\expandafter\x\fmtname xxxxxx\relax \def\y{splain}%
\ifx\x\y   
\gdef\SetFigFont#1#2#3{%
  \ifnum #1<17\tiny\else \ifnum #1<20\small\else
  \ifnum #1<24\normalsize\else \ifnum #1<29\large\else
  \ifnum #1<34\Large\else \ifnum #1<41\LARGE\else
     \huge\fi\fi\fi\fi\fi\fi
  \csname #3\endcsname}%
\else
\gdef\SetFigFont#1#2#3{\begingroup
  \count@#1\relax \ifnum 25<\count@\count@25\fi
  \def\x{\endgroup\@setsize\SetFigFont{#2pt}}%
  \expandafter\x
    \csname \romannumeral\the\count@ pt\expandafter\endcsname
    \csname @\romannumeral\the\count@ pt\endcsname
  \csname #3\endcsname}%
\fi
\fi\endgroup
\mbox{\beginpicture
\setcoordinatesystem units <0.35000cm,0.35000cm>
\unitlength=0.35000cm
\linethickness=1pt
\setplotsymbol ({\makebox(0,0)[l]{\tencirc\symbol{'160}}})
\setshadesymbol ({\thinlinefont .})
\setlinear
%
%
\linethickness= 0.500pt
\setplotsymbol ({\thinlinefont .})
\plot 27.163 20.481 26.782 20.284 /
%
%
\linethickness= 0.500pt
\setplotsymbol ({\thinlinefont .})
\plot 27.163 20.481 26.750 20.680 /
\linethickness= 0.500pt
\setplotsymbol ({\thinlinefont .})
%
%
\plot 25.258 20.481     25.362 20.483
        25.417 20.481
        25.490 20.468
        25.561 20.449
        25.628 20.408
        25.707 20.348
        25.794 20.297
        25.885 20.284
        25.986 20.353
        26.055 20.477
        26.126 20.600
        26.234 20.665
        26.314 20.647
        26.384 20.592
        26.505 20.489
        26.585 20.477
        26.685 20.476
        26.784 20.479
        26.863 20.481
        26.954 20.481
        27.045 20.481
        27.163 20.481
        /
%
%
\linethickness= 0.500pt
\setplotsymbol ({\thinlinefont .})
\put{\makebox(0,0)[l]{\circle*{ 0.288}}} at 10.861 19.844
%
%
\linethickness= 0.500pt
\setplotsymbol ({\thinlinefont .})
\put{\makebox(0,0)[l]{\circle*{ 0.288}}} at 11.273 20.542
%
%
\linethickness= 0.500pt
\setplotsymbol ({\thinlinefont .})
\put{\makebox(0,0)[l]{\circle*{ 0.288}}} at 13.242 21.654
%
%
\linethickness= 0.500pt
\setplotsymbol ({\thinlinefont .})
\put{\makebox(0,0)[l]{\circle*{ 0.288}}} at 14.925 22.733
%
%
\linethickness= 0.500pt
\setplotsymbol ({\thinlinefont .})
\put{\makebox(0,0)[l]{\circle*{ 0.288}}} at 18.131 21.876
%
%
\linethickness= 0.500pt
\setplotsymbol ({\thinlinefont .})
\put{\makebox(0,0)[l]{\circle*{ 0.288}}} at 20.544 20.384
%
%
\linethickness= 0.500pt
\setplotsymbol ({\thinlinefont .})
\put{\makebox(0,0)[l]{\circle*{ 0.288}}} at 22.576 19.907
%
%
\linethickness= 0.500pt
\setplotsymbol ({\thinlinefont .})
\put{\makebox(0,0)[l]{\circle*{ 0.288}}} at 23.116 19.272
%
%
\linethickness= 0.500pt
\setplotsymbol ({\thinlinefont .})
\put{\makebox(0,0)[l]{\circle*{ 0.288}}} at 29.244 19.431
%
%
\linethickness= 0.500pt
\setplotsymbol ({\thinlinefont .})
\put{\makebox(0,0)[l]{\circle*{ 0.288}}} at 29.720 20.256
%
%
\linethickness= 0.500pt
\setplotsymbol ({\thinlinefont .})
\put{\makebox(0,0)[l]{\circle*{ 0.288}}} at 31.593 21.273
%
%
\linethickness= 0.500pt
\setplotsymbol ({\thinlinefont .})
\put{\makebox(0,0)[l]{\circle*{ 0.288}}} at 33.333 22.352
%
%
\linethickness= 0.500pt
\setplotsymbol ({\thinlinefont .})
\put{\makebox(0,0)[l]{\circle*{ 0.288}}} at 35.683 20.923
%
%
\linethickness= 0.500pt
\setplotsymbol ({\thinlinefont .})
\put{\makebox(0,0)[l]{\circle*{ 0.288}}} at 37.715 20.384
%
%
\linethickness= 0.500pt
\setplotsymbol ({\thinlinefont .})
\put{\makebox(0,0)[l]{\circle*{ 0.288}}} at 38.286 19.748
%
%
\linethickness= 0.500pt
\setplotsymbol ({\thinlinefont .})
\put{\makebox(0,0)[l]{\circle*{ 0.449}}} at 14.427 23.226
%
%
\linethickness= 0.500pt
\setplotsymbol ({\thinlinefont .})
\put{\makebox(0,0)[l]{\circle*{ 0.449}}} at 18.633 22.399
\linethickness= 0.500pt
\setplotsymbol ({\thinlinefont .})
%
%
%
\plot   10.183 18.400 10.842 19.860
        10.883 19.950
        10.926 20.037
        10.971 20.122
        11.016 20.205
        11.063 20.286
        11.111 20.364
        11.160 20.439
        11.210 20.513
        11.262 20.583
        11.315 20.652
        11.369 20.718
        11.425 20.781
        11.539 20.901
        11.659 21.012
        11.783 21.113
        11.848 21.160
        11.913 21.204
        11.979 21.246
        12.047 21.286
        12.116 21.323
        12.187 21.358
        12.258 21.390
        12.331 21.420
        12.405 21.448
        12.481 21.473
        12.557 21.496
        12.635 21.516
        12.714 21.534
        12.794 21.550
        12.875 21.565
        12.954 21.581
        13.032 21.598
        13.109 21.617
        13.185 21.637
        13.259 21.658
        13.333 21.680
        13.406 21.704
        13.477 21.729
        13.547 21.754
        13.616 21.782
        13.684 21.810
        13.751 21.839
        13.817 21.870
        13.882 21.902
        13.945 21.935
        14.069 22.005
        14.188 22.080
        14.303 22.160
        14.414 22.244
        14.520 22.333
        14.621 22.428
        14.718 22.527
        14.810 22.630
         /
\plot 14.810 22.630 15.532 23.480 /
\linethickness= 0.500pt
\setplotsymbol ({\thinlinefont .})
%
%
\plot 10.353 19.911     10.258 19.937
        10.211 19.971
        10.176 20.045
        10.160 20.135
        10.160 20.226
        10.173 20.303
        10.199 20.370
        10.241 20.442
        10.294 20.515
        10.354 20.589
        10.415 20.660
        10.475 20.726
        10.569 20.836
        10.632 20.917
        10.711 21.013
        10.801 21.120
        10.896 21.234
        10.989 21.349
        11.075 21.461
        11.148 21.564
        11.201 21.654
        11.237 21.768
        11.184 21.867
        11.275 21.793
        11.375 21.795
        11.465 21.815
        11.573 21.844
        11.692 21.882
        11.816 21.923
        11.940 21.967
        12.056 22.009
        12.160 22.048
        12.245 22.081
        12.323 22.119
        12.415 22.172
        12.516 22.235
        12.623 22.300
        12.734 22.360
        12.845 22.410
        12.953 22.442
        13.056 22.449
        13.164 22.405
        13.248 22.327
        13.277 22.261
        13.299 22.136
        /
\linethickness= 0.500pt
\setplotsymbol ({\thinlinefont .})
%
%
\plot 20.326 21.057     20.304 21.142
        20.309 21.194
        20.426 21.320
        20.504 21.365
        20.576 21.393
        20.647 21.404
        20.730 21.404
        20.820 21.396
        20.913 21.382
        21.005 21.365
        21.092 21.348
        21.170 21.332
        21.234 21.321
        21.335 21.304
        21.458 21.279
        21.525 21.264
        21.595 21.249
        21.667 21.233
        21.740 21.218
        21.814 21.203
        21.886 21.189
        21.957 21.176
        22.025 21.165
        22.151 21.149
        22.257 21.145
        22.371 21.173
        22.432 21.279
        22.418 21.162
        22.462 21.071
        22.521 21.000
        22.598 20.919
        22.688 20.832
        22.784 20.743
        22.883 20.656
        22.977 20.575
        23.061 20.503
        23.131 20.445
        23.202 20.394
        23.293 20.340
        23.397 20.283
        23.506 20.220
        23.612 20.153
        23.710 20.080
        23.790 20.001
        23.846 19.916
        23.861 19.800
        23.834 19.687
        23.789 19.634
        23.690 19.562
        /
\linethickness= 0.500pt
\setplotsymbol ({\thinlinefont .})
%
%
%
\plot   28.567 18.074 29.225 19.534
        29.267 19.624
        29.310 19.711
        29.354 19.797
        29.399 19.879
        29.446 19.960
        29.494 20.038
        29.543 20.113
        29.594 20.187
        29.645 20.257
        29.698 20.326
        29.753 20.392
        29.808 20.455
        29.923 20.575
        30.042 20.686
        30.167 20.787
        30.231 20.834
        30.296 20.878
        30.363 20.920
        30.431 20.960
        30.500 20.997
        30.570 21.032
        30.642 21.064
        30.714 21.094
        30.788 21.122
        30.864 21.147
        30.940 21.170
        31.018 21.190
        31.097 21.208
        31.177 21.224
        31.258 21.239
        31.337 21.255
        31.415 21.273
        31.492 21.291
        31.568 21.311
        31.643 21.332
        31.716 21.354
        31.789 21.378
        31.860 21.403
        31.930 21.428
        32.000 21.456
        32.068 21.484
        32.135 21.513
        32.200 21.544
        32.265 21.576
        32.329 21.609
        32.452 21.679
        32.572 21.754
        32.687 21.834
        32.797 21.918
        32.903 22.007
        33.004 22.102
        33.101 22.201
        33.194 22.304
         /
\plot 33.194 22.304 33.915 23.154 /
\linethickness= 0.500pt
\setplotsymbol ({\thinlinefont .})
%
%
\plot 35.410 21.491     35.387 21.576
        35.393 21.628
        35.509 21.754
        35.587 21.799
        35.659 21.827
        35.730 21.838
        35.813 21.838
        35.903 21.830
        35.996 21.816
        36.089 21.799
        36.176 21.782
        36.254 21.766
        36.318 21.755
        36.419 21.738
        36.541 21.713
        36.609 21.698
        36.679 21.683
        36.751 21.667
        36.824 21.652
        36.897 21.637
        36.969 21.623
        37.040 21.610
        37.109 21.599
        37.235 21.583
        37.340 21.579
        37.454 21.607
        37.516 21.713
        37.501 21.596
        37.545 21.505
        37.605 21.434
        37.682 21.353
        37.771 21.266
        37.868 21.177
        37.966 21.090
        38.060 21.009
        38.145 20.937
        38.214 20.879
        38.285 20.828
        38.376 20.774
        38.480 20.716
        38.589 20.654
        38.696 20.587
        38.793 20.514
        38.873 20.435
        38.930 20.350
        38.944 20.233
        38.917 20.121
        38.873 20.068
        38.773 19.996
        /
\linethickness= 0.500pt
\setplotsymbol ({\thinlinefont .})
%
%
\plot 28.643 19.622     28.548 19.647
        28.501 19.681
        28.466 19.755
        28.450 19.845
        28.450 19.936
        28.463 20.013
        28.489 20.080
        28.531 20.152
        28.584 20.225
        28.644 20.299
        28.705 20.370
        28.765 20.436
        28.859 20.546
        28.922 20.627
        29.002 20.723
        29.092 20.830
        29.186 20.944
        29.279 21.059
        29.365 21.171
        29.438 21.274
        29.492 21.364
        29.527 21.478
        29.475 21.577
        29.566 21.503
        29.665 21.505
        29.755 21.525
        29.863 21.555
        29.982 21.592
        30.106 21.633
        30.230 21.677
        30.347 21.719
        30.450 21.758
        30.535 21.791
        30.614 21.829
        30.705 21.882
        30.806 21.945
        30.913 22.010
        31.024 22.070
        31.135 22.120
        31.244 22.152
        31.346 22.159
        31.454 22.115
        31.538 22.037
        31.567 21.971
        31.589 21.846
        /
\linethickness= 0.500pt
\setplotsymbol ({\thinlinefont .})
%
%
%
\plot   32.300 23.465 33.467 22.178
        33.540 22.099
        33.613 22.023
        33.688 21.949
        33.762 21.877
        33.837 21.808
        33.913 21.742
        33.989 21.678
        34.066 21.616
        34.143 21.558
        34.220 21.501
        34.299 21.448
        34.377 21.397
        34.456 21.348
        34.536 21.302
        34.616 21.258
        34.697 21.217
        34.778 21.179
        34.859 21.143
        34.941 21.110
        35.024 21.079
        35.107 21.051
        35.191 21.025
        35.275 21.002
        35.360 20.981
        35.445 20.963
        35.530 20.948
        35.616 20.935
        35.703 20.925
        35.790 20.917
        35.878 20.911
        35.966 20.909
        36.054 20.908
        36.142 20.909
        36.229 20.907
        36.314 20.903
        36.397 20.897
        36.479 20.889
        36.559 20.880
        36.638 20.868
        36.715 20.855
        36.791 20.839
        36.865 20.822
        36.937 20.803
        37.008 20.781
        37.078 20.758
        37.146 20.733
        37.212 20.706
        37.277 20.677
        37.402 20.613
        37.520 20.541
        37.633 20.462
        37.739 20.375
        37.839 20.280
        37.933 20.177
        38.021 20.067
        38.102 19.949
         /
\plot 38.102 19.949 38.729 18.972 /
\linethickness= 0.500pt
\setplotsymbol ({\thinlinefont .})
%
%
%
\plot   17.168 22.972 18.335 21.685
        18.408 21.606
        18.481 21.529
        18.555 21.455
        18.630 21.384
        18.705 21.315
        18.781 21.248
        18.857 21.185
        18.934 21.123
        19.011 21.064
        19.088 21.008
        19.167 20.955
        19.245 20.903
        19.324 20.855
        19.404 20.809
        19.484 20.765
        19.565 20.724
        19.646 20.686
        19.727 20.650
        19.809 20.617
        19.892 20.586
        19.975 20.558
        20.059 20.532
        20.143 20.509
        20.227 20.488
        20.313 20.470
        20.398 20.455
        20.484 20.442
        20.571 20.431
        20.658 20.424
        20.746 20.418
        20.834 20.415
        20.922 20.415
        21.010 20.415
        21.097 20.413
        21.182 20.410
        21.265 20.404
        21.347 20.396
        21.427 20.387
        21.506 20.375
        21.583 20.362
        21.659 20.346
        21.733 20.329
        21.805 20.309
        21.876 20.288
        21.946 20.265
        22.014 20.240
        22.080 20.213
        22.145 20.184
        22.270 20.120
        22.388 20.048
        22.501 19.969
        22.607 19.882
        22.707 19.787
        22.801 19.684
        22.889 19.574
        22.970 19.455
         /
\plot 22.970 19.455 23.597 18.479 /
%
%
\put{\SetFigFont{9}{10.8}{it}A} [lB] at 10.803 22.043
%
%
\put{\SetFigFont{9}{10.8}{it}B} [lB] at 22.312 21.630
%
%
\put{\SetFigFont{9}{10.8}{it}A} [lB] at 28.884 21.725
%
%
\put{\SetFigFont{9}{10.8}{it}B} [lB] at 37.615 21.963
\linethickness=0pt
\putrectangle corners at 10.139 23.497 and 38.964 18.057
\endpicture}

\end{center}
\vspace{-0pt}

\noindent Finally, the fiber product in the definition (\ref{bddeff}) 
corresponds to the condition that the maps
must take the same value in $X$ on the marked point $\bullet$
in order to be glued. This fiber product is defined
via evaluation maps discussed in the next section.

For $i,j,k,l$ distinct in $[n]$, set
$$D(i,j \mid k,l) = \sum D(A,B; \beta_1, \beta_2).$$
The sum is over all partitions $A \cup B =[n]$ with
 $i,j \in A$ and $\ k,l \in B$ and 
over all classes $\beta_1, \beta_2 \in A_1X$
satisfying $\beta_1+\beta_2=\beta$.
Using the projection $\M_{0,n}(X, \beta) \rarr \M_{0,\{i,j,k,l\}}
\eqq \proj^1$, the fundamental linear equivalences
\begin{equation}
\label{mlinee}
D(i,j\mid k,l)= D(i,k\mid j,l)= D(i,l\mid j,k)
\end{equation}
on $\M_{0,n}(X, \beta)$
are obtained via pull-back of the
the 4-point 
linear equivalences  on $\M_{0, \{i,j,k,l\}}$ as in 
(\ref{divisor}).

\subsection{Gromov-Witten invariants and quantum cohomology}
Let $X$ be a convex variety.
For each marked point $1\leq i\leq n$, there is
a canonical {\em evaluation map} $$\rho_i:\M_{0,n}(X,\beta)\rarr X$$
defined  for  
$[C, p_1, \ldots, p_n, \mu]$ in  $\M_{0,n}(X, \beta)$ by:
$$\rho_i([C, p_1, \ldots, p_n, \mu]) = \mu(p_i).$$
Given classes $\gamma_1, \ldots, \gamma_n$ in $A^*X$,
a product is determined in the ring 
$A^*(\M_{0,n}(X,\beta))$ by:
\begin{equation}
\label{iprod}
\rho_1^*(\gamma_1) \smallcup \cdots \smallcup \rho_n^*(\gamma_n) \in
A^*(\M_{0,n}(X,\beta)).
\end{equation}
If $\sum \text{codim}(\gamma_i)
= \text{dim} (\M_{0,n}(X, \beta))$, the
product (\ref{iprod}) can be evaluated on the
fundamental class of $\M_{0,n}(X, \beta)$. In this case,
the {\em Gromov-Witten invariant} $I_\beta(\gamma_1
\cdots  \gamma_n)$ is defined by:
\begin{equation}
\label{gwitt}
I_\beta(\gamma_1 
\cdots  \gamma_n) = \int_{\M_{0,n}(X, \beta)}
\rho_1^*(\gamma_1) \smallcup \cdots \smallcup \rho_n^*(\gamma_n) .
\end{equation}
The multiplicative notation in the argument
of $I_{\beta}$ is used to indicate
$I_{\beta}$ is a symmetric function of the classes
$\gamma_1, \ldots, \gamma_n$.

Let $X$ be a homogeneous space. 
Poincar\'e duality and Bertini-type transversality arguments
imply a relationship between the Gromov-Witten invariants
and enumerative geometry.
If $\gamma_i=[V_i]$ for a subvariety $V_i \subset X$, the
Gromov-Witten invariant (\ref{gwitt}) equals the number of
marked rational curves in $X$ with $i^{th}$ marked point in 
$V_i$,
suitably counted. For example, when $X=\proj^2$,
$\beta= d[\text{line}]$, $n=3d-1$, and each $V_i$ is a point,
$$N_d= I_{d}(
\underbrace{[p] \cdots  [p]}_{3d-1}).$$
The Gromov-Witten invariants are used to define 
the quantum cohomology ring. Associativity of this ring
is established as a consequence of the 4-point linear equivalences
(\ref{mlinee}). Associativity amounts to many equations
among the Gromov-Witten invariants which often lead
to a determination of all the invariants in terms
of a few basic numbers.

Given $\gamma_1, \ldots, \gamma_n \in H^*X$
(not necessarily of even degrees),
there are more general 
Gromov-Witten invariants in $H^* \M_{0,n}$ defined
by
$$I^{X}_{0,n,\beta}(\gamma_1 \otimes
\cdots  \otimes \gamma_n) = \eta_*(
\rho_1^*(\gamma_1) \smallcup \cdots \smallcup \rho_n^*(\gamma_n))$$
where $\eta: \M_{0,n}(X, \beta) \rarr \M_{0,n}$ is the projection.
The set of multilinear maps 
$$\{ I^{X}_{0,n, \beta}: (H^*X)^{\otimes n} \rarr H^*\M_{0,n}  \}$$ is
called the {\em Tree-Level System of Gromov-Witten
Invariants}. We will not need these generalities here.

The construction and proofs of the basic properties of
$\M_{0,n}(X,\beta)$ are undertaken in sections 1--6.
The theory of Gromov-Witten invariants and
quantum cohomology for homogeneous varieties
is presented in sections 7--10
with the examples of $\proj^2$, $\proj^3$,
and a smooth 
quadric 3-fold ${\mathbf Q}^3$ worked out in detail.
If Theorems 1--3 are taken for granted,
sections 1--6 can be skipped. No originality is
claimed for these notes except for some aspects of the
proofs of Theorems 1--4.
Constructions of Kontsevich's moduli space of stable
maps can also be found in [A], [K], [B-M]. In [A],
a generalization to the case in which the domain is a
surface is analyzed.

\subsection{Calculation of $N_d$}
We end this introduction by sketching how these
moduli spaces of maps can be used to calculate the number
$N_d$ of degree $d$ rational plane curves passing through
$3d-1$ general points in $\proj^2$.
The formula (\ref{rrecc}) will be recovered
in section 9 from the general quantum cohomology
results, but it may be useful now to see a direct proof. 

For $d=1$, $N_1=1$ is the number of lines through 2
points. $N_d$ is determined for $d\geq 2$ by the
{\em recursion formula}:
\begin{equation}
\label{rrecc}
N_d= \sum_{d_1+d_2=d, \ d_1, d_2 > 0}
N_{d_1} N_{d_2} \bigg( d_1^2 d_2^2 \binom{3d-4}{3d_1-2} -
d_1^3 d_2 \binom{3d-4}{3d_1-1} \bigg).
\end{equation}
For example, (\ref{rrecc}) yields$^1$ :
\footnotetext[1]{The number $N_3=12$ is the classically known
number of
singular members in a pencil of cubic curves through
8  given points. The number $N_4=620$ was computed by 
H. Zeuthen in [Z].
Z. Ran reports $N_4=620$ as
well as the higher $N_d$ 's can be 
derived from his formulas in [R1]; see [R2]
for a comparison of the two approaches. Some $N_d$'s are
also computed in [C-H 2].}
$$N_2=1,\ N_3=12, \ N_4=620, \ N_5=87304, \ N_6=26312976, ...$$

The strategy of proof is to utilize the  
fundamental linear relations (\ref{mlinee}) among
boundary components of $\smap_{0,n}(\proj^2,d)$.
Intersection of a curve $Y$ in this moduli space with 
the linear equivalence (\ref{mlinee}) will yield (\ref{rrecc}).
We will take $n=3d$ (not $3d-1$) with $d\geq 2$, so  
$n\geq 6$. Label the marked points by the
set $$\{1,2,3, \ldots, n-4, q, r, s, t\}.$$
The forgetful morphism
$\smap_{0,n}(\proj^2,d) \rarr \barr{M}_{0,\{q,r,s,t\}}$
yields the relations (\ref{mlinee}) on $\M_{0,n}(\proj^2,d)$:
\begin{equation}
\label{mlinel}
D(q,r\mid s,t) = D(q,s\mid r,t).
\end{equation}
Recall from section \ref{kontdef}:
$$D(q,r \mid s,t) = \sum_{q,r \in A, \ s,t \in B, \ d_1+d_2=d} 
D(A,B; d_1,d_2).$$
The curve  $Y \subset \M_{0,n}(\proj^2,d)$ is determined by
a selection of general points and lines in $\proj^2$.
More precisely, let $z_1, \ldots, z_{n-4}, z_s, z_t$ be
$n-2$ general points in $\proj^2$ and let
$l_{q}, l_r$ be general lines. Let the curve $Y$ be
defined by the intersection:
$$Y= \rho_1^{-1}(z_1) \cap \cdots \cap \rho_{n-4}^{-1}(z_{n-4})
     \cap \rho^{-1}_q(l_q) \cap \rho^{-1}_r(l_r) \cap
     \rho^{-1}_s(z_s)  \cap \rho^{-1}_t(z_t).$$
$\M_{0,n}(\proj^2,d)$ is a nonsingular, fine moduli space
on the open set of automorphism-free maps (see section 1.2).
It is not difficult to show the locus of maps
with non-trivial automorphisms in $\M_{0,n}(\proj^2,d)$
is of codimension at least 2 if $(n,d) \neq (0,2)$. Therefore,
by Bertini's theorem applied to each evaluation map
and the generality of the points and lines, we conclude
$Y$ is a nonsingular curve contained
in the automorphism-free locus which intersects all the
boundary divisors transversally at general points
of the boundary. It remains only
to compute the intersection of $Y$ with each
side of the linear equivalence (\ref{mlinel}).

The points of $$Y \cap\  D(A, B; d_1, d_2)$$
correspond
bijectively to maps $\mu:C=C_A\cup C_B \rarr \proj^2$ satisfying:
\begin{enumerate}
\item[(a)] $C_A, C_B \eqq \proj^1$ and meet transversally
           at a point.
\item[(b)] The markings of $A$, $B$ lie on $C_A$, $C_B$ respectively.
\item[(c)] $\mu_*([C_A])= d_1[\text{line}]$, 
$\mu_*([C_B])= d_2[\text{line}].$ 
\item [(d)] $\forall 1\leq i \leq n-4, \ \mu(i)=z_i.$
\item [(e)] $\mu(q)\in l_q$, $\mu(r)\in l_r$, $\mu(s)=z_{s}$, $\mu(t)=z_{t}$.
\end{enumerate}
Let $q, r \in A$ and $s, t \in B$.
$Y\cap\ D(A, B;0,d)$ is nonempty only when $A=\{q,r\}$.
In this case, $C_A$ is required to map to the point
$l_q \cap l_r$. The restriction $\mu:C_B\rarr \proj^2$ must
map the $3d-2$ markings on $C_B$ to the $3d-2$ given
points, and in addition, $\mu$ maps the point $C_A \cap C_B$
to $l_q \cap l_r$.
Therefore,
$$\# \  Y\cap\ D (\{q,r\},\{1,\ldots, n-4,s,t\};0,d) = N_d.$$
For $1\leq d_1 \leq d-1$, $Y\cap D(A,B; d_1, d_2)$ is
nonempty only when $|A|= 3d_1+1$. There are $\binom{3d-4}{3d_1-1}$
partitions satisfying  $q, r\in A$, $s,t\in B$, and $|A|=3d_1+1$.
A simple count of maps satisfying  (a)-(e)
yields
$$\# \ Y\cap \ D(A, B; d_1,d_2) = N_{d_1}N_{d_2} d_1^3 d_2$$
for each partition.
There are $N_{d_1}$ choices for the image of  $C_A$ and
$N_{d_2}$ choices for the image of $C_B$. The points labeled
$q$ and $r$ map to any of the  $d_1$  intersection
points of $\mu(C_A)$ with $l_q$ and $l_r$ respectively. Finally,
there are 
$d_1d_2$ choices for the image of the intersection point
$C_A \cap C_B$ corresponding to the intersection
points of $\mu(C_A) \cap \mu(C_B)\subset \proj^2$. 
The last case is simple:
$Y\cap\ D(A, B;d,0)=\emptyset$. Therefore,
$$\#\  Y\cap \ D(q,r \mid s,t)= N_d + 
\sum_{d_1+d_2=d, \ d_1>0, \ d_2 >0}
N_{d_1} N_{d_2}
d_1^3 d_2 \binom{3d-4}{3d_1-1}.$$

Now consider the other
side of the linear equivalence (\ref{mlinel}).
Let the markings now satisfy  $q,s \in A$ and $r,t \in B$.
$Y\cap \ D(A, B; 0,d)$ and $Y\cap\ (A , B; d,0)$
are both empty.
For $1\leq d_1 \leq d-1$, $Y\cap\  (A\cup B, d_1, d_2)$ is
nonempty only when $|A|=3d_1$. There are $\binom{3d-4}
{3d_1-2}$ such partitions and
$$ \# \ Y\cap\ D(A, B; d_1,d_2)= N_{d_1}N_{d_2} d_1^2 d_2^2$$
for each. Therefore,
$$ \#\ Y\cap\  D(q,s \mid r,t)=  
\sum_{d_1+d_2=d, \ d_1>0, \ d_2 >0}
N_{d_1} N_{d_2}
d_1^2 d_2^2 \binom{3d-4}{3d_1-2}.$$
The linear equivalence (\ref{mlinel}) implies
$$\# \ Y\cap\ D(q,r \mid s,t) = 
 \# \ Y\cap\ D(q,s \mid r,t).$$
The recursion (\ref{rrecc}) follows immediately.

In the general development of quantum cohomology
described in sections 8 and 9,
these numerical relations obtained by intersection
with the basic linear equivalences arise 
as ring associativity relations.

%
%
%

\section{{\bf Stable maps and their moduli spaces}}
\subsection{Definitions}
An $n$-pointed, genus $g$, complex, {\em quasi-stable} curve 
$$(C,\  p_1, \ldots, p_n)$$  is a projective,
connected, reduced,
(at worst) nodal curve of arithmetic genus $g$ with $n$ distinct, 
nonsingular, marked points. Let $S$ be an algebraic scheme over $\com$.
A {\em family} of $n$-pointed, genus
$g$, quasi-stable curves over $S$ is a flat, 
projective map $\pi:\mathcal{C}\rarr S$
with $n$ sections $p_1, \ldots, p_n$ such that each geometric
fiber $(\mathcal{C}_s, \ p_1(s), \ldots, p_n(s))$ is an $n$-pointed,
genus $g$, quasi-stable curve. 
Let $X$ be an algebraic scheme over $\com$. 
A {\em family of maps} over $S$ from $n$-pointed, genus $g$ curves
to $X$ consists of the data $(\pi:\mathcal{C}\rarr S,  
\{ p_i \}_{1\leq i \leq n}\ , \mu:\mathcal{C} \rarr X)$:
\begin{enumerate}
\item[(i)] A family of $n$-pointed, genus $g$,
quasi-stable curves $\pi: \mathcal{C}\rarr S$ 
with $n$ sections $\{p_1, \ldots, p_n \}$. 
\item[(ii)] A morphism $\mu: \mathcal{C} \rarr X$.
\end{enumerate}
Two families of maps over $S$, 
$$(\pi:\mathcal{C}\rarr S,  \{p_i\}, \mu), \ \ 
(\pi':\mathcal{C}' \rarr S,  \{p'_i\}, \mu'),$$ 
are {\em isomorphic}
if there exists a scheme isomorphism 
$\tau: \mathcal{C}\rarr \mathcal{C}'$ satisfying:
$\pi= \pi'\circ \tau$, $p'_i= \tau \circ p_i$, $\mu= \mu' \circ \tau$.
When $\pi:C \rarr \text{Spec}(\com)$ is the structure map,
 $(\pi:C\rarr \text{Spec}(\com), \{p_i\}, \mu)$ is written
as $(C,\{p_i\}, \mu)$.

Let $(C,\{p_i\}, \mu)$ be a
map from an $n$-pointed quasi-stable curve to $X$.
The {\em special points} of an irreducible
component $E\subset C$
are the marked points and the component intersections of $C$ 
that lie on $E$. The  map $(C, \{p_i\}, \mu)$ is
{\em stable} if  the  
following conditions hold for every component $E\subset C$:
\begin{enumerate}
\item[(1)]
If $E\eqq \proj^1$ and $E$ is  mapped to a point by  $\mu$, then  
$E$ must contain at least three special points.
\item[(2)]
If $E$ has arithmetic genus 1 and $E$ is mapped to
a point by $\mu$, then $E$ must contain at least one special point. 
\end{enumerate}
A family of pointed maps 
$(\pi:\mathcal{C}\rarr S, \{p_i\}, \mu)$
is {\em stable} if 
the pointed map on each geometric fiber of $\pi$ is 
stable.

If $X= \proj^r$, stability can be expressed
in the following manner. 
Let $\omega_{\mathcal{C}/S}$ denote the relative dualizing sheaf.
A flat family of maps 
$(\pi:\mathcal{C}\rarr S,\  \{p_i\}, \mu)$ is
{\em stable} if and only if 
 $\omega_{\mathcal{C}/S}(p_1+\ldots+p_n)\otimes
\mu^*(\oh_{\proj^r}(3))$ is $\pi$-relatively ample.

Let $X$ be an algebraic scheme over $\com$. Let $\beta\in A_1 X$.
A map $\mu:C \rarr X$ {\em represents} $\beta$ if
the $\mu$-push-forward of the fundamental class $[C]$ equals $\beta$.
Define a contravariant functor $\map_{g,n}(X,\beta)$ 
from the category of complex algebraic schemes
to sets as follows. Let
$\map_{g,n}(X,\beta)(S)$ be the set of isomorphism classes of
stable families over $S$ of maps from $n$-pointed,
genus $g$ curves to $X$ representing the class $\beta$.

\subsection{Existence} 
Let $X$ be a projective, algebraic scheme over $\com$.
Projective coarse moduli spaces of maps exist for
general $g$. In the genus $0$ case, 
if $X$ is a projective, nonsingular, convex variety, the 
coarse moduli spaces
are normal varieties with finite quotient singularities.
\begin{tm}
\label{rep}
There exists a projective, coarse moduli space $\smap_{g,n}(X,\beta)$.
\end{tm}
$\smap_{g,n}(X,\beta)$ is a scheme together with
a natural transformation of functors $$\phi:
\map_{g,n}(X,\beta) \rarr \mathcal{H}om_{Sch}(*, \smap_{g,n}(X,\beta))$$
satisfying properties:
\begin{enumerate}
\item[(I)] $\phi(\text{Spec}(\com)): \map \st 
(\text{Spec}(\com)) \rarr
\mathcal{H}om(\text{Spec}(\com), \smap \st)$ is a set bijection.
\item[(II)] If $Z$ is a scheme and $\psi: \map\st \rarr
\mathcal{H}om(*,Z)$ is a natural transformation of functors, then there
exists a unique morphism of schemes $$\gamma: \smap\st \rarr
Z$$ such that $\psi= \tilde{\gamma} \circ \phi$. ($\tilde{\gamma}:
\mathcal{H}om(*,\smap\st) \rarr \mathcal{H}om(*,Z)$ is the
natural transformation induced by $\gamma$.)
\end{enumerate}

Let $(C, \{p_i \}, \mu)$ be a 
map of an $n$-pointed, quasi-stable curve to $X$. 
An automorphism of the map is an
automorphism, $\tau$, of the curve $C$ satisfying
$$p_i=\tau(p_i), \ \ \mu=\mu\circ \tau.$$
It is straightforward to check that
$(C,\{p_i \}, \mu)$
is stable if and only if 
 $(C,\{p_i \}, \mu)$
has a finite automorphism group.
Let $\smap^*\st\subset \smap\st$ denote the
open locus of stable maps with no non-trivial automorphisms.

A nonsingular variety $X$ is convex if 
for every map $\mu:\proj^1\rarr X$, $H^1(\proj^1, \mu^*(T_X))=0$
(see section 0.4).
The second and third theorems concern the convex, genus $0$
case.
 
\begin{tm} Let $X$ be a projective, nonsingular, convex variety.
\begin{enumerate}
\item[(i)] $\smap\sto$ is a normal projective
variety of pure dimension $$\text{\em{dim}}
(X)+\int_{\beta} c_1(T_X) +n-3.$$ 
\item[(ii)] $\smap\sto$
is locally a quotient of a nonsingular variety by a finite group.
\item[(iii)] $\smap^*\sto$ is a nonsingular, fine
moduli space (for automorphism-free stable maps)
equipped with a universal family.
\end{enumerate}
\label{t2}
\end{tm}
\noindent 
In part (i), $\smap\sto$ is not claimed in general to be irreducible 
(or even nonempty).

In fact, if the language of stacks is pursued,
it can be seen that the moduli problem of stable maps
from $n$-pointed, genus $0$ curves to a nonsingular,
convex space $X$ determines a complete,
nonsingular, algebraic  stack. 
For simplicity, the stack theoretic
view is not taken in these notes; the experienced reader will see
how to make the required modifications.

The {\em boundary} of $\smap\sto$ is the
locus corresponding to reducible domain curves.
The boundary of the fine moduli space $\M_{0,n}$ is
a divisor with normal crossings. In the coarse
moduli spaces $\M_{g}$ and $\M_{g,n}$, the boundary
is a divisor with normal crossings modulo a finite
group. $\M_{0,n}(X,\beta)$ has the same boundary
singularity type as these moduli spaces of pointed curves.  
\begin{tm} 
\label{t33}
Let $X$ be a nonsingular, projective, convex variety.
The boundary of $\smap_{0,n}(X,\beta)$ is a divisor
with normal crossings (up to a finite group quotient).
\end{tm}

The organization of the construction is as follows. First
$\smap_{g,n}(\proj^r,d)$ is explicitly constructed in sections
2--4.
If $X\subset \proj^r$ is a closed subscheme, it is not difficult
to define a natural, closed 
subscheme $$\smap_{g,n}(X,d)\subset \smap_{g,n}(\proj^r,d)$$
of maps that factor through $X$. $\smap_{g,n}(X,d)$ is a disjoint
union of the spaces $\smap_{g,n}(X,\beta)$ as $\beta$ varies in
$A_1 X$. By the universal property, 
it can be seen that the coarse moduli
spaces $\smap_{g,n}(X,\beta)$ do not depend on the projective
embedding of $X$ (see section  5). 
The deformation arguments required to deduce Theorem 2
from the convexity assumption are covered in section 5.
The boundary
of the space of maps is discussed in section 6.

\subsection{Natural structures}
\label{introint}
The universal property of the moduli space of maps
immediately yields geometric structures on $\smap\st$.
Consider first the marked points.
The $n$ marked points induce
$n$ canonical evaluation maps $\rho_1, \ldots, \rho_n$
on $\smap\st$.
For $1\leq i \leq n$,
define a natural transformation
$$\theta_i: \map \st \rarr \mathcal{H}om(*, X)$$ as follows.
Let $\zeta=(\pi:\mathcal{C} \rarr S, \ \{p_i\},\ \mu)$ be an element of
$\map\st(S)$. Let
$$\theta_i(S)(\zeta)= \mu \circ p_i \in \mathcal{H}om(S,X).$$
$\theta_i$ is easily seen to be a natural transformation.
By Theorem \ref{rep}, $\theta_i$ induces a unique morphism
of schemes
$\rho_i: \smap\st \rarr X.$

By the universal properties of the moduli spaces $\barr{M}_{g,n}$ 
of
$n$-pointed Deligne-Mumford stable genus $g$ curves (in case $2g-2+n>0$), 
each 
element $\zeta\in \map\st(S)$ naturally yields a morphism
$S\rarr \barr{M}_{g,n}$ ([Kn]). 
Therefore, there exist natural forgetful maps
$\eta: \smap\st \rarr \barr{M}_{g,n}$.

%
%
%
%
%

\section{\bf{Boundedness and a quotient approach}}
\subsection{Summary} In this section, 
the case $X=\proj^r$ will  be considered.
The boundedness of the moduli problem of pointed stable
maps is established. The arguments
lead naturally to a quotient approach to
the coarse moduli space. To set up
the quotient approach, a result on equality
loci of families of line bundles is required.

\subsection{Equality of line bundles in families}
Results on scheme theoretic equality loci are 
recalled.
Let $\pi: \cal{C} \rarr S$ be a flat family of quasi-stable curves.
By the theorems of cohomology and base change (cf. 
[H]), there is
a canonical isomorphism 
$\oh_S \eqq \pi_*(\oh_{\cal{C}}) $.
Hence, for any line bundle $\cal{N}$ on $S$, there is a canonical
isomorphism $\cal{N} \eqq \pi_* \pi^*(\cal{N})$.
Suppose $\cal{L}$ and $\cal{M}$ are two line bundles on $\cal{C}$.
The existence of a line bundle $\cal{N}$ on $S$ such that
$\cal{L} \otimes \cal{M}^{-1}\eqq \pi^*(\cal{N})$ is 
equivalent to the joint validity of (a) and (b):
\begin{enumerate}
\item[(a)] $\pi_*(\cal{L} \otimes \cal{M}^{-1})$ is locally free.
\item[(b)] The canonical map 
$\pi^*\pi_*(\cal{L}\otimes \cal{M}^{-1}) \rarr \cal{L}\otimes \cal{M}^{-1}$
is an isomorphism.
\end{enumerate}
Let $\mathcal{L}_s$ be a line bundle on 
the geometric fiber $\cal{C}_s$ of $\pi$. The
{\em multidegree} of $\cal{L}_s$  assigns
to each irreducible component of $\cal{C}_s$ the
degree of the restriction of $\cal{L}_s$ to that component.
\begin{pr}
\label{mumford}
Let $\cal{L}$, $\cal{M}$ be line bundles on $\cal{C}$ such that
the multidegrees of $\cal{L}_s$ and $\cal{M}_s$ coincide on each geometric
fiber $\cal{C}_s$. Then,
 there is a unique closed subscheme
$T\rarr S$ satisfying the following two properties:
\begin{enumerate}
\item[(I)] There is a line bundle $\cal{N}$ on $T$ such that 
$\cal{L}_T \otimes \cal{M}_T^{-1} \eqq \pi^*(\cal{N})$.
\item[(II)] If $(R \rarr S,\  \cal{N})$ is a pair of a morphism from
$R$ to $S$ and a
line bundle on $R$  such that $\cal{L}_R \otimes \cal{M}^{-1}_R \eqq
\pi^*(\cal{N})$, then $R\rarr S$ factors through $T$.
\end{enumerate}
\end{pr}
\bpf
The proof of the Theorem of the Cube (II) in [M1] also establishes
this proposition. The multidegree condition implies
$\cal{L}_s \eqq \cal{M}_s$ if and only if $h^0(\cal{C}_s, \cal{L}_s
\otimes \cal{M}_s^{-1})=1$. The multidegree condition is required
for $T$ to be a {\em closed} subscheme.
\epf

\subsection{Boundedness}
\label{bound}
Let $(C, \{p_i\}, \mu)$ be a 
stable map from an $n$-pointed, genus $g$ curve to $\proj^r$.
Let $$\cal{L}= \omega_{C}(p_1+\ldots +p_n) \otimes \mu^*(\oh
_{\proj^r}(3)).$$
$\cal{L}$ is ample on $C$. A simple argument shows there exists an
$f=f(g,n,r,d)>0$ such that
$\cal{L}^f$ is very ample on $C$ and $h^1(C, \cal{L}^f)=0$, so  
$$\text{degree}(\cal{L}^f)= f\cdot(2g-2+n+3d) =e,$$
$$h^0(C, \cal{L}^f)= e-g+1.$$
Let $W\eqq \com^{e-g+1}$ be a vector space.
An isomorphism
\begin{equation}
\label{choice} 
W^* \stackrel {\sim} {\rarr} H^0(C, \cal{L}^f)
\end{equation}
induces embeddings $\iota:C \hookrightarrow \proj(W)$ and
$\gamma: C \hookrightarrow \proj(W) \times \proj^r$ 
where $\gamma=(\iota, \mu)$.
The $n$ sections $\{p_i\}$ yield $n$ points 
$(\iota \circ p_i, \mu \circ p_i)$ of $\proj(W)\times \proj^r$.
Let $H$ be the Hilbert scheme of genus $g$ curves
in $\proj(W)\times \proj^r$ of multidegree $(e,d)$. 
Let $P_i=\proj(W)\times \proj^r$
be the Hilbert scheme of a point in $\proj(W)\times \proj^r$.
Via the isomorphism (\ref{choice}), a point in
$H\times P_1\times \ldots \times P_n$ is associated to the 
stable map $(C, \{p_i\}, \mu)$.

The locus of points in $H\times P_1\times \ldots \times P_n$ corresponding
to stable maps has a natural quasi-projective scheme
structure. 
There is a natural closed incidence subscheme
$$ I \subset H \times P_1 \times P_2 \times \ldots \times P_n$$
corresponding to the locus where the $n$ points lie on the
curve. There is an open set $U\subset I$ satisfying the following:
\begin{enumerate}
\item[(i)] The curve $C$ is quasi-stable.
\item[(ii)] The natural projection $C \rarr \proj(W)$ is a
non-degenerate embedding.
\item[(iii)] The $n$ points lie in the nonsingular locus of $C$.
\item[(iv)] 
The multidegree of $\oh_{\proj(W)}(1) \otimes \oh_{\proj^r}(1)|_{C}$
equals the multidegree of $$\omega^f_C(fp_1+fp_2+
\ldots + fp_n)\otimes
\oh_{\proj^r}(3f+1)|_{C}.$$
\end{enumerate}
By Proposition \ref{mumford}, there exists a natural closed
subscheme $J\subset U$ 
where the line bundles of condition (iv) above coincide. $J$
corresponds to the locus of stable maps. 
The natural $PGL(W)$-action on $\proj(W)\times \proj^r$
yields $PGL(W)$-actions on $H$,
$P_i$, $I$, $U$, and $J$.
To each stable map from an $n$-pointed, genus $g$
curve to $\proj^r$, we have associated a
$PGL(W)$-orbit in $J$.
If two stable maps are associated to the same orbit,
the two stable maps are isomorphic. The stability condition
implies that a stable map has no 
infinitesimal automorphisms. It follows
that the $PGL(W)$-action on $J$ has finite stabilizers.

\subsection{Quotients}
The moduli space of stable maps is $J/PGL(W)$. 
It may be possible to construct the quotient $J/PGL(W)$ via
Geometric Invariant Theory. Another method will be pursued here.
The quotient will be first constructed as a proper, algebraic
variety by using auxiliary moduli spaces of pointed curves.
Projectivity will then be established via J. Koll\'ar's semipositivity
approach.

%
%
%
%
%
%

\section{{\bf A rigidification of} $\map\stp$}
\subsection{Review of Cartier divisors}
An effective Cartier divisor $D$ on a scheme $Y$ is
a closed subscheme that is locally defined by a
non-zero-divisor. An effective Cartier divisor
determines a line bundle $\cal{L}=\oh(D)$
together with a section $s\in H^0(Y, \cal{L})$
locally not a zero-divisor such that
$D$ is the subscheme defined by $s=0$. (As an
invertible sheaf, $\oh(D)$ can be constructed as
the subsheaf of rational functions with at most
simple poles along $D$ with $s$ equal to the function
$1$, see [M2].)
Conversely, if the pair $(\cal{L}, s)$ satisfies:
\begin{enumerate}
\item[(i)] $\cal{L}$ is line bundle on $Y$.
\item[(ii)] $s\in H^0(Y,\cal{L})$ is a section locally not a zero divisor.
\end{enumerate} 
then the zero
scheme of $s$ is an effective Cartier divisor on $Y$.

\begin{lm}
\label{cart}
Let  the pairs $(\cal{L},s)$ and $(\cal{L}', s')$ satisfy (i) and (ii) above.
If the two pairs yield the same Cartier divisor, then there exists a unique
isomorphism $\cal{L} \rarr \cal{L}'$ taking $s$ to $s'$.
\end{lm}

\subsection{Definitions}
We assume throughout the construction that
$r>0$, $d>0$, and $(g,n,r,d)\neq (0,0,1,1)$.
If 
$r=0$, the functor of stable maps 
to $\proj^0$ is coarsely represented by $\barr{M}_{g,n}$.
If $d=0$, the functor $\map_{g,n}(\proj^r,0)$ 
is coarsely represented by 
$\barr{M}_{g,n}\times \proj^r$ and, $\smap_{0,0}(1,1)$ is easily
seen to be $\text{Spec}(\com)$. For all other values, the construction
of $\smap\st$
will be undertaken.

Let $\proj^r=\proj(V)$. Then, $V^* =H^0(\proj^r, \oh_{\proj^r}(1))$. 
Let $\barr{t}=(t_0, \ldots, t_r)$ span a basis of $V^*$.
A {\em $\barr{t}$-rigid stable} family of 
degree $d$ maps from  $n$-pointed, genus
$g$ curves to $\proj^r$ consists of the 
data 
$$(\pi:\cal{C} \rarr S,\  \{p_i\}_{1\leq i  \leq n}\ ,\ 
 \{q_{i,j}\}_{0\leq i \leq r,\ 1\leq j \leq d}\ ,\  \mu)$$ where:
\begin{enumerate}
\item[(i)] $(\pi:\cal{C} \rarr S,\  \{p_i\}, \mu)$ is a 
stable family of degree $d$ maps from $n$-pointed, genus $g$ curves
to $\proj^r$.
\item[(ii)] $(\pi:\cal{C} \rarr S, \ \{p_i \}, \{ q_{i,j} \})$ is
a flat, projective family of  $n+d(r+1)$-pointed, 
genus $g$,  Deligne-Mumford stable curves
with sections $\{p_i \}$ and $\{ q_{i,j} \}$.
\item[(iii)] For $0\leq i \leq r$, there is an equality of
Cartier divisors $$\mu^*(t_i)= q_{i,1} + q_{i,2} + \ldots + q_{i,d}.$$
\end{enumerate} 
Condition (iii) implies each fibered map of the family intersects each
hyperplane $(t_i)\subset \proj^r$ transversally. Condition (ii)
guarantees these hyperplane intersections are unmarked, nonsingular
points.
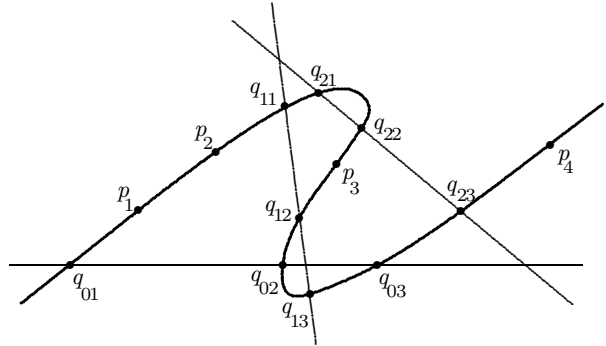
\begin{center}
\font\thinlinefont=cmr5
\mbox{\beginpicture
\setcoordinatesystem units < 0.500cm, 0.500cm>
\unitlength= 0.500cm
\linethickness=1pt
\setplotsymbol ({\makebox(0,0)[l]{\tencirc\symbol{'160}}})
\setshadesymbol ({\thinlinefont .})
\setlinear
%
%
\linethickness= 0.500pt
\setplotsymbol ({\thinlinefont .})
\putrule from  2.095 15.050 to 17.335 15.050
%
%
\put{\SetFigFont{7}{8.4}{rm}
3
} [lB] at 12.383 14.192
%
%
\put{\SetFigFont{9}{10.8}{it}
q
} [lB] at 11.874 14.541
%
%
\put{\SetFigFont{7}{8.4}{rm}
0
} [lB] at 12.129 14.192
%
%
\put{\SetFigFont{7}{8.4}{rm}
2
} [lB] at  8.985 14.287
%
%
\put{\SetFigFont{7}{8.4}{rm}
0
} [lB] at  8.731 14.287
%
%
\put{\SetFigFont{9}{10.8}{it}
q
} [lB] at  8.477 14.637
%
%
\put{\SetFigFont{7}{8.4}{rm}
1
} [lB] at  8.954 19.304
%
%
\put{\SetFigFont{7}{8.4}{rm}
1
} [lB] at  8.731 19.304
%
%
\put{\SetFigFont{9}{10.8}{it}
q
} [lB] at  8.446 19.622
%
%
\linethickness= 0.500pt
\setplotsymbol ({\thinlinefont .})
\put{\makebox(0,0)[l]{\circle*{ 0.191}}} at  9.430 19.272
%
\linethickness= 0.500pt
\setplotsymbol ({\thinlinefont .})
\put{\makebox(0,0)[l]{\circle*{ 0.191}}} at 11.462 18.701
%
\linethickness= 0.500pt
\setplotsymbol ({\thinlinefont .})
\put{\makebox(0,0)[l]{\circle*{ 0.191}}} at  9.811 16.320
%
\linethickness= 0.500pt
\setplotsymbol ({\thinlinefont .})
\put{\makebox(0,0)[l]{\circle*{ 0.191}}} at 14.097 16.478
%
\linethickness= 0.500pt
\setplotsymbol ({\thinlinefont .})
\put{\makebox(0,0)[l]{\circle*{ 0.191}}} at  9.366 15.050
%
\linethickness= 0.500pt
\setplotsymbol ({\thinlinefont .})
\put{\makebox(0,0)[l]{\circle*{ 0.191}}} at 10.097 14.287
%
\linethickness= 0.500pt
\setplotsymbol ({\thinlinefont .})
\put{\makebox(0,0)[l]{\circle*{ 0.191}}} at 11.874 15.050
%
\linethickness= 0.500pt
\setplotsymbol ({\thinlinefont .})
\put{\makebox(0,0)[l]{\circle*{ 0.191}}} at  3.715 15.050
%
\linethickness= 0.500pt
\setplotsymbol ({\thinlinefont .})
\put{\makebox(0,0)[l]{\circle*{ 0.191}}} at  7.588 18.066
%
\linethickness= 0.500pt
\setplotsymbol ({\thinlinefont .})
\put{\makebox(0,0)[l]{\circle*{ 0.191}}} at  5.524 16.510
%
\linethickness= 0.500pt
\setplotsymbol ({\thinlinefont .})
\put{\makebox(0,0)[l]{\circle*{ 0.191}}} at 16.478 18.256
%
\linethickness= 0.500pt
\setplotsymbol ({\thinlinefont .})
\put{\makebox(0,0)[l]{\circle*{ 0.191}}} at 10.795 17.748
%
\linethickness= 0.500pt
\setplotsymbol ({\thinlinefont .})
\put{\makebox(0,0)[l]{\circle*{ 0.191}}} at 10.319 19.622
%
\linethickness= 0.500pt
\setplotsymbol ({\thinlinefont .})
\plot  9.081 22.035 10.255 12.954 /
%
%
\linethickness= 0.500pt
\setplotsymbol ({\thinlinefont .})
\plot  8.033 21.558 17.081 14.002 /
\linethickness=1pt
\setplotsymbol ({\makebox(0,0)[l]{\tencirc\symbol{'160}}})
%
%
%
\plot	 2.413 14.034  5.271 16.304
 	 5.359 16.374
	 5.448 16.444
	 5.535 16.513
	 5.622 16.581
	 5.708 16.649
	 5.793 16.715
	 5.878 16.781
	 5.961 16.847
	 6.044 16.911
	 6.127 16.975
	 6.208 17.038
	 6.289 17.100
	 6.369 17.161
	 6.448 17.222
	 6.527 17.281
	 6.605 17.340
	 6.682 17.399
	 6.759 17.456
	 6.834 17.513
	 6.909 17.569
	 6.983 17.624
	 7.057 17.679
	 7.130 17.732
	 7.202 17.785
	 7.273 17.837
	 7.343 17.889
	 7.413 17.939
	 7.482 17.989
	 7.550 18.038
	 7.618 18.087
	 7.685 18.134
	 7.751 18.181
	 7.816 18.227
	 7.881 18.272
	 7.945 18.317
	 8.008 18.360
	 8.132 18.445
	 8.253 18.527
	 8.371 18.606
	 8.487 18.682
	 8.599 18.755
	 8.708 18.825
	 8.815 18.891
	 8.918 18.955
	 9.019 19.016
	 9.116 19.073
	 9.211 19.127
	 9.303 19.179
	 9.392 19.227
	 9.477 19.272
	 9.561 19.315
	 9.643 19.356
	 9.723 19.394
	 9.802 19.431
	 9.879 19.466
	 9.954 19.498
	10.028 19.529
	10.100 19.558
	10.171 19.585
	10.240 19.610
	10.308 19.632
	10.373 19.653
	10.500 19.689
	10.620 19.717
	10.734 19.737
	10.842 19.748
	10.943 19.752
	11.037 19.748
	11.125 19.737
	11.207 19.717
	11.282 19.689
	11.351 19.653
	11.468 19.571
	11.559 19.482
	11.623 19.386
	11.660 19.284
	11.671 19.176
	11.654 19.061
	11.611 18.940
	11.579 18.877
	11.541 18.812
	11.497 18.744
	11.449 18.671
	11.397 18.593
	11.341 18.509
	11.280 18.421
	11.215 18.327
	11.146 18.229
	11.073 18.125
	10.995 18.017
	10.914 17.903
	10.828 17.784
	10.737 17.660
	10.691 17.596
	10.643 17.531
	10.594 17.464
	10.544 17.397
	10.494 17.328
	10.442 17.258
	10.389 17.186
	10.335 17.113
	10.281 17.040
	10.229 16.967
	10.178 16.895
	10.128 16.824
	10.081 16.752
	10.035 16.682
	 9.990 16.611
	 9.947 16.542
	 9.905 16.473
	 9.865 16.404
	 9.826 16.336
	 9.789 16.268
	 9.754 16.201
	 9.720 16.134
	 9.687 16.068
	 9.656 16.002
	 9.626 15.937
	 9.598 15.872
	 9.572 15.808
	 9.547 15.744
	 9.502 15.618
	 9.462 15.494
	 9.429 15.372
	 9.402 15.252
	 9.381 15.134
	 9.366 15.018
	 9.357 14.907
	 9.354 14.804
	 9.356 14.710
	 9.364 14.624
	 9.378 14.546
	 9.398 14.477
	 9.454 14.363
	 9.532 14.282
	 9.634 14.235
	 9.759 14.221
	 9.829 14.226
	 9.906 14.240
	 9.988 14.259
	10.073 14.283
	10.164 14.309
	10.258 14.339
	10.357 14.373
	10.460 14.410
	10.568 14.450
	10.680 14.494
	10.796 14.541
	10.917 14.592
	11.042 14.646
	11.106 14.675
	11.171 14.704
	11.237 14.734
	11.305 14.765
	11.373 14.797
	11.443 14.830
	11.513 14.864
	11.585 14.898
	11.658 14.934
	11.732 14.970
	11.807 15.008
	11.885 15.049
	11.966 15.092
	12.049 15.137
	12.135 15.185
	12.223 15.236
	12.314 15.289
	12.407 15.345
	12.503 15.403
	12.601 15.464
	12.702 15.528
	12.805 15.594
	12.911 15.662
	13.020 15.733
	13.130 15.807
	13.244 15.883
	13.360 15.962
	13.478 16.043
	13.599 16.127
	13.722 16.213
	13.848 16.302
	13.912 16.347
	13.977 16.393
	14.042 16.440
	14.107 16.487
	14.174 16.535
	14.241 16.583
	14.309 16.633
	14.377 16.682
	14.446 16.733
	14.515 16.784
	14.585 16.836
	14.656 16.888
	14.728 16.941
	14.800 16.995
	14.872 17.049
	14.946 17.104
	15.020 17.159
	15.094 17.216
	15.169 17.272
	15.245 17.330
	15.322 17.388
	15.399 17.447
	 /
\plot 15.399 17.447 17.875 19.336 /
%
%
\put{\SetFigFont{7}{8.4}{rm}
2
} [lB] at 14.002 16.764
%
%
\put{\SetFigFont{7}{8.4}{rm}
3
} [lB] at 14.224 16.764
%
%
\put{\SetFigFont{7}{8.4}{rm}
1
} [lB] at  5.143 16.478
%
%
\put{\SetFigFont{7}{8.4}{rm}
2
} [lB] at  7.239 18.129
%
%
\put{\SetFigFont{7}{8.4}{rm}
1
} [lB] at 10.541 19.812
%
%
\put{\SetFigFont{7}{8.4}{rm}
2
} [lB] at 12.256 18.383
%
%
\put{\SetFigFont{9}{10.8}{it}
p
} [lB] at  4.953 16.859
%
%
\put{\SetFigFont{9}{10.8}{it}
p
} [lB] at  7.017 18.479
%
%
\put{\SetFigFont{9}{10.8}{it}
q
} [lB] at 10.065 20.098
%
%
\put{\SetFigFont{9}{10.8}{it}
q
} [lB] at 11.748 18.669
%
%
\put{\SetFigFont{9}{10.8}{it}
q
} [lB] at 13.684 17.018
%
%
\put{\SetFigFont{9}{10.8}{it}
p
} [lB] at 10.922 17.335
%
%
\put{\SetFigFont{7}{8.4}{rm}
3
} [lB] at 11.144 16.954
%
%
\put{\SetFigFont{9}{10.8}{it}
p
} [lB] at 16.669 17.844
%
%
\put{\SetFigFont{7}{8.4}{rm}
4
} [lB] at 16.828 17.494
%
%
\put{\SetFigFont{7}{8.4}{rm}
2
} [lB] at 10.319 19.812
%
%
\put{\SetFigFont{7}{8.4}{rm}
2
} [lB] at  9.335 16.224
%
%
\put{\SetFigFont{9}{10.8}{it}
q
} [lB] at  8.858 16.510
%
%
\put{\SetFigFont{7}{8.4}{rm}
1
} [lB] at  9.112 16.224
%
%
\put{\SetFigFont{9}{10.8}{it}
q
} [lB] at  9.303 13.780
%
%
\put{\SetFigFont{7}{8.4}{rm}
3
} [lB] at  9.779 13.399
%
%
\put{\SetFigFont{7}{8.4}{rm}
1
} [lB] at  9.557 13.399
%
%
\put{\SetFigFont{7}{8.4}{rm}
2
} [lB] at 12.033 18.383
%
%
\put{\SetFigFont{7}{8.4}{rm}
1
} [lB] at  4.223 14.129
%
%
\put{\SetFigFont{7}{8.4}{rm}
0
} [lB] at  4.000 14.129
%
%
\put{\SetFigFont{9}{10.8}{it}
q
} [lB] at  3.747 14.541
\linethickness=0pt
\putrectangle corners at  2.095 22.035 and 17.875 12.954
\endpicture}

\end{center}

If  $(g,n,r,d)=(0,0,1,1)$, then $n+d(r+1)=2$. There are
no Deligne-Mumford stable $2$-pointed genus $0$ curves. This is
why $(0,0,1,1)$ is avoided.

Define a contravariant functor $\map\stv$ from
the category of complex algebraic schemes  to sets as follows.
Let $\map\stv (S)$ be the set of isomorphism
classes of $\barr{t}$-rigid stable 
families over $S$ of degree $d$ maps from
$n$-pointed, genus $g$ curves to $\proj^r$.
Note that the functor $\map\stv$ depends
only upon the spanning hyperplanes $(t_i)\subset \proj^r$ and
not upon the additional $\com^*$-choices in the
defining equations $t_i$ of the hyperplanes.
Nevertheless, it is natural for the following
constructions to consider the equations of the
hyperplanes $\barr{t}=(t_0, \ldots, t_r)$.
\begin{pr}
\label{bartg}
There exists a quasi-projective coarse moduli space,
$$\smap\stv,$$ and a natural transformation
of functors 
$$\psi: \map\stv\rarr \cal{H}om(*,
\smap\stv)$$
satisfying the analogous conditions (I) and (II) of Theorem \ref{rep}.
\end{pr}
\noindent The genus 0 case is simpler.
\begin{pr}
\label{barto}
$\smap\stvo$ represents the functor $\map\stvo$
and is a  nonsingular algebraic variety.
\end{pr}

\subsection{Proofs}
A complete proof of Proposition \ref{barto} will be given.
The proof of Proposition \ref{bartg} is almost identical. Remarks
indicating the differences will be made. The dependence of the 
coarse and fine moduli property on the genus in Propositions \ref{bartg} and
\ref{barto} is a direct consequence of the fact that $\barr{M}_{g,n}$ is
a coarse moduli space for $g>0$ and a fine moduli space for $g=0$.

The idea behind the construction is the following. 
Let $m=n+d(r+1)$. The data
of the $\barr{t}$-rigid stable family immediately yields
a morphism of the base $S$ to $\barr{M}_{g,m}$. In fact,
the image of $S$ lies in a universal, locally closed subscheme of
$\barr{M}_{g,m}$. This subscheme is denoted by $B$. The first
step of the construction is to identify $B$.
The morphism $S\rarr B$ does not contain all the
data of the $\barr{t}$-rigid stable family. Consider the
case in which the base $S$ is a point. The
corresponding point in
$B$ records the domain curve $C$, the marked
points $\{p_i\}$, and the pull-back divisors under $\mu$ 
of the hyperplanes in
$\proj^r$ determined by $\barr{t}$. The map $\mu$ is determined
by the pull-back divisors
up to the diagonal torus action on $\proj^r$. 
The torus information
is recorded in the total space of $r$ tautological
$\com^*$-bundles over $B$. 
The
$\barr{t}$-rigid moduli space is expressed as the
total space of these $r$ distinct $\com^*$-bundles over $B$.    
To canonically construct the universal family over the
$\barr{t}$-rigid moduli space, the equations $t_i$
of the hyperplanes are needed. This is why the 
equations $t_i$ (rather than the spanning hyperplanes
$(t_i)$) are explicitly chosen. 
 
Proposition \ref{barto} is proved by an explicit
construction of $\smap\stvo$ together with
a universal family of $\barr{t}$-rigid stable maps.
Let $\mk$ be the Mumford-Knudsen compactification of
the moduli space of $m$-pointed, genus 0 curves. Let
$\pi: \mku \rarr \mk$ be the universal curve with $m$ sections
$\{p_i\}_{1\leq i \leq n}$ and $\{q_{i,j} \} _{0 \leq i \leq r, \ 
1 \leq j \leq d}$. Since $\mku$ is nonsingular and the sections
are of codimension $1$, there are canonically defined line bundles:
$$\HH_i = \oh_{\mku}(q_{i,1} + q_{i,2} + \ldots + q_{i,d}),$$
for $0 \leq i \leq r$.
Let $s_i\in H^0(\mku, \HH_i)$ be the canonical section
representing  the Cartier divisor
$(q_{i,1}+ q_{i,2}+ \ldots + q_{i,d})$. 

For any morphism $\gamma:X\rarr \mk$, consider the fiber product:
\begin{equation*}
\begin{CD}
X\times_{\mk} \mku @>{\barr{\gamma}}>> \mku \\
@VV{\pi_X}V @VV{\pi}V \\
X @>{\gamma}>> \mk \\
\end{CD}
\end{equation*}
We call the  morphism $\gamma:X \rarr \mk$  {\em $\HH$-balanced} if
\begin{enumerate}
\item[(a)] For $1\leq i \leq r$, 
$\pi_{X*} \barr{\gamma}^*(\HH_i \otimes \HH_0^{-1})$ is
locally free.
\item[(b)] For $1\leq i \leq r$, the canonical map
$$\pi_X^* \pi_{X*} \barr{\gamma}^*(\HH_i \otimes \HH_0^{-1})
\rarr \barr{\gamma}^*(\HH_i \otimes \HH_0^{-1})$$
is an isomorphism.
\end{enumerate}
If $\gamma$ is $\HH$-balanced, the
line bundles  
$\barr{\gamma}^{*}(\HH_i)$ are isomorphic on the fibers of $\pi_X$. 
Let $B\subset \mk$ be the universal, locally closed subscheme 
satisfying the two
following properties:
\begin{enumerate}
\item[(i)] The inclusion 
$\iota: B \hookrightarrow \mk$ is $\HH$-balanced.
\item[(ii)] Every $\HH$-balanced morphism $\gamma: X \rarr \mk$ factors
(uniquely) through $B$.
\end{enumerate}
By Proposition \ref{mumford},
$B$ exists. In fact, $B\subset \mk$ is
a Zariski open subscheme. In the $g>0$ case, the above constructions
exist over the stacks $\barr{M}_{g,m}$ and $\barr{U}_{g,m}$.
$B_{g,m}$ is a locally closed substack of $\barr{M}_{g,m}$ 
of positive
codimension.

Let $\GG_i=\pi_{B*} \barr{\iota}^*(\HH_i\otimes \HH_0^{-1})$
for $1\leq i \leq r$. Let $\tau_i: Y_i \rarr B$ be the
total space of the canonical $\com^*$-bundle associated to $\GG_i$.
$Y_i$ is the affine bundle associated to $\GG_i$ minus the zero section.
The pull-back $\tau_i^*(\GG_i)$ has a tautological non-vanishing section
and hence is canonically trivial.
Consider the product 
$$Y=Y_1 \times_B \times Y_2 \times_B \ldots \times_B Y_r$$
equipped with projections $\rho_i: Y\rarr Y_i$ and a morphism
$\tau: Y \rarr B$. Form the cartesian square:
\begin{equation*}
\begin{CD}
\cal{U} @>{\barr{\tau}}>> \mku \\
@VV{\pi_Y}V @VV{\pi}V \\
Y @>{\tau}>> B\subset \mk. \\
\end{CD}
\end{equation*}
The line bundles $\barr{\tau}^*(\HH_i)$ for $1\leq i \leq r$ are
canonically isomorphic to 
$\LL=\barr{\tau}^*(\HH_0)$ on $\cal{U}$  since 
$$\barr{\tau}^*(\HH_i\otimes \HH_0^{-1}) 
\eqq \pi_Y^*\rho_i^* \tau_i^*(\GG_i)$$
and $\tau_i^*(\GG_i)$ is canonically trivial.

Via pull-back and
the canonical isomorphisms, $\barr{\tau}^*(s_i)$ canonically
corresponds to a section of $\LL$.
Since these $r+1$ 
sections do not vanish simultaneously, they
define a morphism of $\mu:\cal{U} \rarr \proj^r$. The
canonical method of obtaining $\mu$ is as follows.
Define a vector space map $V^* \rarr H^0(\LL)$ by sending
$t_i$ to $\barr{\tau}^*(s_i)$.
The induced surjection $V^*\otimes \oh \rarr \LL$ canonically
yields a morphism $$\mu: \cal{U}\rarr \proj^r.$$
Note that the equations $t_i$ are used to
define the morphism $\mu$.
The sections $\{p_i \}$, $\{q_{i,j} \}$ pull back to
sections of $\pi_Y$. We claim that the 
family 
\begin{equation}
\label{uufam}
(\pi_Y:\cal{U} \rarr Y,
  \{p_i\},  \{q_{i,j}\}, \mu)
\end{equation}
is a universal family of $\barr{t}$-rigid stable maps, so
$\smap\stvo=Y$. 

The stability of the family of maps
\begin{equation}
\label{ufam}
(\pi_Y: \cal{U} \rarr Y,  \{p_i\},  \mu)
\end{equation}
is straightforward. 
Each fiber $C$ of $\pi_Y$ is an $m$-pointed, genus $0$ 
stable curve with markings $\{p_i\}$ and $\{q_{i,j}\}$. 
Let $E\subset C$ be an irreducible component.
Suppose 
$\text{dim}(\mu(E))=0$. 
By the transversality condition (iii), 
$E$ has no markings from the sections
$\{q_{i,j}\}$. 
Since $C$ is a stable $m$-pointed curve and no $\{q_{i,j}\}$ markings
lie on $E$,  
$\text{deg}_E(\omega_C(p_1+\ldots +p_n))>0$. Hence,
condition (1) in the definition of map stability (section
1.1) holds for $E$.
Therefore (\ref{ufam}) is a stable family of maps. By construction,
it is a $\barr{t}$-rigid stable family.

Finally, it must be shown (\ref{uufam}) is universal.
Let 
\begin{equation}
\label{testf}
(\pi: \cal{C} \rarr S, \{p_i \}, \{ q_{i,j} \},  \nu)
\end{equation} be
a family of $\barr{t}$-rigid stable maps.
Since $(\pi:\cal{C} \rarr S,  \{p_i \},  \{q_{i,j}\})$ is
a flat family of $m$-pointed, genus $0$  stable curves, there
is an induced map
$\lambda : S \rarr \barr{M}_{0,m}$ such that the pull-back family
$S\times_{\mk} \mku$ is canonically isomorphic
to $(\pi:\cal{C} \rarr S, \{p_i\}, \{q_{i,j} \})$.

First we show $\lambda$ is $\HH_i$-balanced. 
The pair 
$(\barr{\lambda}^*(\HH_i), \barr{\lambda}^*(s_i))$
yields the
Cartier divisor 
$q_{i,1} + \ldots + q_{i,d}$ on $\cal{C}$.
The map $\nu$ is induced by a vector space homomorphism
$\psi:V^* \rarr H^0(\cal{C}, \nu^*(\oh_{\proj(V)}(1)))$.
Let $z_i= \psi(t_i)$.  
By condition (iii) of
$\barr{t}$-rigid stability, 
the pair 
$(\nu^*(\oh_{\proj(V)}(1)), z_i)$
yields the Cartier divisor
$q_{i,1} + \ldots + q_{i,d}$ on $\cal{C}$.
By Lemma \ref{cart}, there are {\em canonical} isomorphisms
\begin{equation}
\label{caniso}
\barr{\lambda}^*(\HH_i) \eqq \nu^*(\oh_{\proj(V)}(1))
\end{equation}
for all $i$. Hence $\lambda$ is $\HH_i$-balanced.

By the universal
property of $B$, $\lambda$ factors through $B$: 
$\lambda: S \rarr B$.
There are canonical isomorphisms
\begin{equation}
\label{can2}
\pi_*(\barr{\lambda}^*(\HH_i \otimes \HH_0^{-1})) \eqq \lambda^*(\GG_i).
\end{equation}
The canonical isomorphisms (\ref{caniso}) yield
canonical sections of 
$\barr{\lambda}^*(\HH_i \otimes \HH_0^{-1})$. 
The canonical isomorphisms (\ref{can2}) then
yield
nowhere vanishing sections of $\lambda^*(\GG_i)$ over $S$. Hence there is
a canonical
a map $S \rarr Y$. It is easily checked the
pull-back of the universal family over $Y$
yields a $\barr{t}$-rigid stable family of maps canonically
isomorphic to (\ref{testf}).

%
%
%
%
%

\section{{\bf The construction of} $\smap\stp$}
\subsection{Gluing}
\label{glu}
While a given pointed stable map $\mu: C \rarr \proj^r$ may not
be rigid for a given basis $\barr{t}$ of 
$V^*=H^0(\proj^r,\oh_{\proj^r}(1))$, the map will be
rigid (by Bertini's theorem) for some choice of basis.
The moduli space $\smap\stp$ is obtained by gluing  
together quotients of $\smap\stv$ for different
choices of bases $\barr{t}$.

For notational convenience, set $\M(\barr{t})= \smap\stv$.
We write $(\pi: \cal{U} \rarr \M(\barr{t}), \{p_i\},
\{q_{i,j}\}, \mu)$ for the universal family of
$\barr{t}$-rigid stable maps in the genus 0 case.
If $g>0$, more care is required.

Let $\goth_d$ denote the symmetric group on $d$ letters. 
The group
$$G=G_{d,r} = \goth_d \times \ldots \times \goth_d \ \ \ 
(r+1 \ \text{factors})$$
has a natural action on $\M(\barr{t})$
obtained by
permuting the ordering in each
of the $r+1$ sets of sections $\{q_{i,1}, \ldots, q_{i,d}\}$,
$0 \leq i \leq r$. 
For any $\sigma \in G$, the family
\begin{equation}
\label{permm}
(\pi: \cal{U} \rarr \barr{M}(\barr{t}), \{p_i\}, \{q_{i, \sigma(j)}\}, 
\mu)
\end{equation} 
is also a $\barr{t}$-rigid family over $\M(\barr{t})$.
By the universal property, the permuted family (\ref{permm})
induces an automorphism of $\M(\barr{t})$. Since $\M(\barr{t})$
is quasi-projective and $G$ is finite, 
there is a quasi-projective quotient
scheme $\M(\barr{t})/G$.

Let $\barr{t}$ and $\barr{t}'$ be distinct choices of bases of $V^*$.
Let $\mu: \cal{U} \rarr \proj^r$ be the universal family over
$\smap(\barr{t})$.
Let 
$$\smap(\barr{t}, \barr{t}') \subset \smap(\barr{t})$$
denote the open locus over 
which the divisors $\mu^*(t'_0), \ldots, \mu^*(t'_r)$
are \'etale, disjoint, and disjoint from the sections $\{p_i\}$.
The open set $\smap(\barr{t}, \barr{t}')$ is certainly $G$-invariant. 
Let 
$\M(\barr{t}, \barr{t}')/G$ denote the quasi-projective
quotient.
\begin{pr}
\label{bsugg}
There is a canonical isomorphism
$$\M(\barr{t}, \barr{t}')/G \eqq \M(\barr{t}', \barr{t})/G.$$
\end{pr}
\bpf
The divisors  $\mu^*(t'_i)$ define an \'etale  Galois cover $\cal{E}$ of
$\smap(\barr{t}, \barr{t}')$ with Galois group
$G$ over which a
$\barr{t}'$-rigid stable family is defined.
The fiber of $\cal{E}$ over 
$(C, \{p_i\}, \{q_{i,j}\}, \mu)$ is the set of
orderings $\{q'_{i,j}\}$ of the points
mapped by $\mu$ to the hyperplane $(t_i'=0)$.
Therefore there is a map
\begin{equation}
\label{patcher} 
\cal{E} \rarr \smap(\barr{t}')
\end{equation}
which is easily seen be $G$-equivariant for the
Galois $G$-action on $\cal{E}$ and the $\{q'_{i,j}\}$-permutation
$G$-action on the $\M(\barr{t}')$. 
Moreover (\ref{patcher}) factors through 
$\smap(\barr{t}', \barr{t})$.
Hence there exists a map of quotients
\begin{equation}
\label{patch2}
\smap(\barr{t}, \barr{t}') \eqq 
\cal{E}/ {\text{Galois}} \rarr \M(\barr{t}', \barr{t})/G.
\end{equation}
The map (\ref{patch2}) is $G$-invariant for the 
$\{q_{i,j}\}$-permutation action on $\barr{M}(\barr{t}, \barr{t}')$. 
Therefore (\ref{patch2}) descends to
$\M(\barr{t}, \barr{t}')/G \rarr \M(\barr{t}', \barr{t})/G$.
The inverse is obtained by interchanging $\barr{t}$ and
$\barr{t}'$ in the above construction. In fact,
there is a natural action of $G \times G$ on
$\cal{E}$ and canonical isomorphisms
$\M(\barr{t}, \barr{t}')/G \eqq \cal{E}/(G\times G)
\eqq  \M(\barr{t}', \barr{t})/G$.
\epf

In case $g>0$, the coarse moduli spaces
$\smap\stv$ do not (in general) have universal
families. The permutation action of $G$
can be defined on a Hilbert scheme or a stack
and then descended to $\smap\stv$. The open sets
$\M(\barr{t}', \barr{t})$ and $\M(\barr{t}, \barr{t}')$
are well defined for $g>0$ and still satisfy  
Proposition \ref{bsugg}.

The cocycle conditions on triple intersections are easily established.
Hence, the schemes $\M(\barr{t})/G$ canonically patch together along
the open sets $\M(\barr{t}, \barr{t}')/G$ to form the scheme
$\smap\stp$. The results on boundedness show 
$\smap\stp$ is covered by a finite number  of 
these open sets $\M(\barr{t})/G$.
Hence,
$\smap\stp$ is an algebraic scheme of
finite type over $\com$.
The universal properties of $\smap\stp$ are easily obtained
from the universal properties of the 
moduli spaces of $\barr{t}$-rigid stable maps.

\subsection{Separation and completeness}
Let $(X,x)$ be a nonsingular, pointed curve. Let 
$\iota: X\smallstm \{x\}=U\ \hookrightarrow X$.
Let 
\begin{equation}
\label{first}
(\pi:\cal{C} \rarr X, \ \{p_i \}, \ \mu)
\end{equation}
\begin{equation}
\label{second}
(\pi': \cal{C}' \rarr X, \ \{p_i'\}, \ \mu')
\end{equation} 
be two families over $X$ of
stable maps to $\proj^r=\proj(V)$. 
\begin{pr}
\label{sep}
An isomorphism between the families (\ref{first}) and
(\ref{second}) over $U$ extends to an
isomorphism over $X$.
\end{pr}
\bpf
Choose a basis $\barr{t}=(t_0, \ldots, t_r)$ of $V^*$ that intersects
the maps $\mu:\cal{C}_x \rarr \proj^r$ 
and $\mu':\cal{C}'_x \rarr \proj^r$
transversally at unmarked, nonsingular
points. Since it suffices to prove the isomorphism
extends over a 
local \'etale
cover of $(X,x)$, it can be assumed that the Cartier divisors
$\mu^*(t_i)$ and ${\mu'}^*(t_i)$ split into sections 
$\{q_{i,j}\}$ and $\{q'_{i,j}\}$ of $\pi$ and $\pi'$.
Then $\cal{C}$, $\cal{C}'$ are Deligne-Mumford stable $m=n+d(r+1)$
pointed curves. Therefore, by the separation property of the
functor of Deligne-Mumford stable pointed curves, there
exists an isomorphism (of pointed curves) $\tau: \cal{C} \rarr \cal{C}'$ over
$X$. Since $\tau\circ \mu'$ and $\mu$ agree on an open set, 
$\tau\circ \mu' = \mu$.
\epf

\noindent
Proposition \ref{sep} and the valuative criterion show
$\M_{g,n}(X, \beta)$ is a separated algebraic scheme.

Properness is also established by the valuative
criterion. To complete $1$ dimensional families of
stable maps, semi-stable reduction
techniques for curves are used (as in  [K-K-M] and [Ha]).
\begin{pr}
\label{prop}
Let $\cal{F}=(\pi:\cal{C}\rarr U, \ \{p_i\}, \ \mu)$ be a family
of stable maps to $\proj^r$. There exists a
base change $\gamma:(Y,y) \rarr (X,x)$ satisfying:
\begin{enumerate}
\item[(i)] $\gamma_W: Y\smallstm \{y\}= W \rarr U$ is \'etale.
\item[(ii)] The pull-back family $\gamma_W^*(\cal{F})$ 
extends to a 
stable family over $(Y,y)$.
\end{enumerate}
\end{pr} 
\bpf
First, after restriction to a Zariski open subset
of $U$, it can be assumed that the fibers $\cal{C}_{\xi}$
all have the same number of irreducible components.
There may be non-trivial monodromy around the
point $x \in X$ in the set of 
irreducible components of the fibers $\cal{C}_{\xi}$. After a base
change (possibly ramified at $x$), this monodromy can
be made trivial. It can therefore be
assumed that $\cal{F}$
is a union of 
stable families $\cal{F}_j = (\pi_j:\cal{C}_j \rarr U,
\ \{p^j_i\}, \{p^c_i\}, \ \mu_j )$  where
$\pi_j$ is family of {\em irreducible}, nodal, 
projective curves. The markings
$\{p^j_i\}$ are the markings of $\cal{C}$ that lie on $\cal{C}_j$.
The marking $\{p^c_j\}$ correspond to intersections of
components in $\cal{F}$. It suffices to prove Proposition \ref{prop}
separately for each stable family $\cal{F}_j$.

For technical reasons,  it is convenient to
consider families of nonsingular curves.
After restriction, normalization, and base change of $\cal{F}_j$, a
family  
\begin{equation}
\label{famm}
\tilde{\cal{F}}_j = (\tilde{\pi}_j:\tilde{\cal{C}}_j \rarr U,
 \{p^j_i\}, \{p^c_i\}, \{p^n_i\},  \tilde{\mu}_j )
\end{equation}
can be obtained where  $\tilde{\cal{F}}_j$
is a family of stable maps of 
irreducible, {\em nonsingular}, projective curves.
The additional markings $\{p^n_i\}$ correspond to the
nodes. Consider the nodal locus in $\cal{F}_j$.
This locus consists of curves and isolated points.
Via restriction of $U$ to a Zariski open set, it
can be assumed the nodal locus (if non-empty) 
is of pure dimension
$1$. A normalization now separates the 
sheets along the nodal locus. A base change then 
may be required to make the separated points $\{p^n_i\}$
sections.
If the normalized family $\tilde{\cal{F}}_j$ is completed,
$\cal{F}_j$ can be completed by identifying the nodal
markings on $\tilde{\cal{F}}_j$. This nodal identification
commutes with the map to $\proj^r$.
It therefore suffices to prove Proposition \ref{prop}
for these normalized families (\ref{famm}).

By the above reductions, it suffices to prove Proposition {\ref{prop}}
for a family of stable maps of irreducible, nonsingular,
projective curves.
Let 
\begin{equation}
\label{fammm}
(\pi: \cal{C} \rarr U, 
\{p_i\},  \mu )
\end{equation}
be such a family.
Let $\pi:\cal{E} \rarr X$ be a flat extension of
  $\pi:\cal{C} \rarr U$ over the point $x\in X$.
After blow-ups in the special fiber of $\cal{E}$, it
can be assumed the map $\mu: \cal{C} \rarr \proj^r$
extends to $\mu: \cal{E} \rarr \proj^r$.
By Lemma \ref{ssred} below applied to 
the flat extension  $\pi:\cal{E} \rarr X$,
there exists a base change $\gamma:(Y,y) \rarr (X,x)$
and a 
family of 
pointed curves $$\pi_Y: \cal{C}_Y \rarr (Y,y)$$
satisfying conditions (i)$-$(iii)
of Lemma \ref{ssred}. 
Via $\tau
: \cal{C}_{Y} \rarr \cal{E}$, 
${\mu}$ naturally
induces a map $${\mu}_{Y}:
{\cal{C}}_{Y} \rarr \proj^r.$$

The family $({\pi}_{Y}:{\cal{C}}_{Y} \rarr (Y,y),
\{p_i\}, {\mu}_{Y} )$ 
is certainly an extension of the family over $Y \smallstm \{y\}$
determined
by the $\gamma$ 
pull-back of the stable family (\ref{fammm}).
The special fiber 
is a map of a pointed quasi-stable curve to $\proj^r$.
Unfortunately, the special fiber may not be stable.
A stable family of maps is produced in two steps.
First, unmarked, ${\mu}_{Y}$-collapsed, $-1$-curves in the
special fiber are sequentially blow-down.  
A multiple of the line bundle
\begin{equation}
\label{lbdle}
\omega_{{\pi}_{Y}}( \sum_i p_i) \otimes
{\mu}_{Y}^*(\oh_{\proj^r}(3))
\end{equation}
is then  ${\pi}_{Y}$- relatively basepoint free. Second, 
as in [Kn], the relative
morphism determined by a power of the line bundle 
(\ref{lbdle}) blows-down the remaining destabilizing
$\proj^1$'s to yield a stable extension over
$(Y,y)$. 
\epf

\begin{lm}
\label{ssred}
Let $\pi_X: \cal{S}_X \rarr (X,x)$ be a flat, projective family
of curves with $l$ sections $s_1, \ldots, s_l$ satisfying the
following condition:
$\forall \xi \neq x $, $\pi^{-1}(\xi)=\cal{C}_{\xi}$
is a projective nonsingular curve with $l$ distinct marked
points $s_1(\xi), \ldots, s_l(\xi)$. There exists
a base change $\gamma: (Y,y) \rarr (X,x)$
\'etale except possibly at $y$ with a 
family of 
$l$-pointed curves $\pi_Y: \cal{S}_Y \rarr (Y,y)$ and a 
diagram:
\begin{equation*}
\begin{CD}
\cal{S}_Y @>{\tau}>> \cal{S}_X\\
@VV{\pi_Y}V @VV{\pi_X}V \\
(Y,y) @>{\gamma}>> (X,x) \\
\end{CD}
\end{equation*}
satisfying the following properties:
\begin{enumerate}
\item[(i)] $\cal{S}_Y$ is a nonsingular surface. 
$\pi_Y:\cal{S}_Y \rarr (Y,y)$ is a  
flat, projective family
of $l$-pointed quasi-stable curves.
\item[(ii)] For each marking $1\leq i \leq l$, $\tau\circ s_i=
s_i\circ \gamma$.
\item[(iii)] Over $W=Y\smallstm  \{y\}$, there is isomorphism
$\cal{S}_W \stackrel{\sim}{\rarr} \gamma_W^*(\cal{S}_U)$,
where $U= X \smallstm \{x\}$. 
The morphism $\tau|_{\cal{S}_W}$ is the composition
$$\cal{S}_W \stackrel{\sim}{\rarr} 
\gamma_{W}^*(\cal{S}_U) \rarr \cal{S}_U$$
where the second map is the natural projection. 
\end{enumerate}
\end{lm}
\bpf
The method is by standard semi-stable reduction (cf.
[K-K-M], [Ha]). First, the
singularities of $\cal{S}_X$ are resolved. Note that all singularities
lie in the special fiber. Next, the surface
$\cal{S}_X$ is blown-up sufficiently to ensure the reduced scheme
supported on the special fiber
has normal crossing singularities in $\cal{S}_X$. The required
blow-ups have point
centers in the special fiber. Finally, the resulting surface
is blown-up further (at points in the special
fiber) to ensure the marking sections $s_1, \ldots, s_l$ do not
intersect each other and do not pass through nodes
of the reduced scheme 
supported on the special fiber. 
Let $\hat{\pi}: \hat{\cal{S}}_X \rarr (X,x)$ be the resulting
nonsingular surface. The singularities of the morphism $\hat{\pi}$
are locally of the form $z_1^\alpha  z_2^\beta = t$ where
$z_1$, $z_2$ are coordinates on $\hat{\cal{S}}_X$ and $t$ is a coordinate
on $X$. Let $\{\alpha_j, \beta_j\}$
be the set of exponents that occur at the singularities
of $\hat{\pi}$. Let $\gamma: (Y,y) \rarr (X,x)$ be a
 base change whose ramification index over $x$
is divisible by all $\alpha_j$ and $\beta_j$. 
Let $\cal{S}_Y$ be the normalization of
$\gamma^*(\hat{\cal{S}}_X)$. A straightforward local analysis
shows the family $\pi_Y: \cal{S}_Y\rarr (Y,y)$ 
has an $l$-pointed, reduced, nodal special fiber.
The surface $\cal{S}_Y$ has singularities of the local form $z_1z_2-t^k$ 
in the special fiber.
Blowing-up $\cal{S}_Y$ yields a nonsingular surface with the
required properties (i)$-$(iii).
\epf

By the valuative criterion, Propositions \ref{sep} and \ref{prop}
prove $\smap\stp$ is a separated and proper complex 
algebraic scheme.


\subsection{Projectivity}
The projectivity of the proper schemes $\smap\stp$ is
established here by a method due to J. Koll\'ar ([Ko1]).
Proofs of the projectivity of $\M_{g,n}(\proj^r, \beta)$
can also be found in [A] and [C].
Koll\'ar constructs ample line bundles on proper
spaces via sufficiently nontrivial quotients of 
semipositive vector bundles. 
A vector bundle $E$ on an algebraic scheme $S$ is
{\em semipositive} if for every morphism of a projective curve $f: C \rarr S$,
every quotient line bundle of $f^*(E)$ has nonnegative degree on $C$.

The first step is a semipositivity lemma.
Let 
\begin{equation}
\label{ifam}
\cal{F}=(\pi: \cal{C} \rarr S, \ \{p_i\}, \ \mu)
\end{equation}
be a stable family of maps over $S$ to $\proj^r$. Let
$$E_k(\pi)=\pi_*\bigg(\omega^k_{\pi}
(\sum_{i=1}^{n} k p_i)\otimes \mu^*(\oh(3k)) \bigg).$$ 
\begin{lm}
\label{semip}
$E_k(\pi)$ is a semipositive
vector bundle on $S$ for $k\geq 2$.
\end{lm}
\bpf
A slight perturbation of the arguments in [Ko1] is required.
It suffices to prove semipositivity in case the base is a nonsingular
curve $X$.
Let $\gamma:Y \rarr X$ be a flat base change.
By map stability, Serre duality, and the base change theorems, 
it follows (for $k\geq 2$)
$E_k$ commutes with pull-back:
$$ E_k(\pi_Y)\eqq \gamma^*(E_k(\pi_X)) $$
where $\pi_Y$ is the pull-back family over $Y$.
It therefore suffices to
prove semipositivity after base change.

Using the methods of section 4.2,
it can be assumed (after base change) that $\cal{F}$
is a union of component 
stable families $\cal{F}_j = (\pi_j:\cal{C}_j \rarr X,
\ \{p^j_i\}, \{p^c_i\}, \ \mu_j )$  where
$\pi_j$ is family of stable maps and the generic 
element of $\cal{F}_j$ is a map of an {\em irreducible}, 
projective, nodal
curve. The notation introduced in the proof of
Proposition \ref{prop} is employed. After further base change
and normalization of $\cal{F}_j$, it can be assumed that 
\begin{equation}
\label{famm2}
\tilde{\cal{F}}_j = (\tilde{\pi}_j:\tilde{\cal{C}}_j \rarr X,
\ \{p^j_i\}, \{p^c_i\}, \{p^n_i\}, \ \tilde{\mu}_j )
\end{equation}
is a family of  stable maps where the generic element
is a map of an 
irreducible, projective, {\em nonsingular}  curve.

A semipositivity result for the family $\tilde{\cal{F}}_j$ 
is first established.
Let $H_1$, $H_2$, $H_3 \subset \proj^r$ be general hyperplanes.
After base change,
it can be assumed  
$\tilde{\mu}_j^*(H_l)$ is a union of $d$ reduced
sections for each $l$. 
These $3d$ sections are distinct from the sections $\{p^j_i\}$,
$\{p^c_i\}$, $\{p^n_i\}$.
Therefore,
\begin{equation}
\label{mess}
\omega^k_{\tilde{\pi}_j}
\big( \sum k p^j_i+ \sum (k-1) p^c_i + \sum
(k-1)p^n_i \big)\otimes \tilde{\mu}_j^*(\oh_{\proj^r}(3k)) 
\stackrel{\sim}{\rarr}
\end{equation}  
$$\omega^k_{\tilde{\pi}_j}
(\sum \alpha_q X_q)$$
where $X_q$ are distinct sections of $\tilde{\pi}_j$ and $\alpha_q\leq k$.
The surface $\tilde{\cal{C}}_j$ has finitely many
singularities of the form $z_1z_2-t^{\alpha}$. These
singularities are resolved by blow-up, $$\tau: \cal{S}_j \rarr \tilde
{\cal{C}_j}.$$ Since the relative dualizing sheaf of the
family $\cal{S}_j$ is trivial on the exceptional $\proj^1$'s of
$\tau$,  Lemma \ref{kol} below can be applied to 
deduce the semipositivity of $F_k(\tilde{\pi}_j)$ for 
$k\geq 2$ where 
$$F_k(\tilde{\pi}_j)=\tilde{\pi}_{j*}\bigg(\omega^k_{\tilde{\pi}_j}
\big(\sum kp^j_i+\sum (k-1)p^c_i + \sum
(k-1)p^n_i\big)
 \otimes \mu^*(\oh(3k)) \bigg).$$

For $k\geq 2$, the restriction of the line bundle (\ref{mess}) to
a fiber of $\tilde{\cal{C}}_j$ is equal to
\begin{equation}
\label{prod} \omega \otimes \omega^{k-1}
\big(\sum (k-1)p^j_i+\sum (k-1)p^c_i + \sum
(k-1)p^n_i\big)\otimes 
\end{equation}
$$\tilde{\mu}_j^*(\oh_{\proj^r}(3k-3))
\otimes \mu_j^*(\oh_{\proj^r}(3))(\sum_i p^j_i)$$
where $\omega$ is the dualizing sheaf of the fiber.
By stability for the family $\tilde{\cal{F}}_j$, the
product of the middle two factors 
in (\ref{prod}) is ample for $k\geq 2$. The
last factor in (\ref{prod}) is certainly of non-negative degree.
By Serre duality, for $k\geq 2$,
\begin{equation}
\label{van}
R^1{\tilde{\pi}}_{j*}  \bigg(\omega^k_{\tilde{\pi}_j}
\big(\sum kp^j_i+\sum (k-1)p^c_i + \sum
(k-1)p^n_i\big)
 \otimes \mu^*(\oh(3k)) \bigg)   =0.
\end{equation}
The semipositivity of $E_k(\pi)$ will be obtained
from the semipositivity of $F_k(\tilde{\pi}_j)$.
The $(k-1)$-multiplicities will naturally arise
in considering dualizing sheaves on nodal and reducible curves.

Let $\tilde{{\pi}}_{\cup_j}: \bigcup_{j} \tilde{\cal{C}}_j \rarr 
X$
be the disjoint union of the families $\tilde{\cal{F}}_j$.
There is natural morphism from the
disjoint union to $\cal{C}$
$$\rho: \bigcup_{j} \tilde{\cal{C}}_j \rarr \cal{C}$$
obtained by identifying nodal marked points and gluing components along
intersection marked points. 
Consider the natural sequence of sheaves on $\cal{C}$:
\begin{equation}
\label{seqq}
0 \rarr \rho_*(\omega_{\tilde{\pi}_{\cup_j}}) \rarr \omega_{\pi} \rarr
K \rarr 0.
\end{equation}
The quotient $K$ is easily identified as $\bigoplus_{p^c_i, p^n_i} \oh_p$
where the sum is over all nodal and component intersection sections
of the family $\cal{F}$. Tensoring (\ref{seqq}) with the
line bundle $\omega_{\pi}^{k-1}(\sum kp_i)
\otimes \mu^*(\oh_{\proj(V)}(3k))$
yields the exact sequence:
$$0 \rarr \rho_*\bigg(\omega_{\tilde{\pi}_{\cup_j}}^k
(\sum kp_i+\sum (k-1)p^c_i + \sum
(k-1)p^n_i) \otimes \mu^*(\oh_{\proj(V)}(3k)) \bigg) \rarr $$
$$ \omega^k_{\pi}(\sum_i kp_i) 
\otimes \mu^*(\oh_{\proj(V)}(3k))
\rarr \bigoplus_{p^c_i, p^n_i}
\oh_p\otimes \mu^*(\oh_{\proj(V)}(3k)) \rarr 0.$$
Certainly $\tilde{\pi}_{\cup_{j}*}=
\pi_*\rho_*$. Note the vanishing of $R^1$ determined in (\ref{van}).
These facts imply the $\pi$ direct image of the above
sequence on $\cal{C}$  yields an exact
sequence on $X$:
$$0 \rarr \bigoplus_j F_k(\tilde{\pi}_j) 
\rarr E_k(\pi) \rarr \bigoplus_{p^c_i,
p^n_i} \oh_X \otimes \mu^*(\oh_{\proj(V)}(3k)) \rarr 0.$$
Finally, since an extension of semipositive bundles is semipositive
([Ko1]),
$E_k(\pi)$ is semipositive.
\epf

\begin{lm}
\label{kol}
Let $\pi:\cal{S}\rarr X$ be a map from a nonsingular projective
surface to a nonsingular curve. Assume the general fiber of
$\pi$ is nonsingular. Let $X_q$ be a set of distinct sections
of $\pi$. Then $$\pi_*(\omega_{\cal{S}/X}^k(\sum \alpha_q X_q))$$
is semipositive provided $k\geq 2$ and $\alpha_q\leq k$ for
all $q$.
\end{lm}
\bpf
This is precisely Proposition 4.7 of [Ko1].
\epf

The second step is the construction of a  non-trivial quotient.
Let $\cal{F}$ be the family (\ref{ifam}).
Let $\proj(E_k^*)$ be the projective bundle over $S$ obtained
from the subspace projectivization of $E_k^*$.
The condition of stability implies there is a canonical
$S$-embedding $\iota:\cal{C} \rarr \proj(E_f^*)$ for some
$f=f(d,g,n,r)$ (see section \ref{bound}).
The morphism $\mu$ 
then yields a canonical  $S$-embedding:
$$\gamma: \cal{C} \rarr \proj(E_f^*) \times_{\com} \proj^r.$$
The $n$ sections $\{p_i\}$ yield $n$ sections $\{(\iota \circ p_i,
 \mu \circ p_i) \}$ of
$\proj(E_f^*)\times \proj^r$ over $S$.
Let $S_i$ denote the subscheme of $\proj(E_f^*)\times \proj^r$
defined by the $i^{th}$ section.
Denote the projection of $\proj(E_f^*)\times \proj^r$ to
$S$ also by $\pi$.
Let $\cal{M}=\oh_{\proj(E_f^*)}(1) \otimes \oh_{\proj^r}(1)$.
$\cal{M}$ is an $\pi$-relatively ample line bundle.
Note $\pi_*(\cal{M}^l) \eqq \text{Sym}^l(E_f) 
\otimes \text{Sym}^l(\com^{r+1}).$
By the stability of
semipositivity under symmetric and tensor products ([Ko1]) 
and Lemma \ref
{semip}, $\pi_*(\cal{M}^l)$ is semipositive.
Fix a choice of $l$ 
(depending only on $d$, $g$, $n$, and $r$)
large enough to ensure
\begin{equation}
\label{ntq}
\pi_*(\cal{M}^l) \oplus \ \bigoplus_{i=1}^{n}\pi_*(\cal{M}^l)
\rarr \pi_*(\cal{M}^l\otimes \oh_{\cal{C}}) 
\oplus \ \bigoplus_{i=1}^{n} 
\pi_*(\cal{M}^l\otimes \oh_{S_i}) \rarr 0.
\end{equation}
Such a choice of $l$ is possible by the boundedness 
established in section \ref{bound}. 
Let $Q$ be the quotient in (\ref{ntq}).
By boundedness and the vanishing of higher direct images,
the quotient $Q$ is a vector bundle for large $l$.

The quotient (\ref{ntq}) is nontrivial in the following sense.
Let $G=GL$ be the structure group of the bundle $E_f$.
$G$ is naturally the structure group of $\pi_*(\cal{M}^l)$.
Let $W$ be the $G$-representation inducing the 
bundle $\pi_*(\cal{M}^l) \oplus \ \bigoplus_{1}^{n}\pi_*(\cal{M}^l)$.
Let $q$ be the rank of
the quotient bundle of (\ref{ntq}).
The quotient sequence (\ref{ntq}) yields a
set theoretic classifying map to the Grassmannian: 
$$\rho: S \rarr {\mathbf{Gr}}(q,W ^*)/G.$$
\begin{lm}
\label{a}
There exists a set theoretic injection $$\delta:
\smap\stp \rarr {\mathbf{Gr}} (q,W ^*)/G.$$
Let $\lambda: S \rarr \smap\stp$ be the map
induced by the stable family (\ref{ifam}).
There is a (set theoretic) factorization $\rho= \delta \circ \lambda$. 
\end{lm}
\bpf
For large $l$, the sequence (\ref{ntq}) is equivalent to the
data of a Hilbert point in $J$ (see section \ref{bound}).
Since the $G$ orbits of $J$ are exactly the stable
maps, the lemma follows. 
\epf

\begin{lm}
\label{b}
A stable map has a finite number of automorphisms.
\end{lm}
\bpf
As simple consequence of the definition of stability, there
 are no infinitesimal automorphisms. The total number is
therefore finite.
\epf

Suppose the map to moduli $\lambda:S \rarr \smap\stp$
is a generically finite algebraic morphism.
Then, in the terminology of [Ko1], Lemmas \ref{a} and \ref{b} 
show the classifying map $\rho$ is {\em finite}
on an open set of $S$.

\begin{pr}(Lemma 3.13, [Ko1]) Let the base $S$ of (\ref{ifam})
be a normal projective
variety. Suppose the classifying map is finite on an open 
set of $S$. Then,  the top self-intersection number of 
$\text{Det}(Q)$
on $S$ is positive.
\label{dett}
\end{pr}

If $\smap\stp$ were a fine moduli space equipped
with a universal family, $\text{Det}(Q)$ 
would be well defined
and ample (by Proposition \ref{dett}
 and the Nakai-Moishezon criterion)
on $\smap\stp$.
Since $\smap\stpo$ is expressed locally as a quotient
of a fine moduli space by a finite group, it is easily
seen $\text{Det}(Q)^k$ is a well defined line bundle on 
$\smap\stpo$ for some
sufficiently large $k$. The exponent $k$ is taken to 
trivialize the $\com^*$-representations at the fixed points. 
In the higher genus case,
$\text{Det}(Q)$ is a well defined line bundle
on the Hilbert scheme $J$ or the stack. Since the moduli
problem has finite automorphisms, 
$\text{Det}(Q)^k$ is well defined on the coarse
moduli space for some $k$.

Since the moduli spaces $\smap\stp$ are not fine,
subvarieties are not equipped with stable families.
Proposition (\ref{dett}) and
the Nakai-Moishezon criterion do not directly establish the
ampleness of $\text{Det}(Q)^k$.
An alternative approach (due to J. Koll\'ar) is followed.
Recall the Hilbert scheme $J$ (of section \ref{bound})
is equipped with a universal family and, therefore,
a canonical map
$$J\rarr \smap\stp.$$
Let $X\subset \smap\stp$ be a subvariety.
Using $J$ and the finite automorphism property of a stable
map, a morphism  $Y\rarr X$ of algebraic schemes can be
constructed satisfying
\begin{enumerate}
\item[(i)] $Y\rarr X$ is finite and surjective.
\item[(ii)] $Y$ is equipped with a stable
family of maps such that $Y\rarr X$ is the corresponding morphism
to moduli.
\end{enumerate}
The existence of $Y\rarr X$ is exactly the conclusion of
Proposition 2.7 in [Ko1] under slightly different assumptions.
Nevertheless, the argument is valid in the present setting.
The construction of $Y$ is subtle. First $Y$ is constructed
as an algebraic space. 
Then, a lemma of Artin is used
to find an algebraic scheme $Y$. 
Since $Y$ has a universal family, Proposition \ref{dett}
implies $\text{Det}(Q)^k$ has positive top intersection on $Y$ and
therefore on $X$. 
The Nakai-Moishezon criterion can be applied to conclude
the ampleness of $\text{Det}(Q)^k$ on $\smap\stp$.

\subsection{Automorphisms}
We use the notation of sections 3.2 and 4.1.
In the genus $0$ case, $\smap(\barr{t})$ is nonsingular.
Therefore, the space $\smap\stpo$ is locally a quotient of a nonsingular
variety by a finite group.
\begin{lm}
Let $\xi \in \smap(\barr{t})$ be a point
at which the $G_{d,r}$ action is not free. Then $\xi$ corresponds
to a stable map with nontrivial automorphisms.
\end{lm}
\bpf
$G_{d,r}$ acts by isomorphism on the stable
maps of the universal family over $\smap( \barr{t})$. 
The $G_{d,r}$ action
is not free at $\xi\in \smap(\barr{t})$ 
if and only if there exists a $1\neq \gamma\in 
G_{d,r}$
fixing $\xi$. The element $\gamma$ induces an automorphism
of the map corresponding to $\xi$. The automorphism is nontrivial
on the marked points $\{q_{i,j}\}$.
\epf

Over the automorphism-free locus, the $G_{d,r}$-action on
$\smap(\barr{t})$ (and on 
the universal family over $\smap(
\barr{t})$) is free. 
It follows that the quotient over the automorphism-free
locus is a nonsingular quasi-projective variety denoted by
$\M_{0,n}^*(\proj^r,d)$. A
universal family over
$\M_{0,n}^*(\proj^r,d)$ is obtained by patching.
Theorems \ref{rep} and \ref{t2} have been established
in the case $X\eqq \proj^r$.

%
%
%
%
%

\section{{\bf The construction of} $\smap\st$}
\subsection{Proof of Theorem \ref{rep}} 
Let $X$ be a projective algebraic variety.
Existence of the coarse moduli space $\smap\st$ is
established via a projective embedding
$\iota: X\hookrightarrow \proj^r$. Let
$\iota_*(\beta)$ be $d$ times the class of a line
in $\proj^r$.
\begin{lm}
There exists a natural closed subscheme $$\smap_{g,n}(X,\beta,\barr{t})
\subset \smap_{g,n}(\proj^r,d, \barr{t})$$ satisfying the following
property. Let $(\pi:\cal{C}\rarr{S}, \{p_i\}, \{q_{i,j}\},\mu)$
be a $\barr{t}$-rigid stable family of genus $g$, $n$-pointed,
degree $d$ maps
to $\proj^r$. Then,
the natural morphism $S\rarr \smap_{g,n}(\proj^r,d, \barr{t})$
factors through $\smap_{g,n}(X,\beta,\barr{t})$ if and only
if $\mu$ factors through $\iota$ and each geometric fiber
of $\pi$ is a map to $X$ representing the homology class $\beta\in 
A_1 X$.
\end{lm}
\bpf
The lemma is proved in case $g=0$. If $g>0$, then
$\smap_{g,n}(\proj^r,d, \barr{t})$ is not a fine moduli space
and the argument is more technical.

Let $$(\pi_M: \cal{U} \rarr \smap_{0,n}(\proj^r,d, \barr{t}),
\{p_i\}, \{q_{i,j}\}, \mu)$$ be the universal family
over $\smap_{0,n}(\proj^r,d, \barr{t})$. 
On a genus 0 curve, any vector
bundle generated by global sections has
no higher cohomology. Therefore,
by this cohomology vanishing
and the base change theorems,
$\pi_{M*}\mu^*(\oh_{\proj^r}(k))$ is a vector bundle
for all $k>0$. (This argument must be modified
in the $g>0$ case since
$\pi_{M*}\mu^*(\oh_{\proj^r}(k))$ need not be a vector
bundle even on the Hilbert scheme $J$ or the stack. Nevertheless,
it is not hard to define the closed subscheme determined
by $X$ on the Hilbert scheme $J$ or the stack and then
descend it to the coarse moduli space.) 
Let $\cal{I}_X$ be the ideal sheaf of $X \subset \proj^r$.
Let $I_X(k)=H^0(\proj^r, \cal{I}_X(k))$.
Let $l>>0$ be selected so that $\cal{I}_X(l)$ is generated
by the global sections $I_X(l)$. These sections
$I_X(l)$ yield sections of the vector bundle
$\pi_{M*}\mu^*(\oh_{\proj^r}(l))$. 
Let $Z\subset \smap_{0,n}(\proj^r,d, \barr{t})$ be the 
scheme-theoretic
zero locus of these sections. 
The restriction of $\mu$ to $\pi_{M}^{-1}(Z)$ factors
though $\iota$. 
Since $Z$ is an algebraic scheme, $Z$ is a finite union
of disjoint connected components. The homology class in
$A_1(X)=H_2(X, \mathbb{Z})$ represented by a
map with moduli point in $Z$ is a deformation invariant
of the map. Therefore, the represented homology class
is constant  on each
connected component of $Z$. Let $Z_{\beta}\subset Z$
be the union of components of $Z$ which consist
of maps representing the class $\beta\in A_1 X$.
Let $\smap_{0,n}(X,\beta,\barr{t})=Z_\beta$. The required
properties are easily established.
\epf
 
By the functorial property, $\smap_{g,n}(X,\beta,\barr{t})$
is invariant under the $G_{d,r}\eqq
\goth_d\times \cdots \times 
\goth_d$ action
on $\smap_{g,n}(\proj^r,d,\barr{t})$. The quotient
 $$ \smap_{g,n}(X,\beta,\barr{t})/ G_{d,r}$$ is an open set
of $\smap_{g,n}(X,\beta)$.
A patching argument identical to the patching argument
of section \ref{glu} yields a construction of $\smap\st$
as a closed subscheme of $\smap\stp$. The functorial
property of $\smap\st$ shows the construction is 
independent of the projective embedding of $X$.
Projectivity of $\smap\st$ is obtained from the
projectivity of $\smap\stp$. This completes the
 proof of Theorem \ref{rep}.

\subsection{Proof of Theorem \ref{t2}}
\label{deff}
Let $g=0$. Let $X$ be a projective, nonsingular, convex
variety. Theorem
\ref{t2} is certainly true in case $\beta=0$ since
$\M_{0,n}(X,0)=\M_{0,n} \times X$. In general,
a deformation study is needed to establish Theorem \ref{t2}. 

By the functorial property,
the Zariski tangent space to 
the scheme $\smap_{0,n}(X,\beta,\barr{t})$
at the point $(C, \{p_i\}, \{q_{i,j}\}, \mu:C\rarr X)$ is
canonically isomorphic
to the 
space of first order deformations of the pointed $\barr{t}$-stable map
$(C, \{p_i\}, \{q_{i,j}\}, \mu:C\rarr X)$. 
The later deformation space corresponds bijectively
to the space of first order deformations of the pointed stable map
$(C, \{p_i\}, \mu:C\rarr X)$. 

Let $\text{Def}(\mu)$ denote the
space of first order deformations of the pointed stable map 
$(C, \{p_i\}, \mu:C\rarr X)$.
Consider first the case in which $C\eqq \proj^1$.
Let $\text{Def}_R(\mu)$ be the space of first order
deformations of $(C, \{p_i\}, \mu:C\rarr X)$
with $C$ held rigid. There is an natural exact sequence:
$$0\rarr H^0(C, T_C) \rarr \text{Def}_R(\mu) 
\rarr \text{Def}(\mu)\rarr 0.$$
Stability of $\mu$ implies the
left map is injective.
Let $\text{Hom}(C,X)$ be the quasi-projective scheme of
morphisms from $C$ to $X$ representing the
class $\beta$. $\text{Hom}(C,X)$ is an open subscheme of
the Hilbert scheme of graphs in $C\times X$.
The Zariski tangent space to $\text{Hom}(C,X)$ is
naturally identified: 
$$T_{\text{Hom}(C,X)}([\mu]) \eqq H^0(C, \mu^*T_X)$$
(see [Ko2]). There is an exact sequence:
$$0 \rarr \text{Ker} \rarr \text{Def}_R(\mu) 
\rarr H^0(C, \mu^*T_X)\rarr 0$$
where $\text{Ker}$ corresponds to the deformations of the markings.
Therefore, $\text{dim}_{\com}\text{Ker}=n$.
Since $X$ is convex, the above
sequences suffice to compute the dimension of $\text{Def}(\mu)$:
$$\text{dim}_{\com} \text{Def}(\mu)=  
\text{dim}(X)+ \int_\beta c_1(T_X) + n-3.$$
The dimension of the tangent space to 
$\smap_{0,n}(X,\beta,\barr{t})$ is established in case
$C\eqq \proj^1$.

Before proceeding further, the following deformation
result is needed. A proof can be found in [Ko2].
\begin{lm}
\label{mapdef}
Let $\cal{C}/S$ and $\cal{X}/S$ be flat, projective
schemes over $S$. Let $s\in S$ be a geometric point.
Let $C_s$, $X_s$ be the fibers over $s$ and let
$f:C_s\rarr X_s$ be a morphism. Assume the
following conditions are satisfied:
\begin{enumerate}
\item[(i)] $C_s$ has no
embedded points.
\item[(ii)] $X_s$ is nonsingular.
\item[(iii)] $S$ is equidimensional at $s$.
\end{enumerate}
Then, the dimension of every component of the
quasi-projective variety $\text{Hom}_S(\cal{C}, \cal{X})$
at the point $[f]$ is at least
$$\text{\em{dim}}_{\com}H^0(C_s,f^*T_{X_s})-
  \text{\em{dim}}_{\com} H^1(C_s, f^*T_{X_s}) 
+ \text{\em{dim}}_s S.$$
\end{lm}

Again, let $(C\eqq \proj^1, \{p_i\}, \{q_{i,j}\}, \mu:C\rarr X)$
correspond to a point of the space
$\smap_{0,n}(X,\beta,\barr{t})$.
By Lemma \ref{mapdef}
and the convexity of $X$, every component of
$\text{Hom}(C,X)$ at $[\mu]$
has dimension at least $\text{dim}_{\com} H^0(C, \mu^*T_X)$.
Therefore, every component of  $\smap_{0,n}(X,\beta,\barr{t})$
at $[\mu]$ has dimension at least
$ \text{dim}(X)+\int_\beta c_1(T_X)+ n-3$. By the previous tangent
space computation, it follows
$[\mu]$ is a nonsingular point of
$\smap_{0,n}(X,\beta,\barr{t})$. Before
attacking the reducible case, a lemma is required.

\begin{lm}
\label{redd}
Let $X$ be a nonsingular, projective, convex space.
Let $\mu:C \rarr X$ be a morphism of a 
projective, 
connected, reduced, nodal curve of arithmetic genus 0 to $X$.
Then, 
\begin{equation}
\label{kspc}
H^1(C, \mu^*T_X)=0.
\end{equation}
and $\mu^*T_X$ is generated by global sections on $C$.
\end{lm}
\bpf
Let $E\subset C$ be an irreducible component 
of $C$; $E\eqq \proj^1$. Let
$$\mu^*T_X|_E\eqq \bigoplus \oh_{\proj^1}(\alpha_i).$$
Suppose there exists $\alpha_i<0$. The composition of a rational 
double cover of                                       
$E$ with $\mu$ would then violate the convexity of $X$ . 
It follows that:
\begin{equation}
\label{condy}
\forall i, \ \ \alpha_i\geq 0.
\end{equation}                   

We will prove the following statement by induction on
the number of components of $C$: 
\begin{equation}
\label{ggll}
H^1(C, \mu^* T_X \otimes \oh_C(-p))=0
\end{equation}
for any nonsingular point $p\in C$. Equation (\ref{ggll})
is true by  condition (\ref{condy})
when $C\eqq \proj^1$ is irreducible.
Assume now $C$ is reducible and $p\in E\eqq \proj^1$.
Let $C= C' \cup E$; let
$\{p'_1,\ldots, p'_q\}= C'\cap E$. Since $C$ is a tree,
$C'$ has exactly $q$ connected components
each intersecting $E$ in exactly $1$ point.
There is a component sequence:
 $$ 0 \rarr \mu^* T_X|_{C'} \otimes \oh_{C'}(-\sum_{j=1}^{q}p'_j)
\rarr \mu^* T_X \otimes \oh_C(-p) \rarr 
\mu^* T_X |_{E} \otimes \oh_E(-p)
\rarr 0.$$
Equation (\ref{ggll}) now follows from the
inductive assumptions on $C'$ and $E$. 
The inductive assumption (\ref{ggll}) is applied
to every connected component of $C'$.

We now prove $H^1(C, \mu^* T_X)=0$.
If $C\eqq \proj^1$, then the lemma is established
by condition (\ref{condy}). Assume now
$C= C' \cup E$ where $E\eqq \proj^1$. There is
a component sequence
\begin{equation}
\label{igoth}
0 \rarr \mu^* T_X|_{C'} \otimes \oh_{C'}(-\sum_{j=1}^q p'_j)
\rarr \mu^* T_X  \rarr 
\mu^* T_X |_{E}  
\rarr 0.
\end{equation}
Equation (\ref{kspc}) now follows from (\ref{ggll})
applied to every connected component of $C'$.

Finally, an analysis of sequence (\ref{igoth}) also
yields the global generation result.
$\mu^* T_X|_E$ is generated by global sections by
(\ref{condy}). Sequence (\ref{igoth}) is
exact on global sections by (\ref{ggll}). Hence
$\mu^* T_X$ is generated by global sections on $E$.
But, every point of $C$ lies on some component $E\eqq \proj^1$.
\epf       

In sections 7 and 8, the following related
lemma will be required:
\begin{lm}
\label{george}
Let $\mu: \proj^1 \rarr X$ be a non-constant morphism
to a nonsingular, projective, convex space $X$.
Then $\int_{\mu_*[\proj^1]} c_1(T_X) \geq 2$.
\end{lm}
\bpf
Since $\mu$ is non-constant, the differential
$$d\mu: T_{\proj^1} \rarr \mu^*(T_X)$$
is nonzero. Let $s \in H^0(\proj^1,T_{\proj^1})$
be a vector field with two distinct zeros $p_1, p_2 \in \proj^1$.
Then, $d\mu(s) \in H^0(\proj^1, \mu^*(T_x))\neq 0$
and $d\mu(s)$ vanishes (at least) at $p_1$ and $p_2$.
By the proof of Lemma \ref{redd}, $\mu^*(T_X) \eqq 
\bigoplus \oh_{\proj^1}(\alpha_i)$ where 
$\alpha_i \geq 0$ for all $i$. The existence of $d\mu(s)$ implies
that $\alpha_j \geq 2$ for some $j$.
\epf

Let $C$ now be a reducible curve.
$C$ must be a tree of $\proj^1$'s.
Let  $q$ be the number of nodes of $C$.
Again, let $\text{Def}(\mu)$ be the first order deformation
space of the pointed stable map $\mu$.
The {\em dual graph} of a pointed curve $C$ of
arithmetic genus $0$ consists of vertices and edges
corresponding bijectively to the
irreducible components
 and nodes of $C$ respectively.
The {\em valence} of a vertex in the dual graph is
the numbers of edges incident at that vertex.
Let $\text{Def}_G(\mu)\subset \text{Def}(\mu)$ be the first order
deformation space of the pointed stable map $\mu$ preserving the
dual graph. $\text{Def}_G(\mu)$ is a linear subspace
of codimension at most $q$.
Let $\text{Def}_G(C)$ be the space of first order
deformations of the curve $C$ which preserve the 
dual graph. A simple calculation yields
$$\text{dim}_{\com} \text{Def}_G(C) = \sum_{|\nu|\geq 4} |\nu|-3$$
where the sum is taken over vertices $\nu$ of
the dual graph of valence at least 4.

The natural linear map $\text{Def}_G(\mu) \rarr 
\text{Def}_G(C)$ is
now analyzed. Let $S$ be the nonsingular universal base space of
deformations 
of $C$ preserving the dual graph. Let $\cal{C}$ be
the universal deformation over $S$. Let $\cal{X}=X\times S$.
Let $s_0\in S$ correspond to $C$.
By Lemmas \ref{mapdef} and \ref{redd},
every component of $\text{Hom}_S(\cal{C},\cal{X})$ at $[\mu]$
has dimension
at least $\text{dim}(X)+\int_{\beta} c_1(T_X)+ 
\text{dim}(S)$.
The tangent space to the fiber of $\text{Hom}_S(\cal{C}, \cal{X})$
over $s_o$ at $[\mu]$ is canonically $H^0(C, \mu^*T_X)$. The
latter space has dimension $\text{dim}(X)+\int_{\beta} c_1(T_X)$.
Hence, $\text{Hom}_S(\cal{C},\cal{X})$ is nonsingular at $[\mu]$
of dimension $\text{dim}(X)+ 
\int_{\beta} c_1(T_X)+ \text{dim}(S)$ and the
projection morphism to $S$ is smooth at $[\mu]$. 
Therefore, $\text{Def}_G(\mu) \rarr \text{Def}_G(C)$ is surjective.

The above definitions and results yield a  
natural exact sequence:
$$0\rarr \text{Def}_C(\mu) \rarr \text{Def}_G(\mu) \rarr 
\text{Def}_G(C)\rarr 0$$
where $\text{Def}_C(\mu)$ is the space of first order deformations
of the pointed  stable map $\mu$ which restrict to the trivial
deformation of $C$. As in the case where $C\eqq \proj^1$, 
$\text{Def}_C(\mu)$ differs from  
$\text{Def}_R(\mu)$ only by the
tangent fields obtained from automorphisms:
$$0\rarr H^0(C, T^{auto}_C) \rarr \text{Def}_R(\mu)
\rarr \text{Def}_C(\mu) \rarr 0.$$ 
$H^0(C, T^{auto}_C)$ is the space of tangent
fields on the components of $C$ that vanish
at all the nodes of $C$.
Note $H^0(C, T^{auto}_C)= \sum_{|\nu|\leq 3}3-|\nu|$.
Finally, there is an exact sequence containing 
$\text{Def}_R(\mu)$ and the tangent space to $\text{Hom}(C,X)$:
$$0 \rarr \text{Ker}
\rarr \text{Def}_R(\mu) \rarr H^0(C, \mu^*T_X) \rarr  0.$$
From these exact sequences, Lemma \ref{redd}, and some
arithmetic, it follows
\begin{equation}
\label{ooo}
\text{dim}_{\com}\text{Def}_G(\mu)= 
\text{dim}(X)+ \int_{\beta} c_1(T_X) + n-3 -q.
\end{equation}

Let
$\cal{C}$ be a smoothing of the reducible curve $C$ over
a base $S$ and let $\cal{X}=X\times S$.
 A simple application of Lemma \ref{mapdef} 
shows that $[\mu]\in
\smap_{0,n}(X,\beta,\barr{t})$
lies in the closure of the locus of maps with irreducible
domains. Since the irreducible domain
locus is pure dimensional of dimension  
$\text{dim}(X)+\int_{\beta} c_1(T_X)+ n-3$,
\begin{equation}
\label{tttt}
\text{dim}_{\com}
\text{Def}(\mu) \geq  \text{dim}(X)+\int_{\beta} c_1(T_X)+ n-3.
\end{equation}
It follows from
(\ref{ooo}) and
(\ref{tttt}) that
$\text{Def}_G(\mu)$ is of maximal codimension $q$ in
$\text{Def}(\mu)$ and
that the inequality in (\ref{tttt}) is
an equality.
Since $\text{Def}(\mu)$ is of dimension
$\text{dim}(X)+\int_{\beta} c_1(T_X)+ n-3$, $[\mu]$ is a nonsingular
point of $\smap_{0,n}(X,\beta,\barr{t})$.
Since $\smap_{0,n}(X,\beta,\barr{t})$ is nonsingular of
pure dimension $\text{dim}(X)+\int_{\beta} c_1(T_X)+ n-3$, parts (i) and
(ii) of Theorem \ref{t2} are established. Part (iii) follows
from the corresponding result in the case $X\eqq \proj^r$.

%
%
%
%

\section{{\bf The boundary of} $\M_{0,n}(X,\beta)$}
\subsection{Definitions}
\label{bdry}
Let $X$ be nonsingular, projective, and convex. Let
the genus $g=0$.
The {\em boundary} of $\smap\sto$ is the
locus corresponding to reducible domain curves. 
Boundary properties of the Mumford-Knudsen space $\M_{0,m}$
(where $m=n+d(r+1)$) are
passed to $\M_{0,n}(\proj^r,d)$ by the local quotient construction. 
The boundary locus of  
$\M_{0,m}$ is a divisor with normal
crossings. Since $\M_{0,n}(\proj^r,d, \barr{t})$ is a product of
$\com^*$-bundles over an open set of $\M_{0,m}$,
the boundary locus of
$\M_{0,n}(\proj^r,d, \barr{t})$ is certainly a divisor with
normal crossing. 
$\M_{0,n}(\proj^r,d)$ is locally the $G_{d,r}$-quotient
of $\M_{0,n}(\proj^r,d, \barr{t})$. 
The boundary of $\M_{0,n}(\proj^r,d)$ is therefore a 
union of subvarieties of pure codimension 1. 
Over the automorphism-free locus, the boundary of
$\M_{0,n}(\proj^r,d)$ is a divisor with normal crossings.

Let $X$ be a nonsingular, projective, convex variety.
The corresponding boundary results for
$\smap_{0,n}(X,\beta)$ are consequences of the
deformation analysis of section \ref{deff}. The boundary
locus of $\M_{0,n}(X,\beta, \barr{t})$ is a divisor
with normal crossing singularities. A pointed
map $\mu:C\rarr X$ such that $C$ has $q$ nodes
lies in the intersection of $q$ branches of the 
boundary. The dimension computation
$$\text{dim}_{\com}
\text{Def}_G(\mu)=
\text{dim}\M_{0,n}(X,\beta, \barr{t}) - q $$
shows these branches intersect transversally at $[\mu]$.
This completes the proof of Theorem 3.
In particular   $\smap_{0,n}(X,\beta)$ has the same
boundary singularity type as $\M_{g}$ and $\M_{g,n}$.

A class $\beta\in H_2(X,\mathbb{Z})$ is {\em effective}
if $\beta$ is represented by some genus $0$ stable map
to $X$.
If $n=0$, the boundary
of $\smap_{0,0}(X, \beta)$ decomposes into a union
of divisors which  
are in bijective correspondence
with effective partitions $\beta_1+\beta_2=\beta$. For general $n$,
the boundary decomposes into a union of divisors
in bijective correspondence with  data of weighted partitions
$(A,B; \beta_1, \beta_2)$ where
\begin{enumerate}
\item[(i)] $A \cup B$ is a partition of $[n]=\{1,2, \ldots, n\}$.
\item[(ii)] $\beta_1+\beta_2=\beta$, $\beta_1$ and $\beta_2$ are effective .
\item[(iii)] If $\beta_1=0$ (resp. $\beta_2=0$), 
then $|A|\geq 2$ (resp.
$|B|\geq 2$).
\end{enumerate}
$D(A,B; \beta_1, \beta_2)$, the divisor
corresponding to the data $(A,B; \beta_1, \beta_2)$, is
defined to be the locus of maps
$\mu:C_A\cup C_B \rarr X$ satisfying the following
conditions:
\begin{enumerate}
\item [(a)] $C$ is a union of two quasi-stable curves $C_A$, $C_B$
of genus $0$ meeting in a point.
\item [(b)] The markings of $A$ (resp. $B$) lie on $C_A$ (resp. $C_B$).
\item [(c)] The map $\mu_A=\mu|_{C_A}$ (resp. $\mu_B$) represents 
$\beta_1$
(resp. $\beta_2$).
\end{enumerate}
The deformation results of section 5 show the
locus maps 
satisfying (a)$-$(c) and $C_A\eqq C_B \eqq \proj^1$
is dense in $D(A,B;\beta_1, \beta_2)$.
If $X= \proj^r$, then it is easily seen
that $D(A,B; \beta_1, \beta_2)$
is irreducible. 
In general, we do not claim the divisor 
$D(A,B; \beta_1, \beta_2)$ is  irreducible, although
that is the case in all the examples we have seen.

\subsection{Boundary divisors}
The boundary divisor of $\M_{0,n}$ corresponding to the 
marking partition $A\cup B=[n]$ is 
naturally isomorphic (by gluing) to the product
$$\M_{0,A \cup \{\bullet\}} 
\times \M_{0,B \cup \{\bullet\}}.$$ An  analogous construction exists
for the boundary divisor $D(A, B; \beta_1, 
\beta_2)$ of the space $\M_{0,n}(X,\beta)$.

Let $K=D(A, B; \beta_1, \beta_2)$ 
be a boundary divisor of $\M_{0,n}(X,\beta)$. 
Let $\M_A=\M_{0, A\cup\{\bullet\}}(X,\beta_1)$ and 
$\M_B= \M_{0, B\cup \{\bullet\}}(X, \beta_2)$.
Let $e_A: \M_A \rarr X$ and $e_B: \M_B\rarr X$ be the
evaluation maps obtained from the additional marking $\bullet$.
Let $\tau_A$, $\tau_B$ be the projections of $\M_A \times \M_B$
to the first and second factors respectively.
Let $\tl{K}= \M_A \times_{X} \M_B$ be the fiber
product with respect to the evaluation maps $e_A$, $e_B$.
$\tl{K}\subset \M_A \times_{\com} \M_B$ is the closed
subvariety 
$(e_A\times_{\com} e_B)^{-1} (\Delta)$ 
where $\Delta\subset X \times X$ is
the diagonal. 

Properties of $\tl{K}$ can be deduced from the 
local quotient constructions of $\M_A$ and $\M_B$.
It will be shown that
$\tl{K}$ is a normal
projective variety of pure dimension with finite quotient singularities. 
Let  $\M_A(X,\barr{t}_A)$, $\M_B(X,\barr{t}_B)$ be the
$\barr{t}_A$, $\barr{t}_B$-rigid  moduli spaces.
$\tl{K}$ is the $G_A \times G_B$-quotient of the 
corresponding subvariety 
$$\tl{K}(X,\barr{t}_A, \barr{t}_B)\subset\M_A(X,\barr{t}_A)\times
\M_B(X,\barr{t}_B),$$  
$$\tl{K}(X,\barr{t}_A, \barr{t}_B)=(e_A\times_{\com} e_B)^{-1}
 (\Delta).$$
The differential of $e_A$ at a point $[\mu]$ 
of $\M_A(X,\barr{t}_A)$ is determined in the following manner.
The case in which the domain $C\eqq\proj^1$ is irreducible
is most straightforward. Then, there are  natural
linear maps:
\begin{equation}
\label{hattey}
\text{Def}(\mu) \rarr H^0(\mu^*T_X/ T_C(-p_\bullet)) \rarr
T_X({\mu(p_\bullet)}).
\end{equation}
The first map in (\ref{hattey}) is
the natural surjection of $\text{Def}(\mu)$
onto the deformation space
of the moduli problem obtained by forgetting
all the markings except $\bullet$.
The natural fiber evaluation
$H^0(\mu^* T_X) \rarr T_X({\mu(p_\bullet)})$
is well defined on the space 
$H^0(\mu^*T_X/ T_C(-p_\bullet))$. This is
the second map in (\ref{hattey}).
The composition of maps in (\ref{hattey})
is simply the differential of $e_A$ at $[\mu]$.
Since $\mu^*T_X$ is
generated by global sections by Lemma \ref{redd}, it follows that the
differential of $e_A$ is surjective at $[\mu]$. A similar
argument shows the differential of $e_A$ is surjective
for each $[\mu] \in \M_A(X,\barr{t}_A)$. The differential
of $e_B$ is therefore also surjective. The surjectivity of
the differentials of $e_A$ and $e_B$ imply 
$\tl{K}(X,\barr{t}_A, \barr{t}_B)$ is nonsingular.
Thus
$\tl{K}$ is a normal
projective variety of pure dimension
with finite quotient singularities.

By gluing the universal families over
$\M_A(\barr{t}_A)$ and $\M_B(\barr{t}_B)$ 
along the markings $\bullet$,
a natural family of Kontsevich stable maps 
exists over $\tl{K}(X,\barr{t}_A,
\barr{t}_B)$. The induced map
$$\tl{K}(X,\barr{t}_A, \barr{t}_B) \rarr K$$
is seen to be $G_A\times G_B$ invariant. Therefore, a natural
map $\psi:\tl{K} \rarr K$ is obtained.  

\begin{lm} Results on the morphism $\psi$ :
\label{resy}
\begin{enumerate}
\item[(i)]
If $A\neq \emptyset$ and $B \neq \emptyset$, then $\psi:\tl{K}\rarr K$ 
is an isomorphism. 
\item[(ii)] If $A\neq \emptyset$, or $B\neq \emptyset$, or 
$\beta_A\neq \beta_B$, then
$\psi$ is birational. 
\item[(iii)] If $A=B=\emptyset$ ($n=0$) and $\beta_A=\beta_B= \beta/2$ then
$\psi$ is generically $2$ to $1$.
\end{enumerate}
\end{lm}

\bpf
First part (i) is proven.
Let $q_A \in A$ and $q_B \in B$ be fixed markings
(whose existence is guaranteed by the assumptions of (i)).
Let $\cal{L}$ be a very ample line bundle on $X$ against
which all degrees of maps are computed. Let $d_A$, $d_B$ be
the degrees of $\beta_A$, $\beta_B$ respectively.
Let $K=(A\cup B, \beta_A, \beta_B)$.
Let $\mu: C\rarr \proj^r$ correspond to a
moduli point $[\mu] \in K$. Let $C= \bigcup C_i$ be the
union of irreducible components. Let
$q_A \in C_1$, $q_B\in C_l$ where $1\neq l$ and let
$$C_1, C_2, \ldots, C_l$$
be the unique minimal path from $C_1$ to $C_l$ which
exists since $C$ is a
tree of components. 
For $1\leq i \leq l-1$, let $x_i=C_i\cap C_{i+1}$. 
Each node $x_i$ divides $C$ into two connected curves
$$C= C_{A,i} \cup C_{B,i}$$
labeled by the points $q_A$, $q_B$.
Let $d_i$ be the degree of $\mu$ restricted 
to $C_{A,i}$. 
The degrees $d_i$ increase monotonically. Since $[\mu] \in K$,
$d_i=d_A$ for some $i$. Let $j$ be the minimal value satisfying
$d_j=d_A$. If $d_{j+1} > d_A$, then $\psi^{-1}[\mu]$ is the
unique point determined by cutting at the node $x_j$. 
If $d_{j+1}=d_A$, then the subcurve 
$$C \ \smallstm \ (C_{A,j} \cup C_{B,j+1})$$
must contain (by stability) a nonempty set of marked points $P_{j+1}$.
Let $k$ be maximal index satisfying $d_{j+k}=d_A$.
The analogously defined marked point sets 
$$P_{j+1}, \ldots, P_{j+k}$$
are all nonempty.
There must be a index $t$
satisfying
$P_{j+t'} \subset A$ for $1\leq t' \leq t$ and
$P_{j+t'} \subset B$ for $t < t'\leq k$. $\psi^{-1}[\mu]$ is then
the 
unique point determined by cutting at the node $x_{j+t}$.
Therefore, $\psi$ is bijective in case $A$ and $B$ are
nonempty.

Let $\M_{0,n}(X,\beta, \barr{t})$ be a locally rigidified moduli
space containing
the point $[\mu]\in K$. If $|A|, |B| \geq 1$, a similar argument
shows the boundary components of $\M_{0,n}(X,\beta,\barr{t})$ lying over
$K$ are {\em disjoint}. Therefore, $K$ is normal. In case
$A$ and $B$ are nonempty, $\psi$ is a bijective morphism of 
normal
varieties and  hence an isomorphism.

Note, for example, that the component $K=D(\emptyset, \emptyset; 2,3)$
of $\M_{0,0}(\proj^r,5)$ is not normal. $K$ intersects itself
along the codimension 2 locus of moduli points $[\mu]$
of the form:
$$\mu: C_1\cup C_2 \cup C_3 \rarr \proj^r$$
with restricted degrees $d_1=2$, $d_2=1$, $d_3=2$.
In this case, $\psi: \tl{K} \rarr K$ is a normalization.

Parts (ii) and (iii) follow simply from the
defining properties (a)$-$(c) of $K$.
\epf

The fundamental relations among the Gromov-Witten invariants
will come from the following linear equivalences among 
boundary components in $\overline{M}_{0,n}(X, \beta)$.

\begin{pr}
\label{needit}
For  $i, j, k, l$  distinct in  $[n]$,  set
$$
D(i,j \mid k,l)  =  \sum  D(A,B; \beta_1,\beta_2) ,
$$
the sum over all partitions such that  $i$  and  $j$  are in  $A$,  $k$  
and  $l$  are in  $B$,  and  $\beta_1$  and  $\beta_2$  are 
effective classes in  
$A_1X$  such that  $\beta_1 + \beta_2 = \beta$.  Then, we have 
the linear equivalence of divisors
$$
D(i,j \mid k,l)  \sim  D(i,l \mid j,k) 
$$
on  $\M_{0,n}(X,\beta)$.
\end{pr}
\bpf The proof is obtained by examining the map  
$$
\M_{0,n}(X,\beta)  \rarr  \M_{0,n}  \rarr  \M_{0,\{i,j,k,l\}} \cong 
\proj^1 ,  
$$
and noting that  
the divisor $D(i,j \mid k,l)\subset \M_{0,n}(X,\beta)$ 
is the multiplicity-free inverse 
image of the  point  $D(i,j \mid k,l)\in \M_{0,\{i,j,k,l\}}$.
The deformation methods of section 5 can be used
to prove that the inverse image of the point
$D(i,j \mid k,l) \in \M_{0,\{i,j,k,l\}}$ is multiplicity-free.
Since points are linearly equivalent on $\proj^1$, the linear
equivalence on $\M_{0,n}(X,\beta)$ is established.
\epf

%
%
%
%
%
%

\section{{\bf Gromov-Witten invariants}}
In sections 7--10, unless otherwise stated, $X$ will denote 
a homogeneous variety and the genus $g$ will be zero.
Since the the tangent bundle of $X$ is generated
by global sections, $X$ is convex. The  moduli
spaces $\overline{M}_{0,n}(X, \beta)$ are therefore
available with the properties proved in sections 1--6.
In addition, the cohomology of $X$ has a natural basis of
algebraic cycles (classes of Schubert varieties), so
$A^i X=H^{2i} X$ can be identified with the Chow group
of cycle classes of codimension $i$. The effective
classes $\beta$ in $A_1 X$ (see section 6.1) are non-negative
linear combinations of the Schubert classes of dimension 1.
Each $1$-dimensional Schubert class is represented by
an embedding $\proj^1\subset X$.

The varieties  $\M_{0,n}(X,\beta)$  come equipped with  $n$  
morphisms  $\rho_1, \ldots , \rho_n$  to  $X$,  where  $\rho_i$  
takes the point  $[C, p_1, \ldots , p_n, \mu] 
\in \M_{0,n}(X,\beta)$  to 
the point  $\mu(p_i)$  in  $X$.  Given arbitrary classes  $\gamma_1, 
\ldots , \gamma_n$  in  $A^*X$,  we can construct the cohomology 
class  $${\rho_1}^*(\gamma_1) \smallcup \cdots \smallcup 
{\rho_n}^*(\gamma_n)$$  on  $\M_{0,n}(X,\beta)$,  and we can 
evaluate its homogeneous component of the top codimension on the 
fundamental class, to produce a number, called a {\em Gromov-Witten 
invariant}, that we denote by 
$I_{\beta}(\gamma_1 \cdots  \gamma_n)$:
\begin{equation}
I_{\beta}(\gamma_1 \cdots  \gamma_n) \, =  
\, \int_{\M_{0,n}(X,\beta)} {\rho_1}^*(\gamma_1) \smallcup \cdots 
\smallcup {\rho_n}^*(\gamma_n) .
\end{equation}
If the classes  $\gamma_i$  are homogeneous, this will be a nonzero 
number only if the sum of their codimensions is the dimension of  
$\M_{0,n}(X,\beta)$,  that is,  $$\sum \text{codim} (\gamma_i) = \dim X + 
\int_\beta c_1(T_X) + n - 3.$$  It follows from the definition that  
$I_{\beta}(\gamma_1  \cdots  \gamma_n)$  is 
invariant under permutations of the classes  $\gamma_1, \ldots ,  
\gamma_n$. 

The conventions of [K-M] require $n\geq 3$. However, 
it will be convenient for us to take $n\geq 0$. 
A 0-pointed invariant occurs
when the moduli space $\M_{0,0}(X, \beta)$ is of dimension
0. In this case $I_\beta= \int _{\M_{0,0}(X, \beta)} 1$.
By Lemma \ref{george}, 
$\M_{0,0}(X,\beta)$ is of dimension 0
if and only if $\dim (X)=1 $ and $\int_\beta c_1(X)=2$. Hence,
for homogeneous varieties,
0-pointed invariants only occur on $X\eqq \proj^1$.
In this case, $I_{1}=1$ is the unique 0-pointed invariant.

Let $M_{0,n}^*(X, \beta)=M_{0,n}(X, \beta)\cap
\overline{M}^*_{0,n}(X,\beta)$.
We start with a simple lemma.
\begin{lm}
\label{freddy}
If $n\geq 1$, then
$M_{0,n}^*(X, \beta)\subset \M_{0,n}(X, \beta)$ is
a dense open set. 
\end{lm}
\bpf
If $\beta=0$, then $\overline{M}_{0,n}(X,0)$ is nonempty
only if 
$n\geq 3$. The equality
$\M_{0,n}^*(X, 0)=\M_{0,n}(X,0)$ is deduced from the
corresponding equality for $\M_{0,n}$. Assume $\beta\neq 0$.
By Theorem 3, 
$M_{0,n}(X, \beta)\subset \M_{0,n}(X, \beta)$ is a dense
open set.  
Let $(\proj^1, \{p_i\}, \mu)$ be a point in
$M_{0,n}(X,\beta)$.
It suffices to show that $(\proj^1, \{p_i'\}, \mu)$ is
automorphism-free for general points $p_1', \ldots, p_n'\in \proj^1$.        
The automorphism group $A$ of the unpointed map $\mu: \proj^1 \rarr
X$ is finite since $\beta \neq 0$.
There exists a (nonempty) open set of $\proj^1$
consisting of points with trivial $A$-stabilizers.
If $p_1', \ldots, p_n'$ belong to this open subset,
the pointed map $(\proj^1, \{p_i'\}, \mu)$ is 
automorphism-free.
\epf

Let $X=G/P$, so $G$ acts transitively on $X$.
Let $\Gamma_1, \ldots, \Gamma_n$ be pure dimensional
subvarieties of $X$.
Let  
$[\gamma_i]\in A^*X$  be the corresponding classes
(see our notational conventions in section \ref{nota}).
Assume 
$$\sum_{i=1}^n \text{codim}(\Gamma_i) = 
\text{dim}(X)+\int_{\beta}c_1(T_X) +n-3.$$
Let $g\Gamma_i$ denote the $g$-translate of $\Gamma_i$
for $g\in G$.
\begin{lm}
\label{enuuuu}
Let $n\geq 0$.
Let $g_1, \ldots, g_n \in G$ be general elements.
Then, the scheme theoretic intersection
\begin{equation}
\label{intttt}
\rho_1^{-1}(g_1 \Gamma_1) \cap \cdots \cap
\rho_n^{-1}(g_n \Gamma_n)
\end{equation}
is a finite number of reduced points supported
in $M_{0,n}(X, \beta)$ 
and
$$I_{\beta}(\gamma_1 \cdots \gamma_n) = \# \
\rho_1^{-1}(g_1 \Gamma_1) \cap \cdots \cap
\rho_n^{-1}(g_n \Gamma_n).$$
\end{lm}
\bpf
If $n=0$, $I_1=1$ on $\proj^1$ is the only case and
the lemma holds since $\overline{M}_{0,0}(\proj^1,1)$
is a nonsingular point. Assume $n\geq 1$.
$M_{0,n}^*(X, \beta) \subset
\overline{M}_{0,n}(X, \beta)$ is a dense open set by 
Lemma \ref{freddy}. By simple
transversality arguments (with respect to the $G$-action),
it follows that the intersection (\ref{intttt}) is supported in
${M}_{0,n}^*(X,\beta)$. By Theorem 2,
 ${M}_{0,n}^*(X,\beta)$ is nonsingular.
An application of Kleiman's Bertini theorem ([Kl]) now shows
that the intersection (\ref{intttt}) is a finite set of
reduced points.
To see that the number of
points in (\ref{intttt}) agrees with the intersection number,
consider the fiber diagram:
\begin{equation}
\label{jesdia}
\begin{CD}
\cap_{i=1}^{n} \rho_i^{-1}(g_i \Gamma_i) @>>> \overline{M}\times
\prod_{i=1}^{n} g_i \Gamma_i \\
@VVV @VVV  \\
\overline{M} @>{\iota}>> \overline{M} \times X^n \\
\end{CD}
\end{equation}
where $\overline{M}=\overline{M}_{0,n}(X,\beta)$ and $\iota$
is the graph of the morphism $(\rho_1, \ldots, \rho_n)$.
From (\ref{jesdia}), one sees that
$$\prod_{i=1}^{n}\rho_i^*[g_i \Gamma_i] \cap [\overline{M}]
= \iota^*[\overline{M} \times \prod_{i=1}^{n}g_i \Gamma_i]=
[\cap_{i=1}^{n} \rho_{i}^{-1}(g_i \Gamma_i)]$$
in $A_0(\overline{M})$, which is the required assertion.
\epf 

Lemma \ref{enuuuu} relates the Gromov-Witten
invariants to enumerative geometry. We
see $I_{\beta}(\gamma_1 \cdots \gamma_n)$ equals
the number of pointed maps $\mu$ from $\proj^1$ to $X$
representing the class $\beta \in A_1 X$ 
and satisfying $\mu(p_i) \in g_i \Gamma_i$.
We will need  three basic properties satisfied by
the Gromov-Witten invariants:  

\vspace{+10pt}
 (I)  $\beta = 0$.  
\noindent In this case,  
$\M_{0,n}(X,\beta) = \M_{0,n} \times X$,  and the mappings  $\rho_i$  
are all equal to the projection  $p$  onto the second factor.  
Since  $${\rho_1}^*(\gamma_1) \smallcup \cdots \smallcup  
{\rho_n}^*(\gamma_n) = p^*(\gamma_1  \smallcup \cdots \smallcup 
\gamma_n),$$    
\begin{eqnarray*}       
I_\beta(\gamma_1\cdots  \gamma_n)   
& =  & \int_{\M_{0,n} \times X} 
p^*(\gamma_1 \smallcup \cdots \smallcup \gamma_n)  \\ 
& =  &\int_{p_*[\M_{0,n} \times X]}   
\gamma_1  \smallcup \cdots \smallcup \gamma_n  .
\end{eqnarray*}
Note that $\overline{M}_{0,n}$ is empty if $0\leq n \leq 2$.
If  $n > 3$,  $p_*[\M_{0,n} \times X] = 0$,  since the fibers of  $p$  
have positive dimension.  The only way the number  
$I_\beta(\gamma_1 \cdots  \gamma_n)$  
can be nonzero is when  $n = 3$,  so that  $\M_{0,n}$  is just a point.  
In this case,  $I_\beta(\gamma_1 {\cdot} \gamma_2 
{\cdot} \gamma_3)$  is the classical intersection number  
$\int_X  \gamma_1 \smallcup \gamma_2 \smallcup \gamma_3$. 
 
\vspace{+10pt}  
(II) $\gamma_1 =1 \in A^0X$.    
\noindent If 
$\beta\neq 0$, then the 
product  ${\rho_1}^*(\gamma_1) \smallcup \cdots \smallcup 
{\rho_n}^*(\gamma_n)$  is the pullback of a class on  
$\M_{0,n-1}(X,\beta)$  by the map from  $\M_{0,n}(X,\beta)$  to  
$\M_{0,n-1}(X,\beta)$  that forgets the first point.  Since the fibers 
of this map have positive dimension, the evaluation
$I_\beta(\gamma_1 {\cdots}\gamma_n)$ must vanish.
Therefore, by (I), 
$I_\beta(\gamma_1\cdots  
\gamma_n)$  vanishes unless  $\beta = 0$  and  $n = 3$. In this case,
$I_0( 1{\cdot} \gamma_2 {\cdot} \gamma_3) =  
\int_X \gamma_2 \smallcup \gamma_3$.

\vspace{+10pt}
(III)   $\gamma_1 \in A^1X$ and $\beta\neq 0$.  
\noindent In this case,
\begin{equation}
\label{axxer}
I_\beta(\gamma_1\cdots  \gamma_n)  =  
\left(\int_\beta \gamma_1\right) \cdot I_\beta(\gamma_2 
\cdots  \gamma_n).
\end{equation}
For a map  $\mu: C \rarr X$  with  
$\mu_*[C] = 
\beta$,  there are  $\left( \int_\beta \gamma_1 \right)$  choices 
for the point  $p_1$  in  $C$  to map to a point in  $\Gamma_1$,  
where  $\Gamma_1$  is a variety representing  $\gamma_1$.  
Equation (\ref{axxer}) 
is therefore a consequence of Lemma \ref{enuuuu}.

For a 
formal intersection-theoretic
proof of (\ref{axxer}), consider the mapping
$$
\psi:  \M_{0,n}(X,\beta)   \rarr   X   \times  
\M_{0,n-1}(X,\beta)
$$
which is the product of $\rho_1$  and the map that forgets the first 
point.  
By the K\"unneth formula, we can write
$\psi_*[\M_{0,n}(X,\beta)] = \beta' \times 
[\M_{0,n-1}(X,\beta)] + \alpha$,  where  $\beta'$  is a class in  
$A_1 X$,  and $\alpha$ is some homology
class that is supported over a proper 
closed subset of $\M_{0,n-1}(X,\beta)$.  The class  $\beta'$  can 
be calculated by restricting to what happens over a generic point of 
$\M_{0,n-1}(X,\beta)$.  Representing such a point by  
$(C, p_2, \ldots, p_n, \mu)$ with $C\eqq \proj^1$,  one sees that
the fiber over this point is isomorphic to $C$ and   
$\beta' = \mu_*[C] = \beta$. 
Using the projection formula as in (I) and (II), it 
follows that
\begin{eqnarray*}
I_\beta(\gamma_1 \cdots  \gamma_n)  
& =  & \int_{\beta \times [\M_{0,n-1}(X,\beta)]} \,  \gamma_1  \, 
\times  \, {\rho_2}^*(\gamma_2) \smallcup \cdots \smallcup 
{\rho_n}^*(\gamma_n) \\
& = & \int_\beta \gamma_1   \cdot  
\int_{\M_{0,n-1}(X,\beta)}
{\rho_2}^*(\gamma_2)  \smallcup \cdots \smallcup 
{\rho_n}^*(\gamma_n) ,
\end{eqnarray*}
as asserted.

\vspace{+10pt}
It should be noted that the generic element of
$\overline{M}_{0,0}(X, \beta)$ may not be 
a {\em birational} map of $\proj^1$ to $X$.
This is seen immediately for $X\eqq \proj^1$
where  the generic
element of $\overline{M}_{0,0}(\proj^1, d)$ is
a $d$-fold branched covering of $\proj^1$.
This phenomenon occurs in higher dimensions.
For example, let $X$ be the complete
flag variety $\mathbf{Fl}(\com^3)$ 
(the space of 
pairs $(p,l)$ satisfying $p\in l$ where
$p$ and $l$ are a point and a line in 
$\proj^2$).
Let $\beta \in A_1 \mathbf{Fl}(\com^3)$ be the
class of the curve $\proj^1 \subset \mathbf{Fl}(\com^3)$
determined by all pairs $(p,l)$ for a fixed line $l$.
One computes 
$\int_{\beta} c_1(T_{\mathbf{Fl}(\com^3)})=2$, so
the dimension of $\M_{0,0}(\mathbf{Fl}(\com^3), \beta)$
is $3+2-3=2$ 
by  Theorem 2. Directly, one sees that
$\M_{0,0}(\mathbf{Fl}(\com^3), \beta)$ is isomorphic to the
space of lines in $\proj^2$. In particular, 
$\M_{0,0}(\mathbf{Fl}(\com^3), \beta)$
has no boundary.
As in the case of $\proj^1$,
it is seen that every element of
$M_{0,0}(\mathbf{Fl}(\com^3), 2\beta)$ corresponds 
to a double cover of an element of
$\M_{0,0}(\mathbf{Fl}(\com^3), \beta)$. The
boundary of $\M_{0,0}(\mathbf{Fl}(\com^3), 2\beta)$
consists of degenerate double covers.
Note also that every 
element of $\M_{0,0}(\mathbf{Fl}(\com^3), 2\beta)$
has a nontrivial automorphism. 
Since the space of image curves
of maps in $\M_{0,0}(\mathbf{Fl}(\com^3), 2\beta)$
is only $2$-dimensional, it follows that all Gromov-Witten
invariants of $\mathbf{Fl(\com^3)}$ of the form
$I_{2\beta}(\gamma_1\cdots \gamma_n)$
vanish. 
%
%
%
%

\section{{\bf Quantum cohomology }}
We keep the notation of section 7.
Let  $T_0 = 1 \in A^0X$,  let  $T_1, \ldots , T_p$  be a basis of  
$A^1X$,  and let  $T_{p+1}, \ldots , T_m$  be a basis for the other 
cohomology groups.  The classes of Schubert varieties form the natural
basis for homogeneous varieties.
The fundamental numbers counted by the Gromov-Witten  
invariants are the numbers
\begin{equation}
\label{idad}
N(n_{p+1}, \ldots , n_m; \beta)  =  I_\beta({T_{p+1}}^{n_{p+1}} 
\cdots   {T_m}^{n_m})
\end{equation}
for $n_i\geq 0$.
The invariant (\ref{idad}) 
is nonzero only when $\sum n_i \left(\text{codim}
(T_i) - 1 
\right) = \dim X + \int_\beta c_1(T_X) \, - \, 3$.  
In this case, it is the number of pointed rational maps  meeting  $n_i$  
general representatives of  $T_i$  for each  $i$,  $p+1 \leq i 
\leq m$.

Define the numbers  $g_{i j},  0 \leq i, j \leq m$,  by the 
equations
\begin{equation}
g_{i j}  =  \int_X  T_i \smallcup T_j .
\end{equation}
(If the $T_i$ are the Schubert classes, then for
each $i$ there is a unique $j$ such that $g_{ij}\neq 0$.
For this $j$, $g_{ij}=1$.)

Define  $\left( g^{i j}\right)$  to be the inverse matrix to the 
matrix $\left( g_{i  j}\right)$.  Equivalently, the class of the 
diagonal  $\Delta$  in  $X \times X$  is given by the formula
\begin{equation}
[\Delta]  =  \sum_{e\, f} g^{e f} \, T_e \otimes T_f  
\end{equation}
in $A^*(X \times X) = A^*X \otimes A^*X$.
The following equations hold:
\begin{equation}
T_i \smallcup T_j  =  \sum_{e, \, f}  \left( \int_X T_i \smallcup T_j 
\smallcup T_e \right) g^{e  f} T_f  =  \sum_{e, \, f} I_0(T_i 
{\cdot} T_j {\cdot} T_e ) 
g^{e  f} T_f . \label{2.3}
\end{equation}
        
The idea is to define a ``quantum  deformation'' of 
the cup multiplication of (\ref{2.3}) by 
allowing nonzero classes  $\beta$.  Here enters a key idea from 
physics -- to write down a ``potential function'' that carries all the 
enumerative information.  

Define, for a class  $\gamma$  in  $A^*X$,  
\begin{equation}
\label{2.4}
\Phi(\gamma)  =  \sum_{n \geq 3} \sum_\beta \frac{1}{n!} 
I_{\beta}(\gamma^n) ,
\end{equation}
where  $\gamma^n$  denotes  $\gamma  \cdots 
\gamma$  ($n$  times).  
\begin{lm}
\label{uunet}
For a given  integer $n$,  
there are only finitely many effective classes  $\beta\in A_1 X$  
such that  
$I_\beta(\gamma^n)$  is not zero. 
\end{lm}
\bpf
Since $X$ is a homogeneous space, the effective
classes in $A_1 X$ are the non-negative
linear combination of finitely many (nonzero) effective
classes $\beta_1, \ldots, \beta_p$. By Lemma
\ref{george}, $\int_{\beta_i} c_1(T_X)\geq 2$.
Hence, for a given
integer $N$, there are only a finite 
number of  effective $\beta$  for which $\int_\beta c_1(T_X) \leq N$.
If $I_\beta(\gamma^n)$ is nonzero, then
$$\text{dim} \\M_{0,n}(X, \beta)  \leq  n \cdot \text{dim}\ X$$
which implies
that  $\int_\beta c_1(T_X) \leq (n-1) \cdot \text{dim}\ X +3 -n$.
\epf

\noindent  Let
$\gamma = \sum y_i \, T_i$. By
Lemma \ref{uunet},  $\Phi(\gamma) = \Phi(y_0, 
\ldots , y_m)$  becomes a formal power series in  
$\mathbb Q[[y]] = \mathbb Q[[y_0, \ldots , y_m]]$:
\begin{equation}
\Phi(y_0, \ldots , y_m)  =  \sum_{n_0 + \ldots + n_m \geq 3}
\sum_\beta I_\beta({T_0}^{n_0} \cdots  
{T_m}^{n_m} ) \frac{{y_0}^{n_0}}{n_0!} \cdots 
\frac{{y_m}^{n_m}}{n_m!}  \label{2.5} .
\end{equation}

Define   $\Phi_{i  j  k}$  to be the partial derivative:
\begin{equation}
\Phi_{i  j  k} = \frac{\partial^3\Phi}{\partial y_i \, \partial y_j
\partial y_k} \; , \; \;  0 \leq i, j, k \leq m .
\end{equation}
A simple formal calculation, using (\ref{2.5}), gives the following 
equivalent formula:
\begin{equation}
\Phi_{i  j  k} = \sum_{n \geq 0}  \sum_\beta   \frac{1}{n!} \,
I_\beta(\gamma^n \cdot {T_i {\cdot} T_j {\cdot} T_k}) . 
\label{2.7}
\end{equation}  

 Now we define a new ``quantum'' product  $*$  by the rule:
\begin{equation}
T_i \, * \,  T_j  = \sum_{e, \, f}\Phi_{i  j  e} \, g^{e  f} \,  T_f . 
\label {2.8}
\end{equation}
The product in (\ref{2.8}) is extended 
$\mathbb Q[[y]]$-linearly to the 
$\mathbb Q[[y]]$-module  
$A^*X\otimes_{\mathbb Z}\mathbb Q[[y]]$,  thus 
making it a  $\mathbb Q[[y]]$-algebra.  One thing is evident from this 
remarkable definition: this product is commutative, since the partial 
derivatives are symmetric in the subscripts.  

It is less obvious, but not difficult, to see   $T_0 = 1$  is a unit 
for the $*$-product.  In fact, it follows from property (I) 
of section 7, together with (\ref {2.7}), that  
$$\Phi_{0  j  k}  \, = \,  I_0(T_0 {\cdot} T_j {\cdot} T_k)  
\, = \,  \int_X  T_j \smallcup T_k  \, = \, g_{j k} ,$$
and from this we see that  $T_0 * T_j = \sum g_{j e} \, g^{e f} \, 
T_f = T_j$.

The essential point, however, is the associativity:

\begin{tm} 
\label{assss}
This definition makes   $A^*X\otimes\mathbb Q[[y]]$  into a 
commutative, associative  $\mathbb Q[[y]]$-algebra, with unit  $T_0$.  
\end{tm}

We start the proof by writing down what associativity says:
$$
(T_i \, * \, T_j) \, * \, T_k  \, = \,  \sum_{e, \, f} \, \Phi_{i  j  e} 
\,  g^{e f} \, T_f \, * \, T_k  \, = \,  \sum_{e, \, f} \, \sum_{c, \, 
d}\Phi_{i  j  e} \, g^{e\, f} \, \Phi_{f  k  c} \, g^{c d} \, T_d , 
$$
$$
T_i \, * \, (T_j \, * \, T_k)  \, = \,  \sum_{e, \, f} \, \Phi_{j  k  e} 
\,  g^{e f} \, T_i \, * \, T_f  \, = \,  \sum_{e, \, f} \, \sum_{c, \, 
d}\Phi_{j k e} \, g^{e f} \, \Phi_{i  f  c} \, g^{c d} \, T_d .
$$ 
Since the matrix  $\left( g^{c\, d} \right)$  is nonsingular, the 
equality of  $(T_i \, * \, T_j) \, * \, T_k$  and  $T_i \, * \, (T_j \, * 
\, T_k)$  is equivalent to the equation
$$
\sum_{e, \, f}\Phi_{i j e} \, g^{e f} \, \Phi_{f k  l} \, = \, 
\sum_{e, \, f}\Phi_{j k e} \, g^{e  f} \, \Phi_{i  f  l}
$$
for all  $l$.  If we set  
\begin{equation}
\label{rrrr}
F(i,j  \mid  k,l) \,  = \, \sum_{e, \, f}\Phi_{i  j  e} \, g^{e  f} \, 
\Phi_{f k  l} ,  
\end{equation}
and use the symmetry $\Phi_{i  f  l}  = \Phi_{f  i  l}$,  we see 
that the associativity is equivalent to the equation
\begin{equation}
F(i,j  \mid  k,l) \,  =  \, F(j,k  \mid  i,l) . \label {2.9}
\end{equation}
        
It follows from (\ref {2.7}) that
\begin{equation}
F(i,j  \mid  k,l)  \,  = \,   \sum \frac{1}{n_1! \, n_2!}
I_{\beta_1}(\gamma^{n_1} {\cdot} T_i {\cdot} T_j 
{\cdot} T_e) \, g^{e f} \,  I_{\beta_2}(\gamma^{n_2} 
{\cdot} T_k {\cdot} T_l {\cdot} T_f) , \label {2.10}
\end{equation}
where the sum is over all nonnegative  $n_1$  and  $n_2$,  over all  
$\beta_1$  and  $\beta_2$  in  $A_1X$,  and over all  $e$  and  $f$  
from  $0$  to  $m$.  We need the following lemma.  Recall from 
section 6, the divisor $D(A,B;\beta_1, \beta_2)$.
In case $A$ and $B$ are nonempty,
$$
D(A,B;\beta_1,\beta_2)  =  \M_{0,A\cup\{\bullet\}}(X,\beta_1) \, 
\times_X \, \M_{0,B\cup\{\bullet\}}(X,\beta_2) .
$$
\begin{lm}  
Let  $\iota$  denote the natural inclusion of  
$D(A,B;\beta_1\beta_2)$  in the Cartesian product  
$\M_{0,A\cup\{\bullet\}}(X,\beta_1) \, \times \, 
\M_{0,B\cup\{\bullet\}}(X,\beta_2)$,  and let  $\alpha$  be the 
embedding of  $D(A,B;\beta_1, \beta_2)$  as a divisor in  
$\M_{0,n}(X,\beta)$, with
$\beta=\beta_1+\beta_2$.  Then for any classes  $\gamma_1, \ldots , 
\gamma_n$  in  $A^*X$,  
$$ 
\iota_* \circ \alpha^* ({\rho_1}^*(\gamma_1) \smallcup \cdots 
\smallcup {\rho_n}^*(\gamma_n)) = $$
$$ \sum_{e, \, f}  g^{e  f} 
\left( \prod_{a \in A} {\rho_a}^*(\gamma_a)  {\cdot} 
{\rho_\bullet}^*(T_e) \right) \, \times \, \left( \prod_{b \in B} 
{\rho_b}^*(\gamma_b) {\cdot} {\rho_\bullet}^*(T_f) \right)  .
$$
\label{toast}
\end{lm}
\bpf
Let  $M_1 = \M_{0,A\cup\{\bullet\}}(X,\beta_1)$,  $M_2 
= M_{0,B\cup\{\bullet\}}(X,\beta_2)$,  $M = \M_{0,n}(X,\beta)$,  and  
$D = D(A,B;\beta_1\beta_2)$.  From the identification of  $D$  with  
$M_1 \, \times_X \, M_2$,  we have a commutative diagram, with 
the right square a fiber square:
\begin{equation}
\begin{CD}
M @<{\alpha}<< D @>{\iota}>> M_1 \times M_2 \\
@V{\rho}VV @V{\eta}VV @VV{\rho'}V \\
X^n @<<{p}< X^{n+1} @>>{\delta}> X^{n+2} \\
\end{CD}
\end{equation}

Here  $\rho$  is the product of the 
evaluation maps denoted  $\rho_i$,  $\rho'$  
is the product of maps  $\rho_i$  and the two others denoted  
$\rho_\bullet$,  $\delta$  is the diagonal embedding that repeats the 
last factor, and  $p$  is the projection that forgets the last factor.  
Then we have
\begin{eqnarray*}
\iota_* \circ \alpha^* ({\rho_1}^*(\gamma_1) \smallcup \cdots 
\smallcup {\rho_n}^*(\gamma_n)) & = &
\iota_* \circ \alpha^* \circ \rho^* (\gamma_1 \times \cdots \times 
\gamma_n) \\
& = &  \iota_* \circ \eta^* \circ p^*(\gamma_1 \times \ldots 
\times \gamma_n)\\
& = & \iota_* \circ \eta^*(\gamma_1 \times 
\ldots \times \gamma_n \times [X]) \\
& = & \rho'{}^* \circ \delta_*(\gamma_1 \times \ldots \times 
\gamma_n \times [X]) \\
& = &  \rho'{}^* (\gamma_1 \times \ldots 
\times \gamma_n \times [\Delta]) \\
& = & \sum_{e, \, f} \, g^{e  f}\rho'{}^*(\gamma_1 \times \ldots 
\times \gamma_n \times T_e \times T_f)
\end{eqnarray*}
$$=  \sum_{e, \, f} \, g^{e  f} \left( \prod_{a \in A} 
{\rho_a}^*(\gamma_a) {\cdot} {\rho_\bullet}^*(T_e) \right)
\, \times \, \left( \prod_{b \in B} {\rho_b}^*(\gamma_b) {\cdot} 
{\rho_\bullet}^*(T_f) \right) .$$
\epf

Fix  $\beta\in A_1 X$ and  $\gamma_1, \ldots , \gamma_n \in A^*X$,  
and 
fix four distinct integers  $q, r, s$,  and  $t$  in  $[n]$.  Set
\begin{equation}
G(q,r  \mid  s,t)  \, = \, \sum g^{e  f} I_{\beta_1} \left( \prod_{a \in 
A} \gamma_a {\cdot} T_e \right) \cdot I_{\beta_2}
\left( \prod_{b \in B} \gamma_b {\cdot} T_f \right) ,
 \label{2.11}
\end{equation}
where the sum is over all partitions of $[n]$  into two sets  $A$  and  
$B$  such that  $q$  and  $r$  are in  $A$  and  $s$  and  $t$  are in  
$B$,  and over all  $\beta_1$  and  $\beta_2$  that sum to  $\beta$,  
and over  $e$  and  $f$   between  $0$  and  $m$.
It follows from Lemma \ref{toast} that  
$$
G(q,r  \mid  s,t) =  \sum  \int_ {D(A,B;\beta_1,\beta_2)}  
{\rho_1}^*\gamma_1 \smallcup \cdots \smallcup 
{\rho_n}^*\gamma_n ,
$$
the sum over  $A$  and  $B$  and  $\beta_1$  and  $\beta_2$  as 
above.  Now Proposition \ref{needit} from section 7 implies 
\begin{equation}
G(q,r  \mid  s,t)  \, = \,  G(r,s  \mid  q,t) .  \label{2.12}
\end{equation}
Apply (\ref{2.12}) 
in the following case :
$$\gamma_i = \gamma, \ \ \ \text{for} \ \  1 \leq i \leq n - 4, $$   
$$\gamma_{n-3} = T_i, \ \ \gamma_{n-2} = T_j, \ \ \gamma_{n-1} = 
T_k, \ \ \gamma_n = T_l,$$  
$$q = n-3, \  \   r = n-2, \ \  s = 
n-1, \  \ t = n.$$  Then 
(\ref{2.11}) becomes
$$
G(q,r  \mid  s,t)  =  \sum \binom {n-4}{n_1 - 2} g^{e f}   
I_{\beta_1}(\gamma^{n_1-2} {\cdot} T_i {\cdot} T_j 
{\cdot} T_e) \cdot  I_{\beta_2}(\gamma^{n_2-2} {\cdot} 
T_k {\cdot} T_l {\cdot} T_f) ,
$$
the sum over  $n_1$  and  $n_2$,  each at least  $2$,  adding to  $n$,  
and  $\beta_1$  and  $\beta_2$  adding to  $\beta$;  the binomial 
coefficient is the number of partitions  $A$  and $B$  for which  $A$  
has  $n_1$  elements, and  $B$  has  $n_2$  elements.  This can be 
rewritten
\begin{equation}
G(q,r  \mid  s,t)  =  n! \sum \frac{1}{n_1! \, n_2!} 
g^{ef}I_{\beta_1}(\gamma^{n_1} {\cdot} T_i {\cdot} T_j 
{\cdot} T_e) \cdot  I_{\beta_2}(\gamma^{n_2} {\cdot} T_k 
{\cdot} T_l {\cdot} T_f) ,   \label{2.13}
\end{equation}
the sum over nonnegative  $n_1$  and  $n_2$  adding to  $n-4$,  and  
$\beta_1$  and  $\beta_2$  adding to  $\beta$.
        
The required equality (\ref{2.9}) then follows immediately from 
(\ref{2.12}) and (\ref{2.13}), together with (\ref{2.10}).   This 
completes the proof of Theorem \ref{assss}.

While
the definition of the quantum cohomology ring depends upon
a choice of basis $T_0, \ldots, T_m$ of $A^* X$,
the rings obtained from different basis choices
are canonically isomorphic.
The variables $y_0, \ldots, y_m$ should be identified
with the dual basis to $T_0, \ldots, T_m$. If
$T_0', \ldots, T_m'$ is another basis of $A^* X$ and
$T_i'= \sum a_{ij} T_j$ is the change of coordinates, let
\begin{equation}
\label{duell}
y_i = \sum a_{ji}y'_j
\end{equation}
be the dual coordinate change.
Relation (\ref{duell})
yields an isomorphism of $\Q$-vector spaces
$$A^* X \otimes \Q[[y]] \eqq A^*X \otimes \Q[[y']].$$
It is easy to check that the quantum products
defined respectively on the left and right 
by the $T$ and $T'$ bases
agree with this identification. 

Let $V$ denote the underlying free abelian
group of $A^* X$.
Let $\Q[[V^*]]$ be the completion of the
graded polynomial ring 
$\bigoplus_{i=0}^{\infty} Sym^i(V^*)\otimes \Q$
at the unique maximal graded ideal.
The quantum product defines
a canonical ring structure on the free $\Q[[V^*]]$-module
$V \otimes_\Z \Q[[V^*]]$. 
Let $QH^* X = (V \otimes_\Z \Q[[V^*]], *)$ denote
the quantum cohomology ring.
There is a canonical injection of abelian groups
$$\iota: A^* X \hookrightarrow QH^* X$$ 
determined by $\iota(v)=v \otimes 1$ for $v\in V$.
The injection $\iota$ is {\em not} compatible with
the $\smallcup$ and $*$ products.

It is worth noting that the quantum cohomology ring $QH^* X$
is {\em not} in general a formal deformation of  $A^* X$ over the
local ring $\Q[[V^*]]$. 
It can be seen directly from the definitions that
the $*$-product does not specialize to the $\smallcup$-product
when the formal parameters are set to 0.
At the end of section 9, a presentation
of $QH^* \proj^2$ shows explicitly the difference
between  $A^* \proj^2$ and the specialization
of $QH^* \proj^2$.
In section 10, a ring deformation of $A^* X$ will
be constructed via a smaller quantum cohomology ring.
%
%
%
%

\section{\bf Applications to enumerative geometry}
We write the potential function as a sum:
$$
\Phi(y_0, \ldots , y_m)   =    \Phi_{\text{classical}}(y)  +  
\Phi_{\text{quantum}}(y).
$$
The classical part has the terms for $\beta=0$:
$$
 \Phi_{\text{classical}}(y) =  \sum_{n_0 + \ldots + n_m = 3}
\int_X \left({T_0}^{n_0} \smallcup \cdots \smallcup 
{T_m}^{n_m}\right) \frac{{y_0}^{n_0}}{n_0!} \cdots 
\frac{{y_m}^{n_m}}{n_m!} .
$$ 
Since the associativity equations involve only third
derivatives, we can modify $\Phi$ by any terms of degree
at most 2. Using properties (I)--(III) of section 7, we see
that $\Phi_{\text{quantum}}(y)$  can be replaced by
$\Gamma(y)$: 
$$
\Gamma(y)  = \sum_{n_{p+1} + \ldots + n_m \geq 0} 
\sum_{\beta \neq 0} 
N(n_{p+1}, \ldots , n_m; \beta) \prod_{i = 1}^p e^{\left(\int _\beta 
T_i\right) y_i} \prod_{i = p+1}^m \frac{{y_i}^{n_i}}{n_i!} , 
$$
where  $N(n_{p+1}, \ldots , n_m; \beta)  =  
I_\beta({T_{p+1}}^{n_{p+1}} \cdots 
{T_m}^{n_m})$.  The partial derivatives of  $\Phi_{\text{classical}}$  
involve only the numbers  
$\int_X T_i \smallcup T_j \smallcup T_k$,  while  $\Gamma$  
involves the interesting enumerative geometry numbers.  From this 
form of  $\Gamma$,  it is easy to read off its partial derivatives.

Let us look again at the projective plane from this point of view.  
Take the obvious basis:  $T_0 = 1$,  $T_1$  the class of a line, and  
$T_2$  the class of a point.  Note that  $g_{i  j}$  is  $1$  if  $i + j 
= 2$,  and  $0$  otherwise, so the same is true for  
$g^{i  j}$.  Therefore,
$$
T_i \, * \,  T_j  =  \Phi_{i  j  0} T_2 + \Phi_{i  j  1} T_1 +  
\Phi_{i  j  2} T_0 .
$$
For example,
\begin{eqnarray*}
 T_1 \, * \,  T_1 & = & T_2 + \Gamma_{111} T_1 + 
\Gamma_{112} T_0 , \\ 
T_1 \, * \,  T_2  &= & \Gamma_{1 21} T_1 + \Gamma_{122} T_0 ,\\
T_2 \, * \,  T_2 & = & \Gamma_{221} T_1 + \Gamma_{222} T_0 .
\end{eqnarray*}
Therefore,
$$
(T_1 \, * \,  T_1)  \, * \,   T_2  =  (\Gamma_{221} T_1 + 
\Gamma_{222} T_0) + \Gamma_{111}(\Gamma_{121} T_1 
+ \Gamma_{122} T_0) + \Gamma_{112} T_2 , 
$$
$$
 T_1 \, * \,  (T_1  \, * \,   T_2)  = \Gamma_{121} (T_2 + 
\Gamma_{111} T_1 + \Gamma_{112} T_0) + 
\Gamma_{122} T_1 .
$$
The fact that the coefficients of  $T_0$  must be equal in these last 
two expressions gives the equation:
\begin{equation}
\Gamma_{222}  =  {\Gamma_{112}}^2 - \Gamma_{111} \, 
\Gamma_{122} . \label{2.14}
\end{equation}

If  $\beta = d[\text{line}]$,  
the number  $N(n,\beta)$  is nonzero only when  
$n = 3d-1$,  when it is the number  $N_d$  of plane rational curves 
of degree  $d$  passing through  $3d-1$  general points.  So,
$$
\Gamma(y)  =  \sum_{d \geq 1} N_d    e^{dy_1}  
\frac{{y_2}^{3d-1}}{(3d-1)!} .
$$
From this we read off the partial derivatives:
\begin{eqnarray*}
\Gamma_{222} &  = & \sum_{d \geq 2} N_d     e^{dy_1}   
\frac{{y_2}^{3d-4}}{(3d-4)!}  \\
\Gamma_{112} &  = & \sum_{d \geq 1} d^2  N_d    e^{dy_1} 
\  \frac{{y_2}^{3d-2}}{(3d-2)!}  \\
\Gamma_{111} &  = & \sum_{d \geq 1} d^3  N_d    e^{dy_1} 
\  \frac{{y_2}^{3d-1}}{(3d-1)!} \\
\Gamma_{122} & = & \sum_{d \geq 1} d  N_d    e^{dy_1} 
\  \frac{{y_2}^{3d-3}}{(3d-3)!}  .
\end{eqnarray*}
Therefore,
$$
{\Gamma_{112}}^2  =  \sum_{d \geq 2} \sum_{d_1 + d_2 = 
d}{d_1}^2  N_{d_1}  {d_2}^2  N_{d_2}  e^{dy_1} 
\frac{{y_2}^{3d-4}}{(3d_1-2)! \ (3d_2-2)!} ,$$ 
$$\Gamma_{111} \, \Gamma_{122}  =  \sum_{d \geq 2} 
\sum_{d_1 + d_2 = d}{d_1}^3  N_{d_1}  d_2  N_{d_2} 
e^{dy_1} \frac{{y_2}^{3d-4}}{(3d_1-1)! \ (3d_2-3)!} .
$$
In all these sums,  $d_1$  and  $d_2$  are positive.  Equating the 
coefficients of   $$e^{dy_1} {y_2}^{3d-4}/(3d-4)!,$$  we get the 
identity $(d\geq 2)$:
\begin{equation}
N_d  =  \sum_{d_1 + d_2 = d} N_{d_1}   N_{d_2}   \left[ {d_1}^2 
 {d_2}^2 \ \binom {3d-4}{3d_1-2} - {d_1}^3   d_2 
\ \binom {3d-4}{3d_1-1} \right] . \label{2.15}
\end{equation}
Here a binomial coefficient $\binom n m$ is defined to be zero if any 
of  $n$,  $m$,  or  $n - m$  is negative. 
This is the recursion formula discussed in the introduction.

Note that the quantum formalism has removed any necessity to be 
clever.  One simply writes down the associativity equations, and 
reads off enumerative information.
One can organize the information in these associativity
equations more systematically as follows (see [DF-I]).
Let $F(i,j \mid k,l)$ be defined by (\ref {rrrr}). For
$0 \leq i,j,k,l \leq m$, define:
\begin{equation*}
A(i,j,k,l)  =  F(i,j\mid k,l)- F(j,k \mid i,l) 
\end{equation*}
$$ =  \sum_{e,f} \Phi_{ije}g^{ef}\Phi_{fkl}
                  - \Phi_{jke} g^{ef}\Phi_{fil}.$$
Associativity (Theorem 4) amounts to the
equations $A(i,j,k,l)=0$ for all $i,j,k,l$.
The symmetry of $\Phi_{ijk}$ in the
subscripts and $g^{ef}$ in the superscripts and
the basic facts about $\Phi_{0jk}$ imply:
\begin{enumerate}
\item[(i)]  $A(k,j,i,l) = -A(i,j,k,l)$,
\item[(ii)] $A(l,k,j,i) =  A(i,j,k,l)$,
\item[(iii)] $A(i,j,k,l)=0$ if $i=k$ or $j=l$ or
             if any of the indices $i,j,k,l$ equals 0. 
\end{enumerate}
We consider equations equivalent if they differ by sign.
For distinct $i,j,k,l$, the 24 possible equations
divide into 3 groups of 8. The equation $A(i,j,k,l)=0$
that
says $F(i,j\mid k,l)= F (j,k\mid i,l)$ can
be labelled by  a duality diagram from topological
field theory (see [DF-I]):

\vspace{-0pt}
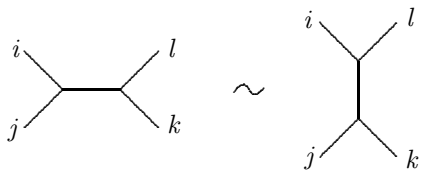
\begin{center}
 
\font\thinlinefont=cmr5
\begingroup\makeatletter\ifx\SetFigFont\undefined
\def\x#1#2#3#4#5#6#7\relax{\def\x{#1#2#3#4#5#6}}%
\expandafter\x\fmtname xxxxxx\relax \def\y{splain}%
\ifx\x\y   
\gdef\SetFigFont#1#2#3{%
  \ifnum #1<17\tiny\else \ifnum #1<20\small\else
  \ifnum #1<24\normalsize\else \ifnum #1<29\large\else
  \ifnum #1<34\Large\else \ifnum #1<41\LARGE\else
     \huge\fi\fi\fi\fi\fi\fi
  \csname #3\endcsname}%
\else
\gdef\SetFigFont#1#2#3{\begingroup
  \count@#1\relax \ifnum 25<\count@\count@25\fi
  \def\x{\endgroup\@setsize\SetFigFont{#2pt}}%
  \expandafter\x
    \csname \romannumeral\the\count@ pt\expandafter\endcsname
    \csname @\romannumeral\the\count@ pt\endcsname
  \csname #3\endcsname}%
\fi
\fi\endgroup
\mbox{\beginpicture
\setcoordinatesystem units <0.40000cm,0.40000cm>
\unitlength=0.40000cm
\linethickness=1pt
\setplotsymbol ({\makebox(0,0)[l]{\tencirc\symbol{'160}}})
\setshadesymbol ({\thinlinefont .})
\setlinear
%
%
\linethickness= 0.500pt
\setplotsymbol ({\thinlinefont .})
\plot  1.905 24.765  3.175 23.495 /
\putrule from  3.175 23.495 to  5.080 23.495
\plot  5.080 23.495  6.350 24.765 /
%
%
\linethickness= 0.500pt
\setplotsymbol ({\thinlinefont .})
\plot  5.080 23.495  6.350 22.225 /
%
%
\linethickness= 0.500pt
\setplotsymbol ({\thinlinefont .})
\plot  3.175 23.495  1.905 22.225 /
%
%
\linethickness= 0.500pt
\setplotsymbol ({\thinlinefont .})
\plot 13.003 22.528 11.733 21.258 /
%
%
\linethickness= 0.500pt
\setplotsymbol ({\thinlinefont .})
\plot 14.268 25.711 12.998 24.441 /
\putrule from 12.998 24.441 to 12.998 22.536
\plot 12.998 22.536 14.268 21.266 /
%
%
\linethickness= 0.500pt
\setplotsymbol ({\thinlinefont .})
\plot 12.990 24.450 11.720 25.720 /
\linethickness= 0.500pt
\setplotsymbol ({\thinlinefont .})
%
%
\plot  8.871 23.417      8.968 23.534
         9.052 23.609
         9.129 23.647
         9.207 23.654
         9.303 23.601
         9.368 23.496
         9.432 23.391
         9.525 23.336
         9.599 23.340
         9.673 23.373
         9.756 23.439
         9.855 23.544
        /
%
%
\put{\SetFigFont{10}{12.0}{it}j} [lB] at 11.225 21.004
%
%
\put{\SetFigFont{10}{12.0}{it}i} [lB] at  1.446 24.591
%
%
\put{\SetFigFont{10}{12.0}{it}j} [lB] at  1.367 21.939
%
%
\put{\SetFigFont{10}{12.0}{it}l} [lB] at  6.653 24.575
%
%
\put{\SetFigFont{10}{12.0}{it}k} [lB] at  6.636 22.035
%
%
\put{\SetFigFont{10}{12.0}{it}i} [lB] at 11.159 25.542
%
%
\put{\SetFigFont{10}{12.0}{it}l} [lB] at 14.590 25.624
%
%
\put{\SetFigFont{10}{12.0}{it}k} [lB] at 14.556 20.894
\linethickness=0pt
\putrectangle corners at  1.367 26.056 and 14.590 20.786
\endpicture}

\end{center}
\vspace{-0pt}

\noindent This  diagram corresponds to the equations:
$$A(i,j,k,l)  =  A(j,i,l,k) =  A(k,l,i,j)  =  A(l,k,j,i)  =  0$$
$$-A(i,l,k,j) = -A(k,j,i,l)  =  -A(l,i,j,k)  = -A(j,k,l,i) = 0.$$
To obtain the equations, read the labels around
the left or right diagram (either clockwise or counterclockwise,
but always reading two grouped together at an end first).
The other sixteen equations correspond similarly to
the diagrams:

\vspace{-0pt}
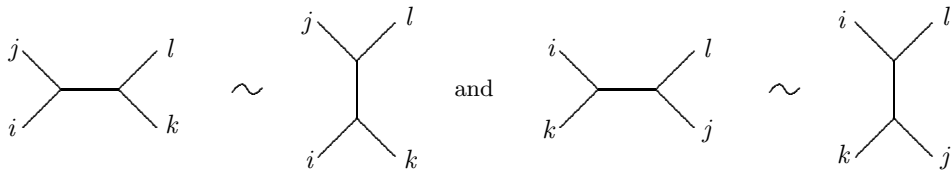
\begin{center}
  
\font\thinlinefont=cmr5
\begingroup\makeatletter\ifx\SetFigFont\undefined
\def\x#1#2#3#4#5#6#7\relax{\def\x{#1#2#3#4#5#6}}%
\expandafter\x\fmtname xxxxxx\relax \def\y{splain}%
\ifx\x\y   
\gdef\SetFigFont#1#2#3{%
  \ifnum #1<17\tiny\else \ifnum #1<20\small\else
  \ifnum #1<24\normalsize\else \ifnum #1<29\large\else
  \ifnum #1<34\Large\else \ifnum #1<41\LARGE\else
     \huge\fi\fi\fi\fi\fi\fi
  \csname #3\endcsname}%
\else
\gdef\SetFigFont#1#2#3{\begingroup
  \count@#1\relax \ifnum 25<\count@\count@25\fi
  \def\x{\endgroup\@setsize\SetFigFont{#2pt}}%
  \expandafter\x
    \csname \romannumeral\the\count@ pt\expandafter\endcsname
    \csname @\romannumeral\the\count@ pt\endcsname
  \csname #3\endcsname}%
\fi
\fi\endgroup
\mbox{\beginpicture
\setcoordinatesystem units <0.40000cm,0.40000cm>
\unitlength=0.40000cm
\linethickness=1pt
\setplotsymbol ({\makebox(0,0)[l]{\tencirc\symbol{'160}}})
\setshadesymbol ({\thinlinefont .})
\setlinear
%
%
\linethickness= 0.500pt
\setplotsymbol ({\thinlinefont .})
\plot 20.858 24.748 22.128 23.478 /
\putrule from 22.128 23.478 to 24.033 23.478
\plot 24.033 23.478 25.303 24.748 /
%
%
\linethickness= 0.500pt
\setplotsymbol ({\thinlinefont .})
\plot 24.033 23.478 25.303 22.208 /
%
%
\linethickness= 0.500pt
\setplotsymbol ({\thinlinefont .})
\plot 22.128 23.478 20.858 22.208 /
%
%
\linethickness= 0.500pt
\setplotsymbol ({\thinlinefont .})
\plot 31.955 22.511 30.685 21.241 /
%
%
\linethickness= 0.500pt
\setplotsymbol ({\thinlinefont .})
\plot 33.221 25.694 31.951 24.424 /
\putrule from 31.951 24.424 to 31.951 22.519
\plot 31.951 22.519 33.221 21.249 /
%
%
\linethickness= 0.500pt
\setplotsymbol ({\thinlinefont .})
\plot 31.943 24.433 30.673 25.703 /
%
%
\linethickness= 0.500pt
\setplotsymbol ({\thinlinefont .})
\plot  3.001 24.748  4.271 23.478 /
\putrule from  4.271 23.478 to  6.176 23.478
\plot  6.176 23.478  7.446 24.748 /
%
%
\linethickness= 0.500pt
\setplotsymbol ({\thinlinefont .})
\plot  6.176 23.478  7.446 22.208 /
%
%
\linethickness= 0.500pt
\setplotsymbol ({\thinlinefont .})
\plot  4.271 23.478  3.001 22.208 /
%
%
\linethickness= 0.500pt
\setplotsymbol ({\thinlinefont .})
\plot 14.099 22.511 12.829 21.241 /
%
%
\linethickness= 0.500pt
\setplotsymbol ({\thinlinefont .})
\plot 15.365 25.694 14.095 24.424 /
\putrule from 14.095 24.424 to 14.095 22.519
\plot 14.095 22.519 15.365 21.249 /
%
%
\linethickness= 0.500pt
\setplotsymbol ({\thinlinefont .})
\plot 14.086 24.433 12.816 25.703 /
\linethickness= 0.500pt
\setplotsymbol ({\thinlinefont .})
%
%
\plot 27.824 23.400     27.921 23.517
        28.004 23.592
        28.081 23.630
        28.160 23.637
        28.255 23.584
        28.321 23.479
        28.385 23.374
        28.478 23.319
        28.551 23.323
        28.626 23.356
        28.709 23.422
        28.808 23.527
        /
\linethickness= 0.500pt
\setplotsymbol ({\thinlinefont .})
%
%
\plot  9.967 23.400     10.065 23.517
        10.148 23.592
        10.225 23.630
        10.304 23.637
        10.399 23.584
        10.464 23.479
        10.529 23.374
        10.621 23.319
        10.695 23.323
        10.769 23.356
        10.852 23.422
        10.952 23.527
        /
%
%
\put{\SetFigFont{9}{10.8}{rm}and } [lB] at 17.340 23.273
%
%
\put{\SetFigFont{10}{12.0}{it}i} [lB] at 20.398 24.575
%
%
\put{\SetFigFont{10}{12.0}{it}l} [lB] at 25.605 24.558
%
%
\put{\SetFigFont{10}{12.0}{it}l} [lB] at 33.543 25.607
%
%
\put{\SetFigFont{10}{12.0}{it}k} [lB] at 20.254 21.859
%
%
\put{\SetFigFont{10}{12.0}{it}j} [lB] at 25.590 21.905
%
%
\put{\SetFigFont{10}{12.0}{it}i} [lB] at 30.048 25.510
%
%
\put{\SetFigFont{10}{12.0}{it}k} [lB] at 30.048 21.035
%
%
\put{\SetFigFont{10}{12.0}{it}j} [lB] at 33.513 21.006
%
%
\put{\SetFigFont{10}{12.0}{it}j} [lB] at  2.542 24.575
%
%
\put{\SetFigFont{10}{12.0}{it}i} [lB] at  2.464 21.922
%
%
\put{\SetFigFont{10}{12.0}{it}l} [lB] at  7.749 24.558
%
%
\put{\SetFigFont{10}{12.0}{it}k} [lB] at  7.732 22.018
%
%
\put{\SetFigFont{10}{12.0}{it}l} [lB] at 15.687 25.607
%
%
\put{\SetFigFont{10}{12.0}{it}k} [lB] at 15.653 20.877
%
%
\put{\SetFigFont{10}{12.0}{it}i} [lB] at 12.368 20.860
%
%
\put{\SetFigFont{10}{12.0}{it}j} [lB] at 12.291 25.480
\linethickness=0pt
\putrectangle corners at  2.464 26.039 and 33.543 20.758
\endpicture}

\end{center}
\vspace{-0pt}

\noindent In practice, one only needs to write down one
equation for each such diagram.

When 3 of the 4 labels are distinct, say
$i,i,j,k$, there is only 1 equation up to sign
(which occurs 8 times). It corresponds to:

\vspace{-0pt}
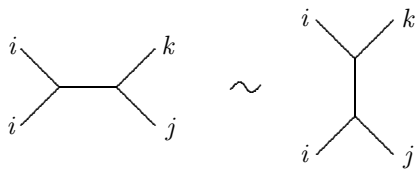
\begin{center}
  
\font\thinlinefont=cmr5
\begingroup\makeatletter\ifx\SetFigFont\undefined
\def\x#1#2#3#4#5#6#7\relax{\def\x{#1#2#3#4#5#6}}%
\expandafter\x\fmtname xxxxxx\relax \def\y{splain}%
\ifx\x\y   
\gdef\SetFigFont#1#2#3{%
  \ifnum #1<17\tiny\else \ifnum #1<20\small\else
  \ifnum #1<24\normalsize\else \ifnum #1<29\large\else
  \ifnum #1<34\Large\else \ifnum #1<41\LARGE\else
     \huge\fi\fi\fi\fi\fi\fi
  \csname #3\endcsname}%
\else
\gdef\SetFigFont#1#2#3{\begingroup
  \count@#1\relax \ifnum 25<\count@\count@25\fi
  \def\x{\endgroup\@setsize\SetFigFont{#2pt}}%
  \expandafter\x
    \csname \romannumeral\the\count@ pt\expandafter\endcsname
    \csname @\romannumeral\the\count@ pt\endcsname
  \csname #3\endcsname}%
\fi
\fi\endgroup
\mbox{\beginpicture
\setcoordinatesystem units <0.40000cm,0.40000cm>
\unitlength=0.40000cm
\linethickness=1pt
\setplotsymbol ({\makebox(0,0)[l]{\tencirc\symbol{'160}}})
\setshadesymbol ({\thinlinefont .})
\setlinear
%
%
\linethickness= 0.500pt
\setplotsymbol ({\thinlinefont .})
\plot  1.905 24.765  3.175 23.495 /
\putrule from  3.175 23.495 to  5.080 23.495
\plot  5.080 23.495  6.350 24.765 /
%
%
\linethickness= 0.500pt
\setplotsymbol ({\thinlinefont .})
\plot  5.080 23.495  6.350 22.225 /
%
%
\linethickness= 0.500pt
\setplotsymbol ({\thinlinefont .})
\plot  3.175 23.495  1.905 22.225 /
%
%
\linethickness= 0.500pt
\setplotsymbol ({\thinlinefont .})
\plot 13.003 22.528 11.733 21.258 /
%
%
\linethickness= 0.500pt
\setplotsymbol ({\thinlinefont .})
\plot 14.268 25.711 12.998 24.441 /
\putrule from 12.998 24.441 to 12.998 22.536
\plot 12.998 22.536 14.268 21.266 /
%
%
\linethickness= 0.500pt
\setplotsymbol ({\thinlinefont .})
\plot 12.990 24.450 11.720 25.720 /
\linethickness= 0.500pt
\setplotsymbol ({\thinlinefont .})
%
%
\plot  8.871 23.417      8.968 23.534
         9.052 23.609
         9.129 23.647
         9.207 23.654
         9.303 23.601
         9.368 23.496
         9.432 23.391
         9.525 23.336
         9.599 23.340
         9.673 23.373
         9.756 23.439
         9.855 23.544
        /
%
%
\put{\SetFigFont{10}{12.0}{it}i} [lB] at  1.446 24.591
%
%
\put{\SetFigFont{10}{12.0}{it}i} [lB] at 11.159 25.542
%
%
\put{\SetFigFont{10}{12.0}{it}i} [lB] at  1.414 21.924
%
%
\put{\SetFigFont{10}{12.0}{it}k} [lB] at  6.572 24.606
%
%
\put{\SetFigFont{10}{12.0}{it}j} [lB] at  6.716 21.893
%
%
\put{\SetFigFont{10}{12.0}{it}i} [lB] at 11.159 21.002
%
%
\put{\SetFigFont{10}{12.0}{it}k} [lB] at 14.541 25.512
%
%
\put{\SetFigFont{10}{12.0}{it}j} [lB] at 14.624 20.974
\linethickness=0pt
\putrectangle corners at  1.414 25.948 and 14.624 20.841
\endpicture}

\end{center}
\vspace{-0pt}

\noindent When two labels are distinct, there
is again only 1 equation up to sign (occurring 4 times):

\vspace{-0pt}
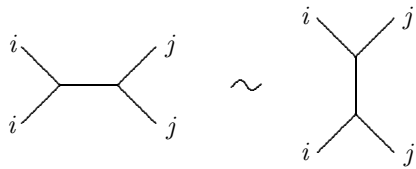
\begin{center}
  
\font\thinlinefont=cmr5
\begingroup\makeatletter\ifx\SetFigFont\undefined
\def\x#1#2#3#4#5#6#7\relax{\def\x{#1#2#3#4#5#6}}%
\expandafter\x\fmtname xxxxxx\relax \def\y{splain}%
\ifx\x\y   
\gdef\SetFigFont#1#2#3{%
  \ifnum #1<17\tiny\else \ifnum #1<20\small\else
  \ifnum #1<24\normalsize\else \ifnum #1<29\large\else
  \ifnum #1<34\Large\else \ifnum #1<41\LARGE\else
     \huge\fi\fi\fi\fi\fi\fi
  \csname #3\endcsname}%
\else
\gdef\SetFigFont#1#2#3{\begingroup
  \count@#1\relax \ifnum 25<\count@\count@25\fi
  \def\x{\endgroup\@setsize\SetFigFont{#2pt}}%
  \expandafter\x
    \csname \romannumeral\the\count@ pt\expandafter\endcsname
    \csname @\romannumeral\the\count@ pt\endcsname
  \csname #3\endcsname}%
\fi
\fi\endgroup
\mbox{\beginpicture
\setcoordinatesystem units <0.40000cm,0.40000cm>
\unitlength=0.40000cm
\linethickness=1pt
\setplotsymbol ({\makebox(0,0)[l]{\tencirc\symbol{'160}}})
\setshadesymbol ({\thinlinefont .})
\setlinear
%
%
\linethickness= 0.500pt
\setplotsymbol ({\thinlinefont .})
\plot  1.905 24.765  3.175 23.495 /
\putrule from  3.175 23.495 to  5.080 23.495
\plot  5.080 23.495  6.350 24.765 /
%
%
\linethickness= 0.500pt
\setplotsymbol ({\thinlinefont .})
\plot  5.080 23.495  6.350 22.225 /
%
%
\linethickness= 0.500pt
\setplotsymbol ({\thinlinefont .})
\plot  3.175 23.495  1.905 22.225 /
%
%
\linethickness= 0.500pt
\setplotsymbol ({\thinlinefont .})
\plot 13.003 22.528 11.733 21.258 /
%
%
\linethickness= 0.500pt
\setplotsymbol ({\thinlinefont .})
\plot 14.268 25.711 12.998 24.441 /
\putrule from 12.998 24.441 to 12.998 22.536
\plot 12.998 22.536 14.268 21.266 /
%
%
\linethickness= 0.500pt
\setplotsymbol ({\thinlinefont .})
\plot 12.990 24.450 11.720 25.720 /
\linethickness= 0.500pt
\setplotsymbol ({\thinlinefont .})
%
%
\plot  8.871 23.417      8.968 23.534
         9.052 23.609
         9.129 23.647
         9.207 23.654
         9.303 23.601
         9.368 23.496
         9.432 23.391
         9.525 23.336
         9.599 23.340
         9.673 23.373
         9.756 23.439
         9.855 23.544
        /
%
%
\put{\SetFigFont{10}{12.0}{it}i} [lB] at  1.446 24.591
%
%
\put{\SetFigFont{10}{12.0}{it}i} [lB] at 11.159 25.542
%
%
\put{\SetFigFont{10}{12.0}{it}i} [lB] at  1.414 21.924
%
%
\put{\SetFigFont{10}{12.0}{it}j} [lB] at  6.716 21.893
%
%
\put{\SetFigFont{10}{12.0}{it}i} [lB] at 11.159 21.002
%
%
\put{\SetFigFont{10}{12.0}{it}j} [lB] at 14.624 20.974
%
%
\put{\SetFigFont{10}{12.0}{it}j} [lB] at  6.684 24.591
%
%
\put{\SetFigFont{10}{12.0}{it}j} [lB] at 14.592 25.514
\linethickness=0pt
\putrectangle corners at  1.414 26.048 and 14.624 20.841
\endpicture}

\end{center}
\vspace{-0pt}

\noindent The symmetry in these diagrams reflects
the symmetry in the equations. Taking just one
equation for each diagram, one sees that the
number $N(m)$ of equations for $\text{rank} (A^*X)=m+1$
is
$$
N(m)= 3 \binom {m}{4} + m \binom {m-1}{2}+ \binom{m}{2} =
\frac
{m(m-1)(m^2-m+2)}{8},$$
so 
$N(2)=1$, $ N(3)=6$,  $N(4)=21$, $N(5)=55$,  $N(6)=120$, 
and $N(7)=
231$.
For the complete flag manifold $\mathbf{Fl}(\com^n)$, $m=n!-1$.
The number of equations for $\mathbf{Fl}(\com^4)$ is 
$N(23)=30861$.

Let us work this out for the two varieties $X= \proj^3$ and
$X= \mathbf{Q}^3$ (a smooth quadric 3-fold), which have
very similar classical cohomology rings.
Each has a basis :
\begin{eqnarray*}
T_0 & = & 1,\\
T_1 & = & \text{hyperplane class},\\
T_2 & = & \text{line class}, \\
T_3 & = & \text{point class}. 
\end{eqnarray*}
The difference in the classical product is that
$T_1 \smallcup T_1= T_2$ for $\proj^3$ but
$T_1 \smallcup T_1= 2T_2$ for $\mathbf{Q}^3$.
Let $c=1$ for $\proj^3$ and $c=2$ for $\mathbf{Q}^3$.
The $N(3)=6$ equations are:
\begin{align*}
\font\thinlinefont=cmr5
\begingroup\makeatletter\ifx\SetFigFont\undefined
\def\x#1#2#3#4#5#6#7\relax{\def\x{#1#2#3#4#5#6}}%
\expandafter\x\fmtname xxxxxx\relax \def\y{splain}%
\ifx\x\y   
\gdef\SetFigFont#1#2#3{%
  \ifnum #1<17\tiny\else \ifnum #1<20\small\else
  \ifnum #1<24\normalsize\else \ifnum #1<29\large\else
  \ifnum #1<34\Large\else \ifnum #1<41\LARGE\else
     \huge\fi\fi\fi\fi\fi\fi
  \csname #3\endcsname}%
\else
\gdef\SetFigFont#1#2#3{\begingroup
  \count@#1\relax \ifnum 25<\count@\count@25\fi
  \def\x{\endgroup\@setsize\SetFigFont{#2pt}}%
  \expandafter\x
    \csname \romannumeral\the\count@ pt\expandafter\endcsname
    \csname @\romannumeral\the\count@ pt\endcsname
  \csname #3\endcsname}%
\fi
\fi\endgroup
\mbox{\beginpicture
\setcoordinatesystem units <0.25000cm,0.25000cm>
\unitlength=0.25000cm
\linethickness=1pt
\setplotsymbol ({\makebox(0,0)[l]{\tencirc\symbol{'160}}})
\setshadesymbol ({\thinlinefont .})
\setlinear
%
%
\put{\SetFigFont{6}{7.2}{rm}1} [lB] at  1.486 24.384
%
%
\put{\SetFigFont{6}{7.2}{rm}1} [lB] at  1.507 23.004
%
%
\linethickness= 0.500pt
\setplotsymbol ({\thinlinefont .})
\plot  2.857 23.971  2.223 23.336 /
%
%
\linethickness= 0.500pt
\setplotsymbol ({\thinlinefont .})
\plot  2.223 24.606  2.857 23.971 /
\putrule from  2.857 23.971 to  3.810 23.971
\plot  3.810 23.971  4.445 24.606 /
%
%
\linethickness= 0.500pt
\setplotsymbol ({\thinlinefont .})
\plot  3.810 23.971  4.445 23.336 /
%
%
\put{\SetFigFont{6}{7.2}{rm}2} [lB] at  4.610 24.420
%
%
\put{\SetFigFont{6}{7.2}{rm}2} [lB] at  4.646 22.989
\linethickness=0pt
\putrectangle corners at  1.486 24.801 and  4.646 22.894
\endpicture}
&\hspace{.4in}& 2 \Gamma_{123}-c\Gamma_{222} & =  \Gamma_{111}\Gamma_{222} -
                                 \Gamma_{112}\Gamma_{122} \\
 \font\thinlinefont=cmr5
\begingroup\makeatletter\ifx\SetFigFont\undefined
\def\x#1#2#3#4#5#6#7\relax{\def\x{#1#2#3#4#5#6}}%
\expandafter\x\fmtname xxxxxx\relax \def\y{splain}%
\ifx\x\y   
\gdef\SetFigFont#1#2#3{%
  \ifnum #1<17\tiny\else \ifnum #1<20\small\else
  \ifnum #1<24\normalsize\else \ifnum #1<29\large\else
  \ifnum #1<34\Large\else \ifnum #1<41\LARGE\else
     \huge\fi\fi\fi\fi\fi\fi
  \csname #3\endcsname}%
\else
\gdef\SetFigFont#1#2#3{\begingroup
  \count@#1\relax \ifnum 25<\count@\count@25\fi
  \def\x{\endgroup\@setsize\SetFigFont{#2pt}}%
  \expandafter\x
    \csname \romannumeral\the\count@ pt\expandafter\endcsname
    \csname @\romannumeral\the\count@ pt\endcsname
  \csname #3\endcsname}%
\fi
\fi\endgroup
\mbox{\beginpicture
\setcoordinatesystem units <0.25000cm,0.25000cm>
\unitlength=0.25000cm
\linethickness=1pt
\setplotsymbol ({\makebox(0,0)[l]{\tencirc\symbol{'160}}})
\setshadesymbol ({\thinlinefont .})
\setlinear
%
%
\put{\SetFigFont{6}{7.2}{rm}1} [lB] at  1.507 24.418
%
%
\put{\SetFigFont{6}{7.2}{rm}1} [lB] at  1.496 23.036
%
%
\linethickness= 0.500pt
\setplotsymbol ({\thinlinefont .})
\plot  2.857 23.971  2.223 23.336 /
%
%
\linethickness= 0.500pt
\setplotsymbol ({\thinlinefont .})
\plot  2.223 24.606  2.857 23.971 /
\putrule from  2.857 23.971 to  3.810 23.971
\plot  3.810 23.971  4.445 24.606 /
%
%
\linethickness= 0.500pt
\setplotsymbol ({\thinlinefont .})
\plot  3.810 23.971  4.445 23.336 /
%
%
\put{\SetFigFont{6}{7.2}{rm}2} [lB] at  4.646 22.989
%
%
\put{\SetFigFont{6}{7.2}{rm}3} [lB] at  4.650 24.462
\linethickness=0pt
\putrectangle corners at  1.496 24.843 and  4.650 22.894
\endpicture}
  &\hspace{.4in}&  \Gamma_{133}-c\Gamma_{223} & =  \Gamma_{111}\Gamma_{223} -
                                 \Gamma_{113}\Gamma_{122} \\
 \font\thinlinefont=cmr5
\begingroup\makeatletter\ifx\SetFigFont\undefined
\def\x#1#2#3#4#5#6#7\relax{\def\x{#1#2#3#4#5#6}}%
\expandafter\x\fmtname xxxxxx\relax \def\y{splain}%
\ifx\x\y   
\gdef\SetFigFont#1#2#3{%
  \ifnum #1<17\tiny\else \ifnum #1<20\small\else
  \ifnum #1<24\normalsize\else \ifnum #1<29\large\else
  \ifnum #1<34\Large\else \ifnum #1<41\LARGE\else
     \huge\fi\fi\fi\fi\fi\fi
  \csname #3\endcsname}%
\else
\gdef\SetFigFont#1#2#3{\begingroup
  \count@#1\relax \ifnum 25<\count@\count@25\fi
  \def\x{\endgroup\@setsize\SetFigFont{#2pt}}%
  \expandafter\x
    \csname \romannumeral\the\count@ pt\expandafter\endcsname
    \csname @\romannumeral\the\count@ pt\endcsname
  \csname #3\endcsname}%
\fi
\fi\endgroup
\mbox{\beginpicture
\setcoordinatesystem units <0.25000cm,0.25000cm>
\unitlength=0.25000cm
\linethickness=1pt
\setplotsymbol ({\makebox(0,0)[l]{\tencirc\symbol{'160}}})
\setshadesymbol ({\thinlinefont .})
\setlinear
%
%
\put{\SetFigFont{6}{7.2}{rm}1} [lB] at  1.492 24.416
%
%
\put{\SetFigFont{6}{7.2}{rm}1} [lB] at  1.482 23.036
%
%
\linethickness= 0.500pt
\setplotsymbol ({\thinlinefont .})
\plot  2.857 23.971  2.223 23.336 /
%
%
\linethickness= 0.500pt
\setplotsymbol ({\thinlinefont .})
\plot  2.223 24.606  2.857 23.971 /
\putrule from  2.857 23.971 to  3.810 23.971
\plot  3.810 23.971  4.445 24.606 /
%
%
\linethickness= 0.500pt
\setplotsymbol ({\thinlinefont .})
\plot  3.810 23.971  4.445 23.336 /
%
%
\put{\SetFigFont{6}{7.2}{rm}3} [lB] at  4.635 23.036
%
%
\put{\SetFigFont{6}{7.2}{rm}3} [lB] at  4.652 24.448
\linethickness=0pt
\putrectangle corners at  1.482 24.829 and  4.652 22.940
\endpicture}
 &\hspace{.4in}& c\Gamma_{233} & =  2\Gamma_{113}\Gamma_{123} -
                                 \Gamma_{112}\Gamma_{133}
                               -\Gamma_{111}\Gamma_{233} \\
\font\thinlinefont=cmr5
\begingroup\makeatletter\ifx\SetFigFont\undefined
\def\x#1#2#3#4#5#6#7\relax{\def\x{#1#2#3#4#5#6}}%
\expandafter\x\fmtname xxxxxx\relax \def\y{splain}%
\ifx\x\y   
\gdef\SetFigFont#1#2#3{%
  \ifnum #1<17\tiny\else \ifnum #1<20\small\else
  \ifnum #1<24\normalsize\else \ifnum #1<29\large\else
  \ifnum #1<34\Large\else \ifnum #1<41\LARGE\else
     \huge\fi\fi\fi\fi\fi\fi
  \csname #3\endcsname}%
\else
\gdef\SetFigFont#1#2#3{\begingroup
  \count@#1\relax \ifnum 25<\count@\count@25\fi
  \def\x{\endgroup\@setsize\SetFigFont{#2pt}}%
  \expandafter\x
    \csname \romannumeral\the\count@ pt\expandafter\endcsname
    \csname @\romannumeral\the\count@ pt\endcsname
  \csname #3\endcsname}%
\fi
\fi\endgroup
\mbox{\beginpicture
\setcoordinatesystem units <0.25000cm,0.25000cm>
\unitlength=0.25000cm
\linethickness=1pt
\setplotsymbol ({\makebox(0,0)[l]{\tencirc\symbol{'160}}})
\setshadesymbol ({\thinlinefont .})
\setlinear
%
%
\put{\SetFigFont{6}{7.2}{rm}3} [lB] at  1.513 23.019
%
%
\put{\SetFigFont{6}{7.2}{rm}1} [lB] at  1.556 24.443
%
%
\linethickness= 0.500pt
\setplotsymbol ({\thinlinefont .})
\plot  2.857 23.971  2.223 23.336 /
%
%
\linethickness= 0.500pt
\setplotsymbol ({\thinlinefont .})
\plot  2.223 24.606  2.857 23.971 /
\putrule from  2.857 23.971 to  3.810 23.971
\plot  3.810 23.971  4.445 24.606 /
%
%
\linethickness= 0.500pt
\setplotsymbol ({\thinlinefont .})
\plot  3.810 23.971  4.445 23.336 /
%
%
\put{\SetFigFont{6}{7.2}{rm}2} [lB] at  4.610 24.420
%
%
\put{\SetFigFont{6}{7.2}{rm}2} [lB] at  4.646 22.989
\linethickness=0pt
\putrectangle corners at  1.513 24.824 and  4.646 22.894
\endpicture}
 &\hspace{.4in}& \Gamma_{233}& = \Gamma_{113}\Gamma_{222} -
                                 \Gamma_{112}\Gamma_{223} \\
\font\thinlinefont=cmr5
\begingroup\makeatletter\ifx\SetFigFont\undefined
\def\x#1#2#3#4#5#6#7\relax{\def\x{#1#2#3#4#5#6}}%
\expandafter\x\fmtname xxxxxx\relax \def\y{splain}%
\ifx\x\y   
\gdef\SetFigFont#1#2#3{%
  \ifnum #1<17\tiny\else \ifnum #1<20\small\else
  \ifnum #1<24\normalsize\else \ifnum #1<29\large\else
  \ifnum #1<34\Large\else \ifnum #1<41\LARGE\else
     \huge\fi\fi\fi\fi\fi\fi
  \csname #3\endcsname}%
\else
\gdef\SetFigFont#1#2#3{\begingroup
  \count@#1\relax \ifnum 25<\count@\count@25\fi
  \def\x{\endgroup\@setsize\SetFigFont{#2pt}}%
  \expandafter\x
    \csname \romannumeral\the\count@ pt\expandafter\endcsname
    \csname @\romannumeral\the\count@ pt\endcsname
  \csname #3\endcsname}%
\fi
\fi\endgroup
\mbox{\beginpicture
\setcoordinatesystem units <0.25000cm,0.25000cm>
\unitlength=0.25000cm
\linethickness=1pt
\setplotsymbol ({\makebox(0,0)[l]{\tencirc\symbol{'160}}})
\setshadesymbol ({\thinlinefont .})
\setlinear
%
%
\put{\SetFigFont{6}{7.2}{rm}3} [lB] at  1.524 22.987
%
%
\put{\SetFigFont{6}{7.2}{rm}1} [lB] at  1.501 24.428
%
%
\linethickness= 0.500pt
\setplotsymbol ({\thinlinefont .})
\plot  2.857 23.971  2.223 23.336 /
%
%
\linethickness= 0.500pt
\setplotsymbol ({\thinlinefont .})
\plot  2.223 24.606  2.857 23.971 /
\putrule from  2.857 23.971 to  3.810 23.971
\plot  3.810 23.971  4.445 24.606 /
%
%
\linethickness= 0.500pt
\setplotsymbol ({\thinlinefont .})
\plot  3.810 23.971  4.445 23.336 /
%
%
\put{\SetFigFont{6}{7.2}{rm}2} [lB] at  4.610 24.420
%
%
\put{\SetFigFont{6}{7.2}{rm}3} [lB] at  4.604 22.989
\linethickness=0pt
\putrectangle corners at  1.501 24.809 and  4.610 22.892
\endpicture}
 &\hspace{.4in}&     \Gamma_{333}& = \Gamma_{123}^2 -
                                 \Gamma_{122}\Gamma_{133}
                        +\Gamma_{113}\Gamma_{223}-\Gamma_{112}\Gamma_{233}\\
\font\thinlinefont=cmr5
\begingroup\makeatletter\ifx\SetFigFont\undefined
\def\x#1#2#3#4#5#6#7\relax{\def\x{#1#2#3#4#5#6}}%
\expandafter\x\fmtname xxxxxx\relax \def\y{splain}%
\ifx\x\y   
\gdef\SetFigFont#1#2#3{%
  \ifnum #1<17\tiny\else \ifnum #1<20\small\else
  \ifnum #1<24\normalsize\else \ifnum #1<29\large\else
  \ifnum #1<34\Large\else \ifnum #1<41\LARGE\else
     \huge\fi\fi\fi\fi\fi\fi
  \csname #3\endcsname}%
\else
\gdef\SetFigFont#1#2#3{\begingroup
  \count@#1\relax \ifnum 25<\count@\count@25\fi
  \def\x{\endgroup\@setsize\SetFigFont{#2pt}}%
  \expandafter\x
    \csname \romannumeral\the\count@ pt\expandafter\endcsname
    \csname @\romannumeral\the\count@ pt\endcsname
  \csname #3\endcsname}%
\fi
\fi\endgroup
\mbox{\beginpicture
\setcoordinatesystem units <0.25000cm,0.25000cm>
\unitlength=0.25000cm
\linethickness=1pt
\setplotsymbol ({\makebox(0,0)[l]{\tencirc\symbol{'160}}})
\setshadesymbol ({\thinlinefont .})
\setlinear
%
%
\put{\SetFigFont{6}{7.2}{rm}3} [lB] at  1.507 24.431
%
%
\put{\SetFigFont{6}{7.2}{rm}3} [lB] at  1.492 23.019
%
%
\linethickness= 0.500pt
\setplotsymbol ({\thinlinefont .})
\plot  2.857 23.971  2.223 23.336 /
%
%
\linethickness= 0.500pt
\setplotsymbol ({\thinlinefont .})
\plot  2.223 24.606  2.857 23.971 /
\putrule from  2.857 23.971 to  3.810 23.971
\plot  3.810 23.971  4.445 24.606 /
%
%
\linethickness= 0.500pt
\setplotsymbol ({\thinlinefont .})
\plot  3.810 23.971  4.445 23.336 /
%
%
\put{\SetFigFont{6}{7.2}{rm}2} [lB] at  4.610 24.420
%
%
\put{\SetFigFont{6}{7.2}{rm}2} [lB] at  4.646 22.989
\linethickness=0pt
\putrectangle corners at  1.492 24.812 and  4.646 22.894
\endpicture}
 &\hspace{.4in}&0 &= \Gamma_{133}\Gamma_{222}-2\Gamma_{123}\Gamma_{223}
        +\Gamma_{122}\Gamma_{233}
\end{align*}
The function $\Gamma$ has the form:
\begin{equation}
\label{dddef}
\Gamma= \sum N_{a,b}  e^{dy_1}  \frac{y_2^a}{a!} 
                                     \frac{y_3^b}{b!} .
\end{equation}
For $\proj^3$ the sum in (\ref{dddef})
is over non-negative $a,b$ satisfying
$a+2b=4d$, $d\geq 1$. A crucial difference is that for
$\mathbf{Q}^3$, the sum in (\ref{dddef}) is over
$a+2b=3d$, $d\geq 1$  reflecting the fact
that $c_1(T_{\proj^3})=4T_1$ while
$c_1(T_{\mathbf{Q}^3})=3T_1$. In each case, $N_{a,b}$
is the number of degree $d$ rational curves in $X$
meeting $a$ general lines and $b$ general points of $X$.

Each of the six differential equations above yields
a recursion among the $N_{a,b}$:
\vspace{+10pt}

\noindent (1)
For   $ a \geq 3,\ b \geq 0$, $\ \ 
2d N_{a-2,b+1}-c N_{a,b}=$
$$\sum N_{a_1,b_1}N_{a_2,b_2} \binom{b}{b_1}
\bigg( d_1^3  \binom{a-3}{a_1} -
d_1^2 d_2 \binom{a-3}{a_1-1} \bigg)$$
\noindent (2)
For   $ a \geq 2,\ b \geq 1$, $\ \ 
d N_{a-2,b+1}-c N_{a,b}=$
$$\sum N_{a_1,b_1}N_{a_2,b_2} \binom{a-2}{a_1}
\bigg( d_1^3  \binom{b-1}{b_1} -
d_1^2 d_2 \binom{b-1}{b_1-1} \bigg)$$
\noindent (3)
For $ a \geq 1,\ b \geq 2$, $\ \ 
c N_{a,b}=$
$$\sum N_{a_1,b_1}N_{a_2,b_2} 
\bigg(2 d_1^2 d_2  \binom{a-1}{a_1} \binom{b-2}{b_1-1} -
d_1^2 d_2 \binom{a-1}{a_1-1} \binom{b-2}{b_1}$$ $$-
d_1^3\binom{a-1}{a_1}\binom{b-2}{b_1} \bigg)$$
\noindent (4)
For  $ a \geq 3,\ b \geq 1$, $\ \ 
N_{a-2,b+1}=$
$$\sum N_{a_1,b_1}N_{a_2,b_2} d_1^2 
\bigg(  \binom{a-3}{a_1} \binom{b-1}{b_1-1} -
 \binom{a-3}{a_1-1} \binom{b-1}{b_1} \bigg)$$
\noindent (5)
For    $ a \geq 2,\ b \geq 2$, $\ \ 
N_{a-2,b+1}=$
$$\sum N_{a_1,b_1}N_{a_2,b_2} 
\bigg( d_1d_2  \binom{a-2}{a_1-1} \binom{b-2}{b_1-1} -
d_1 d_2 \binom{a-2}{a_1-2} \binom{b-2}{b_1}$$ $$+
d_1^2\binom{a-2}{a_1}\binom{b-2}{b_1-1}
- d_1^2\binom{a-2}{a_1-1}\binom{b-2}{b_1} \bigg)$$
\noindent (6)
For $ a \geq 3,\ b \geq 2$, $\ \ 
 0 =$
$$\sum N_{a_1,b_1}N_{a_2,b_2} d_1
\bigg(   \binom{a-3}{a_1} \binom{b-2}{b_1-2} -
2 \binom{a-3}{a_1-1} \binom{b-2}{b_1-1}$$ $$+
\binom{a-3}{a_1-2}\binom{b-2}{b_1} \bigg)$$
In these formulas, the sum is over
non-negative $a_1, a_2, b_1, b_2$ satisfying
\begin{enumerate}
\item[(i)] $a_1+a_2=a$, $b_1+b_2=b$, 
\item[(ii)] $a+2b=4d$, $a_i+2b_i=4d_i$, $d_i>0$  for $\proj^3$,
\item[{}] $a+2b=3d$, $a_i+2b_i=3d_i$, $d_i>0$ for $\mathbf{Q}^3$.
\end{enumerate}

For $\proj^3$, one starts with the $N_{0,2}=1$
for the number of lines through two points.
For $\mathbf{Q}^3$,  $N_{1,1}=1$ is not hard
to compute directly. In each case, the six recursions
are more than enough to solve for all the other
$N_{a,b}$. These numbers for $\proj^3$
include the classical results: there are
{$N_{4,0}=2$} lines meeting 4 general lines, 
{$N_{8,0}=92$} conics meeting 8 general lines, and 
{$N_{12,0}=80160$} twisted cubics meeting
12 general lines.
See [DF-I] for more of these numbers$^2$.
\footnotetext[2]{The numbers $N_{a,b}$ given in
[DF-I] are correct, although their version of
equation (6) has a misprint.}
For $\mathbf{Q}^3$, computations yield:

\begin{tabular}{ll}
$(d=1)$ & $ N_{1,1}=1$,   $N_{3,0}=1$  \\
$(d=2)$ & $N_{0,3}=1$, $N_{2,2}=1$, $N_{4,1}=2$ ,
             $N_{6,0}=5$ \\
$(d=3)$ & $N_{1,4}=2$, $N_{3,3}=5$,  $N_{5,2}=16$, 
             $N_{7,1}=59$, $N_{9,0}=242$ \\
$(d=4)$  & $N_{0,6}=6$, $N_{2,5}=20$, $N_{4,4}=74$,  
           $N_{6,3}=320$, $N_{8,2}=1546$, \\ 
 {} &  $N_{10,1}=8148$,
            $N_{12,0}=46230$ \\
$(d=5)$  & $N_{1,7}=106$, $N_{3,6}=448$, $N_{5,5}=2180$,
         $N_{7,4}=11910$, \\ &
$N_{9,3}=71178$, $N_{11,2}=457788$,
         $N_{13,1}=3136284$, \\ & $N_{15,0}=22731810$. 
\end{tabular}

The reader is invited to work out the equations for
some other simple homogeneous
spaces such as $\proj^4$, $\proj^1 \times \proj^1$,
$\mathbf{Gr}(2,4)$, or the  incidence variety
$\mathbf{Fl}(\com^3)$ of points on lines in the plane.
For  very pleasant excursions along these paths, see
[DF-I].

There is a simple method of obtaining a presentation
of $QH^* X$ from $\Phi$ and
a presentation of $A^*X$. It will be convenient
to consider 
$A^* X_\Q= H^*(X, \Q)$, the cohomology ring 
of $X$ with rational coefficients. Following the
notation of section 8,
let $QH^* X=(V\otimes _Z \Q[[V^*]],*)$.  
There is a
canonical embedding:
$$\iota_\Q: A^* X_{\Q} \hookrightarrow  QH^* X$$
of $\Q$-vector spaces. In the discussion below,
$A^* X_{\Q}$ is viewed as a $\Q$-subspace of $QH^* X$
via $\iota_\Q$. The results
relating presentations of $A^* X_\Q$ and $QH^* X$ are
established in Propositions \ref{ett} and \ref{tvo}.
\begin{pr}
\label{ett}
Let $z_1,..., z_r$ be homogeneous elements of positive
codimension that generate
$A^* X_{\Q}$ as a $\Q$-algebra. Then, $z_1, \ldots,
z_r$ generate $QH^* X$ as a $\Q[[V^*]]$-algebra.
\end{pr}

The proof requires a lemma. Note that for
$\gamma\in \Q[[V^*]]$ there is a well-defined
constant term $\gamma(0)\in \Q$.
\begin{lm} 
\label{mygid}
Let $T_0, \ldots, T_m$ be any homogeneous $\Q$-basis
of $A^* X_{\Q}$. Let $w_1, w_2\in A^* X_{\Q}$ be
homogeneous elements. Let
$$
w_1 \smallcup w_2  =  \sum_{k=0}^{m} c_{k}T_k,  \ \ c_{k}\in \Q, $$
$$
w_1 * w_2  =  \sum_{k=0}^{m} \gamma_{k} T_k,  \ \ 
\gamma_{k}\in \Q[[V^*]],$$
be the unique expansions in $A^* X_\Q$ and $QH^* X$ respectively. 
\begin{enumerate}
\item[(i)] If $\text{\em codim}(T_k)> 
\text{\em codim}(w_1)+\text{\em codim}(w_2)$,
           then $\gamma_{k}(0)=0$.
\item[(ii)] If $\text{\em codim}(T_k)= \text{\em codim}(w_1)+
                \text{\em codim}(w_2)$,
           then $\gamma_{k}(0)=c_{k}$.
\end{enumerate}
\end{lm}
\bpf 
By linearity of the $*$-product, it can be assumed
that $w_1$ and $w_2$ are basis elements $T_i$ and $T_j$
respectively.
In the basis $T_0, \ldots, T_m$ of $A^*_\Q X$, the
$*$-product is determined by:
$$T_i * T_j= T_i \smallcup T_j 
+ \sum_{i=1}^{m} \Gamma_{ijl} g^{lk}T_k $$
where the dual coordinates $y_0, \ldots, y_m$ are
taken in $V^*\otimes \Q$. $\Gamma_{ijl}(0)=\sum_{\beta\neq 0}
I_\beta(T_i\cdot T_j \cdot T_l)$.
Therefore, if
$\Gamma_{ijl}(0) \neq 0$, 
there must exist a nonzero effective class $\beta\in A_1 X$
such that  
$$\text{dim} \overline{M}_{0,3}(X, \beta) = \text{codim}(T_i)+ 
\text{codim}(T_j) +\text{codim}(T_l).$$
Since $X$ is homogeneous, $\int_\beta c_1(X)\geq 2$
by Lemma \ref{george}. By the dimension formula,
\begin{equation}
\label{arithh}
\text{codim}(T_i)+ 
\text{codim}(T_j) + \text{codim}(T_l)
\geq \text{dim}(X)+ 2.
\end{equation}
Equation (\ref{arithh}) yields
$\text{codim}(T_l)  \geq \text{dim} (X) - \text{codim}(T_i)
-\text{codim}(T_j) +2$. For $g^{lk}$
to be nonzero, it follows that
$\text{codim}(T_k) \leq \text{codim}(T_i)+\text{codim}(T_j)-2$.  
The lemma is proven.
\epf

We will apply Lemma \ref{mygid} to 
products in
a basis of $A^* X_\Q$ consisting of
monomials $z^I= z_1^{i_1} \smallcup \cdots \smallcup z_r^{i_r}$.
Let 
$$z^{*I}= 
\underbrace{z_1* \cdots *z_1}_{i_1}
*\underbrace{z_2*\cdots *z_2}_{i_2}* \cdots 
*\underbrace{z_r* \cdots *z_r}_{i_r}$$
denote the corresponding monomial in $QH^* X$.
Let 
\begin{equation}
\label{monny}
\{z^I \ | \ I \in \mathcal{S}\}
\end{equation}
be a monomial $\Q$-basis 
of $A^* X_\Q$. Choose an ordering of the set $\mathcal{S}$ so that
$\text{codim}(z^I) \leq \text{codim} (z^J)$ for $I<J$.
Let 
$$z^{*I}= \sum_{J\in \mathcal{S}} \gamma_{IJ} z^J, \ \ \gamma_{IJ}\in 
\Q[[V^*]]$$
be the unique expansion in $QH^* X$.
An inductive application of Lemma \ref{mygid} yields:
\begin{enumerate}
\item[(i)] If $J>I$, then $\gamma_{IJ}(0)=0$.
\item[(ii)] $\gamma_{II}(0)=1$.
\end{enumerate}
Therefore, 
the matrix $(\gamma_{IJ}(0))$ is invertible
over $\Q$. It follows that the matrix $(\gamma_{IJ})$
is invertible over $\Q[[V^*]]$. In particular,
$\{ z^{*I} \ | \ I \in \mathcal{S}\}$ is a $\Q[[V^*]]$-basis
of $QH^* X$. Proposition \ref{ett} is proved.

Let $K$ be the kernel of the
surjection 
$$\phi:\Q[Z]=\Q[Z_1, \ldots, Z_r] \rarr A^* X_\Q$$
determined by $\phi(Z_i)=z_i$. Let $K'$ be the
kernel of the corresponding surjection
$$\phi':\Q[[V^*]][Z] \rarr QH^* X$$
determined by $\phi'(Z_i)=z_i$
Using our choice (\ref{monny}) of monomial basis, there is
a method of constructing elements of $K'$
from elements of $K$. 
Let $f\in K$. The polynomial $f$ is also an
element of $\Q[[V^*]][Z]$. There is a unique expansion:
$$\phi'(f)= \sum_{I \in \mathcal{S}} \xi_I z^{*I}, \ \
\xi_I \in \Q[[V^*]].$$
Then, $f'= f(Z_1, \ldots, ,Z_r)- \sum_{I \in \mathcal{S}} \xi_I
Z^I$ is in $K'$.

The ideal $K$ is homogeneous provided the degree of $Z_i$
is taken to be the codimension of $z_i$. We need the
following fact.
\begin{lm}
\label{ghead}
Let $f\in K$ be homogeneous of degree $d$ and let $I\in \mathcal{S}$. 
If $\text{\em deg}(Z^I)\geq d$, then $\xi_I(0)=0$.
\end{lm}
\bpf
If $d> \text{dim}(X)$, the statement is vacuous.
Assume $d\leq \text{dim}(X)$.
Let $\phi'(f)= \sum_{I \in \mathcal{S}} \tilde{\xi}_I z^{I}, \ \
\tilde{\xi}_I \in \Q[[V^*]]$
be the unique expansion.
Apply Lemma \ref{mygid} repeatedly to 
the monomials of $f$ in
the basis
$\{ z^I\ | \ I \in \mathcal{S}\}$ of $A^* X_\Q$.
It follows that if
$\text{deg}(Z^I)\geq d$, then $\tilde{\xi}_I(0)=0$.
The
change of basis relations (i) and (ii) for the $\Q[[V^*]]$-basis
$\{ z^{*I} \ | \ I \in \mathcal{S}\}$ now imply the lemma.
\epf

Now suppose the elements $f_1, \ldots, f_s$ are homogeneous
generators of $K$, so
$$A^*X_\Q = \Q[Z]/(f_1, \ldots,f_s)$$
is a presentation of the cohomology ring.
\begin{pr}
\label{tvo}
The ideal $K'$ is generated by $f_1', \ldots, f_s'$, so
$$QH^*X = \Q[[V^*]][Z]/(f'_1, \ldots, f'_s)$$
is a presentation of the quantum cohomology ring.
\end{pr}

\bpf 
Since we have a surjection
$$\Q[[V^*]][Z]/ (f_1' \ldots, f_s') \rarr QH^*X$$
and $QH^*X$ is a free $\Q[[V^*]]$-module with basis
$\{ z^{*I} \ | \ I \in \mathcal{S}\}$, it suffices
to show that the monomials  
$\{ Z^{I} \ | \ I \in \mathcal{S}\}$
span the $\Q[[V^*]]$-module on the left.
By Nakayama's lemma, it suffices to show that these
monomials generate the $\Q$-vector space
\begin{equation}
\label{frod}
\Q[[V^*]][Z]/(f_1', \ldots, f_s', \frak{m}),
\end{equation}
where $\frak{m}\subset \Q[[V^*]]$ is the maximal ideal.
Let $f'_i= f_i- \sum\xi_{iI} Z^I$. Define
$\barr{f}'_i \in \Q[Z]$ by $\barr{f}'_i= f_i- \sum\xi_{iI}(0)
Z^I$.
The $\Q$-algebra (\ref{frod}) can be identified with
$$\Q[Z]/(\barr{f}'_1, \ldots, \barr{f}'_s).$$
By Lemma \ref{ghead}, all the terms $\xi_{iI}(0)Z^I$
have strictly lower degree than $f_i$. It is then
a simple induction on the degree to see that
the same monomials $\{Z^I\}$ that span modulo
$(f_1, \ldots, f_s)$ will also span modulo
$(\barr{f}'_1, \ldots, \barr{f}'_s)$.
\epf

For example, let $X=\mathbf{P}^2$. Let $Z=Z_1$ and let
$A^*_\Q \mathbf{P}^2=\Q[Z]/Z^3$ be the standard
presentation with the monomial basis $1, Z, Z^2$.
A presentation of $QH^* \mathbf{P}^2$ is obtained:
\begin{equation}
\label{pezzx}
QH^*\mathbf{P}^2 \stackrel{\sim}{\rightarrow} \Q[[y_0,y_1, y_2]][Z]
/(Z^3- \Gamma_{111}Z^2 -2 \Gamma_{112} Z - \Gamma_{122})
\end{equation}
where $\Gamma$ is the quantum potential of $\mathbf{P}^2$.
By (\ref{pezzx}) and the determination of $\Gamma$,
$$QH^*\proj^2 \otimes_{\Q[[V^*]]} \Q[[V^*]]/\frak{m}
= \Q[Z]/(Z^3-1).$$
Note that $QH^* \proj^2$ does not specialize to $A^* \proj^2$.
%
%
%

\section{\bf{Variations}}
The algebra  $QH^* X=
A^*X \otimes \mathbb Q [[V^*]]$  may be regarded as the ``big'' 
quantum cohomology ring.  There is also a ``small'' quantum 
cohomology ring, $QH^*_{s} X$, that incorporates only the $3$-point
Gromov-Witten invariants in its product.
$QH^*_{s} X$ is obtained by
restricting
the $*$-product
to the formal deformation parameters of the divisor
classes. Most computations of quantum cohomology rings
have been of this small ring, which is often easier to
describe; the small ring is often denoted $QH^* X$.

It is simplest to define $QH^*_{s} X$
in the Schubert basis
$T_0, \ldots, T_m$. Let
\begin{equation}
\overline\Phi_{i  j  k}  =  \Phi_{i  j  k}(y_0, y_1,  \ldots , 
y_p,   0,  \ldots  , 0) 
= \int_X  T_i \smallcup T_j \smallcup T_k  
\ + \ \overline\Gamma_{i j k}  \label{2.16}.
\end{equation}
The modified quantum potential $\barr{\Gamma}_{i j k}$
is determined by
$$\overline\Gamma_{i j k}  = \sum_{n \geq 0} \frac{1}{n!} 
\sum_{\beta \neq 0} I_\beta(\gamma^ n {\cdot} T_i 
{\cdot} T_j {\cdot} T_k)$$
where  $\gamma = y_1T_1 + \ldots + 
y_pT_p$.  By the divisor property (III) of section 7, 
\begin{equation}
\overline\Gamma_{i  j  k} =  \sum_{\beta \neq 0}   I_\beta(T_i 
{\cdot} T_j {\cdot} T_k)  {q_1}^{\int_\beta T_1} 
 \cdots {q_p}^{\int_\beta T_p}  , \label{2.17}
\end{equation}
where  $q_i = e^{y_i}$. Note that only $3$-point
invariants occur. 
Let ${\mathbb Z}[q]=[q_1, \ldots, q_p]$.
By Theorem 4, the product
$$
T_i \, * \,  T_j  =  
\sum_{e, \, f} \overline\Phi_{i  j  e}   g^{e  f}   
T_f  \, = \,  T_i  \smallcup T_j \ + \ \sum_{e, \, f} 
\overline\Gamma_{i  j  e}   g^{e  f}   T_f 
$$
then makes  the ${\mathbb Z}[q]$-module
$A^*X \otimes_{\mathbb Z} {\mathbb Z}[q]$  into a 
commutative, associative  ${\mathbb Z}[q]$-algebra with 
unit  $T_0$.  
From equation (\ref{2.17}), it easily follows that
the small quantum cohomology is a deformation
of $A^* X$ is the usual sense: $A^* X$ is
recovered by setting the variables $q_i=0$.

For example, let  $X = \proj^r$.
Then, $q=q_1$.
If  $T_i$  
is the class of a linear subspace of codimension  $i$ and  $\beta$  
is  $d$  times the class of a line, then the number  
$I_\beta(T_i {\cdot} T_j {\cdot} T_k)$  can be nonzero 
only if  $i + j + k = r + (r \! + \! 1)d$;  this can happen only for  $d = 
0$  or  $d = 1$,  and in each case the number is  $1$.  It follows that,
\begin{enumerate}
\item[(i)] if $i+j \leq r$, then $T_i\, * \, T_j =T_{i+j}$;
\item[(ii)] if $r+1 \leq i+j \leq 2r$, then $T_i\, * \, T_j =qT_{i+j-r-1}$.
\end{enumerate} 
Therefore the small quantum cohomology ring is:  
$$
QH^*_s \proj^r =\mathbb{Z} [ T,\, q ] / (T^{r+1} - q) ,  
$$
where  $T = T_1$  is the class of a hyperplane.

The following variation of Proposition \ref{tvo} is 
valid for the small quantum cohomology ring (cf. [S-T]).
As before let $z_1, \ldots, z_r$ be homogenous
elements of positive codimension that generate $A^* X$.
(We use integer coefficients but rational coefficients
could be used as well). Let ${\mathbb Z}[Z]=
{\mathbb Z}[Z_1, \ldots, Z_r]$, and let 
$$A^*X = {\mathbb Z}[Z]/ (f_1, \ldots, f_s)$$
be a presentation with arbitrary homogeneous
generators $f_1, \ldots, f_s$ for the ideal
of relations. Let ${\mathbb Z}[q,Z]={\mathbb Z}
[q_1, \ldots, q_p, Z_1, \ldots, Z_r]$. The variables
$q_i, Z_j$ are graded by the following degrees:
$\text{deg}(q_i)= \int_{\beta_i} c_1(T_X)$
where $\beta_i$ is the class of the Schubert variety
dual to $T_i$ and $\text{deg} (Z_j)= \text{codim}(z_j)$.
Let $QH_s^* X= A^*X \otimes  {\mathbb Z}[q]$ with
the quantum product.
\begin{pr}
\label{tva}
Let $f'_1, \ldots, f'_s$ be any homogeneous
elements in ${\mathbb Z}[q,Z]$ such that:
\begin{enumerate}
\item[(i)] $f_i'(0,\ldots,0,Z_1,\ldots, Z_r)= f_i(Z_1, \ldots,Z_r)$
in ${\mathbb Z}[q,Z]$,
\item[(ii)] $f'_i(q_1, \ldots, q_p, Z_1, \ldots, Z_r)=0$
in $QH^*_s X$.
\end{enumerate}
Then, the canonical map
\begin{equation}
\label{fffff}
{\mathbb Z}[q,Z]/ (f'_1, \ldots, f'_s) \rarr QH^*_s X
\end{equation}
is an isomorphism.
\end{pr}
\bpf
The proof is by a Nakayama-type induction. As the arguments
are similiar to the proof of Proposition \ref{tvo}, we
will be brief. The fact that each $q_i$ has positive degree
implies the following statement. If
$\psi:M\rarr N$ is a homogeneous map of finitely
generated ${\mathbb Z}[q,Z]$-modules that
is surjective modulo the ideal $(q)=(q_1, \ldots, q_p)$,
then $\psi$ is surjective.
Hence, by (i), the map (\ref{fffff}) is surjective.
Similarly, if $\tilde{T}_0, \ldots, \tilde{T}_m$
are homogeneous lifts to ${\mathbb Z}[Z]$ of a basis
of $A^* X$, an easy induction shows that their
images in ${\mathbb Z}[q,Z]/(f_1', \ldots, f_s')$ generate
this ${\mathbb Z}[q]$-module. Since
$QH^*_s X$ is free over ${\mathbb Z}$ of rank
$m+1$, the map (\ref{fffff}) must be
an isomorphism.
\epf

A similar calculation, as in [S-T], yields
the small quantum cohomology ring of 
the Grassmannian  $X = 
\mathbf{Gr}(p,n)$  of  $p$-dimensional subspaces of  
$\mathbb C^n$.  Let  $k = n - p$,  let  
$0 \rarr S \rarr {{\mathbb C}^n}_X 
\rarr Q \rarr 0$  be the universal exact sequence of bundles on  $X$,  
and let  $\sigma_i = c_i(Q)$.  Set  $S_r(\sigma) = 
\text{det}
\left(\sigma_{1+j-i}\right)_{1 \leq i, j \leq r}$,  and let  $q = 
q_1$.   

\begin{pr}  The small quantum 
cohomology ring of  $\mathbf{Gr}(p,n)$  is
$$
{\mathbb Z}[\sigma_1, \ldots ,  \sigma_k,  q]  /  \left(S_{p+1}(\sigma),   
S_{p+2}(\sigma), \ldots ,   S_{n-1}(\sigma),   S_n(\sigma)  +  (-1)^kq  
\right) .$$
\end{pr}
\bpf
We use some standard facts about the Grassmannian.  In 
particular, the cohomology has an additive basis of Schubert classes  
$\sigma_\lambda$,  as $\lambda$  varies over partitions with  $k 
\geq \lambda_1 \geq \ldots \geq \lambda_p \geq 0$;  
$\sigma_\lambda = [\Omega_\lambda]$  is the class of a Schubert 
variety  
$$
\Omega_\lambda  =  \{  L \in X   :  \dim   L \cap V_{k+i-\lambda_i} 
\geq i\  \text {for}\  1 \leq i \leq p  \} ,
$$
where  $V_1 \subset V_2 \subset \ldots \subset V_n = \mathbb C^n$  is 
a given flag of subspaces.  In  $A^*(X)$,  $S_r(\sigma)$  represents 
the  $r^{\text th}$  Chern class of  $S^\vee$,  from which we have  
$$
A^*(X)  =  {\mathbb Z}[\sigma_1, \ldots , \sigma_k]  /  \left( 
S_{p+1}(\sigma), \ldots , S_n(\sigma)\right)  .  
$$
By Proposition \ref{tva},
it suffices to show that the relations displayed in the 
proposition are valid in  $QH^*_s X$.
        
Since  $c_1(T_X) = n \sigma_1$,  a number  $I_\beta(\gamma_1 
{\cdot} \gamma_2 {\cdot} \gamma_3)$  can be nonzero 
only if the sum of the codimensions of the  $\gamma_i$  is equal to  
$\dim  X + nd$,  where  $\beta$  is  $d$  times the class of a 
line.  If  $d \geq 1$,  all such numbers vanish when  
$\text{codim}(\gamma_1) + \text{codim}
(\gamma_2) < n$.  In particular, the 
relations  $S_i(\sigma) = 0$  for  $p < i < n$  remain valid in  
$QH^*_s X$.  
From the formal identity  
$$
S_n(\sigma) - \sigma_1\, S_{n-1}(\sigma) + \sigma_2 \, S_{n-
2}(\sigma) - \ldots + (-1)^k\sigma_k \, S_{n-k}(\sigma)  =  0 ,
$$
we therefore have  $S_n(\sigma) = (-1)^{k-1}\sigma_k \, S_{n-
k}(\sigma)$  in  $QH^*_s X$.  Since  $S_{n-k}(\sigma) =  \sigma_{(1^{n-k})}$,  
the proof will be completed by verifying that  $\sigma_k * 
\sigma_{(1^{n-k})} = q$.  Equivalently, when  $\beta$  is the class of 
a line, we must show that
$$
I _\beta (\sigma_k , \sigma_{(1^p)},  \sigma_{(k^p)})  =  1 .
 $$  
This is a straightforward calculation.  First we have
$$
\sigma_k = [\{ L : L \supset A \}] ,\    \sigma_{(1^p)} = [\{ L : L \subset 
B \}] , \    \sigma_{(k^p)} =  [\{ L : L = C \}] ,
$$
where  $A$,  $B$,  and  $C$  are linear subspaces of  ${\mathbb C}^n$  of 
dimensions  $1$,  $n-1$,  and  $p$  respectively.  It is not hard to 
verify that any line in  $X$  is a Schubert variety of the form  
$\{ L : U \subset L \subset V\}$,  where  $U \subset V$  are 
subspaces of  ${\mathbb C}^n$  of dimensions  $p-1$  and  $p+1$.  Such a 
line will meet the three displayed Schubert varieties only if  $V$  
contains  $A$  and  $C$,  and  $U$  is contained in  $B$  and  $C$.  For  
$A$,  $B$,  and  $C$  general, there is only one such line, with  $U = B 
\cap C$  and  $V$  spanned by  $A$  and  $C$.  \epf

This proposition was proved in another way by Bertram [Ber], where 
the beginnings of some ``quantum Schubert calculus'' can be found.  
For the small quantum cohomology ring of a flag manifold, following 
ideas of Bertram, Givental, and Kim, see [CF]$^3$.
\footnotetext[3]{This Schubert calculus is extended to
flag manifolds in [F-G-P].}

As with the big quantum cohomology ring, the
small ring has a basis independent
description. 
Let $\Z[A_1 X]$ be the group algebra.
The small $*$-product is naturally
defined on the free $\Z[A_1 X]$-module $A^* X \otimes_\Z \Z[A_1 X]$.
If $\beta_1, \ldots, \beta_p$ is a basis
of $A_1 X$ consisting of Schubert classes, then
the dual Schubert classes $T_1, \ldots, T_p$ satisfy
$\int_\beta T_i \geq 0$ for every effective class $\beta$.
In this case, the small $*$-product on
$A^* X \otimes_\Z \Z[A_1 X]$ preserves the 
$\Z[q_1, \ldots, q_p]$-submodule:
$$A^*X \otimes _\Z \Z[q_1, \ldots, q_p]
\subset A^* X \otimes_\Z \Z[A_1 X].$$
Hence, in the Schubert basis, the small
quatum cohomology ring can be taken to be
$QH^*_s X= (A^* X \otimes_\Z \Z[q_1, \ldots, q_p],*).$

The numbers $I_{\beta}(\gamma_1 \cdots \gamma_n)$
should not be confused with the numbers denoted
by the expression $\langle\gamma_1, \ldots, \gamma_n \rangle_{\beta}$
which often occur in discussions of small quantum
cohomology rings ([B-D-W], [Ber], [CF]). To define the latter,
one {\em fixes} $n$ distinct points 
$p_1, \ldots, p_n$ in $\proj^1$.
Then, $\langle\gamma_1, \ldots, \gamma_n\rangle_{\beta}$ is the
number of maps $\mu:\proj^1 \rarr X$ satisfying:
$\mu_*[\proj^1]=\beta$ and
$\mu(p_i)\in \Gamma_i$ for $1 \leq i \leq n$
(where $\Gamma_i$ is a subvariety in general position
representing the class $\gamma_i$). For $n=3$, the
numbers agree: $I_{\beta}(\gamma_1 \cdot
\gamma_2 \cdot \gamma_3) = \langle\gamma_1, \gamma_2, \gamma_3
\rangle_{\beta}$.
For $n>3$, the numbers $\langle\gamma_1, \ldots, \gamma_n
\rangle_{\beta}$
and $I_{\beta}(\gamma_1 \cdots \gamma_n)$
are solutions to different enumerative problems.
In fact, $\langle\gamma_1, \ldots, \gamma_n\rangle_{\beta}$
can  be expressed in terms of the
$3$-points numbers while  $I_{\beta}(\gamma_1 \cdots \gamma_n)$
cannot.
For $1 < k < n-1$, $\langle\gamma_1, \ldots, \gamma_n\rangle_{\beta} =$
\begin{equation}
\label{ggpp}
\sum_{\beta_1+\beta_2=\beta} \sum_{e,f} 
\langle\gamma_1, \ldots, \gamma_k, T_e\rangle_{\beta_1}g^{ef}
\langle T_f,\gamma_{k+1}, \ldots, \gamma_n\rangle_{\beta_2}
\end{equation}
Equation (\ref{ggpp}) can be seen geometrically by
deforming $\proj^1$ to a union of two $\proj^1$'s meeting
at a point with $p_1, \ldots, p_k$ going
to fixed points on the first line and
$p_{k+1}, \ldots, p_n$ going to fixed points on the
second. Algebraically, in the
small quantum cohomology ring,
$$\gamma_1 * \cdots * \gamma_n =
\sum_{\beta} \sum_{ e,f} q^{\beta}
\langle\gamma_1, \ldots, \gamma_n, T_e\rangle_{\beta} g^{ef} T_f.$$
Equation (\ref{ggpp}) amounts to the associativity of this
product.

We conclude with a few general remarks to relate the discussion and 
notation here to that in [K-M 1]. 

The numbers that we have denoted  $I_\beta(\gamma_1  
\cdots  \gamma_n)$  are part of a more general story.  
Let  $\eta$  denote the forgetful map from  $\M_{0,n}(X,\beta)$  to  
$\M_{0,n}$.  For any cohomology classes  $\gamma_1, \ldots ,  
\gamma_n$  on  $X$,  one can construct a class
\begin{equation}
I^X_{0, \,  n ,  \,  \beta}   (\gamma_1 \otimes  \cdots  \otimes 
\gamma_n)  =  \eta_* \left({\rho_1}^*(\gamma_1)  \smallcup \cdots 
\smallcup {\rho_n}^*(\gamma_n)\right) \label{2.18}
\end{equation}
in the cohomology ring  $H^*(\M_{0,n})$.   These are called (tree-level, 
or genus zero) Gromov-Witten classes.  The number we 
denoted  $I_\beta(\gamma_1 \cdots  
\gamma_n)$  is the degree of the zero-dimensional component of 
this class, which they denote by $ \langle I^X_{0, \,  n ,  \,  \beta} 
\rangle  (\gamma_1 \otimes  \cdots  \otimes \gamma_n)$.  The 
intersections with divisors that we have carried out on  
$\M_{0,n}(X,\beta)$  can be carried out with the corresponding 
divisors on  $\M_{0,n}$;  this has the advantage that the 
intersections take place on a nonsingular variety.

One of the main goals of [K-M 1] and especially [K-M 2] is to show 
how Gromov-Witten classes can be reconstructed from the numbers 
obtained by evaluating them on the fundamental classes.  The idea is 
that a cohomology class in  $H^*(\M_{0,n})$  is known by evaluating 
it on the classes of the closures of the strata determined by the 
combinatorial types of the labeled trees.  As we saw and exploited 
for divisors, these numbers can be expressed in terms of the 
numbers  $I_\beta$  for the pieces making up the tree.

Kontsevich and Manin also allow cohomology classes of odd degrees, 
in which case one has to be careful with signs and the ordering of 
the terms.  For an interesting application to some Fano varieties, 
see [Bea]. 
 
Since the space  $H = H^*(X,{\mathbb Q})$  can be identified with its dual 
by Poincar\'e duality, the maps  $I^X_{0, \,  n ,  \,  \beta}$  can be 
regarded as maps
\begin{equation}
H_*(\M_{0,n+1})  \rarr  \text{Hom}
 (H^{\otimes n} \, , \, H ) .\label{2.19}
\end{equation}
Both of these, for varying  $n$,  have a natural operad structure, that 
on the first coming from all the ways to glue together labeled trees 
of  $\proj^1$'s  to form new ones, and the second from all the ways 
to compose homomorphisms.  Remarkably, the associativity 
(Theorem 4), is equivalent to the assertion that (\ref{2.19}) 
is a morphism of operads.  

The structure constants  $g^{i  j}$  put a metric on the cohomology 
space  $H^*(X,{\mathbb C})$;  with coordinates given by the basis for the 
cohomology, there is a (formal) connection given by the formula  
$A_{i  j}^k = \sum \Phi_{i j e} g^{e  k}$.  In this formalism of 
Dubrovin, the associativity translates to 
the assertion that this is a 
flat connection.

The numbers calculated here are part of a much more ambitious 
program described in [K-M 1] and [K].  The hope is to extend the story 
to varieties without the positivity assumptions made here, with 
some other construction of what should be the fundamental class of  
$\M_{g,n}(X,\beta)$.  (For varieties whose tangent bundles are not as 
positive as those considered here, the definition of the potential 
function  $\Phi$  is modified by multiplying the summands in 
(\ref{2.4}) by  $e^{-\int_\beta \omega}$,  for a K\"ahler class  
$\omega$,  in the hopes of making the power series converge on 
some open set of the cohomology space  $H$.)

Even if this program is carried out, however -- and associativity 
has been proved by symplectic methods [R-T] in some cases beyond 
those mentioned here$^4$ -- the interpretation cannot always be in 
enumerative terms as simple as those we have discussed, cf. [C-M].  
\footnotetext[4]{At this conference, J. Li lectured
on his work with G. Tian (see [L-T 2]) using cones
in vector bundles to construct a ``virtual fundamental
class'' to use in place of $[\M_{g,n}(X, \beta)]$ in
case $\M_{g,n}(X,\beta)$ has the wrong dimension. This
approach has been clarified and extended by K. Behrend and
B. Fantechi using the language of stacks ([B-F], [B]).
Algebraic computations in the non-convex case can
be found, for example, in [G], [G-P], [K].}
On the other hand, these ideas from quantum cohomology have 
inspired some recent work in enumerative geometry, even in cases 
where the associativity formalism does not apply directly, cf. [C-H 1] 
and [P]$^5$.
\footnotetext[5]{Also, [E-K], [K-Q-R].}

\end{document}